\documentclass[10pt]{smfart}

\usepackage{pict2e}

\usepackage{enumerate}
\def\tnorm#1{\left|\mskip -2 mu\left|\mskip -2
mu\left|#1\right|\mskip -2mu \right|\mskip -2mu\right|}
\usepackage{stmaryrd}
\usepackage{graphicx}
\usepackage{amssymb, amsmath, amsfonts}
\usepackage[latin1]{inputenc}
\usepackage{amsthm}
\theoremstyle{theorem}
\newtheorem{theo}{\textbf{Théorème}}[section]
\newtheorem{prop}[theo]{\textbf{Proposition}}
\newtheorem{defi}[theo]{\textbf{Définition}}
\newtheorem{coro}[theo]{\textbf{Corollaire}}
\newtheorem{lemm}[theo]{\textbf{Lemme}}
\newtheorem{rema}[theo]{\textbf{Remarque}}
\newenvironment{demo}{\noindent\textsc{Preuve.}~}{\hfill$\square$\bigbreak}

\def\bema{\begin{displaymath}}
\def\enma{\end{displaymath}}

\def\bemar{\begin{displaymath}\begin{array}{rcl}}
\def\enmar{\end{array}\end{displaymath}}

\def\bess{Bessaga-Pe\l czy\'nski}
\def\horm{Avakumovi\v c-Levitan-Hörmander }
\def\R{\mathbb{R}}
\def\C{\mathbb{C}}

\def\E{\mathbf{E}}
\def\P{\mathbf{P}}
\def\N{\mathbb{N}}
\def\Z{\mathbb{Z}}

\def\CC{\mathcal{C}}

\def\HH{\mathcal{H}}
\def\S{\mathbb{S}}
\def\hl{\pl^p(\R^d,\oplus E_{d,n})}

\def\tt{{}^t}
\def\tr{\textnormal{\textbf{tr}}}
\def\ep{\varepsilon}

\newcommand\comb[2]{ \left( \begin{array}{c} #1 \\ #2 \end{array}   \right) }
\newcommand\Sum{ \displaystyle\sum}
\newcommand\Int{ \displaystyle\int}
\newcommand\Prod{\displaystyle\prod}
\newcommand\Frac{ \displaystyle\frac}

\def\pun{\mathbf{1}}
\def\pl{\mathbf{PL}}

\title{Concentration et randomisation universelle de sous-espaces propres}
\author{Rafik Imekraz}
\address{Institut de Math\'ematiques de Bordeaux, UMR 5251 du CNRS, Universit\'e de Bordeaux 1, 351, cours de la Lib\'eration F33405 Talence Cedex, France     }
\email{rafik.imekraz@math.u-bordeaux1.fr}

\date{}
\begin{document}
\maketitle

\begin{abstract}
Nous développons une théorie de randomisation multidimensionnelle dans les espaces de Lebesgue $L^p$ à partir des inégalités de Kahane-Khintchine-Marcus-Pisier.
De façon précise, nous obtenons un théorème de type Maurey-Pisier en termes de matrices aléatoires, à savoir que la randomisation multidimensionnelle est universelle dans $L^p$.
Ensuite, se pose la question de déterminer des conditions nécessaires et suffisantes de convergence presque sûre des combinaisons aléatoires de fonctions propres propres du Laplacien dans les espaces de Lebesgue $L^p$ d'une variété riemannienne compacte (le célèbre théorème de Paley-Zygmund répond à la question pour les tores).
Nous proposons une nouvelle méthode qui permet de répondre à cette question lorsque les fonctions considérées ont des propriétés de concentration.
Comme application, nous étudions la théorie $L^p$ presque sûre des variétés compactes, harmoniques sphériques et oscillateur harmonique.
Nous écrirons aussi les injections de Sobolev probabilistes associées au sens de Burq et Lebeau et obtiendrons des solutions globales en régime sur-critique à l'équation cubique des ondes sur une variété compacte sans bord de dimension $3$.
\end{abstract}

\begin{altabstract}
We develop a theory of multidimensional randomization in Lebesgue spaces $L^p$ with the aid of Kahane-Khintchine-Marcus-Pisier inequalities. 
More precisely, we obtain a result in the spirit of Maurey-Pisier's theorem which involves random matrices and proves that the multidimensional randomization is universal in $L^p$.
Then, we deal with the question of studying necessary and sufficient conditions to get the almost sure convergence in $L^p$ of random linear combinations of  eigenfunctions in a compact Riemannian manifold (the famous Paley-Zygmund theorem gives the answer for tori).
We introduce a new method which solves the problem if the considered eigenfunctions have some concentration property.
We give several applications like the almost sure $L^p$ convergence for compact manifolds, spherical harmonics and the harmonic oscillator.
We also get probabilistic Sobolev embeddings in the sense of Burq and Lebeau and 
we obtain global solutions for the supercritical cubic wave equation on a boundaryless compact $3$-manifold.
\end{altabstract}

\setcounter{section}{-1}

\newpage
\begin{center}
\textbf{Ordre de lecture}
\end{center}

\setlength{\unitlength}{1cm}
\begin{center}
\begin{picture}(16,5)
\multiput(6.5,0.125)(0,0.5){10}{\line(0,1){0.25}}
\put(2,3.8){\makebox(2,1)[c]{Théorie} }
\put(8,3.8){\makebox(2,1)[c]{Applications} }
\put(1,2.5){\framebox(2,1)[c]{Partie \ref{theo-gene}}}
\put(3,3){\vector(1,-1){1}}
\put(3,3){\vector(4,0){4}}
\put(4,1.5){\framebox(2,1)[c]{Partie \ref{enonc-interpo}}}
\put(7,0.5){\framebox(2,1)[c]{Partie \ref{Lp-oscillo}}}
\put(6,2){\vector(4,0){4}}
\put(6,2){\vector(1,-1){1}}
\put(7,2.5){\framebox(2,1)[c]{Partie \ref{parti-comp}}}
\put(13,2.5){\framebox(2,1)[c]{Partie \ref{pzoll}}}
\put(9,3){\vector(4,0){4}}
\put(9,3){\vector(1,-1){1}}
\put(10,1.5){\framebox(2,1)[c]{Partie \ref{harmo-sph}}}
\end{picture}
\end{center}

\tableofcontents

\section{Introduction}

L'objet de ce travail est l'étude de la convergence dans les espaces $L^p(X)$, avec $p\in [1,+\infty[$, des combinaisons aléatoires de modes propres des opérateurs de type Laplace-Beltrami sur une variété riemannienne compacte $X$ (ou oscillateur harmonique $-\Delta+|x|^2$ sur $X=\R^d$).
 Ce type de problème prend sa source dans le théorème de Paley-Zygmund \cite{paley1930,paley1932} : si l'on considère une suite $(a_n)_{n\in \Z}$ appartenant à $\ell^2(\Z)$ et une suite $(\ep_n)_{n\in\Z}$ de variables aléatoires indépendantes suivant une loi $\frac{1}{2}$-Bernoulli à valeurs dans $\{-1,+1\}$ alors la série de Fourier aléatoire $\sum_{n\in \Z} \ep_n a_n e^{inx}$ définit presque sûrement un élément de $L^p(\mathbb{T})$ pour tout réel $p\in [1,+\infty[$.
La randomisation permet ainsi un gain d'intégrabilité alors qu'il n'y a évidemment aucun gain de régularité dans l'échelle des espaces de Sobolev $H^s(\mathbb{T})$.
Le théorème de Paley-Zygmund doit être comparé à l'injection de Sobolev $H^{\frac{1}{2}-\frac{1}{p}}(\mathbb{T})\subset L^p(\mathbb{T})$ valide pour tout réel $p\geq 2$.
  Les démonstrations du théorème de Paley-Zygmund utilisent généralement l'inégalité de Khintchine.
L'ouvrage \cite{kahane} de Kahane contient de nombreux résultats importants dans ce thème et introduit notamment des versions banachiques de l'inégalité de Khintchine, à savoir les inégalités de Kahane-Khintchine (voir plus loin \eqref{KK}).
Citons maintenant trois résultats connus concernant la randomisation dans les espaces de Lebesgue.

Le premier résultat, qui peut être vu comme une conséquence des inégalités de Kahane-Khintchine, énonce ceci : pour toute suite $
(u_n)_{n\in \N}$ de $L^p(X)$ où $X$ est un espace mesuré $\sigma$-fini, on a l'équivalence 
\begin{equation}\label{carac-uni}
\sum \ep_n u_n \quad \mbox{converge presque sûrement dans }L^p(X) \quad \Leftrightarrow \quad 
\sqrt{\sum_{n\in \N} |u_n|^2}\in L^p(X).
\end{equation}
Un théorème de Maurey permet d'ailleurs de remplacer $L^p(X)$ par n'importe quel treillis de Banach qui dispose de la propriété de cotype fini (\cite[pages  21-22]{maurey1973} et \cite[Theorem 1.d.6, Corollary 1.f.9]{linden2}).

Le deuxième résultat est ce que nous appellerons le théorème de randomisation universelle de Maurey-Pisier \cite[corollaire 1.3]{maurey-pisier76} (pour le distinguer des autres théorèmes de \cite{maurey-pisier76}). Ce dernier assure, avec les mêmes notations, que les convergences des séries aléatoires $\sum  \ep_n u_n$ et $\sum g_n u_n$ dans $L^p(X)$ sont équivalentes où $(g_n)_{n\in \N}$ désigne une suite i.i.d. de gaussiennes suivant une loi $\mathcal{N}_\C(0,1)$.
En fait, le théorème de randomisation universelle de Maurey-Pisier permet de remplacer $L^p(X)$ par un espace de Banach complexe de cotype fini et les gaussiennes $g_n$ par des variables aléatoires centrées plus générales (voir l'appendice \ref{cotype} pour des rappels et énoncés précis).

Le troisième résultat a été obtenu par 
Fig\`a-Talamanca et Rider qui ont pu remplacer le tore $\mathbb{T}$ par un groupe compact quelconque $G$ 
dans le théorème de Paley-Zygmund (\cite[Theorem 4]{figa1966} et \cite[Corollary 4]{figa1967}). On pourra aussi consulter \cite{figa1970}.
Dans \cite{pisier1981}, Marcus et Pisier ont résolu le problème de la convergence presque sûre des séries de Fourier dans $L^\infty(G)$ et leur analyse permet de retrouver les résultats dans l'échelle des espaces $L^p(G)$ avec $p\in [1,+\infty[$ à l'aide d'une version multidimensionnelle des inégalités de Kahane-Khintchine.
On qualifiera ces dernières d'inégalités de Kahane-Khintchine-Marcus-Pisier (voir plus loin \eqref{KKMP}).

Récemment, ces problèmes ont ressurgi en considérant les fonctions trigonométriques $x\mapsto e^{inx}$ non pas comme les caractères 
du groupe abélien compact $\mathbb{T}$ mais plutôt comme les modes propres de l'opérateur de Laplace-Beltrami $\frac{d^2}{dx^2}$.
Deux types de résultat ont notamment motivé cette résurgence :  

\begin{enumerate}[i) ]
\item les bases hilbertiennes aléatoires de modes propres ont des propriétés non triviales comme l'ergodicité quantique \cite{zeldi1992,PRT3} ou des estimations de normes $L^p$ bien meilleures que celles des bases hilbertiennes ``canoniques'' \cite{PRT1,burq-lebeau}.
\item l'étude des équations non-linéaires de type Schrödinger ou ondes sur une variété riemannienne compacte $X$ avec des conditions initiales à faible régularité dans les espaces de Sobolev $H^s(X)$ (on parle de régime sur-critique) est un problème difficile quand $\dim(X)$ est grand.
 La randomisation donne un gain d'intégrabilité qui permet de construire des solutions qui sont hors d'atteinte avec les méthodes déterministes actuelles.
Ces travaux concernent les constructions de mesures de Gibbs (voir les articles \cite{Bo94,Bo96,BTT,deng2012,tzvetkov08,bourg-bulutI,bourg-bulutII,bourg-bulut3d} et leurs références) ou les conditions initiales aléatoires plus générales \cite{burq-tz-inv2008-1,burq-tz-inv2008-2,burq-tz-jems2011,PRT2,BTT-global}.
Concernant ces derniers articles, on pourra consulter le séminaire Bourbaki \cite{bourbaki-anne-bou}.
\end{enumerate}

Les deux points précédents ont même été combinés dans \cite{suzzo-S3} : de Suzzoni a étudié l'équation cubique des ondes sur la sphère $\S^3$ à l'aide d'une base hilbertienne aléatoire de $L^2(\S^3)$.

Il est donc légitime d'étudier la randomisation non plus sur un groupe compact mais sur une variété riemannienne compacte $X$ (que l'on supposera toujours lisse, sans bord, munie de sa mesure volume $\mu$ et de dimension $d\geq 2$).
\'Etant données une famille orthonormée $(\phi_n)_{n\geq 0}$ de $L^2(X)$ constituée de fonctions propres de l'opérateur de Laplace-Beltrami $\Delta$ sur $X$ et une fonction 
\bema   \sum_{n\geq 0} a_n \phi_n \in L^2(X), \quad (a_n)\in \ell^2(\N), \enma
on s'intéresse à la convergence presque sûre dans $L^p(X)$ de la série aléatoire 
$  \sum \ep_n a_n \phi_n $.
Bien entendu, le critère \eqref{carac-uni} répond à la question de façon théorique mais il n'est pas du tout évident de le traduire en un comportement asymptotique sur la suite des coefficients $(a_n)_{n\in \N}$.
Sans surprise, cette étude est intimement liée à la suite des normes $\left\Vert \phi_n\right\Vert_{L^p(X)}$. 
En fait, nous verrons que c'est l'éventuelle concentration des fonctions $\phi_n$ qui rentre en jeu.
Rappelons les résultats connus.
Tzvetkov montre ce que l'on appellera par la suite une injection de Sobolev probabiliste (selon la terminologie introduite par Burq et Lebeau dans \cite{burq-lebeau}) : pour tout réel $p\geq 2$ il existe un réel explicite $\delta(d,p)>0$ tel que
\begin{equation}\label{sobo-tzv}  \sum_{n\geq 0} a_n \phi_n\in  H^{ \delta(d,p) }(X)\quad \Rightarrow \quad \Sum_{n\geq 0} \ep_n a_n \phi_n \quad \mbox{converge p.s. dans } L^p(X).\end{equation}
Avant de donner l'expression de $\delta(d,p)$, signalons que l'on a l'inégalité $\delta(d,p)\leq d\left(\frac{1}{2}-\frac{1}{p}\right)$ et donc \eqref{sobo-tzv} améliore l'injection de Sobolev déterministe $H^{d\left(\frac{1}{2}-\frac{1}{p}\right)}(X)\subset L^p(X)$.
Par construction (voir \cite[Theorem 4]{Tzvetkov4}), le nombre $\delta(d,p)$ vient des inégalités de Sogge \cite{sogge88} que nous rappelons : si l'on a $\Delta \phi_n=-\lambda_n^2 \phi_n$, avec $\lambda_n>0$ alors on a 
\begin{equation}\label{sogge}\forall n\in \N^\star \quad \left\Vert\phi_n \right\Vert_{L^p(X)}\leq C(X,p) \lambda_n^{\delta(d,p)} ,\quad \delta(d,p):=\left\{\begin{array}{lcl} \frac{d-1}{2}\left(\frac{1}{2}-\frac{1}{p}\right) & \mbox{ si } & 2\leq p\leq \frac{2(d+1)}{d-1} , \\
\frac{d-1}{2}-\frac{d}{p}& \mbox{ si } & \frac{2(d+1)}{d-1} \leq p \leq +\infty, \end{array}  \right.   \end{equation}
où $C(X,p)>1$ ne dépend que de $X$ et $p$.
Le cas $p=\infty$ est dû à Avakumovi\v c, Levitan et Hörmander \cite{horm68} et les inégalités \eqref{sogge} sont optimales pour $X=\S^d$.
Dans \cite[Theorem 1]{tzvetkov-ay}, Ayache et Tzvetkov obtiennent un éclairage gaussien du théorème de Paley-Zygmund, à savoir l'équivalence pour tout $p\in [2,+\infty[$ des deux assertions suivantes :
\begin{enumerate}[i) ]
\item les fonctions propres $\phi_n$ sont uniformément bornées dans $L^p(X)$,
\item pour toute suite complexe $(a_n)_{n\in \N}$, la fonction gaussienne aléatoire $\sum_{n\geq 0} g_n a_n \phi_n$ appartient presque sûrement à $L^p(X)$ si et seulement si $(a_n)_{n\in \N}$ appartient à $\ell^2(\N)$.
\end{enumerate}
Lorsque les fonctions $\phi_n$ ne sont pas uniformément bornées dans $L^p(X)$, Ayache et Tzvetkov posent la question de déterminer en fonction des coefficients $a_n$, la borne supérieure des réels $p\in [2,+\infty[$ tels que la série aléatoire $\sum g_n a_n \phi_n$ converge presque sûrement dans $L^p(X)$.
Sans aucune information sur les fonctions $\phi_n$, cette question est trop générale et
Ayache et Tzvetkov examinent les modes propres radiaux $\psi_n$ de l'opérateur Laplacien $\Delta$ avec condition de Dirichlet au bord sur la boule ouverte unité $\mathbb{B}_d(0,1)$ de $\R^d$.
Il s'avère que les fonctions $\psi_n$ ne sont pas uniformément bornées dans $L^p(\mathbb{B}_d(0,1))$ pour $p\gg 1$, et l'on en déduit l'existence de suites $(a_n)_{n\in \N}\in \ell^2(\N)$ telles que la série aléatoire gaussienne $\sum g_n a_n \psi_n$ diverge presque sûrement dans $L^p(\mathbb{B}_d(0,1))$.
Il s'agit là d'une différence majeure avec le théorème de Paley-Zygmund sur le tore $\mathbb{T}$. Bien que la boule $\mathbb{B}_d(0,1)$ ne soit pas compacte, l'analyse précédente se transfère sans problème aux fonctions propres zonales de l'opérateur de Laplace-Beltrami sur la sphère $\S^d$.
Dans le cas des fonctions $\psi_n$, Ayache et Tzvetkov obtiennent une réponse partielle \cite[Theorem 4]{tzvetkov-ay}.
Dans \cite{Grivaux} apparaissent deux idées qui vont jouer un rôle important dans notre travail : 
\begin{enumerate}[i) ]
\item d'une part, Grivaux remarque que le théorème de randomisation universelle de Maurey-Pisier assure l'équivalence des convergences des séries aléatoires $\sum \ep_n a_n \psi_n$ et $\sum g_n a_n \psi_n$.
\item d'autre part, Grivaux répond à la question de Ayache et Tzvetkov en utilisant la concentration des fonctions $\psi_n$ en l'origine $0\in \R^d$.
\end{enumerate}

Les deux précédentes idées ont été exploitées dans \cite{randomh} par Robert, Thomann et l'auteur pour l'oscillateur harmonique $-\Delta+|x|^2$ sur $\R^d$ qui est un exemple fondamental d'opérateur à spectre discret.
La nouveauté avec l'oscillateur harmonique est que l'on peut atteindre presque sûrement des fonctions dans $L^p(\R^d)$ à partir de distributions qui ne sont même pas des fonctions (on précisera ce point plus loin).

Ajoutons que, dans les papiers précédents, apparaît souvent l'astuce d'examiner des modes propres invariants par symétrie afin de réduire un problème multidimensionnel à un problème unidimensionnel. 
De façon indépendante, des travaux font intervenir la notion de randomisation multidimensionnelle sur une variété riemannienne compacte $X$ en tenant compte de la décomposition spectrale de l'opérateur de Laplace-Beltrami sur $L^2(X)$.
Cette notion apparaît sous des formes en apparence différente, dans les travaux de Shiffman et Zelditch \cite{ShZe,zeldi1992}, dans celui de Burq et Lebeau \cite{burq-lebeau} ainsi que dans celui de Robert et Thomann \cite{PRT1} (pour l'oscillateur harmonique sur $\R^d$) avec des arguments de ``concentration de la mesure'' et de ``grandes déviations".
C'est dans l'article \cite{burq-lebeau} que le terme ``injection de Sobolev probabiliste" apparaît pour la première fois pour exprimer rigoureusement le gain d'intégrabilité obtenu par la randomisation.
Même si cela n'est pas explicitement écrit, il nous semble qu'un des intérêts du papier \cite{burq-lebeau} est précisément de s'émanciper de l'astuce, citée plus haut, de considérer des modes propres invariants par symétrie.
C'est ainsi que Burq et Lebeau obtiennent un résultat d'existence locale pour l'équation semi-linéaire des ondes sur une variété riemannienne compacte en régime sur-critique et  leurs solutions ne sont pas spectralement supportées par des sous-suites particulières de fonctions propres.

Nous présenterons un cadre unifié qui nous permettra d'atteindre les trois buts suivants : 

\textbf{But 1 :} écrire une théorie de randomisation multidimensionnelle dans les espaces $L^p$ avec des hypothèses optimales de moment en termes de matrices aléatoires. Nous considérerons des espaces mesurés $\sigma$-finis quelconques.

\textbf{But 2 :} proposer une nouvelle méthode d'interpolation et de dualité afin de résoudre une version plus forte du problème de Ayache et Tzvetkov lorsque les fonctions $\phi_n$ ont des propriétés de concentration.
En l'occurrence, nos arguments permettent d'obtenir des conditions nécessaires et suffisantes de convergence presque sûre de la série aléatoire $\sum \ep_n a_n \phi_n$ dans $L^p(X)$ en fonction des coefficients $a_n$.
Ces conditions, quoique plus délicates à obtenir que le calcul de la borne supérieure des réels $p\geq 1$ convenables, s'exprimeront très facilement en assimilant chaque fonction $\phi_n$ à sa restriction sur sa zone de concentration.
 En application, nous obtiendrons des ``injections de Sobolev probabilistes" au sens de Burq et Lebeau qui seront optimales dans diverses configurations.

\textbf{But 3 :} démontrer l'existence de ``beaucoup'' de solutions globales en régime sur-critique 
à l'équation cubique des ondes sur une variété riemannienne compacte sans bord de dimension $3$ dans l'esprit de \cite{burq-tz-jems2011,burq-lebeau,suzzo-S3} à l'aide de nos nouvelles estimations.

La trame logique de nos contributions est décrite dans les parties \ref{theo-gene} et \ref{enonc-interpo} (modèle théorique) et les parties \ref{parti-comp}, \ref{pzoll}, \ref{harmo-sph} et \ref{Lp-oscillo} (applications).
Afin de faciliter la lecture de nos applications, nous avons fait le choix d'écrire la quasi-totalité des preuves après la partie \ref{Lp-oscillo}.
Nous avons par ailleurs ajouté quelques appendices.

La partie \ref{theo-gene} est dévolue à l'étude théorique d'un procédé de randomisation multidimensionnelle associé à une suite $(E_n)_{n\in\N}$ de sous-espaces de \textbf{dimension finie} de $L^2(X)$,
où $X$ est un espace mesuré $\sigma$-fini muni d'une mesure $\mu$.
Pour tout $p\in [1,+\infty[$, nous définirons un espace de Lebesgue probabiliste, noté $\mathbf{PL}^p(X,\oplus E_n)$, qui sera un sous-espace vectoriel de $\prod_n E_n$.
Le cas $\dim(E_n)=1$, pour tout $n\in \N$, est le plus simple à décrire : une suite $(u_n)_{n\in \N}$ appartient à l'espace $\pl^p(X,\oplus E_n)$ si et seulement si la série aléatoire $\sum \ep_n u_n$ converge presque sûrement dans $L^p(X)$.
Ces suites $(u_n)_{n\in \N}$ sont caractérisées par le critère déterministe \eqref{carac-uni}.
Dans la théorie unidimensionnelle, les inégalités de Kahane-Khintchine, que nous rappelons, sont très importantes :
pour tous réels $p>q\geq 1$, il existe une constante $K_{p,q}\geq 1$ telle que, pour tout espace de Banach $B$ et tous éléments $u_0,\dots,u_N\in B$, 
les moments des variables aléatoires $\left\Vert\sum_{n=0}^N \ep_n u_n\right\Vert_B$ sont du même ordre de grandeur au sens suivant
\begin{equation}\label{KK}   \E\left[\left\Vert \sum_{n=0}^N \ep_n u_n \right\Vert_B^q \right]^{\frac{1}{q}}\leq   \E\left[\left\Vert \sum_{n=0}^N \ep_n u_n \right\Vert_B^p \right]^{\frac{1}{p}} \leq K_{p,q}     \E\left[\left\Vert \sum_{n=0}^N \ep_n u_n \right\Vert_B^q \right]^{\frac{1}{q}}.\end{equation}
D'après Borell et Latala-Oleszkiewicz, on peut estimer précisément les nombres $K_{p,q}$ mais nous n'aurons besoin que des inégalités $K_{p,q}\leq K_{p,1}\leq K\sqrt{p}$ valides pour une constante universelle $K>0$ (\cite[Part 3.III]{queff}).
Pour étudier la randomisation dans $L^p(X)$ dans un cadre multidimensionnel, nous devrons remplacer de façon adéquate les variables de Bernoulli $\ep_n$ par des matrices aléatoires qui suivent une loi uniforme dans un groupe orthogonal ou unitaire.
L'idée de base est tout simplement de faire tourner chaque $u_n$ de façon uniforme dans la sphère de $E_n$ de centre $0$ et de rayon $\left\Vert u_n \right\Vert_{L^2(X)}$.
Cela apparaît dans les travaux de 
Shiffman-Zelditch, Burq-Lebeau et Robert-Thomann.
Mais notre approche est très différente et repose sur la version multidimensionnelle des inégalités \eqref{KK}, à savoir les inégalités de Kahane-Khintchine-Marcus-Pisier que nous avons mentionnées plus haut.
L'idée maîtresse de notre article est que ces inégalités sont précisément celles qui permettent d'aborder de façon satisfaisante la randomisation multidimensionnelle avec des hypothèses optimales de moments sur un espace mesuré $\sigma$-fini quelconque.
Un des intérêts esthétiques de notre approche est l'unification de la randomisation unidimensionnelle classique et la randomisation multidimensionnelle abordée par les précédents auteurs.
Pour exprimer un critère déterministe analogue à \eqref{carac-uni}, on doit définir la fonction spectrale $e(n,\cdot)$ du sous-espace $E_n$ : 
\bema  \forall x\in X \quad e(n,x) := \sup\{ |u(x)|^2 , \quad u\in E_n \quad  \left\Vert u\right\Vert_{L^2(X)}=1\}. \enma
Si $\phi_{n,1},\dots,\phi_{n,d_n}$ est une base hilbertienne de $E_n$ alors on vérifie aisément l'égalité suivante (voir \eqref{defi-en}) 
\bema e(n,x)=|\phi_{n,1}(x)|^2+\dots+|\phi_{n,d_n}(x)|^2. \enma
En particulier, on a 
\bema
\int_{X} e(n,x)d\mu(x)= \dim(E_n).
\enma
On verra que toute information liée à l'aléa est codée dans le comportement des fonctions $\frac{e(n,\cdot)}{\dim(E_n)}$.
Par exemple, le théorème \ref{mplp} donnera la caractérisation suivante des espaces de Lebesgue probabilistes :
\bema
(u_n)_{n\in \N}\in \pl^p(X,\oplus E_n) \quad \Leftrightarrow \quad \sqrt{\sum_{n\in \N} \left\Vert u_n\right\Vert_{L^2(X)}^2 \frac{e(n,x)}{\dim(E_n)}}\in L_x^p(X).
\enma
Plus important, le théorème \ref{mplp} est l'analogue multidimensionnel du théorème de randomisation universelle de Maurey-Pisier sur $L^p(X)$ : nous serons capables d'utiliser des matrices aléatoires beaucoup plus générales que celles qui suivent une loi uniforme à valeurs dans un groupe orthogonal ou unitaire. En précisant une inégalité de matrices aléatoires due à Latala, nous donnerons des exemples naturels de matrices aléatoires convenables.
Expliquons maintenant pourquoi les arguments de concentration de la mesure et les estimées de grandes déviations sont inutilisables dans notre article : ces derniers nécessitent d'utiliser des variables aléatoires sous-gaussiennes (qui ont donc des moments de tout ordre), or nous travaillons avec un nombre fini de moments.
Ces arguments sont utilisés par les précédents auteurs (sauf Marcus et Pisier dans \cite{pisier1981}) pour estimer des probabilités mais nous préférerons estimer des intégrales et des espérances (ce qui ne nécessitera que le théorème de Fubini et les inégalités de Kahane-Khintchine-Marcus-Pisier).
Nos arguments sont inspirés du livre \cite{pisier1981} mais diffèrent en ceci que la théorie $L^\infty$ est essentiellement métrique tandis que la théorie $L^p$, avec $p<\infty$, autorise à travailler dans des espaces mesurés. On verra même dans l'annexe \ref{treillis} que notre analyse n'utilise en réalité que la structure de treillis de Banach de cotype fini de $L^p(X)$.

Dans la partie \ref{enonc-interpo}, on s'intéresse aux propriétés de dualité et d'interpolation des espaces de Lebesgue probabilistes $\pl^p(X,\oplus E_n)$.
A priori, on s'attend à ce que le dual de $\pl^p(X,\oplus E_n)$ soit $\pl^{\frac{p}{p-1}}(X,\oplus E_n)$ et l'on espère que les espaces $\pl^p(X,\oplus E_n)$ soient stables par interpolation complexe et réelle.
En fait, cela s'avérera faux pour la randomisation des fonctions zonales sur la sphère $\S^d$ (voir la partie \ref{harmo-sph}) : on n'a pas de dualité et l'interpolation se réalise seulement si $p$ parcourt $]\frac{2d}{d-1},+\infty[$.
Nous avons cependant des résultats positifs : les théorèmes \ref{theo-interpo} et \ref{theo-dual} donnent des hypothèses pratiques sur les fonctions spectrales $e(n,\cdot)$ qui assurent les propriétés attendues de dualité et d'interpolation.
La preuve montrera que $\pl^p (X,\oplus E_n)$ est un rétracte de l'espace de Bochner-Lebesgue $L^p(X,\ell^2(\N))$, ce qui nous ramènera à ces derniers espaces.
Mais pour y arriver, nous serons obligés de justifier la continuité de certains projecteurs de $L^p(X,\ell^2(\N))$ à l'aide d'un nouveau critère de $\mathcal{R}$-bornitude sur $L^p(X)$.
Les notions de rétracte et de $\mathcal{R}$-bornitude, que nous rappellerons plus loin, peuvent paraître abstraites mais elles sont incontournables sous une forme ou une autre (cela est formalisé par la proposition \ref{Lamba-surj}).
C'est à ce moment qu'interviendra l'hypothèse abstraite
\begin{equation}\label{lpin}
\sup\limits_{n\in  \N} \frac{\sqrt{e(n,x)}}{\left\Vert \sqrt{e(n,\cdot)}\right\Vert_{L^p(X)}}\in L_x^{p,\infty}(X),
\end{equation}
où $L_x^{p,\infty}(X)$ est l'usuel espace $L^p$ faible.
La contribution la plus importante ce travail est la mise en évidence que \eqref{lpin} est vérifié dans des exemples importants issus de la physique mathématique où les fonctions $e(n,\cdot)$ se concentrent sur une même région de $X$.
Revenons au problème de Ayache et Tzevtkov.
Sous sa forme originelle, il se reformule ainsi : 
si l'on considère une suite $(a_n)_{n\in \N}\in \ell^2(\N)$, peut-on déterminer la borne supérieure des réels $p\geq 1$ tels que $(a_n\phi_n )\in \pl^p(X,\oplus \C \phi_n )$ ?
Nous allons étudier la version plus forte : peut-on décrire explicitement $\pl^p(X,\oplus \C \phi_n )$ ?
 Si $|\phi_n|$ se concentre sur une partie $A_n\subset X$ avec une amplitude $\alpha_n>0$, alors il est légitime d'espérer l'équivalence suivante pour toute suite complexe $(a_n)_{n\in \N}$ 
\bema
\sqrt{\sum_{n\in\N} |a_n \phi_n|^2  }\in L^p(X) \quad \Leftrightarrow \quad
\sqrt{\sum_{n\in\N} |a_n \alpha_n \pun_{A_n}|^2  }\in L^p(X).
\enma
En partitionnant les parties $A_n$ en sous-parties disjointes, la condition du membre droit est explicite et résoudrait le problème grâce au critère \eqref{carac-uni}.
Si l'approximation de $|\phi_n|$
par $\alpha_n \pun_{A_n}$ n'est pas bonne, alors l'équivalence précédente est difficile à prouver. L'hypothèse \eqref{lpin}, avec $e(n,x)=|\phi_n(x)|^2$, a vocation à mesurer le reste de cette approximation et les arguments de dualité et l'interpolation sont alors utiles pour simplifier les calculs.

La partie \ref{parti-comp} est basique et prépare les parties \ref{pzoll} et \ref{harmo-sph}.
On étudie le cas où $X$ est une variété riemannienne compacte sans bord munie de sa mesure volume $\mu$.
 Les sous-espaces $E_n$ seront, par exemple, des sous-espaces propres de l'opérateur de Laplace-Beltrami $\Delta$.
Les espaces $\pl^p (X,\oplus E_n)$ s'identifieront à des sous-espaces de distributions sur $X$.
Ensuite, la proposition \ref{inje-mani} donnera une condition nécessaire et suffisante sur les fonctions spectrales $e(n,\cdot)$ afin d'obtenir un résultat de type Paley-Zygmund, à savoir que 
les espaces $\pl^p (X,\oplus E_n)$ coïncident avec $L^2(X)$.
Cela se produira si les densités de probabilité $\frac{1}{\dim(E_n)} e(n,\cdot)$ se localisent uniformément sur la variété compacte $X$ ou si les isométries de $X$ agissent de façon transitive (par exemple si $X$ est un produit fini de sphères ou de groupes de Lie).

Dans la partie \ref{pzoll}, on continue l'étude sur une variété riemannienne compacte $X$.
Suivant une idée de Burq et Lebeau, on montrera qu'il existe un choix naturel de sous-espaces $E_n\subset L^2(X)$ pour lesquels le théorème de Paley-Zygmund est vérifié.
En caricaturant, on peut dire que presque toute fonction dans $L^2(X)$ appartient à $L^p(X)$ pour tout $p\in [1,+\infty[$.
Nous retrouverons ainsi un résultat obtenu dans \cite{burq-lebeau} mais nous obtiendrons aussi l'information duale : l'espace des distributions sur $X$ qui appartiennent presque sûrement à $L^1(X)$ est précisément $L^2(X)$.
Il s'agit du théorème de Paley-Zygmund poussé dans ces derniers retranchements sur une variété riemannienne compacte sans bord quelconque.
Avec à peine un peu plus de régularité, presque toute fonction appartenant à $\cup_{s>0} H^s(X)$ est continue (cela généralise un cas connu sur $\mathbb{T}^d$ et améliore de $\frac{\dim X}{2}$ dérivées l'injection de Sobolev).
On utilisera ce dernier fait et la randomisation multidimensionnelle pour l'équation cubique des ondes en régime sur-critique sur une variété riemannienne compacte.
Comparons notre résultat aux travaux cités plus haut :
\vspace*{0.2 cm}
\begin{center}
\begin{tabular}{|c|c|c|c|c|} \hline 
variété & existence &  nature de la  & hypothèse de  & référence  \\ 
  riemannienne & temporelle  & randomisation & moments &  \\\hline
sans bord, $\dim 3$ & locale & unidimensionnelle & moments pairs & \cite{burq-tz-inv2008-1}\\
sans bord, $\dim 3$ & locale & multidimensionnelle & grandes déviations & \cite[Théorème 14]{burq-lebeau}\\
$\mathbb{T}^3$ & globale & unidimensionnelle & grandes déviations & \cite{burq-tz-jems2011} \\
$\mathbb{S}^3$ & globale & unidimensionnelle & grandes déviations & \cite{suzzo-S3}  \\ \hline  
sans bord, dim $3$ & globale & multidimensionnelle & moment d'ordre $3$ & notre résultat \\ \hline
\end{tabular}
\end{center}
\vspace*{0.2 cm}

Notre partie déterministe est la même que dans \cite{burq-tz-jems2011} et \cite{suzzo-S3}.
La randomisation multidimensionnelle nous libérera de la recherche d'une base de fonctions propres uniformément bornées dans tout espace $L^p(X)$.
En fait, toute base hilbertienne (non forcément constituée de modes propres) adaptée à la décomposition $L^2(X)=\bigoplus_{n\geq 0} E_n$ sera convenable.
Cette idée apparaît sous une forme différente dans \cite[Théorème 14]{burq-lebeau} mais 
notre randomisation de conditions initiales est plus générale et fait intervenir des matrices aléatoires dont les moments d'ordre $3$ sont uniformément bornés.
\`A notre connaissance, \cite{burq-tz-inv2008-1} est l'un des seuls articles où apparaît l'idée d'utiliser des variables aléatoires avec une hypothèse de finitude de moments (au lieu d'invoquer des estimations de grandes déviations).
Mais le lemme 4.2 de \cite{burq-tz-inv2008-1} oblige à considérer des moments qui sont nécessairement des entiers pairs.
Nos arguments probabilistes sont différents et prouvent que le bon moment est précisément l'ordre de la non-linéarité, c'est-à-dire $3$.

Le lecteur devra attendre la partie \ref{harmo-sph} pour avoir des exemples concrets où les espaces $\pl^p(X,\oplus E_n)$ sont différents de $L^2(X)$.
Nous étudierons en détail les deux exemples les plus importants d'harmoniques sphériques sur la sphère $\S^d$, $d\geq 2$, à savoir la suite $(Y_n)_{n\geq 1}$ des fonctions ``de plus haut poids" qui se concentrent sur une géodésique et la suite $(Z_n)_{n\geq 1}$ des fonctions zonales qui se concentrent sur deux pôles opposés.
Les fonctions $Y_n$ ne vérifient pas le théorème de Paley-Zygmund :
la fonction $\sum_{n\geq 2} \frac{1}{\sqrt{n}\ln(n)} Y_n$ appartient à $L^2(\S^d)$ tandis que la série aléatoire $\sum_{n\geq 2} \frac{\ep_n}{\sqrt{n}\ln(n)} Y_n$ diverge presque sûrement dans $L^p(\S^d)$ pour tout $p>2$.
Quant aux fonctions $Z_n$, elles sont uniformément bornées dans $L^p(\S^d)$ pour tout $p\in [1,\frac{2d}{d-1}[$, et l'on a donc d'après Ayache-Tzvetkov la convergence presque sûre de la série aléatoire $\sum \ep_n a_n Z_n$ dans $L^p(\S^d)$ pour toute suite 
$(a_n)_{n\geq 1}\in \ell^2(\N^\star)$.
Dans le cas $p>\frac{2d}{d-1}$, le meilleur résultat atteignable avec les méthodes connues serait la formule de Grivaux (\cite{Grivaux} et \eqref{form-gene}) qui exprime de façon explicite le nombre
 \bema
p_c:=\sup\left\{p\geq\frac{2d}{d-1},  \quad 
\sum_{n\geq 1} a_n Z_n \in \pl^p(\S^d,\oplus \C Z_{n})  \right\},
\enma
en fonction des coefficients $(a_n)_{n\geq 1}$. Cela nous donne le schéma
\setlength{\unitlength}{1cm}
\begin{center}
\begin{picture}(10,3.8)
\put(0,1){\vector(1,0){10}}
\put(0,1){\line(0,1){2.5}}
\put(1,1){\line(0,1){0.5}}
\put(3,1){\line(0,1){0.5}}
\put(5,1){\line(0,1){2.5}}
\put(0,1.5){\makebox(5,2){\shortstack{$\sum \ep_n a_n Z_n $ converge \\ p.s. dans $L^p(\S^d)$}}}
\put(5,1.5){\makebox(5,2){\shortstack{$\sum \ep_n a_n Z_n $ diverge \\ p.s. dans $L^p(\S^d)$}}}
\put(10,0.5){\makebox(1,1){$p$}}
\put(-0.5,0){\makebox(1,1){$1$}}
\put(0.5,0){\makebox(1,1){$2$}}
\put(2.5,0){\makebox(1,1){$\frac{2d}{d-1}$}}
\put(4.5,0){\makebox(1,1){$p_c$}}
\put(9.5,0){\makebox(1,1){$+\infty$}}
\end{picture}
\end{center}

Il reste à décider ce qui se passe en $p_c$, c'est-à-dire si $p_c$ est un maximum ou non.
Nous  justifierons rigoureusement que l'on peut assimiler $|Z_n|$ à sa restriction $\widetilde{Z}_n$ sur la boule de rayon $\frac{1}{n}$ centrée en l'un des pôles de symétrie afin d'obtenir l'encadrement suivant pour tout $p>\frac{2d}{d-1}$ : 
\bema\quad  \left\Vert \sqrt{\sum_{n\geq 1} |a_n \widetilde{Z}_n|^2}\right\Vert_{L^p(\S^d)}\leq \left\Vert \sqrt{\sum_{n\geq 1} |a_n Z_n|^2}\right\Vert_{L^p(\S^d)}\leq C(d,p)\left\Vert \sqrt{\sum_{n\geq 1} |a_n \widetilde{Z}_n|^2}\right\Vert_{L^p(\S^d)},
\enma
où $C(d,p)\geq 1$ ne dépend que de $d$ et $p$.
Comme les normes dans $L^{\frac{p}{2}}(\S^d)$ des fonctions $\sum_{n\geq 1} |a_n \widetilde{Z}_n|^2$ sont calculables explicitement, 
le critère \eqref{carac-uni} résoudra complètement la question si l'on a $p_c>\frac{2d}{d-1}$.
Nos preuves utilisent les propriétés de concentration des fonctions $Y_n$ et $Z_n$, les résultats abstraits d'interpolation et de dualité développés dans la partie \ref{enonc-interpo} et en particulier le théorème de Marcinkiewicz afin de corriger des facteurs logarithmiques : une bornitude $L^p\rightarrow L^p$ équivaut à une bornitude $L^p\rightarrow L^{p,\infty}$ dès lors que $p$ parcourt un intervalle ouvert.
Grâce à leurs descriptions, les espaces de Lebesgue probabilistes des fonctions $Y_n$ et $Z_n$ pourront être comparés aux espaces de Sobolev $H^s(\S^d)$ (c'est-à-dire que l'on écrira des injections de Sobolev probabilistes).
Il apparaîtra que l'exposant $\delta(d,p)$ de l'injection de Sobolev probabiliste \eqref{sobo-tzv} de Tzvetkov est optimal.

Dans la partie \ref{Lp-oscillo}, nous étudions la randomisation multidimensionnelle de la décomposition hilbertienne 
$ L^2(\R^d)=\bigoplus\limits_{n\in \N} E_{d,n}$
 en sous-espaces propres
de l'oscillateur harmonique multidimensionnel $-\Delta+|x|^2$.
Au niveau de l'analyse du spectre, l'oscillateur harmonique est souvent comparé aux sphères et surtout aux variétés de Zoll.
Cependant, l'étude des fonctions propres est très différente.
L'analyse $L^p$ probabiliste de l'oscillateur harmonique \textbf{multidimensionnel} a été entamée dans \cite{randomh} avec le procédé de randomisation \textbf{unidimensionnelle} en utilisant seulement des estimées dans $L^p(\R^d)$ de la fonction spectrale de $E_{d,n}$.
Pour remédier à la nature multidimensionnelle de l'oscillateur harmonique et obtenir des résultats indépendants d'une base hilbertienne spécifique, une condition de contrôle (``squeezing condition'') a été imposée sur les coefficients des séries aléatoires dans les papiers \cite{PRT2,PRT1,PRT3} et a été améliorée dans \cite{randomh} en tenant compte du principe de contraction.
 C'est cette condition de contrôle qui fait apparaître la fonction spectrale dans le cadre de la randomisation unidimensionnelle.
Notre approche nous permet d'obtenir les nouveaux points : 

$\bullet$ la proposition \ref{ed-maj} dit, en première approximation, que la fonction spectrale du sous-espace propre $E_{d,n}$ se concentre uniformément sur la boule $\mathbb{B}_d(0,\sqrt{2n+1})$. Ce résultat, essentiellement optimal, nous permettra de décrire les espaces $\pl^p(\R^d,\oplus E_{d,n})$.
Comme fait frappant, les inclusions de ces espaces sont inversées par rapport à celles attendues pour les variétés compactes : si une distribution tempérée appartient presque sûrement à $L^{p_1}(\R^d)$ alors elle appartient presque sûrement à tous les espaces $L^{p_2}(\R^d)$ avec $p_2\in [p_1,+\infty[$. Le théorème de Paley-Zygmund est donc vérifié et l'on peut même parler d'un effet régularisant de l'oscillateur harmonique.

$\bullet$ nous écrirons les injections de Sobolev probabilistes de la famille $(\mathcal{H}^s(\R^d))_{s\in \R}$ des espaces de Sobolev 
de l'oscillateur harmonique.
Pour tout réel $p\in ]2,+\infty[$, rappelons que l'on a l'injection de Sobolev déterministe $\mathcal{H}^{d\left(\frac{1}{2}-\frac{1}{p}\right)}(\R^d)\subset L^p(\R^d)$.
Le théorème 2.2 de \cite{randomh} montre, sous les conditions de contrôle évoquées ci-dessus, que presque toute distribution appartenant à un sous-espace strict, noté   
$\mathcal{Z}_{\varphi}^{-d\left(\frac{1}{2}-\frac{1}{p}\right)}(\R^d)$, de $\mathcal{H}^{-d\left(\frac{1}{2}-\frac{1}{p}\right)}(\R^d)$ appartient à $L^p(\R^d)$.
Pour la randomisation multidimensionnelle, les conditions techniques de contrôle sont inutiles et le théorème \ref{inje-sobo} donnera l'inclusion $\mathcal{H}^{-d\left( \frac{1}{2}-\frac{1}{p}\right)}(\R^d)\subset \pl^p(\R^d,\oplus E_{d,n})$, c'est-à-dire que presque toutes les distributions dans $\mathcal{H}^{-d\left(\frac{1}{2}-\frac{1}{p}\right)}(\R^d)$ appartiennent à $L^p(\R^d)$.
Nous invoquerons aussi le théorème \ref{theo-dual} de dualité développé dans la partie \ref{enonc-interpo} afin d'obtenir gratuitement l'inclusion duale $\pl^{\frac{p}{p-1}}(\R^d,\oplus E_{d,n})\subset \mathcal{H}^{d\left(\frac{1}{2}-\frac{1}{p}\right)}(\R^d)$.

Les appendices contiennent des rappels sur des notions classiques dans la théorie des espaces de Banach qui permettent de donner un éclairage sur certains de nos résultats : type, cotype, théorème de Maurey-Pisier, treillis complet de cotype fini et suites basiques inconditionnelles.

\vspace*{1 cm}

En résumé, la randomisation multidimensionnelle avec des matrices aléatoires est un point de vue intrinsèque qui évite la recherche de bases hilbertiennes spécifiques ou les conditions de contrôle sur des coefficients.
En outre, cela permet d'unifier plusieurs cadres différents avec des hypothèses optimales de moment.
Enfin, la méthode d'interpolation et de dualité des espaces de Lebesgue probabilistes permet d'obtenir plus ou moins facilement des résultats que l'on pourrait probablement obtenir directement mais en faisant usage d'astuces techniques et spécifiques aux situations rencontrées.

\vspace*{0.3 cm}

\textbf{Conventions.} Dans cet article, on fera les conventions suivantes 
\begin{enumerate}[$\bullet$]
\item la lettre $d$ sera réservée pour désigner la dimension des espaces en jeu ($\R^d$, $\S^d$ ou plus généralement une variété) et sera supérieure ou égale à $2$.
\item on notera $C$ une constante universelle supérieure ou égale à $1$ qui peut changer d'une ligne à l'autre.
\item on notera $C(d,p)$ une constante supérieure ou égale à $1$ qui ne dépend que des paramètres $d$ et $p$ et qui peut changer d'une ligne à l'autre.
\item si $A$ et $B$ sont deux nombres réels, la formule $A\lesssim B$ signifiera qu'il existe une constante universelle $C\geq 1$ telle que $A\leq CB$.
De même, l'assertion $A\gtrsim B$ signifie $B\lesssim A$.
\item on écrira $A\simeq B$ si l'on a $A\lesssim B$ et $B\lesssim A$.
\item si l'on fait intervenir des paramètres $d$ et $p$, alors on notera $A\lesssim_{d,p} B$ pour signifier qu'il existe une constante $C(d,p)\geq 1$ qui ne dépend que de $d$ et $p$ telle que $A\leq C(d,p)B$. On définit de même les symboles $\gtrsim_{d,p}$ et $A\simeq_{d,p}$.
\end{enumerate}

Par exemple, une assertion de la forme 
\bema \forall n\in \N \quad  A_n \simeq_{d,p} B_n ,\enma
signifiera 
\bema
\exists C(d,p) \geq 1\quad \forall n\in \N \quad  \frac{A_n}{C(d,p)} \leq B_n \leq C(d,p)A_n.
\enma

\section{Universalité de la randomisation multidimensionnelle dans $L^p$}\label{theo-gene}

Dans la suite, $(\Omega,\mathcal{F},\P)$ désignera un espace probabilisé de référence (par exemple $\Omega=[0,1]$ muni de la mesure de Lebesgue sur la tribu borélienne) et l'on notera $\omega$ les éléments de $\Omega$ ainsi que $\E$ l'opérateur d'espérance.
Sauf précision ultérieure, $X$ désignera un espace mesuré $\sigma$-fini muni d'une mesure $\mu$. Les espaces de Banach $L^p(X)$, avec $p\geq 1$, seront considérés sur le corps $\C$.
On pourrait aussi raisonner avec le cas réel mais le cas complexe nous paraît plus général et plus naturel.
Le moyen le plus basique de randomiser une suite $(u_n)_{n\geq 0}$ d'éléments de $L^p(X)$, est de considérer la série aléatoire
\bema     \sum_{n\geq 0} \ep_n u_n  , \enma
où $(\ep_n)_{n\geq 0}$ est une suite de variables aléatoires indépendantes de loi $\frac{1}{2}$-Bernoulli à valeurs dans $\{-1,1\}$. En fait, on peut remplacer cette loi par beaucoup d'autres sans changer la nature de la série aléatoire.
Afin de situer le résultat de randomisation que nous allons présenter plus loin, commençons par énoncer dans le cas unidimensionnel un théorème plus ou moins bien connu vis-à-vis de l'optimalité de ses hypothèses.

\begin{theo}\label{maurey-pi}
Considérons un réel $p\in [1,+\infty[$ et une suite de variables aléatoires $X_n:\Omega\rightarrow \R$ indépendantes, centrées et vérifiant
\begin{equation}\label{hypo-minim}
 \sup\limits_{n\in \N} \E [|X_n|^{\max(2,p)}]<+\infty   \quad  \mbox{et} \quad \inf\limits_{n\in \N} \E[|X_n|]>0.
\end{equation}
Pour toute suite $(u_n)_{n\geq 0}$ de $L^p(X)$, les propriétés suivantes sont équivalentes : 
\begin{enumerate}[i) ]
\item la fonction $x\mapsto \sqrt{\Sum_{n\geq 0} |u_n(x)|^2}$ appartient à $L^p(X)$,
\item la série aléatoire $\sum \ep_n u_n $ converge presque sûrement dans $L^p(X)$,
\item la série aléatoire $\sum X_n u_n $ converge presque sûrement dans $L^p(X)$.\end{enumerate}
L'exposant $\max(2,p)$ est optimal en toute généralité dans \eqref{hypo-minim} (par exemple si $L^p(X)=\ell^p(\N)$).
Si l'on considère des variables aléatoires $X_n$ à valeurs complexes, alors il faut remplacer l'hypothèse ``centrée'' par ``symétrique au sens complexe'' : pour tout réel $\alpha\in \R$ les variables $e^{i\alpha} X_n$ et $X_n$ ont la même loi.
\end{theo}

La considération de variables aléatoires $X_n$ non nécessairement identiquement distribuées prendra tout son sens lorsque nous examinerons plus loin des matrices aléatoires de tailles différentes (qui ne peuvent donc pas avoir la même loi).
Au niveau de la théorie des espaces de Banach, l'équivalence i) $\Leftrightarrow$ ii) est classique et peut être montrée grâce aux inégalités de Kahane-Khintchine (dans $L^p(X)$ et $\C$) et le théorème de Fubini.
L'équivalence ii) $\Leftrightarrow$ iii) prend certainement sa source dans \cite[corollaire 1.3]{maurey-pisier76} où $L^p(X)$ est généralement remplacé par un espace de Banach de cotype fini. L'exposant $\max(2,p)$ vient d'ailleurs du fait que $L^p(X)$ est un espace de Banach de cotype $\max(2,p)$.
On trouvera dans l'appendice \ref{cotype} les définitions et résultats classiques sur les espaces de cotype fini, notamment le théorème de randomisation universelle de Maurey-Pisier (Théorème \ref{theo-maupi}) qui éclaire le théorème \ref{maurey-pi} au prix d'un léger affaiblissement de la condition \eqref{hypo-minim}.
Comme nous travaillons sur $L^p(X)$ et non sur un espace de Banach de cotype fini quelconque, nous nous sommes rendu compte que nous pouvions écrire une preuve du théorème \ref{maurey-pi} sans faire intervenir la notion de cotype et qui n'utilise essentiellement que les inégalités de Kahane-Khintchine.
Nous ne connaissons pas de référence écrite mais la preuve en question n'est certainement pas nouvelle.
Cela dit, c'est cette preuve que nous allons adapter dans un cadre multidimensionnel en invoquant les inégalités de Kahane-Khintchine-Marcus-Pisier.
Le lecteur intéressé par une preuve simple du théorème \ref{maurey-pi} pourra lire les arguments du théorème \ref{mplp} avec $d_n=1$ et les appendices \ref{preuve-mp1} et \ref{preuve-mp2}.
Rappelons aussi que l'optimalité des moments nous interdit d'invoquer des arguments de grandes déviations (ces dernières nécessitent d'utiliser des variables aléatoires sous-gaussiennes, voir plus loin \eqref{sous-gaus}).

Lorsque chaque fonction $u_n$ appartient à un sous-espace $E_n$ de $L^2(X)$ de dimension finie, il existe un procédé de randomisation bien mieux adapté.
L'exemple que l'on doit avoir à l'esprit est celui d'une décomposition d'un sous-espace fermé de $L^2(X)$ en somme directe orthogonale de sous-espaces $E_n$.
On notera dans la suite 
\bema
d_n:=\dim_\C (E_n)\geq 1.
\enma
De manière informelle, on va faire tourner chaque fonction $u_n$ 
de façon aléatoire et uniforme dans la sphère de $E_n$ centrée en l'origine et de rayon $\left\Vert u_n\right\Vert_{L^2(X)}$.
Définissons ce procédé de façon rigoureuse.
Pour chaque $n\in \N$, fixons une fois pour toutes une base hilbertienne $\phi_{n,1},\dots,\phi_{n,d_n}$ de $E_n$ en tant que $\C$-espace vectoriel, c'est-à-dire que l'on a la somme directe
\bema E_n=\bigoplus\limits_{1\leq j\leq d_n} \C \phi_{n,j} \subset L^2(X) \enma
et considérons de plus une suite $(W_n)_{n\geq 0}$ de matrices aléatoires indépendantes de lois respectivement uniformes dans les groupes unitaires $(U_{d_n}(\C))_{n\geq 0}$, c'est-à-dire que la loi de $W_n$ est la mesure de Haar normalisée du groupe $U_{d_n}(\C)$.
Faire tourner la suite $(u_n)_{n\geq 0}$ revient à étudier la série aléatoire 
\begin{equation}\label{defi-alea}
\sum_{n\geq 0} \sum_{i=1}^{d_n}\left( \sum_{j=1}^{d_n} W_{n,i,j}^\omega \langle u_n,\phi_{n,j}\rangle_{L^2(X)} \right)\phi_{n,i},
\end{equation}
où $W_n^\omega \in U_{d_n}(\C)$ est la valeur de la matrice aléatoire $W_n$ en $\omega \in \Omega$.
Faisons deux remarques : 
\begin{enumerate}[$\bullet$ ]
\item lorsque l'égalité $d_n=1$ a lieu pour tout $n\in \N$, le terme général de la série aléatoire \eqref{defi-alea} 
prend la forme $S_n(\omega)u_n$ où $S_n$ suit une loi uniforme sur le cercle unité. Dans l'échelle des espaces $L^p(X)$, cela nous ramène au théorème \ref{maurey-pi}.
\item changer de base hilbertienne $\phi_{n,1},\dots,\phi_{n,d_n}$ de $E_n$ revient à remplacer $W_{n}$ par $P^{-1} W_n P$ pour une certaine matrice unitaire $P\in U_{d_n}(\C)$.
Mais cela n'a bien entendu aucun effet sur la loi matricielle de $W_n$.
En d'autres termes, l'étude de la série aléatoire \eqref{defi-alea}
est indépendante de la base hilbertienne $\phi_{n,1},\dots,\phi_{n,d_n}$ de $E_n$. Ce point, à lui tout seul, justifie la pertinence de la définition \eqref{defi-alea}.
\end{enumerate}

Afin de recouvrir les variables aléatoires $X_n$ du théorème \ref{maurey-pi}, nous allons considérer des matrices aléatoires plus générales.
Introduisons la notion suivante qui généralise la notion de symétrie pour les variables aléatoires scalaires \cite[page 82, line (2.3)]{pisier1981}.
\begin{defi}
Une matrice aléatoire $M_n:\Omega \rightarrow \mathcal{M}_{d_n}(\R)$ est \textbf{orthogonalement invariante} si pour tout $P\in O_{d_n}(\R)$ les matrices aléatoires $M_n$ et $PM_n$ ont la même loi. On définit de même l'invariance unitaire pour une matrice aléatoire $M_n:\Omega\rightarrow \mathcal{M}_{d_n}(\C)$ en remplaçant le groupe $O_{d_n}(\R)$ par $U_{d_n}(\C)$.
\end{defi}

Fixons maintenant quelques notations usuelles d'algèbre linéaire : 
\begin{eqnarray} \nonumber
\forall y \in \mathcal{M}_{d_n,1}(\C)=\C^{d_n} \quad | y|&:=& \sqrt{|y_1|^2+\dots+|y_{d_n}|^2},\\ \nonumber
 \forall A\in \mathcal{M}_{d_n}(\C) \quad \left\Vert A \right\Vert_{op}& :=& \sup\limits_{y\in \C^{d_n}\backslash \{0\}} \Frac{ | Ay|}{|y|},\\  \label{val-abs}
  |A|& := &  \sqrt{\tt \overline{A} A }\in \mathcal{M}_n(\C),\\ \nonumber
  \sigma(A) &:= & \inf\limits_{y\neq 0} \Frac{|Ay |}{| y|}.     
\end{eqnarray}

La matrice $|A|$ est une matrice hermitienne positive issue d'une décomposition polaire de $A=P|A|$ avec $P\in U_{n}(\C)$ (une telle décomposition est unique dès lors que $\det(A)\neq 0$). 
Enfin, $\sigma(A)$ est la plus petite valeur singulière de $A$, il s'agit aussi de la plus petite valeur propre de $|A|$.

Introduisons aussi la fonction spectrale du sous-espace $E_n\subset L^2(X)$ : 
\bema   \forall x\in X \quad e(n,x):=\sup\{ |u_n(x)|^2, \quad  u_n\in E_n \quad \left\Vert u_n\right\Vert_{L^2(X)}=1 \}.  \enma
Dans une base hilbertienne $(\phi_{n,1},\dots,\phi_{n,d_n})$ de $E_n$, la fonction spectrale prend la forme plus courante :

\begin{equation}\label{defi-en}
\begin{array}{rcl}   e(n,x) & =&   \sup\left\{ |a_1 \phi_{n,1}(x)+\dots+a_{d_n} \phi_{n,d_n}(x)|^2 ,\quad |a_1|^2+\dots+|a_{d_n}|^2= 1\right\}  \\
& = & \Sum_{k=1}^{d_n} |\phi_{n,k}(x)|^2.  
 \end{array} \end{equation}
On remarque ainsi que $x\mapsto \frac{1}{d_n} e(n,x)$ est une densité de probabilité sur $X$. Le théorème suivant est dans le même état d'esprit que le théorème \ref{maurey-pi} et énonce un phénomène d'universalité pour la randomisation multidimensionnelle dans $L^p(X)$ de séries aléatoires plus générales que \eqref{defi-alea}.

\begin{theo}\label{mplp}
Fixons un réel $p\in [1,+\infty[$ et considérons une suite de matrices aléatoires $M_n:\Omega\rightarrow \mathcal{M}_{d_n}(\R)$ indépendantes, orthogonalement invariantes et vérifiant
\begin{eqnarray}\label{hyp1}
  \sup\limits_{n\in \N} \E \left[\left\Vert M_n\right\Vert_{op}^{\max(2,p)}\right]<+\infty   \qquad \mbox{et} \qquad 
\inf\limits_{n\in \N} \sigma(\E[|M_n|]) >0.
\end{eqnarray}
On suppose que les sous-espaces $E_n$ sont inclus dans $L^p(X)$.
Pour toute suite $(u_n)_{n\geq 0}$, avec $u_n\in E_n$, on définit la série aléatoire à valeurs dans $L^p(X)$
\begin{equation}\label{seri-al}
\Sum_n \sum_{i=1}^{d_n}\left( \sum_{j=1}^{d_n} M_{n,i,j} \langle u_n,\phi_{n,j}\rangle \right)\phi_{n,i}.
\end{equation}
Les propriétés suivantes sont équivalentes : 
\begin{enumerate}[i) ]
\item la fonction $x\mapsto \sqrt{\Sum_{n\geq 0} \left\Vert u_n \right\Vert_{L^2(X)}^2 \frac{e(n,x)}{d_n}}$ appartient à $L^p(X)$,
\item la série aléatoire \eqref{seri-al} converge presque sûrement dans $L^p(X)$,
\item la série aléatoire \eqref{seri-al} converge dans l'espace de Bochner-Lebesgue $L^{\max(2,p)}(\Omega,L^p(X))$.
\end{enumerate}
Le même énoncé est valide si les matrices aléatoires $M_n:\Omega \rightarrow \mathcal{M}_{d_n}(\C)$ sont supposées unitairement invariantes.
\end{theo}

Ce théorème appelle à plusieurs commentaires : 
\begin{enumerate}[$\bullet$ ]
\item on insiste sur le point crucial que la condition i) du théorème \ref{mplp} est \textbf{indépendante} de la suite des lois des matrices aléatoires $M_n$ et des bases hilbertiennes $\phi_{n,1},\dots,\phi_{n,d_n}$ de $E_n$. 
Cela traduit précisément un phénomène d'universalité pour la randomisation dans $L^p(X)$.
\item on ne doit pas être surpris de demander l'invariance orthogonale alors que les espaces vectoriels sont complexes. En effet, ne serait-ce que pour la dimension $1$, nous avons bien la formule suivante
\begin{displaymath} \forall (a,b)\in \C^2 \quad \E[ |\ep_1 a +\ep_2 b|^2  ]= |a|^2+|b|^2  .\end{displaymath}
\item dans le cas où chaque dimension $d_n$ vaut $1$, les séries aléatoires considérées sont de la forme 
\bema \sum X_n \langle u_n,\phi_{n,1} \rangle \phi_{n,1}=\sum X_n u_n, \enma
 où $X_n$ est une variable aléatoire symétrique (au sens réel ou complexe). On se retrouve dans la situation unidimensionnelle du théorème \ref{maurey-pi}. Cela dit, nous ne voyons pas comment le théorème de randomisation universelle de Maurey-Pisier (Théorème \ref{theo-maupi}), avec la définition usuelle d'un espace de Banach de cotype fini, peut impliquer le théorème \ref{mplp} (même en affaiblissant la condition de moments \eqref{hyp1}) lorsque les dimensions $d_n$ sont différentes de $1$.
Puisque l'exposant $\max(2,p)$ dans \eqref{hyp1} est précisément le meilleur cotype pour l'espace de Banach $L^p(X)$, il doit y avoir une interprétation similaire de la randomisation multidimensionnelle dans la catégorie des espaces de Banach (voir par exemple l'inégalité \eqref{pietch2} qui apparaît au cours de la preuve).
Dans l'annexe \ref{treillis}, nous expliquerons comment une telle interprétation est possible si l'on considère $L^p(X)$ non pas comme un espace de Banach de cotype fini mais comme un treillis de Banach de cotype fini.

\item les conditions \eqref{hyp1} apparaissent sous une forme analogue dans l'étude de la convergence presque sûre des séries de Fourier dans l'espace fonctionnel $L^\infty(G)$ sur un groupe compact $G$ (voir \cite[page 97, Theorem 3.5]{pisier1981}) mais $\left\Vert M_n\right\Vert_{op}^{\max(2,p)}$ est remplacé par $\left\Vert M_n\right\Vert_{op}^2$. Dans le cas où le groupe $G$ est abélien et localement compact, Marcus et Pisier expliquent que l'exposant $2$ est dû au fait que l'espace des fonctions ``presque sûrement continues" sur $G$ est curieusement de cotype $2$ (alors que $L^\infty(G)$ n'est même pas de cotype fini dès que $G$ est infini).
\end{enumerate}

Il convient dès maintenant d'examiner des exemples concrets de matrices aléatoires orthogonalement (ou unitairement) invariantes $M_n$ qui satisfont les hypothèses \eqref{hyp1} :
\begin{enumerate}[$\bullet$ ]
\item un exemple évident est une suite $(\mathcal{E}_n)_{n\geq 0}$ de matrices aléatoires indépendantes de lois respectivement uniformes dans les groupes orthogonaux $(O_{d_n}(\R))_{n\geq 0}$. On a évidemment $\left\Vert \mathcal{E}_n\right\Vert_{op}=1$ et $|\mathcal{E}_n|=I_{d_n}$.
\`A juste titre, on peut considérer que le procédé de randomisation multidimensionnelle proposé par Burq et Lebeau dans \cite[Appendice C]{burq-lebeau} correspond au cas $M_n(\omega)=m_n(\omega)\mathcal{E}_n(\omega)$ où $(m_n)_{n\geq 0}$ est une suite indépendante de variables aléatoires sous-gaussiennes (voir \eqref{sous-gaus}).
\item la suite $(W_n)_{n\geq 0}$ qui a motivé la définition \eqref{defi-alea} fonctionne de même trivialement.
\item voici un autre exemple très important : on considère une famille $(g_{n,i,j})$, avec $n\in\N$ et $1\leq i,j\leq d_n$, de variables aléatoires i.i.d. suivant une loi $\mathcal{N}_{\R}(0,1)$). On pose alors 
\begin{equation}\label{defi-gau} G_n=\frac{1}{\sqrt{d_n}}[g_{n,i,j}].\end{equation}
Les propriétés d'invariance orthogonale des variables gaussiennes assurent que $G_n$ est orthogonalement invariante.
Quant aux propriétés \eqref{hyp1}, elles sont démontrées dans \cite[page 78, Proposition 1.5]{pisier1981}.
En fait, le cas gaussien est un cas particulier de la proposition suivante.
\end{enumerate}

\begin{prop}\label{ex-mat}
Considérons des variables aléatoires $X_{ij}:\Omega \rightarrow \R$ centrées, i.i.d. et où $(i,j)$ parcourt $\N^{\star 2}$.
Pour tout $n\in \N^\star$, on note $M_n$ la matrice aléatoire $\frac{1}{\sqrt{n}}[X_{ij}]:\Omega\rightarrow \mathcal{M}_n(\R)$.
Alors, il existe une constante universelle $C>1$ telle que si l'on a $0<\E[|X_{1,1}|^4]<+\infty$ alors 
\begin{equation}\label{nucl}   \frac{\E[|X_{11}|]^2}{C\E[|X_{11}|^4]^{\frac{1}{4}}}
\leq \sigma\left(\E\left[|M_n|\right]\right).\end{equation}
En outre, si $X_{1,1}$ admet un moment d'ordre $p\geq 4$ alors nous avons l'inégalité de Latala précisée 
\begin{equation}\label{super-lata}  \E\left[ \left\Vert M_n\right\Vert_{op}^p\right]  \leq C(p)\E[|X_{1,1}|^{p}].  \end{equation}
Par conséquent, on a 
\bema
\sup\limits_{n\geq 1}  \E\left[ \left\Vert M_n\right\Vert_{op}^p\right]<+\infty \quad \mbox{et} \quad \inf\limits_{n\geq 1} \sigma\left(\E\left[|M_n|\right]\right)>0.
\enma
\end{prop}
\begin{rema}
Si les variables aléatoires $X_{ij}$ sont à valeurs complexes et d'espérance nulle, alors on peut aussi montrer \eqref{super-lata} en considérant partie réelle et partie imaginaire des variables aléatoires $X_{ij}$.
Quant à \eqref{nucl}, une preuve similaire est valide mais l'on doit remplacer l'hypothèse ``centrée'' par ``symétrique au sens complexe''
: pour tout réel $\alpha\in \R$ les variables $e^{i\alpha} X_{ij}$ et $X_{ij}$ ont la même loi.
\end{rema}
\begin{rema}\label{rema-lata}
Dans le cadre du théorème \ref{mplp}, nous avons besoin de matrices aléatoires orthogonalement invariantes.
Or les matrices $\frac{1}{\sqrt{n}}[X_{ij}]$ de la proposition \ref{ex-mat} ne sont pas orthogonalement invariantes en général (le cas gaussien est une exception notable).
Il suffit alors d'examiner des matrices aléatoires de la forme $\frac{1}{\sqrt{n}}\mathcal{E}_n [X_{ij}]$ où $\mathcal{E}_n:\Omega \rightarrow O_{d_n}(\R)$ est indépendante de $[X_{ij}]$ et de remarquer que $\frac{1}{\sqrt{n}} \mathcal{E}_n [X_{ij}]$ et $\frac{1}{\sqrt{n}} [X_{ij}]$ ont la même norme d'opérateur et évidemment que l'on a $|\mathcal{E}_n [X_{ij}]|=|[X_{ij}]|$ (au sens de \eqref{val-abs}).
\end{rema}
On montrera l'inégalité \eqref{nucl} à l'aide d'un argument d'interpolation  simplifiant celui du cas gaussien \cite[page 78, Proposition 1.5]{pisier1981} et 
de l'inégalité de Latala \cite{latala} :
\begin{equation}\label{lata} \E\left[\left\Vert M_n\right\Vert_{op}\right]\leq C  \E[|X_{1,1}|^4]^{\frac{1}{4}}. \end{equation}
Nous n'avons pas trouvé l'estimation \eqref{super-lata} dans des références déjà publiées mais nous pensons qu'il est possible de l'obtenir avec la méthode des moments (par exemple expliquée dans \cite[Part 2.3]{tao-rand-matrix}).
Cependant, cette méthode est assez technique à mettre en place.
Un cas particulier très bien compris dans la littérature des matrices aléatoires est le cas ``sous-gaussien''.
Parmi plusieurs définitions équivalentes, $X_{1,1}$ est sous-gaussienne si l'on a $\E[X_{1,1}]=0$ et 
\begin{equation}\label{sous-gaus} \exists K>0 \qquad \forall p\in \N^\star  \qquad \E[|X_{1,1}|^p]^{\frac{1}{p}} \leq K\sqrt{p}.  \end{equation}
Cela équivaut par exemple à la condition
\bema
\exists K'>0 \quad \forall t\in \R \quad \E[\exp(tX_{11})]\leq \exp(K't^2).
\enma
D'après \cite[Proposition 2.3]{rudelson2009} (voir aussi la preuve de \cite[Fact 2.4]{litvak}), on sait que l'on a
\bema
\forall p\in \N^\star \quad \exists C(p,K)>0 \quad \forall n\in \N^\star \quad \E\left[\left\Vert  M_n\right\Vert_{op}^p\right]\leq C(p,K).
\enma
L'inégalité \eqref{super-lata} est plus générale et notre argument combine la preuve de l'inégalité \eqref{lata}  et 
 les inégalités de Kahane-Khintchine \eqref{KK} dans l'espace de Banach $(\mathcal{M}_n(\R),\left\Vert \cdot\right\Vert_{op})$.
Cela explique pourquoi \eqref{super-lata} implique \eqref{lata} avec $p=4$.
La renormalisation de $[X_{ij}]$ par le facteur $\frac{1}{\sqrt{n}}$ et le fait que \eqref{lata} nécessite au moins un moment d'ordre $4$ peuvent surprendre a priori mais il s'agit de faits bien connus dans la théorie des matrices aléatoires (voir la remarque \ref{lata-optim} ainsi que \cite{bai88b} et \cite[Theorem 3.1]{bai88a}).
Enfin, l'inégalité triviale $|X_{1,1}|\leq \sqrt{n} \left\Vert M_n\right\Vert_{op}$ nous fait comprendre que l'inégalité \eqref{super-lata} est optimale par rapport à la condition de moment d'ordre $p$.

En résumé, la randomisation multidimensionnelle dans $L^p(X)$ est universelle et la proposition \ref{ex-mat} donne des exemples naturels et non triviaux de procédés de randomisation.
Le théorème \ref{mplp} nous amène à poser la définition suivante.

\begin{defi}\label{defi-cv}
Pour tout $p\in [1,+\infty[$, on note $\pl^p(X,\oplus E_n)$ l'espace vectoriel des suites $(u_n)_{n\geq 0}$, avec $u_n \in E_n$ pour tout $n\geq 0$, qui satisfont les assertions équivalentes i), ii), iii) du théorème \ref{mplp}. 
On appelle $\pl^p(X,\oplus E_n)$ l'espace de Lebesgue probabiliste associé à la suite de sous-espaces $(E_n)_{n\geq 0}$ et on le munit de la norme 
\bema    \left\Vert (u_n)\right\Vert_{\pl^p(X,\oplus E_n)} := \left\Vert \sqrt{\Sum_{n\geq 0}  \left\Vert u_n\right\Vert_{L^2(X)}^2 \frac{e(n,x)}{\dim(E_n)} }\right\Vert_{L_x^p(X)}  .  \enma
\end{defi}
Faisons quelques remarques sur cette définition :
\begin{enumerate}[$\bullet$ ]
\item dans le cas où les sous-espaces $E_n$ orthogonaux deux à deux, l'espace $\pl^p(X,\oplus E_n)$ est linéairement isométrique à un sous-espace fermé de $L^p(X,\oplus E_n)$ et donc est un espace de Banach.
Cela justifie au passage pourquoi l'on a choisi la notation $\pl^p(X,\oplus E_n)$ pour désigner l'espace de Lebesgue probabiliste.
En fait, $\pl^p(X,\oplus E_n)$ est un espace de Banach même si les sous-espaces $E_n$ ne sont pas en somme directe dans $L^2(X)$. Pour le voir, il suffirait de raisonner avec l'espace de Hilbert abstrait $\bigoplus E_n$ (voir plus loin \eqref{defi-S}).
\item nous verrons que dans le cas où $X$ est une variété riemannienne compacte, si les sous-espaces $E_n$ sont ``linéairement localisés en fréquence'', alors l'espace $\pl^p(X,\oplus E_n)$ pourra être considéré comme un espace de distributions sur $X$. Il en sera de même pour l'oscillateur harmonique $-\Delta+|x|^2$ sur $\R^d$.
\item rappelons que dans le cas unidimensionnel $d_n=\dim(E_n)=1$, l'espace $\pl^p(X,\oplus E_n)$ est exactement l'espace des suites $(u_n)_{n\in \N}$, avec $u_n\in E_n$ telles que la série aléatoire $\sum \ep_n u_n$ converge presque sûrement dans $L^p(X)$ (voir le théorème \ref{maurey-pi}).
\item si une suite $(u_n)$, vérifiant $u_n\in E_n$ pour tout $n\in \N$, n'appartient pas à $\pl^p(X,\oplus E_n)$
alors les séries aléatoires du théorème \ref{mplp} divergent dans $L^p(X)$ avec une probabilité strictement positive et donc divergent presque sûrement d'après la loi du tout ou rien.
\end{enumerate}

\section{Interpolation et dualité des espaces de Lebesgue probabilistes $\pl^p$}\label{enonc-interpo}

Puisque les fonctions $\frac{1}{d_n}e(n,\cdot)$ sont des densités de probabilité sur $X$, on a l'égalité
\bema
\pl^2(X,\oplus E_n)=\left\{(u_n)_{n\in \N}, \quad u_n \in E_n,\quad \sum_{n\in \N} \left\Vert u_n\right\Vert_{L^2(X)}^2<+\infty \right\}.
\enma
Lorsque $p$ est différent de $2$, il semble nécessaire d'exploiter des propriétés des fonctions spectrales $e(n,\cdot)$ afin de mieux comprendre l'espace $\pl^p(X,\oplus E_n)$.
Pour aborder cette étude, nous avons besoin de rappeler la définition de l'espace de Lorentz $L^{p,\infty}(X)$ pour tout réel $p\in [1,+\infty[$. Il s'agit de l'espace vectoriel des fonctions mesurables $f:X\rightarrow \C\cup \{\infty\}$ telles que 
\bema \exists c>0\quad \forall t>0 \quad \mu\left\{ x\in X, |f(x)|>t   \right\}\leq \frac{c^p}{t^p}, \enma
où $\mu$ est la mesure de l'espace mesuré $X$.
L'inclusion $L^p(X)\subset L^{p,\infty}(X)$ est toujours vraie et est généralement stricte. Le contre-exemple typique est $|x|^{-\frac{d}{p}}\in L_x^{p,\infty}(\R^d)\backslash L_x^p(\R^d)$.
Les deux résultats suivants permettent de prévoir les résultats d'interpolation et de dualité des espaces $\pl^p(X,\oplus E_n)$.
Rappelons que les espaces interpolés des espaces $\pl^p(X,\oplus E_n)$ sont définis de façon abstraite mais peuvent être vus comme des sous-espaces vectoriels du même espace ambiant, à savoir $\prod_{n\in \N} E_n$.
En première lecture, on suggère au lecteur de considérer l'énoncé suivant dans le cas unidimensionnel $E_n=\C \phi_n$, c'est-à-dire si l'on a $d_n=1$ et $\sqrt{e(n,x)}=|\phi_n(x)|$.

\begin{theo}\label{theo-interpo}
Considérons $p_1<p_2$ deux nombres appartenant à $]1,+\infty[$ et vérifiant la condition $\frac{1}{p_1}+\frac{1}{p_2}\leq 1$.
On suppose que l'on a l'inclusion $E_n\subset L^{p_1}(X)\cap L^{p_2}(X)$, pour tout $n\in \N$, et les assertions
\begin{eqnarray} \label{loren1}
\forall p\in \{p_1,p_2\} \qquad \sup\limits_{n\in \N} \frac{\sqrt{e(n,x)}}{\left\Vert \sqrt{e(n,x')}\right\Vert_{L_{x'}^{p}(X)}}& \in &L_x^{p,\infty}(X), 
 \\[3mm]
\label{hypo-interpo}
\exists p\in ]p_1,p_2[ \quad \sup\limits_{n\in \N} \Frac{\left\Vert \sqrt{e(n,x')}\right\Vert_{L_{x'}^{p_1}(X)}^{1-\theta}\left\Vert \sqrt{e(n,x')}\right\Vert_{L_{x'}^{p_2}(X)}^{\theta}}{\left\Vert \sqrt{e(n,x')}\right\Vert_{L_{x'}^{p}(X)}} &< & +\infty,
\end{eqnarray}
où $\theta\in [0,1]$ est déterminé par la relation $\frac{1}{p}=\frac{1-\theta}{p_1}+\frac{\theta}{p_2}$.
Alors les espaces $\pl^p(X,\oplus E_n)$,
pour $p$ parcourant $]p_1,p_2[$, sont stables par interpolation \textbf{complexe} : en d'autres termes, si l'on a 
\bema p_1<p_1'<p<p_2'<p_2, \quad \theta'\in ]0,1[,\quad  \frac{1}{p}=\frac{1-\theta'}{p_1'}+\frac{\theta'}{p_2'} ,   \enma
alors on a l'égalité $\pl^p(X,\oplus E_n)=[\pl^{p_1'}(X,\oplus E_n), \pl^{p_2'}(X,\oplus E_n)]_{\theta'}$ avec équivalence des normes.

De même, les espaces $\pl^p(X,\oplus E_n)$,
pour $p$ parcourant $]p_1,p_2[$, sont stables par d'interpolation \textbf{réelle} : 
on a l'égalité $\pl^p(X,\oplus E_n)= [\pl^{p_1'}(X,\oplus E_n), \pl^{p_2'}(X,\oplus E_n)]_{\theta',p}$
avec équivalences des normes.
\end{theo}

La condition $\frac{1}{p_1}+\frac{1}{p_2}\leq 1$ est de nature technique et est peut-être inutile en toute généralité. Elle est toujours vérifiée si l'on a $2\leq p_1 <p_2$.
Une autre situation intéressante se produit si $p_1$ et $p_2$ sont deux exposants conjugués, i.e. vérifient $\frac{1}{p_1}+\frac{1}{p_2}=1$.
Dans ce cas, il sera commode de considérer \eqref{hypo-interpo} avec $p=2$ et $\theta=\frac{1}{2}$ pour avoir l'assertion suivante qui est en apparence plus faible (mais équivalente comme le lemme \ref{majo-Q} le montrera plus loin) :
\begin{equation}\label{hol1}
\sup\limits_{n\in \N} \frac{\left\Vert\sqrt{ e(n,x')}\right\Vert_{L_{x'}^{p_1}(X)}\left\Vert\sqrt{ e(n,x') }\right\Vert_{L_{x'}^{p_2}(X)}  }{d_n} <+\infty.
\end{equation}
En prime de l'interpolation, on a une propriété de dualité.
\begin{theo}\label{theo-dual}
Fixons deux réels $p_1<p_2$ appartenant à $]1,+\infty[$ et vérifiant l'égalité $\frac{1}{p_1}+\frac{1}{p_2}=1$.
Pour tout $n\in \N$, on suppose que l'on a l'inclusion $E_n\subset L^{p_1}(X)\cap L^{p_2}(X)$.
Pour tous $p\in ]p_1,p_2[$ et $(u,w)\in \pl^{p}(X,\oplus E_n)\times \pl^q(X,\oplus E_n)$, avec $q:=\frac{p}{p-1}$, on a
\begin{equation}\label{dual-triv} \sum_{n\geq 0} |\langle u_n,w_n\rangle_{L^2(X)}| \leq \left\Vert u\right\Vert_{\pl^p(X,\oplus E_n)}\left\Vert w\right\Vert_{\pl^q(X,\oplus E_n)}.\end{equation}
En outre, si l'on suppose \eqref{loren1} et \eqref{hol1} alors les espaces $\pl^p(X,\oplus E_n)$, pour $p$ parcourant $]p_1,p_2[$, sont stables par dualité au sens suivant : l'injection canonique $\Lambda_p: \pl^q(X,\oplus E_n) \rightarrow \pl^p(X,\oplus E_n)'$ qui à un élément $w\in \pl^q(X,\oplus E_n)$ associe la forme linéaire 
\bemar  \pl^p(X,\oplus E_n) &\mapsto & \C\\
u & \mapsto &   \Sum_{n\geq 0} \langle u_n,w_n\rangle_{L^2(X)}   \enmar
est un isomorphisme d'espaces de Banach. En d'autres termes, le dual de $\pl^p(X,\oplus E_n)$ est canoniquement isomorphe à $\pl^q(X,\oplus E_n)$ pour tout $p\in ]p_1,p_2[$.
\end{theo}
Les deux théorèmes précédents appellent aussi à quelques remarques :
\begin{enumerate}[$\bullet$ ]
\item on introduira dans la partie \ref{partV} la notion de ``défaut d'interpolation'' qui montrera que l'hypothèse \eqref{hypo-interpo} est inévitable si l'on veut interpoler les espaces $\pl^p(X,\oplus E_n)$ (voir \eqref{defi-Q} et la proposition \ref{neces-interpo}). De même, l'hypothèse \eqref{hol1} est nécessaire pour montrer un théorème de dualité (voir la proposition \ref{Lamba-surj}).
En pratique, les hypothèses \eqref{hypo-interpo} et \eqref{hol1} se réalisent lorsque les densités $\frac{1}{d_n} e(n,\cdot)$ ont tendance à se concentrer uniformément sur un borélien $A_n\subset X$ de mesure finie et strictement positive, c'est-à-dire si l'on peut assimiler la densité de probabilité $\frac{1}{d_n} e(n,\cdot)$ à une fonction de la forme $x\mapsto \frac{1}{\mu(A_n)} \pun_{A_n}(x)$.
Dans les exemples issus de la physique mathématique qui font intervenir des polynômes orthogonaux, on a très souvent des propriétés de concentration qui forcent 
$\left\Vert \sqrt{e(n,\cdot)}\right\Vert_{L^p(X)}$ à être équivalent à $n^{\left(a +\frac{b}{p}\right)}$, pour un certain couple $(a,b)\in \R^2$, lorsque $p$ parcourt un intervalle donné $]p_1,p_2[$. Ainsi, \eqref{hypo-interpo} sera vérifié.
De même \eqref{hol1} sera vérifié si  $\left\Vert\sqrt{e(n,\cdot)}\right\Vert_{L^p(X)}$ est équivalent à $n^{a\left(\frac{1}{2}-\frac{1}{p}\right)} \sqrt{d_n}$ pour un certain réel $a$.
\item l'hypothèse \eqref{loren1} est très importante dans nos preuves.
Bien qu'il ne semble pas facile de l'interpréter, elle est toujours satisfaite si les fonctions spectrales $e(n,\cdot)$ se concentrent sur des parties de $X$ qui ne sont pas trop éloignées les unes des autres.
 Il est sans doute possible de relaxer cette hypothèse (et peut-être même de la supprimer), mais permettons-nous d'expliquer l'intervention les espaces de Lorentz.
D'une part, nos démonstrations utilisent l'interpolation réelle, théorie dans laquelle les espaces de Lorentz jouent un rôle clé.
D'autre part, nous verrons des exemples bien concrets, à savoir les fonctions propres qui se concentrent sur une géodésique de la sphère $\S^d\subset \R^{d+1}$, où l'hypothèse \eqref{loren1} ne se réalise pas si l'on remplace l'espace de Lorentz $L^{p,\infty}(\S^d)$ par l'espace de Lebesgue $L^{p}(\S^d)$ (voir la remarque \ref{loren-perti}).
\end{enumerate}

\section{Théorie $L^p$ presque sûre des variétés compactes, propriétés générales}\label{parti-comp}

Dans cette partie, $X$ est une variété riemannienne \textbf{compacte} lisse, de dimension $d\geq 2$, sans bord et munie de sa mesure volume $\mu$. L'opérateur de Laplace-Beltrami $\Delta$ de $X$ est diagonalisable dans une base orthonormée de $L^2(X)$ et l'on peut définir le spectre d'une fonction $u\in L^2(X)$ :
\bema  \mbox{Sp}(u):=\{\lambda \geq 0,\quad \exists \phi \in L^2(X) \quad \sqrt{-\Delta}\phi=\lambda \phi \quad \mbox{et} \quad \langle u,\phi\rangle \neq 0 \}.\enma
Par souci de cohérence avec les parties \ref{pzoll} et \ref{harmo-sph}, on considère de façon générale une suite de sous-espaces non nuls $(E_n)_{n\in \N}$ de $L^2(X)$ vérifiant : 
\begin{enumerate}[i) ]
\item les sous-espaces $E_n$ sont orthogonaux deux à deux,
\item chaque sous-espace $E_n$ admet une base finie constituée de modes propres de $\Delta$,
\item il existe un réel $\delta\geq 1$ et une suite de réels positifs $(\lambda_n^\star)_{n\in \N}$ tels que 
\begin{equation}\label{hypo-spec}
\forall n\in \N \quad \forall u\in E_n \quad \mbox{Sp}(u)\subset \left[\frac{\lambda_n^\star}{\delta},\delta \lambda_n^\star\right].
\end{equation}
\end{enumerate}

En notant $\Pi_n$ le projecteur orthogonal de $L^2(X)$ sur $E_n$, on déduit de ce qui précède que pour toute distribution $u$ sur $X$ supportée par les sous-espaces $E_n$ et tout réel $s\in \R$, la norme dans l'espace de Sobolev $H^s(X)$ s'écrit 
\bema
\left\Vert u \right\Vert_{H^s(X)}=\left\Vert (I-\Delta)^{\frac{s}{2}} u\right\Vert_{L^2(X)} \simeq_{s,\delta} \left(\sum_{n\in \N} \lambda_n^{\star 2s} \left\Vert \Pi_n(u) \right\Vert_{L^2(X)}^2\right)^{\frac{1}{2}}.
\enma

Le lemme élémentaire suivant montre que l'on peut identifier $\pl^p(X,\oplus E_n)$ à un sous-espace de l'espace vectoriel $\mathcal{D}'(X)$ des distributions sur $X$.

\begin{lemm}\label{identif}
Pour tous $p\in [1,+\infty[$ et $(u_n)_{n\in \N}\in \pl^p(X,\oplus E_n)$, la série $\sum u_n$ converge dans $H^{-\frac{(d-1)}{4}}(X)$.
On identifie alors $\pl^p(X,\oplus E_n)$ au sous-espace des distributions $u\in \mathcal{D}'(X)$ vérifiant
\begin{equation}\label{def-dist}
\left\Vert u\right\Vert_{\pl^p(X,\oplus E_n)}:=\left\Vert \sqrt{\sum_{n\in \N} \left\Vert \Pi_n(u) \right\Vert_{L^2(X)}^2 \frac{e(n,x)}{\dim(E_n)} }\right\Vert_{L_x^p(X)}<+\infty.
\end{equation}
\end{lemm}
\begin{demo}
L'inégalité de Hölder et l'inégalité triangulaire (entre la norme euclidienne et la norme $\left\Vert \cdot\right\Vert_{L^1(X)}$) nous donnent la minoration
\bemar
\left\Vert (u_n)_{n\in \N} \right\Vert_{\pl^p(X,\oplus E_n)}  & \geq & \Frac{1}{\mu(X)^{1-\frac{1}{p}}}  \left\Vert (u_n)_{n\in \N} \right\Vert_{\pl^1(X,\oplus E_n)}   \\
& \geq &  \Frac{1}{\mu(X)^{1-\frac{1}{p}}} \Int_{X} \sqrt{\Sum_{n\in \N} \frac{1}{d_n} \left\Vert u_n \right\Vert_{L^2(X)}^2 e(n,x)}d\mu(x) \\
& \geq &   \Frac{1}{\mu(X)^{1-\frac{1}{p}}}\sqrt{\Sum_{n\in \N} \frac{1}{d_n}\left\Vert u_n\right\Vert_{L^2(X)}^2 \left(\int_{X} \sqrt{e(n,x)}d\mu(x) \right)^2 } .   \enmar
Nous pourrons conclure si l'on arrive à prouver l'inégalité
\begin{equation}\label{minor-L1}     \left\Vert \sqrt{e(n,\cdot)} \right\Vert _{L^1(X)}\geq \frac{1}{C(X,\delta)} \left(\lambda_n^{\star}\right)^{ -\left(\frac{d-1}{4}\right)}\sqrt{d_n}     .  \end{equation}
Nous allons utiliser l'argument de la preuve de la proposition 2 de \cite{sogge2010lower}.
Fixons $r\in ]2,\frac{2(d+1)}{d-1}[$ et invoquons l'inégalité de Hölder :
\begin{equation}\label{hol-sogge}
\sqrt{d_n}=\left\Vert \sqrt{e(n,\cdot)} \right\Vert _{L^2}\leq \left\Vert \sqrt{e(n,\cdot)} \right\Vert _{L^1(X)}^{1-\theta}\left\Vert \sqrt{e(n,\cdot)} \right\Vert _{L^r(X)}^{\theta}.
\end{equation}
où $\theta\in [0,1]$ est défini par la relation $\frac{1}{2}=\frac{1-\theta}{1}+\frac{\theta}{r} $.
On peut facilement majorer $\left\Vert \sqrt{e(n,\cdot)} \right\Vert _{L^r(X)}$ en considérant une base hilbertienne $\phi_{n,1},\dots,\phi_{n,d_n}$ de $E_n$ constituée de modes propres, les inégalités de Sogge \eqref{sogge} et l'hypothèse \eqref{hypo-spec} :
\bema
\left\Vert \sqrt{e(n,\cdot)} \right\Vert _{L^r(X)}^2 = \left\Vert \Sum_{i=1}^{d_n} | \phi_{n,i}|^2\right\Vert_{L^{\frac{r}{2}}(X)} \leq  \sum_{i=1}^{d_n}  \left\Vert \phi_{n,i}\right\Vert_{L^r(X)}^2 \leq C(X,r,\delta) (\lambda_n^{\star})^{(d-1)\left(\frac{1}{2}-\frac{1}{r} \right)} d_n.
\enma
En injectant cette inégalité dans \eqref{hol-sogge}, on obtient \eqref{minor-L1}.
\end{demo}

Donnons quelques exemples pour situer l'identification autorisée par le lemme précédent : 
\begin{enumerate}[$\bullet$ ]
\item si $(\phi_n)_{n\in \N}$ est une base hilbertienne de $L^2(X)$ constituée de modes propres, alors $\pl^p(X,\oplus \C \phi_n)$ est un espace de distributions sur $X$ pour tout $p\in [1,+\infty[$.
\item l'espace $\pl^2(X,\oplus E_n)$ coïncide avec la somme directe orthogonale $\bigoplus\limits_{n\in \N} E_n\subset L^2(X)$. Par conséquent, pour tout réel $p\in [2,+\infty[$, l'espace $\pl^p(X,\oplus E_n)$ est un espace de fonctions.
\item si les sous-espaces $E_n$ sont tous unidimensionnels, alors $\pl^p(X,\oplus E_n)$ est exactement le sous-espace des distributions $u=\sum_{n\geq 0} \Pi_n u$ telles que la série $\sum \ep_n \Pi_n(u)$ converge presque sûrement dans $L^p(X)$ (voir le théorème \ref{maurey-pi}).
Par exemple, le théorème de Paley-Zygmund, énoncé au début de ce travail, dit que l'on a, pour tout $p\in[1,+\infty[$, l'inclusion $L^2(\mathbb{T})\subset \pl^p(\mathbb{T},\oplus \C e^{inx})$ où $n$ parcourt $\Z$.
On va retrouver immédiatement que l'inclusion précédente est en fait une égalité.
\end{enumerate}

On peut maintenant exprimer simplement une condition nécessaire et suffisante pour obtenir un théorème de Paley-Zygmund (le cas $p\geq 2$ a essentiellement été écrit dans \cite{tzvetkov-ay} même si le formalisme n'est pas le même).
\begin{prop}\label{bornLp}
Les deux assertions suivantes sont vraies
\begin{enumerate}[i) ]
\item pour tout $p\in [2,+\infty[$, on a l'équivalence
\bema \pl^p(X,\oplus E_n)=\bigoplus_{n\in \N} E_n \quad \Leftrightarrow \quad \sup\limits_{n\in \N}      \left\Vert \frac{\sqrt{e(n,\cdot)}}{\sqrt{\dim(E_n)}}\right\Vert_{L^p(X)}<+\infty.   \enma
\item pour tout $q\in [1,2]$, on a l'équivalence
\bema \pl^q(X,\oplus E_n)=\bigoplus_{n\in \N} E_n \quad \Leftrightarrow \quad \inf\limits_{n\in \N}      \left\Vert\frac{\sqrt{e(n,\cdot)}}{\sqrt{\dim(E_n)}}\right\Vert_{L^q(X)}>0.\enma
\end{enumerate}
\end{prop}
\begin{demo}
i) Commençons par le sens $\Leftarrow$.
On a déjà l'inclusion $\pl^p(X,\oplus E_n)\subset \pl^2(X,\oplus E_n)=\bigoplus\limits_{n\in \N} E_n$.
L'autre inclusion découle de l'inégalité triangulaire dans $L^{\frac{p}{2}}(X)$, on a pour tout $u\in \bigoplus\limits_{n\in \N} E_n$ :
\bema
\left\Vert u\right\Vert_{\pl^p(X,\oplus E_n)}=\left\Vert\sum_{n\in \N} \left\Vert \Pi_n(u) \right\Vert_{L^2(X)}^2 \frac{e(n,x)}{d_n} \right\Vert_{L_x^{\frac{p}{2}}(X)}^{\frac{1}{2}}\leq \left( \sum_{n\in \N} \left\Vert \Pi_n(u) \right\Vert_{L^2(X)}^2 \right)^{\frac{1}{2}} \sup\limits_{n\in \N} \left\Vert \frac{\sqrt{e(n,\cdot)}}{\sqrt{d_n}}\right\Vert_{L^p(X)}.
\enma
Traitons le sens $\Rightarrow$ et supposons par l'absurde que l'on a 
\bema 
\sup\limits_{n\in \N} \left\Vert \frac{e(n,\cdot)}{d_n}\right\Vert_{L^{\frac{p}{2}}(X)} = \sup\limits_{n\in \N} \left( \left\Vert \frac{\sqrt{e(n,\cdot)}}{\sqrt{d_n}}\right\Vert_{L^{p}(X)}\right)^2 =+\infty.
\enma
Il existe donc un sous-ensemble infini $J\subset \N$ et une fonction $u\in \bigoplus\limits_{n\in \N} E_n$ qui vérifie 
\bemar
n\in J & \Rightarrow & \left\Vert \Pi_n(u) \right\Vert_{L^2(X)}=\left\Vert \Frac{e(n,\cdot)}{d_n} \right\Vert_{L^{\frac{p}{2}}(X)}^{-\frac{1}{2}}, \\
n\not\in J & \Rightarrow & \left\Vert \Pi_n(u) \right\Vert_{L^2(X)}=0.
\enmar
Cette fonction $u$ ne peut pas appartenir à $\pl^p(X,\oplus E_n)$ en vertu de 
\bema \sum_{n\in J} \left\Vert \frac{e(n,\cdot)}{d_n} \right\Vert_{L^{\frac{p}{2}}(X)}^{-1} \frac{e(n,x)}{d_n} \not\in L_x^{\frac{p}{2}}(X).     \enma

ii) Pour le sens $\Leftarrow$, on a de même l'inclusion 
$\bigoplus\limits_{n\in \N} E_n=\pl^2(X,\oplus E_n)\subset \pl^q(X,\oplus E_n)$.
L'autre inclusion découle de l'inégalité élémentaire suivante (voir \eqref{Lp-type})  qui dit que $L^q(X)$ est un espace de Banach de cotype $\max(2,q)=2$  
\bemar
\left(\Int_{X} \left(\Sum_{n\in \N} \frac{\left\Vert \Pi_n(u) \right\Vert_{L^2(X)}^2}{d_n}e(n,x)  \right)^{\frac{q}{2}}d\mu(x)\right)^{\frac{1}{q}} & 	 \geq & \left( \Sum_{n\in \N} \left\Vert \Pi_n(u) \right\Vert_{L^2(X)}^2 \left\Vert \frac{\sqrt{e(n,\cdot)}}{\sqrt{d_n}}\right\Vert_{L^q(X)}^2 \right)^{\frac{1}{2}} \\
& \geq & \left\Vert u\right\Vert_{L^2(X)} \inf\limits_{n\in \N}\left\Vert \Frac{\sqrt{e(n,\cdot)}}{\sqrt{d_n}}\right\Vert_{L^q(X)} .
\enmar
Traitons le sens $\Rightarrow$ par l'absurde.
Il existe un sous-ensemble infini $J\subset \N$ tel que 
\bema  \sum_{n\in J} \left\Vert\frac{\sqrt{e(n,\cdot)}}{\sqrt{d_n}} \right\Vert_{L^q(X)} <+\infty.  \enma
Considérons alors une distribution $u\in \mathcal{D}'(X)$ telle que $\left\Vert \Pi_n(u) \right\Vert_{L^2(X)}$ vaut $1$ pour $n\in J$ et $0$ sinon. Il est clair que $u$ n'appartient pas à $L^2(X)$ mais appartient à $\pl^q(X,\oplus E_n)$ :
\bema   \left\Vert \sqrt{\sum_{n\in J} \frac{e(n,\cdot)}{d_n} }\right\Vert_{L^q(X)}\leq 
 \left\Vert \sum_{n\in J}  \frac{\sqrt{e(n,\cdot)}}{\sqrt{d_n}} \right\Vert_{L^q(X)} <+\infty     .\enma
\end{demo}
\begin{rema}
La proposition précédente suggère déjà que, sans hypothèses sur les fonctions $e(n,\cdot)$, le dual de $\pl^p(X,\oplus E_n)$ ne peut pas être canoniquement isomorphe à $\pl^q(X,\oplus E_n)$, avec $q=\frac{p}{p-1}$. En effet, l'inégalité de Hölder donne 
\bema
1 = \left\Vert \frac{e(n,\cdot)}{d_n} \right\Vert_{L^1(X)} \leq \left\Vert \frac{\sqrt{e(n,\cdot)}}{\sqrt{d_n}}\right\Vert_{L^p(X)} \left\Vert\frac{\sqrt{e(n,\cdot)}}{\sqrt{d_n}}\right\Vert_{L^q(X)}.
\enma
Si l'on suppose par exemple que l'on a $\bigoplus\limits_{n\in \N} E_n=L^2(X)$, 
on en déduit que l'on a seulement les implications (et non les équivalences) pour tout $p\geq 2$ :
\bema\begin{array}{ccc}
\sup\limits_{n\in \N}     \left\Vert\frac{\sqrt{e(n,\cdot)}}{\sqrt{d_n}}\right\Vert_{L^p(X)}<+\infty & \Rightarrow & \inf\limits_{n\in \N}    \left\Vert\frac{\sqrt{e(n,\cdot)}}{\sqrt{d_n}}\right\Vert_{L^q(X)}>0 \\[6mm]
\Updownarrow &  & \Updownarrow \\[3mm]
\pl^p(X,\oplus E_n)=L^2(X) & \Rightarrow & \pl^q(X,\oplus E_n)=L^2(X).
\end{array}
\enma 
Cette intuition sera confirmée par la randomisation des fonctions zonales sur $\S^d$ dans la partie \ref{harmo-sph}.
\end{rema}

La proposition précédente nous permet de prouver le résultat suivant.

\begin{prop}\label{inje-mani}
Nous supposons qu'il existe un groupe $G$ d'isométries de $X$ qui vérifie :
\begin{enumerate}[i) ]
\item l'action de $G$ est transitive : pour tout $(x,y)\in X^2$ il existe $\chi\in G$ tel que $y=\chi(x)$,
\item pour tout $n\in \N$, le sous-espace $E_n$ est invariant par $G$ : pour tout $u\in E_n$ on a $u\circ \chi \in E_n$.
\end{enumerate}
Alors pour tout réel $p\in [1,+\infty[$ on a $ \pl^p(X,\oplus E_n )=\bigoplus\limits_{n\in \N} E_n        $.
\end{prop}
\begin{demo}
Quelle que soit la base hilbertienne $\phi_{n,1},\dots,\phi_{n,d_n}$ de $E_n$ et quelle que soit l'isométrie $\chi\in G$, la famille $(\phi_{n,1}\circ\chi,\dots,\phi_{n,d_n}\circ\chi)$ est aussi une base hilbertienne de $E_n$.
La formule \eqref{defi-en} implique que l'on a $e(n,x)=e(n,\chi(x))$ pour tout $x\in X$.
Par transitivité, la densité de probabilité $\frac{1}{d_n} e(n,\cdot)$ est constante (en l'occurrence égale à $\frac{1}{\mu(X)}$ où $\mu(X)$ est le volume de $X$).
La conclusion découle de la proposition \ref{bornLp} ou directement de \eqref{def-dist}.
\end{demo}

La proposition précédente s'applique immédiatement si les sous-espaces $E_n$ sont des sous-espaces propres de $\Delta$ et si le groupe des isométries de $X$ agit transitivement : par exemple si $X$ est un produit fini de sphères (dont les rayons ne sont pas nécessairement identiques).

Décrivons un autre cas auquel la proposition \ref{inje-mani} s'applique et où les sous-espaces $E_n$ ne sont pas spectralement saturés.
Supposons que $X$ soit un groupe de Lie réel compact.
On munit $X$ de son unique mesure de Haar bi-invariante de masse $1$ et par suite d'un produit scalaire sur son algèbre de Lie $\mbox{ad}(X)$-invariant.
La multiplication à droite permet de transporter ce produit scalaire afin d'obtenir une structure de variété riemannienne sur $X$. 
Par construction, la multiplication à droite par un élément de $X$ est isométrique.
Ensuite, l'analyse spectrale de l'opérateur de Laplace-Beltrami $\Delta$ s'exprime en fonction d'un système complet $(\rho_n)_{n\in \N}$ de représentations unitaires irréductibles de $X$.
Pour tout $n\in \N$, $\rho_n:X\rightarrow U(d_n)$ est un morphisme de groupes à valeurs dans un groupe unitaire et le théorème de Peter-Weyl énonce que la famille des coefficients normalisés $(\sqrt{d_n} \rho_{n,i,j})_{n,i,j}$ est une base orthonormée de $L^2(X)$ avec $n\in \N$ et $1\leq i,j\leq d_n$.
On définit alors le sous-espace $E_{n,i}\subset L^2(X)$ engendré par les coefficients de la $i$-ème ligne de $\rho_n$ :
\bema
E_{n,i}:=\mbox{Vect}\left\{\sqrt{d_n} \rho_{n,i,j},\quad 1\leq j\leq d_n \right\}.
\enma 
On peut démontrer qu'il existe un nombre $\lambda_n\geq 0$ tel que pour $\Delta \rho_{n,i,j}=-\lambda_n^2 \rho_{n,i,j}$ pour tout $i,j$.
Ainsi, l'hypothèse de localisation spectrale \eqref{hypo-spec} est vérifiée et les $d_n$ sous-espaces $E_{n,1},\dots,E_{n,d_n}$ sont associés à la même valeur propre.
L'action de $X$ sur lui-même par multiplication à droite  est transitive et laisse invariants les sous-espaces $E_{n,i}$ :
\bema
\forall x\in X \quad \forall \chi \in X \quad \rho_{n,i,j}(x\chi)=\sum_{k=1}^{d_n} \rho_{n,i,k}(x)\rho_{n,k,j}(\chi).\enma
La somme directe orthogonale $L^2(X)=\bigoplus E_{n,i}$ est en fait la décomposition en représentations irréductibles de la représentation régulière à droite de $X$ sur $L^2(X)$.
Finalement, la proposition \ref{inje-mani} avec $G=X$ nous fournit un résultat de type Paley-Zygmund : 
\bema 
\forall p\in [1,+\infty[ \quad \pl^p(X,\oplus E_{n,i}) =L^2(X).
\enma

\section{Application à l'équation cubique des ondes sur une variété de dimension $3$}\label{pzoll}

Comme dans la partie \ref{parti-comp}, $X$ continue à désigner une variété riemannienne compacte lisse sans bord que l'on suppose pour le moment de dimension $d\geq 2$. Si l'on fixe $\kappa>0$, alors on dispose d'un choix naturel de sous-espaces $E_n$ : 
\bemar
& E_0:= &   \{u\in L^2(X),\quad \mbox{Sp}(u)\subset [0,\kappa]   \}    ,\\
\forall n\in \N^\star \quad & E_n:=& \{u\in L^2(X),\quad \mbox{Sp}(u)\subset ]\kappa n,\kappa n+\kappa] \}.
\enmar
La somme directe orthogonale des sous-espaces $E_n$ est égale à $L^2(X)$.
Pour certaines variétés, il se peut très bien qu'une infinité de sous-espaces $E_n$ soient réduits à $\{0\}$ (penser à $X=\S^d$ avec $\kappa>0$ suffisamment petit).
Si $\kappa$ est suffisamment grand, alors non seulement les sous-espaces $E_n$ ne sont pas triviaux mais les densités de probabilité $\frac{e(n,\cdot)}{\dim(E_n)}$ sont essentiellement constantes.

\begin{lemm}\label{en-var}
Il existe une constante $\kappa=\kappa(X)>0$ telle que pour tout $(n,x)\in \N \times X$ on a 
\bema
\begin{array}{rcccrl}
\Frac{1}{C(X)} (1+n)^{d-1}& \leq & e(n,x) & \leq & C(X) (1+n)^{d-1}, \\[4mm]
\Frac{1}{C(X)}(1+n)^{d-1} & \leq & \dim(E_n) & \leq & C(X) (1+n)^{d-1}.
\end{array}
\enma
\end{lemm}
\begin{demo}
Remarquons que la fonction $L^2(X)$-normalisée $\frac{1}{\sqrt{\mbox{Vol}(X)}}$ appartient à $E_0$.
Il s'ensuit que l'on a $e(0,x)\geq \frac{1}{\mbox{Vol}(X)}$ et $\dim(E_0)\geq 1$.
La continuité de $x\mapsto e(0,x)$ et la compacité de $X$ rendent trivial le cas $n=0$.
Le cas $n\geq 1$ se traite avec les estimées de la fonction spectrale sur une variété riemannienne compacte sans bord \cite{horm68} et la formule de Weyl (qui s'obtient par intégration) :
\bemar
e(0,x)+\dots+e(n,x) & = &  (2\pi)^{-d}\mbox{Vol}(\mathbb{B}_d(0,1))(\kappa n+\kappa)^{d}    + (\kappa n+\kappa)^{d-1}\mathcal{O}\left(1\right),\\
\dim(E_0)+\dots+\dim(E_n) & = & (2\pi)^{-d}\mbox{Vol}(\mathbb{B}_d(0,1))\mbox{Vol}(X)(\kappa n+\kappa)^{d}   +(\kappa n+\kappa)^{d-1}\mathcal{O}\left(1\right).
\enmar 
où le reste $\mathcal{O}\left(1\right)$ est uniforme en $x\in X$ et indépendant de $\kappa>0$. Il vient donc
\bemar
e(n,x) & = &  (2\pi)^{-d}\mbox{Vol}(\mathbb{B}_d(0,1)) \kappa^d \left[(n+1)^{d}-n^d \right]    + \kappa^{d-1}n^{d-1} \mathcal{O}\left(1\right), \\
\dim(E_n) & = & (2\pi)^{-d}\mbox{Vol}(\mathbb{B}_d(0,1))\mbox{Vol}(X) \kappa^d\left[(n+1)^{d}- n^d \right]  +\kappa^{d-1} n^{d-1} \mathcal{O}\left(1\right).
\enmar 
La conclusion découle en choisissant $\kappa=\frac{\kappa^d}{\kappa^{d-1}}$ suffisamment grand.
\end{demo}

Jusqu'à la fin de cette partie, le nombre $\kappa>0$ est choisi de sorte que la conclusion du lemme précédent soit vérifiée.
La proposition \ref{bornLp} donne immédiatement un théorème de Paley-Zygmund sur la variété $X$.

\begin{theo}\label{zoll-Lp}
Pour tout réel $p\in [1,+\infty[$, l'espace $\pl^p(X,\oplus E_n)$ est égal à $L^2(X)$.
\end{theo}

Le cas $p\geq 2$ du résultat précédent est prouvé dans \cite{burq-lebeau} mais le cas $p\in [1,2[$ nous paraît aussi digne d'intérêt et nécessite une minoration de la fonction spectrale comme le montre la proposition \ref{bornLp}.
Ajoutons que le théorème \ref{zoll-Lp} est précisément l'analogue multidimensionnel de la randomisation unidimensionnelle d'une base hilbertienne de fonctions propres uniformément bornées dans tout espace $L^p(X)$.

D'après \cite{burq-tz-jems2011}, le théorème \ref{zoll-Lp} et l'injection de Sobolev déterministe $W^{s,p}(X)\subset L^\infty(X)$ (pour $s>0$ et $p\gg 1$) sont suffisants pour étudier l'équation cubique des ondes mais nécessitent beaucoup de moments de notre procédé de randomisation.
Une version $L^\infty$ du théorème de Paley-Zygmund va nous permettre de pallier ce défaut.
Rappelons d'abord ce qu'il en est sur le tore $\mathbb{T}^d$ : 
on vérifie facilement grâce à 
\cite[Chapter 6 and 7]{kahane} que l'on a pour toute suite $(a_n)_{n\in \Z^d}$ l'implication
\begin{equation}\label{pz-to} \exists s>0\quad \sum_{n\in \Z} a_n e^{i\langle n,x\rangle} \in H_x^s(\mathbb{T}^d) \quad  \Rightarrow \quad \sum \ep_n a_n e^{i\langle n,x\rangle} \quad \mbox{converge p.s. dans }L_x^\infty(\mathbb{T}^d). \end{equation}
Bien que l'espace de Banach $L^\infty(\mathbb{T}^d)$ ne soit pas un espace de Banach de cotype fini (voir l'annexe \ref{cotype}), Marcus et Pisier ont démontré un phénomène d'universalité pour la randomisation des fonctions trigonométriques $x\mapsto e^{i\langle n,x\rangle}$ (voir \cite{pisier1981} ou \cite[Page 527, théorèmes III.5 et III.6]{queff} sous une forme intégrale).
En l'occurrence, si $(X_n)_{n\in \Z^d}$ est une suite de variables aléatoires indépendantes symétriques dont les moments $\E[|X_n|^2]$ sont uniformément bornés, alors on a l'implication 
\begin{equation}\label{marcu-pi}
\sum \ep_n a_n e^{i\langle n,x\rangle} \quad \mbox{converge p.s. dans }L_x^\infty(\mathbb{T}^d) \quad \Rightarrow \quad \sum X_n a_n e^{i\langle n,x\rangle} \quad \mbox{converge p.s. dans }L_x^\infty(\mathbb{T}^d).
\end{equation}
L'implication précédente est très surprenante car la condition de moments $\sup\limits_{n\in \N} \E[|X_n|^2]<+\infty$ est très faible par comparaison avec le théorème \ref{maurey-pi}. De ce qui précède, on déduit évidemment l'implication 
\begin{equation}\label{univ-linf}
\exists s>0\quad \sum_{n\in \Z^d} a_n e^{i\langle n,x\rangle} \in H_x^s(\mathbb{T}^d) \quad \Rightarrow \quad \sum X_n a_n e^{i\langle n,x\rangle} \quad \mbox{converge p.s. dans }L_x^\infty(\mathbb{T}^d).
\end{equation}
On peut donc se demander comment prolonger ces idées dans le cas d'une variété riemannienne compacte et sans bord.
Un obstacle majeur est la preuve de l'implication \eqref{marcu-pi}, elle est en effet assez délicate, si l'on n'impose aucune condition sur la suite $(a_n)_{n\in \Z^d}$, car elle est intimement liée au fait que l'espace de Banach des fonctions ``presque sûrement continues" sur $\mathbb{T}^d$ est de cotype $2$.
\'Etrangement, on peut démontrer l'implication \eqref{univ-linf}, sans invoquer \eqref{marcu-pi}, à l'aide d'une version intégrale de \eqref{pz-to} et de l'astuce qui consiste à remarquer que $X_n$ et $\ep_n X_n$ ont la même loi (où $\ep_n$ est indépendante de $X_n$).
Pour toute suite $(a_n)_{n\in \Z^d}$ à support fini, nous écrivons l'argument en abrégé 
\bemar
\E\left[\left\Vert \Sum_{n\in \Z^d} X_n a_n e^{i\langle n,x\rangle }  \right\Vert_{L_x^\infty(\mathbb{T}^d)}^2 \right] & = & \E_{\omega'}\E_{\omega}\left[\left\Vert \Sum_{n\in \Z ^d} \ep_n(\omega)X_n(\omega') a_n e^{i\langle n,x\rangle }  \right\Vert_{L_x^\infty(\mathbb{T}^d)}^2 \right] \\
& \leq & C(s)\E_{\omega'} \left[\Sum_{n\in \Z^d} (1+|n|)^{2s} |X_{n}(\omega')a_n|^2 \right] \\
& \leq & C(s) \left( \sup\limits_{n\in \Z^d} \E[|X_n|^2]\right) \Sum_{n\in \Z^d} (1+|n|)^{2s} |a_n|^2 .
\enmar
L'argument précédent justifie pourquoi il est plus intéressant d'estimer des espérances (qui font intervenir des moments des variables aléatoires) plutôt que des probabilités (qui font intervenir des estimations de grandes déviations).
Cela nous amène au théorème suivant qui doit être comparé à l'injection de Sobolev $H^{s}(X)\subset \CC^0(X)$ valide pour tout $s>\frac{d}{2}$ et dit moralement que presque toute fonction appartenant à $\bigcup_{s>0} H^s(X)$ est continue (on examinera le cas particulier $M_n=\mathcal{E}_n$ des isométries aléatoires).
On conviendra que $\phi_{n,1},\dots,\phi_{n,d_n}$ est une base hilbertienne \textbf{quelconque} de $E_n$.

\begin{theo}\label{pz-linf}
Considérons un réel $p\in [2,+\infty[$ et une suite de matrices aléatoires $M_n:\Omega\rightarrow \mathcal{M}_{n}(\R)$ indépendantes, orthogonalement invariantes et vérifiant la condition $\sup\limits_{n\in \N} \left[\left\Vert M_n \right\Vert_{op}^p \right]<+\infty$.
Pour tout réel $s>0$, pour tout $N\in \N$ et pour tout $u\in H^s(X)$, on a 
\begin{equation}\label{zoll-pal}
\E\left[ \left\Vert \sum_{n=0}^N \sum_{i=1}^{d_n} \left(\sum_{j=1}^{d_n} M_{n,i,j} \langle u,\phi_{n,j}\rangle \right)\phi_{n,i}(x) \right\Vert_{L_x^\infty(X)}^p\right]\leq C(X,p,s) \left(\sup\limits_{n\in \N} \E\left[\left\Vert M_n\right\Vert_{op}^p\right] \right)\left\Vert u\right\Vert_{H^s(X)}^p.
\end{equation}
En particulier, la série aléatoire 
\bema \sum_{n=0}^N \sum_{i=1}^{d_n}\left( \sum_{j=1}^{d_n} M_{n,i,j}(\omega) \langle u,\phi_{n,j}\rangle \right)\phi_{n,i} \enma
converge dans $L^p(\Omega,L^\infty(X))$ et presque sûrement dans $L^\infty(X)$.
Une conclusion similaire est valide si l'on suppose que les matrices $M_n:\Omega\rightarrow \mathcal{M}_{d_n}(\C)$ sont unitairement invariantes.
\end{theo}

La théorie développée par Burq et Tzvetkov dans \cite{burq-tz-jems2011} nous permet de donner une application du théorème de Paley-Zygmund à l'équation cubique des ondes : 
\bema
(\partial_t^2-\Delta)v+v^3=0,\quad (t,x)\in \R\times X.
\enma
Cette équation est $\frac{1}{2}$-critique en dimension $3$, c'est-à-dire que l'équation posée sur $\R^3$ et la norme $\left\Vert v\right\Vert_{\dot{H}^{\frac{1}{2}}(\R^3)}+
\left\Vert \partial_t v \right\Vert_{\dot{H}^{\frac{1}{2}-1}(\R^3)}$, où $\dot{H}^s(\R^3)$ désigne l'espace de Sobolev homogène, sont stables par le changement d'échelle qui consiste à remplacer $v$ par $\frac{1}{\delta} v(\frac{t}{\delta},\frac{x}{\delta})$ quel que soit $\delta>0$.
Afin d'alléger la rédaction, nous choisissons d'utiliser le vocabulaire de la théorie des probabilités et non celui des mesures sur $H^s(X)\times H^{s-1}(X)$.

\begin{theo}\label{edp-zoll}
Supposons $\dim(X)=3$, considérons $s>0$ et deux suites de matrices aléatoires indépendantes $M_n:\Omega\rightarrow \mathcal{M}_{d_n}(\R)$ et $M_n':\Omega\rightarrow \mathcal{M}_{d_n}(\R)$ orthogonalement invariantes et vérifiant
\begin{equation}\label{zoll-mom}
\sup\limits_{n\in \N}\E\left[\left\Vert M_n \right\Vert_{op}^{3} +\left\Vert M_n' \right\Vert_{op}^{3} \right]<+\infty.
\end{equation}
Pour tout $(v_0,v_1)\in H^s(X)\times H^{s-1}(X)$, on définit les conditions initiales aléatoires  
\begin{eqnarray} \label{zoll-ini1} \forall \omega \in \Omega \qquad  v_0^\omega & :=& \sum_{n\in \N} \sum_{i=1}^{d_n}  \left(\sum_{j=1}^{d_n}M_{n,i,j}(\omega) \langle v_0,\phi_{n,j}\rangle \right)\phi_{n,i}      ,\\ \label{zoll-ini2}
  v_1^\omega & :=& \sum_{n\in \N} \sum_{i=1}^{d_n}  \left(\sum_{j=1}^{d_n}M_{n,i,j}'(\omega) \langle v_1,\phi_{n,j}\rangle \right)\phi_{n,i}         .\end{eqnarray}
Alors pour presque tout $\omega \in \Omega$, la fonction aléatoire $(v_0^\omega,v_1^\omega)$ appartient à $H^s(X)\times H^{s-1}(X)$ et l'équation cubique des ondes
\bema  (\partial_t^2-\Delta)v+v^3=0   , \qquad    v(0,\cdot)=v_0^\omega,\qquad \dot{v}(0,\cdot)=v_1^\omega, \enma
admet une unique solution globale $v$ qui vérifie 
\bema v(t)-\cos(t\sqrt{-\Delta})v_0^\omega-\frac{\sin(t\sqrt{-\Delta})}{\sqrt{-\Delta}}v_1^\omega \in \CC_t^0(\R,H^1(X))\cap \CC_t^1(\R,L^2(X)).
\enma
Enfin, si l'on fait les deux hypothèses supplémentaires 
\bema  \min\left(\inf\limits_{n\in \N}\sigma\left(\E[|M_n|]\right),\inf\limits_{n\in \N}
\sigma\left(\E[|M_n'|]\right)\right)>0 \quad \mbox{et} \quad (v_0,v_1)\not \in \bigcup\limits_{\rho>0} H^{s+\rho}(X)\times H^{s-1+\rho}(X),\enma 
alors presque sûrement on a 
\bema (v_0^\omega,v_1^\omega)\not \in \bigcup\limits_{\rho>0} H^{s+\rho}(X)\times H^{s-1+\rho}(X) . \enma
\end{theo}

Faisons quelques commentaires sur le résultat précédent : 
\begin{enumerate}[$\bullet$ ]
\item la fin de l'énoncé précédent, la proposition \ref{ex-mat} et la remarque \ref{rema-lata} assurent l'existence de beaucoup de matrices aléatoires $M_n$ et $M_n'$ pour lesquelles les conditions initiales aléatoires \eqref{zoll-ini1} et \eqref{zoll-ini2} ont exactement la même régularité dans les espaces de Sobolev que $v_0$ et $v_1$.
Cela signifie que l'on peut bien résoudre l'équation cubique des ondes en régime sur-critique.

\item on insiste de nouveau sur le fait que la base hilbertienne $(\phi_{n,i})$ de $E_n$ est \textbf{quelconque} dans les conditions initiales aléatoires \eqref{zoll-ini1} et \eqref{zoll-ini2}.

\item de même que dans \cite{burq-tz-jems2011}, il est certainement possible de prouver un résultat analogue au théorème \ref{edp-zoll} pour $s=0$ mais le théorème \ref{pz-linf} serait inutilisable et il faudrait affaiblir sensiblement la condition de moments \eqref{zoll-mom}. On devrait alors utiliser la version $L^p$ du théorème de Paley-Zygmund.
\end{enumerate}

\section{Théorie $L^p$ presque sûre pour les harmoniques sphériques de $\S^d$}\label{harmo-sph}
Nous allons voir des résultats très différents sur une même variété selon le choix des sous-espaces $E_n\subset L^2(X)$.
L'exemple modèle est la sphère $\S^d$ sur laquelle deux exemples importants de concentration se produisent : la concentration gaussienne autour d'une géodésique et la concentration polaire.
Par cohérence avec ce qui précède, on notera $\mu_d$ la mesure volume de $\S^d$.
Faisons quelques rappels sur les harmoniques sphériques \cite[Chapter IV.2]{stein71}.
Le spectre de l'opérateur de Laplace-Beltrami sur 
\bema\S^d=\{(x_1,x_2,\dots,x_{d+1})\in \R^{d+1},\quad x_1^2+\dots+x_{d+1}^2=1\}
\enma
est ordonné par la suite $(n(n+d-1))_{n\in \N}$.
En l'occurrence, nous serons intéressés par la suite $(Y_{n})_{n\geq 1}$ de modes propres définis par
\begin{equation}\label{circul} \forall n\in \N^\star \quad Y_{n}(x):=c_{d,n}(x_1+i x_2)^n,\quad \left\Vert Y_n\right\Vert_{L^2(\S^d)}=1, \quad c_{d,n}>0.   \end{equation}
On vérifie que l'on a $\Delta Y_{n}=-n(n+d-1)Y_{n}$ et $c_{d,n}\simeq_d n^{\frac{d-1}{4}}$.
Les fonctions $Y_{n}$ sont connues pour avoir des estimations de normes dans $L^p(\S^d)$ maximales pour $2<p\leq \frac{2(d+1)}{d-1}$ et minimales pour $p\in [1,2[$ parmi les modes propres $L^2(\S^d)$-normalisés associés à la valeur propre $-n(n+d-1)$ (et cela est même optimal pour toutes les variétés riemanniennes compactes d'après les inégalités de Sogge \eqref{sogge} et \cite[Proposition 2]{sogge2010lower}).
De manière précise, il est connu et facile à vérifier que \eqref{circul} implique les estimations suivantes
\begin{equation}\label{circul-normlp}  \forall p\in [1,+\infty[\cup \{+\infty\} \qquad \forall n\in \N^\star \qquad  \left\Vert Y_{n}\right\Vert_{L^p(\S^d)}\simeq_{d,p} n^{\frac{d-1}{2}\left( \frac{1}{2}-\frac{1}{p}\right)}    .\end{equation}
La chose la plus facile que nous en déduisons, à partir de la définition \ref{defi-cv}, est que les espaces $\pl^p(\S^d,\oplus \C Y_n)$ sont distincts deux à deux.
En effet, on a $\left\Vert Y_n\right\Vert_{\pl^p(\S^d,\oplus \C Y_{n})}=\left\Vert Y_n\right\Vert_{L^p(\S^d)}$ et il est donc clair que les normes $\left\Vert\cdot\right\Vert_{\pl^p(\S^d,\oplus \C Y_n)}$ ne sont pas équivalentes pour deux valeurs distinctes de $p\geq 1$.
Comme expliqué dans la partie \ref{parti-comp}, les espaces $\pl^p(\S^d,\oplus \C Y_n)$ s'identifient à des sous-espaces de distributions sur $\S^d$.
La proposition suivante résout complètement la question de la description de ces espaces pour $p>1$.

\begin{prop}\label{beam-gaus}
On considère une suite complexe $(a_n)_{n\geq 1}$ à croissance polynomiale.
Pour tout réel $p\in ]1,+\infty[$, la distribution 
$\Sum_{n\geq 1} a_n Y_{n}$ appartient à $\pl^p(\S^d,\oplus \C Y_n)$ si et seulement si 
l'on a la condition
\begin{equation}\label{circu-cond}\Sum_{n\geq 1} \frac{1}{n^{\frac{d+1}{2}}}\left(\sum_{k=1}^n k^{\frac{d-1}{2}}|a_k|^2 \right)^{\frac{p}{2}}<+\infty.  \end{equation}
En outre, les espaces $\pl^p(\S^d,\oplus \C Y_n)$ sont stables par dualité et interpolation réelle et complexe pour $p$ parcourant $]1,+\infty[$ au sens du théorème \ref{theo-interpo}.
\end{prop}

Comme application immédiate, on obtient le fait suivant qui contraste avec le théorème de Paley-Zygmund sur $\mathbb{T}$ :  la fonction $\sum_{n\geq 2}\frac{1}{\sqrt{n}\ln(n)}Y_n$ appartient à $L^2(\S^d)$ mais la série aléatoire $\sum \frac{\ep_n}{\sqrt{n} \ln(n)}Y_n$ diverge presque sûrement dans $L^p(\S^d)$ pour tout réel $p>2$.
Expliquons brièvement la preuve de la proposition \ref{beam-gaus}.
L'idée est d'estimer les nombres
\bema
\left\Vert \sqrt{\sum_{n\geq 1} |a_n Y_n|^2} \right\Vert_{L^p(\S^d)}
\enma
Cela sera très facile pour $p\in 2\N$ en utilisant des estimées optimales vérifiées par les intégrales 
\begin{equation}\label{Yn-esti-mul} \int_{\S^d} |Y_{n_{(1)}}(x)\dots Y_{n_{(p/2)}}(x)|^2 d\mu_d(x). \end{equation}
Ensuite, on utilisera la concentration gaussienne de $Y_n$ autour de la géodésique $\{x_1^2+x_2^2=1\}\subset \S^d$ afin de vérifier l'hypothèse \eqref{loren1} : 
\bema 
\forall p>1 \quad \sup\limits_{n\geq 1} \frac{|Y_n|}{\left\Vert Y_n \right\Vert_{L^p(\S^d)}} \in L^{p,\infty}(\S^d).
\enma
On conclura alors par interpolation et dualité à l'aide des théorèmes 
\ref{theo-interpo} et \ref{theo-dual}.

\begin{rema}\label{rema-conc}
On verra que la condition \eqref{circu-cond} est exactement celle que l'on obtient avec le critère \eqref{carac-uni} en remplaçant $Y_n$ par la fonction $x\mapsto n^{\frac{d-1}{4}} \pun_{V_n}(x)$ où $V_n\subset \S^d$ est une bande de largeur $\frac{1}{\sqrt{n}}$ autour de la géodésique d'équation $x_1^2+x_2^2=1$.
Afin d'alléger la rédaction, nous n'avons traité que le cas des géodésiques sur $\S^d$ mais notre démarche est en fait plus générale.
En effet, on peut estimer 
les intégrales \eqref{Yn-esti-mul} en utilisant seulement la concentration gaussienne de $Y_n$ autour de la géodésique $\{x_1^2+x_2^2=1\}$.
Ce type de concentration se réalise dans un autre cas important, à savoir celui d'une surface $X$ qui admet une géodésique $\Gamma\subset X$ fermée elliptique et non-dégénérée. Il est alors connu, mais assez délicat à rédiger, que l'on peut construire des quasi-modes qui se concentrent autour de la géodésique $\Gamma$ avec des estimations gaussiennes (ce sont les travaux de Ralston et Babich).
\end{rema}

Par dualité, on obtient gratuitement la moitié des injections de Sobolev probabilistes des fonctions $Y_n$.

\begin{coro}\label{Yn-injesobo}
Considérons $p\in]2,+\infty[$ et $q=\frac{p}{p-1}\in ]1,2[$. Nous avons les inclusions strictes
\bema
\begin{array}{rcccl}
H^{\frac{d-1}{2}\left(\frac{1}{2}-\frac{1}{p}\right)}(\S^d) & \subsetneq & \pl^p(\S^d,\oplus \C Y_n) & \subsetneq & \bigcap\limits_{\ep>0} H^{\frac{d-1}{2}\left(\frac{1}{2}-\frac{1}{p}\right)-\ep}(\S^d),\\
\bigcup\limits_{\ep>0} H^{\frac{-(d-1)}{2}\left( \frac{1}{q}-\frac{1}{2}\right)+\ep}(\S^d) & \subsetneq & \pl^q(\S^d, \oplus \C Y_n) & \subsetneq & H^{\frac{-(d-1)}{2}\left(\frac{1}{q}-\frac{1}{2}\right)}(\S^d).
\end{array}
\enma
\end{coro}
Faisant tendre $q$ vers $1^+$ et $\ep$ vers $0^+$ dans le corollaire précédent, on constate que l'exposant $-\left(\frac{d-1}{4}\right)$ de l'énoncé du lemme \ref{identif} est optimal.

La condition \eqref{circu-cond} doit être comparée à des formules, déjà présentes dans la littérature pour d'autres familles de modes propres $(\phi_n)_{n\geq 1}$, de l'exposant critique de convergence défini comme suit
\bema \forall (a_n)\in \ell^2(\N^\star) \quad p_c(a_n):=\sup\left\{ p>2, \quad \Sum_{n\geq 1} \ep_n a_n \phi_n \mbox{ converge p.s. dans } L^p(X)\right\}.\enma
Pour les modes propres $\psi_n$ à symétrie radiale de l'opérateur Laplacien sur la boule unité $\mathbb{B}_d(0,1)\subset \R^d$, un examen de $p_c(a_n)$ a été entamé par Ayache et Tzvetkov 
\cite[Theorem 4, Proposition 2.8]{tzvetkov-ay}. Puis, Grivaux en a donné une formule très satisfaisante en fonction de la suite $(a_n)$ dans \cite{Grivaux}.
Cette analyse a été reprise pour les modes propres $\phi_n$ à symétrie radiale de l'oscillateur harmonique $-\Delta+|x|^2$ par Robert, Thomann et l'auteur \cite[Theorem 2.3]{randomh}.
Dans les deux cas précédents, dû au fait que la $n$-ième fonction propre $\phi_n$ se concentre essentiellement sur une boule $\mathbb{B}_d(0,n^{-\delta})\subset \R^d$ avec $\delta>0$, l'analyse des preuves montrerait que la formule de Grivaux prend la forme unifiée :
\begin{equation}\label{form-gene}
p_c(a_n)= \sup\limits\left\{ p>2, \quad \Sum_{k=1}^{n} \left\Vert\phi_k\right\Vert_{L^\infty}^2 |a_k|^2=\mathcal{O}\left(n^{\frac{2d\delta}{p}} \right)\right\}.
\end{equation}

On notera cependant que, contrairement à la formule \eqref{circu-cond}, la formule précédente ne permet pas de décider la convergence
de la série aléatoire $\sum \ep_n a_n \phi_n$ dans l'espace de Lebesgue critique $L^{p_c(a_n)}(X)$ (cette question a été soulevée dans \cite[remark page 16]{randomh}).

Comme le remarquent Ayache et Tzvetkov \cite[Theorem 4, remark (d)]{tzvetkov-ay}, les fonctions propres zonales de $\Delta$ sur $\S^d$ ont des propriétés de concentration similaires aux fonctions propres radiales de l'opérateur Laplacien sur $\mathbb{B}_d(0,1)$ et il est donc possible d'adapter les résultats d'un cadre à l'autre.
Faisons cette comparaison et convenons qu'une fonction sur $\S^d$ est zonale si elle ne dépend que de la première coordonnée $x_1$ d'un point $x\in \S^d$. 
D'après \cite[Chapter IV.2, Theorem 2.14]{stein71}, on sait que l'opérateur de Laplace-Beltrami $\Delta$ admet une suite de fonctions propres zonales $(Z_{n})_{n\geq 1}$ vérifiant  :
\begin{equation}\label{defi-Zn}
Z_n(x) =n^{\frac{1}{2}} P_n^{\left(\frac{d-2}{2},\frac{d-2}{2} \right)}(x_1)  ,\quad \Delta Z_n=-n(n+d-1) Z_n,\quad \left\Vert Z_n\right\Vert_{L^2(\S^d)}\simeq_d 1,
\end{equation}
où $P_n^{\left(\frac{d-2}{2},\frac{d-2}{2} \right)}$ est le $n$-ième polynôme de Jacobi, c'est-à-dire le $n$-ième polynôme orthogonal pour le poids $w\in [-1,1]\mapsto (1-w^2)^{\frac{d-2}{2}}$ et normalisé de sorte que 
\bema \forall n\in \N^\star \quad P_n^{\left(\frac{d-2}{2},\frac{d-2}{2} \right)}(1)=\comb{n+\frac{d-2}{2}}{n}\simeq_d n^{\frac{d-2}{2}}. \enma
Parmi les modes propres de $\Delta$, les fonctions $Z_n$ sont connues pour maximiser la croissance des quotients $\frac{\left\Vert Z_n\right\Vert_{L^p(\S^d)}}{\left\Vert Z_n\right\Vert_{L^2(\S^d)}}$ si $n$ tend vers $+\infty$ avec $p\geq  \frac{2(d+1)}{d-1}$ (et cela est même optimal pour toutes les variétés d'après les inégalités de Sogge \eqref{sogge}).
En l'occurrence, une formule de changement de variables (voir \eqref{chg-var}) donne
\bema
  \left\Vert Z_n\right\Vert_{L^p(\S^d)}^p =\int_{\S^d} |Z_n(x)|^p d\mu_d(x)\simeq_{d,p} n^{\frac{p}{2}} 
\int_{-1}^1     |P_n^{\left(\frac{d-2}{2},\frac{d-2}{2} \right)}(x_1)|^p (1-x_1^2)^{\frac{d-2}{2}}dx_1, 
\enma 
puis \cite[Page 391]{szeg} nous fournit les estimations des normes dans $L^p(\S^d)$ des fonctions $Z_n$ :
\begin{equation}\label{born-Lp-Yn0}\begin{array}{cccl}
 1\leq p < \frac{2d}{d-1}  & \Rightarrow \quad \left\Vert Z_{n}\right\Vert_{L^p(\S^d)} &\simeq_{d,p} & 1, \\
p=\frac{2d}{d-1} & \Rightarrow \quad  \left\Vert Z_{n}\right\Vert_{L^p(\S^d)} & \simeq_{d,p} & \sqrt[p]{\ln(n+1)},\\
\frac{2d}{d-1}<p < \infty  & \Rightarrow \quad  \left\Vert Z_n\right\Vert_{L^p(\S^d)} & \simeq_{d,p} & n^{\frac{d-1}{2}-\frac{d}{p}}.
\end{array}
\end{equation}

La proposition \ref{bornLp} et l'argumentation tenue après \eqref{circul-normlp} montrent alors que l'on a 
\begin{equation}\label{non-dual}\begin{array}{cccl} 
 \forall p\in \left[1,\frac{2d}{d-1}\right[ \quad & \pl^p(\S^d,\oplus \C Z_n) & =&\bigoplus\limits_{n\geq 1} \C Z_n,\\[4mm]
\forall p_{2}>p_1\geq \frac{2d}{d-1} \quad & \pl^{p_2}(\S^d,\oplus \C Z_n) &\subsetneq & \pl^{p_1}(\S^d,\oplus \C Z_n) . \end{array}\end{equation}

Rappelons que \eqref{non-dual} s'interprète comme une sorte de théorème de Paley-Zygmund pour les fonctions zonales : 
pour toute fonction zonale $f=\sum_{n\geq 1} a_n Z_n \in L^2(\S^d)$, la fonction aléatoire $f^\omega=\sum_{n\geq 1} \ep_n a_n Z_n$ appartient presque sûrement à $L^p(\S^d)$ pour tout réel $p\in [1,\frac{2d}{d-1}[$.
Par contre, il existe des fonctions zonales $f\in L^2(\S^d)$ telles que 
la série aléatoire définissant $f^\omega$ diverge presque sûrement dans $L^{\frac{2d}{d-1}}(\S^d)$.

D'après \eqref{non-dual}, on voit que les espaces $\pl^p(\S^d,\oplus \C Z_n)$ ne peuvent raisonnablement pas être stables par dualité. 
Des arguments rigoureux sont donnés plus loin aussi bien pour l'absence de dualité (voir la proposition \ref{Lamba-surj}) que pour l'absence d'interpolation si $p$ parcourt un intervalle ouvert contenant $\frac{2d}{d-1}$ (voir la proposition \ref{neces-interpo}).

Pour analyser de façon plus fine les espaces $\pl^p(\S^d,\oplus  \C Z_n)$ avec $p>\frac{2d}{d-1}$, on doit rappeler que les fonctions $Z_{n}$ se concentrent autour de deux boules centrées aux pôles $(\pm1,0,\dots,0)$ de rayon d'ordre $\frac{1}{n}$ avec une amplitude d'ordre $n^{\frac{d-1}{2}}$.
D'une part, les formules \eqref{born-Lp-Yn0} expliquent que cette concentration polaire est significative dans $L^p(\S^d)$ pour $p>\frac{2d}{d-1}$.
D'autre part, on pourrait appliquer la formule de Grivaux \eqref{form-gene} : pour toute suite $(a_n)\in \ell^2(\N^\star)$ et toute fonction zonale $u=\Sum_{n\geq 1} a_n Z_{n}\in L^2(\S^d)$ on aurait
\begin{equation}\label{griv-Zn}
\sup\left\{p\geq\frac{2d}{d-1}, u\in \pl^p(\S^d,\oplus \C Z_{n})  \right\} = \sup\left\{p\geq \frac{2d}{d-1},\quad \Frac{1}{n^{d+1}}\left(\Sum_{k=1}^{n} k^{d-1} |a_k|^2\right)^{\frac{p}{2}}=\mathcal{O}\left( \frac{1}{n}\right) \right\}.
\end{equation}
Le prochain énoncé rendra légitime la mise en évidence du terme $\mathcal{O}(\frac{1}{n})$.
Le théorème d'interpolation \ref{theo-interpo}, malgré son apparence abstraite, va nous permettre de répondre complètement à la question.

\begin{prop}\label{prop-zon}
On considère une suite complexe $(a_n)_{n\geq 1}$ à croissance polynomiale.
Pour tout réel $p\in ]\frac{2d}{d-1},+\infty[$, la distribution $\Sum_{n\geq 1} a_n Z_n$ appartient à $\pl^p(\S^d,\oplus \C Z_{n})$ si et seulement si l'on a 
\begin{equation}\label{zon-cond}  \sum_{n\geq 1} \frac{1}{n^{d+1}} \left(\sum_{k=1}^n k^{d-1} |a_k|^2 \right)^{\frac{p}{2}}<+\infty.  \end{equation}
Les espaces $\pl^p(\S^d,\oplus \C Z_{n})$ sont stables par interpolation réelle et complexe pour $p$ parcourant $]\frac{2d}{d-1},+\infty[$ au sens du théorème \ref{theo-interpo}.
Enfin, les injections de Sobolev probabilistes des fonctions $Z_n$ sont données par les inclusions strictes 
\begin{equation}\label{Zn-inclus}
H^{\frac{d-1}{2}-\frac{d}{p}}(\S^d) \subsetneq \pl^p(\S^d,\oplus\C Z_n) \subsetneq \bigcap\limits_{\ep>0} H^{\frac{d-1}{2}-\frac{d}{p}-\ep}(\S^d).
\end{equation}
\end{prop}

La proposition \ref{prop-zon} nous permet de préciser la formule \eqref{griv-Zn}
par comparaison logarithmique.
Examinons par exemple des coefficients de la forme $a_n=\frac{1}{n^{d\left(\frac{1}{2}-\frac{1}{p}\right)}\ln^{\frac{\beta}{p}}(n)}
$ avec $\beta\geq 0$, $p>\frac{2d}{d-1}$ et $n\geq 2$, il vient
\bema
 \frac{1}{n^{d+1}} \left(\sum_{k=2}^n k^{d-1} |a_k|^2\right)^{\frac{p}{2}}\simeq_{d,\beta} \frac{1}{n\ln^{\beta}(n)}.
\enma
La borne supérieure \eqref{griv-Zn} vaut donc $p$, mais la série aléatoire $\sum_{n\geq 2} \ep_n a_n Z_n$ converge presque sûrement dans $L^p(\S^d)$ si et seulement si $\beta>1$.

De nouveau, nous verrons que la condition \eqref{zon-cond} est exactement celle que l'on obtient avec 
le critère \eqref{carac-uni} en remplaçant $Z_n$ par $x\mapsto n^{\frac{d-1}{2}}\pun_{B(N,\frac{1}{n})}(x)$, où $N$ désigne le pôle $(1,0,\dots,0)$.
En d'autres termes, les propositions \ref{beam-gaus} et \ref{prop-zon} affirment que la concentration sur une zone déterminée (voisinage d'une géodésique ou d'un pôle) est la seule pertinente dans la théorie $L^p$ presque sûre dès lors qu'elle est significative dans $L^p(\S^d)$.
Malgré la ressemblance de leurs énoncés, les preuves des propositions \ref{beam-gaus} et \ref{prop-zon} vont être différentes car la concentration des fonctions $Y_{n}$ autour de la géodésique $\{x_1^2+x_2^2=1\}\subset \S^d$ est bien meilleure que la concentration des fonctions zonales $Z_{n}$ autour des pôles $(\pm 1,0,\dots,0)$ (voir la discussion concernant les estimées multilinéaires au début de la preuve de la proposition \ref{prop-zon}).

Pour conclure cette partie, remarquons que l'injection de Sobolev probabiliste de Tzvetkov \eqref{sobo-tzv} se reformule pour toute base hilbertienne $(\phi_n)_{n\in \N}$ de fonctions propres de $L^2(X)$ :
\bema \forall p\geq 2 \qquad H^{\delta(d,p)}(X)\subset \pl^p(X,\oplus \C \phi_n). \enma
Les injections de Sobolev probabilistes des harmoniques sphériques $Y_n$ et $Z_n$ (Corollaire \ref{Yn-injesobo} et \eqref{Zn-inclus}) montrent que l'on ne peut pas choisir un exposant plus petit que $\delta(d,p)$ en toute généralité.
La proposition \ref{bornLp} montre cependant qu'un moyen d'abaisser $\delta(d,p)$ à $0$ serait de randomiser une base hilbertienne spécifique de $L^2(\S^d)$ constituée de modes propres uniformément bornés dans tout espace $L^p(\S^d)$ avec $p\in [2,+\infty[$ (une  telle base existe par des arguments de concentration de la mesure \cite{burq-lebeau,suzzo-S3}).

\section{Théorie $L^p$ presque sûre pour l'oscillateur harmonique sur $\R^d$}\label{Lp-oscillo}

Nous allons maintenant étudier l'oscillateur harmonique multidimensionnel $-\Delta+|x|^2$ sur $L^2(\R^d)$.
Pour tout entier $n\in \N$, le $n$-ième polynôme de Hermite $H_n\in \R[X]$ est défini par la formule de Rodrigues 
\bema      H_n(x)=(-1)^n e^{x^2} \frac{d^n}{dx^n}e^{-x^2}.       \enma
On appelle $n$-ième fonction de Hermite la fonction $h_n$ normalisée dans $L^2(\R)$ définie par la formule 
\bema        h_n(x)=\frac{H_n(x)}{\sqrt{n!2^n \sqrt{\pi}}}e^{-\frac{1}{2} x^2} .\enma
En particulier, on a $h_n(-x)=(-1)^n h_n(x)$.
Pour tout $n\in \N$, on définit $E_{d,n}$ le sous-espace de $L^2(\R^d)$ dont une base hilbertienne est constituée des fonctions de Hermite multidimensionnelles
\bema  x\in \R^d \mapsto  h_{i_1}(x_1)\dots h_{i_d}(x_n) ,\quad (i_1,\dots,i_d)\in \N^d \quad i_1+\dots+i_d=n  .     \enma
Il s'ensuit que $\dim(E_{d,n})$ est un polynôme en $n$ de degré $d-1$ : 
\begin{equation}\label{herm-dn} d_n:=\dim(E_{d,n})=\frac{(n+1)\dots(n+d-1)}{(d-1)!} \sim \frac{n^{d-1}}{(d-1)!},\quad n\rightarrow +\infty. \end{equation}
Il faut savoir de plus que $E_{d,n}$ est le sous-espace propre de $-\Delta+|x|^2$ associé à la valeur propre $d+2n$ et que l'on a la somme directe orthogonale 
\bema   L^2(\R^d)=\bigoplus_{n\in \N} E_{d,n}. \enma
L'étude de l'espace $\pl^p(\R^d,\oplus E_{d,n})$ passe par la compréhension de la localisation des fonctions spectrales
\begin{eqnarray}\nonumber
   e_d(n,x)& :=& \sup\{ |u_n(x)|^2,\quad u_n\in E_{d,n}, \quad \left\Vert u_n\right\Vert_{L^2(\R^d)}=1\}      \\ \label{h-sumcar} & = & \sum_{ \substack{(i_1,\dots,i_d)\in \N^d \\i_1+\dots+i_d=n}} h_{i_1}(x_1)^2\dots h_{i_d}(x_d)^2 .\end{eqnarray}
Comme mentionné dans l'introduction, cet exemple a déjà été étudié dans \cite{randomh} par randomisation unidimensionnelle.
Le point clé qui a permis d'obtenir de meilleurs résultats que ceux des sphères est le bon comportement dans les espaces $L^p(\R^d)$ de la fonction spectrale $e_d(n,\cdot)$.
Comme l'oscillateur harmonique est par essence multidimensionnel, l'usage des espaces $\pl^p(\R^d,\oplus E_{d,n})$ est plus naturel et permet d'éviter une condition technique de contrôle des coefficients (appelée ``squeezing condition'').
Suivant l'idée selon laquelle seule la concentration de $e_d(n,\cdot)$ devrait être significative, on démontre le résultat suivant.

\begin{prop}\label{ed-maj}
Il existe une constante universelle $\gamma>0$ et 
un réel $C(d)>1$ tels que pour tout entier $n\in \N^\star$ et tout vecteur $x\in \R^d$ on a 
\begin{equation}\label{majo} \begin{array}{rcl}
|x|\leq \sqrt{2n+1}& \Rightarrow & e_d(n,x)\leq C(d) n^{\frac{d}{2}-1}\\
|x|\geq \sqrt{2n+1}& \Rightarrow & e_d(n,x)\leq C(d) n^{\frac{d}{2}-1}e^{-\gamma |x|^2} \\
\end{array}  \end{equation}
Il existe aussi une constante universelle $\alpha\in ]0,1[$ et un entier $n(d)\in \N^\star$ tels que pour tout entier $n\geq n(d)$ on a, quitte à augmenter $C(d)>1$, l'implication
\begin{equation}\label{mino} \frac{C(d)}{\sqrt{2n+1}}\leq |x|\leq \alpha \sqrt{2n+1} \quad \Rightarrow \quad  \frac{n^{\frac{d}{2}-1} }{C(d)} \leq e_{d}(n,x)\leq C(d)n^{\frac{d}{2}-1}.   \end{equation}
\end{prop}

On justifiera plus loin que la fonction $x\mapsto e_d(n,x)$ est invariante par rotation autour de $0$, si bien que l'on a $e_d(n,x)=e_d(n,(|x|,0,\dots,0))$.
\`A titre d'exemple, on examine les graphes pour $d=2$ et $n\in\{5,10,50\}$ de $e_d(n,x)$ en fonction de $|x|\in \left[0,\frac{3}{2}\sqrt{2n+1}\right]$.
\begin{center}
\begin{tabular}{ccc}
   \includegraphics[scale=0.24]{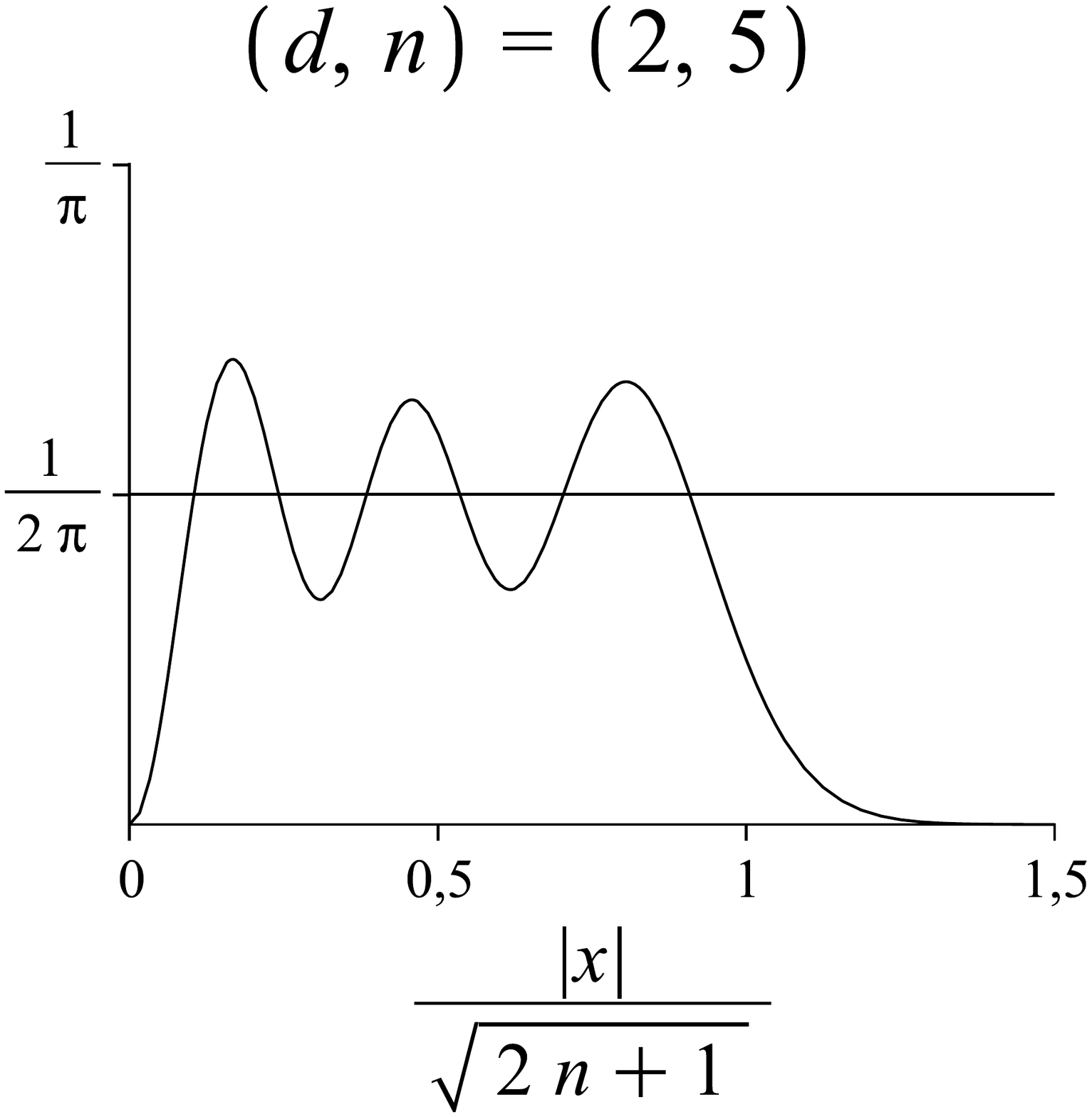} &
   \includegraphics[scale=0.24]{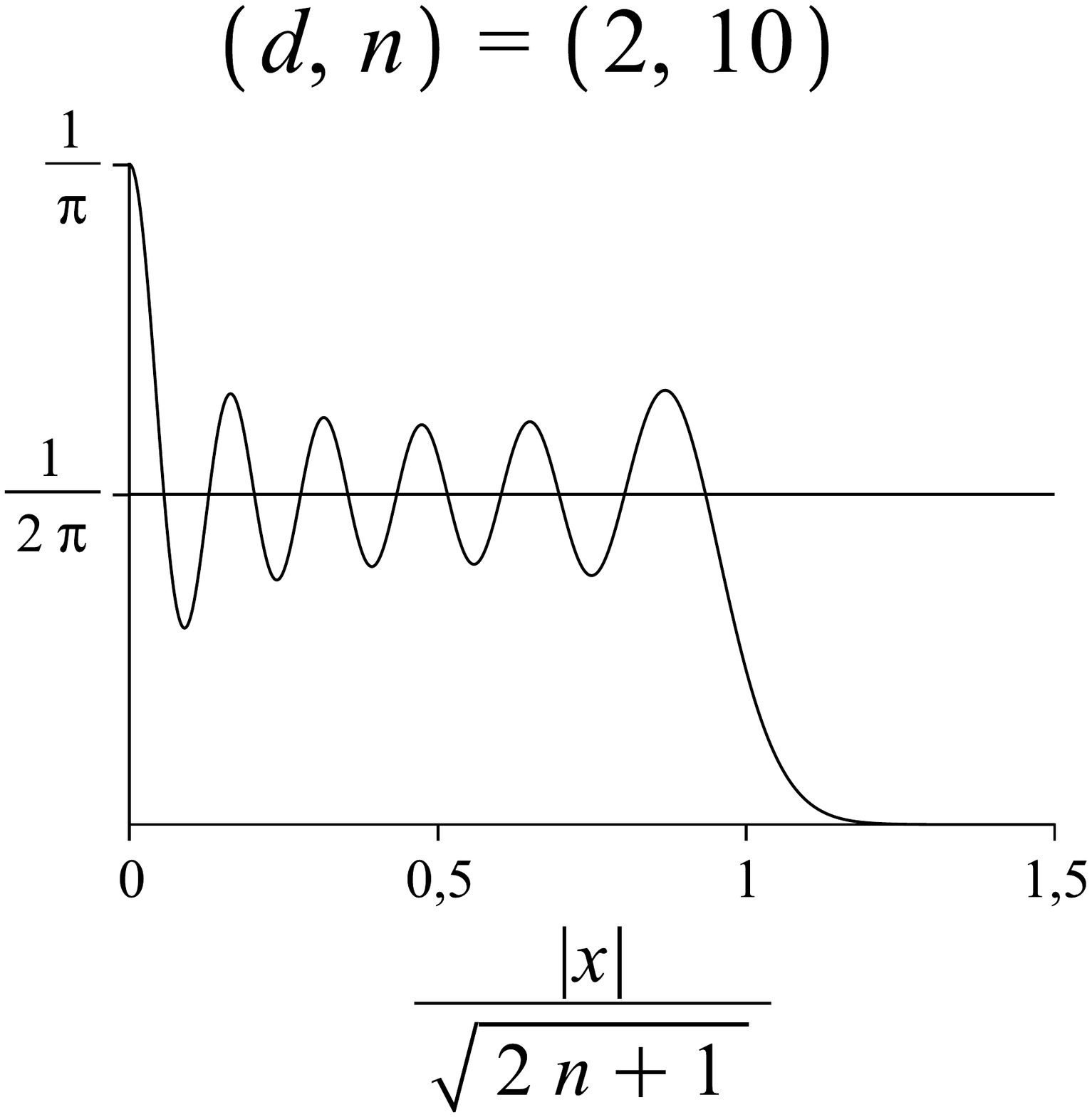} & \includegraphics[scale=0.24]{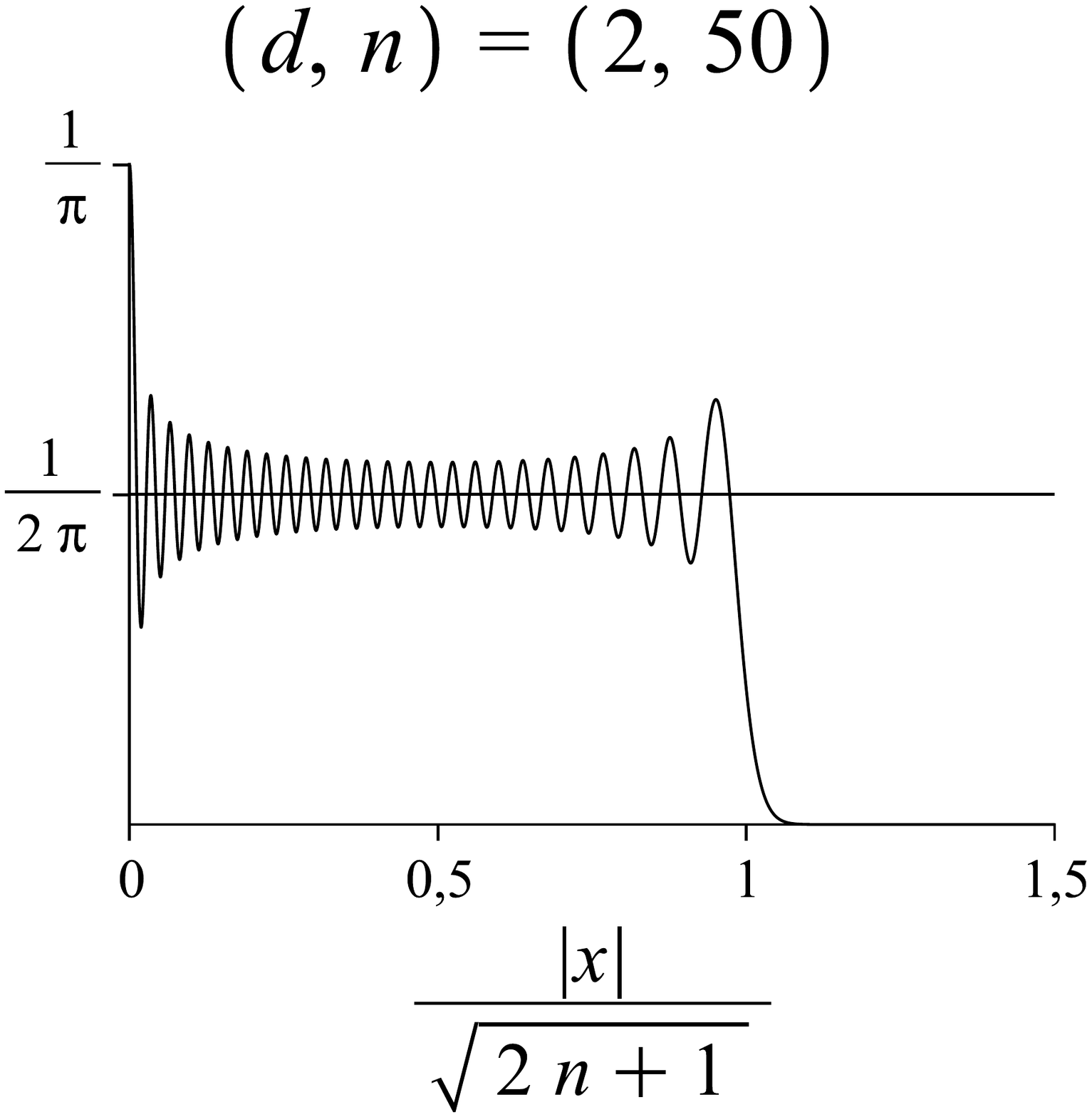} 
\end{tabular}
\end{center}

La proposition \ref{ed-maj} et les graphes ci-dessus suggèrent que l'on peut considérer, en première approximation, que la fonction spectrale $e_d(n,\cdot)$ se localise uniformément sur le compact $\overline{\mathbb{B}_d(0,\sqrt{2n+1})}$.
C'est-à-dire que l'on peut assimiler $e_d(n,x)$ à $c(d)n^{\frac{d}{2}-1} \pun_{\mathbb{B}_d(0,\sqrt{2n+1})}(x)$ où $c(d)$ ne dépend que de $d$.
Dans ce cas, l'équivalent \eqref{herm-dn} et le fait que $\frac{1}{d_n} e_d(n,\cdot)$ soit une densité permettent de prévoir que l'on devrait avoir
\bema
c(d)=\frac{1}{(d-1)!  2^{\frac{d}{2}} \mbox{Vol}(\mathbb{B}_d(0,1))}.
\enma
Par exemple, on devrait comparer $e_2(n,x)$ à $\frac{1}{2\pi}\pun_{\mathbb{B}_2(0,\sqrt{2n+1})}(x)$.
Par contre, nous sommes obligés d'éviter, en toute rigueur, un voisinage de l'origine pour effectuer une minoration uniforme de $e_d(n,\cdot)$.
En effet, pour tout entier $n$ impair la fonction $h_n$ est impaire et donc $e_d(n,0)$ est nul (voir \eqref{h-sumcar}).

La majoration \eqref{majo} est obtenue en suivant les constantes 
dans la démonstration de \cite[Lemma 3.2.2]{thanga}.
Quant à la minoration de \eqref{mino}, elle semble nouvelle. Sa démonstration est plus subtile et utilise des approximations essentiellement optimales des fonctions de Hermite dues à Muckhenhoupt.
Une application immédiate de la proposition \ref{ed-maj} est donnée par les estimations
\begin{equation}\label{hermLp}
\forall p\in [1,+\infty[ \cup \{+\infty\} \quad \forall n\in \N^\star \quad  
\left\Vert\sqrt{e_d(n,\cdot)}\right\Vert_{L^{p}(\R^d)} \simeq_{d,p} n^{\frac{d-1}{2}-\frac{d}{2}\left( \frac{1}{2}-\frac{1}{p}\right)}=n^{\frac{1}{2}\left[\frac{d}{2}-1+\frac{d}{p}\right]}.
\end{equation}

De nouveau, la majoration $\left\Vert e_d(n,\cdot)\right\Vert_{L^{\frac{p}{2}}(\R^d)} \lesssim_{d,p} n^{\frac{d}{2}-1+\frac{d}{p}}$ est connue pour $p\geq 2$ et est généralement traitée par interpolation entre $p=2$ et $p=+\infty$, mais nous ne connaissons pas de référence où l'optimalité est prouvée (voir \cite[Lemma 3.5]{PRT1} et les références indiquées).
Les hypothèses du théorème \ref{theo-dual} de dualité sont alors très faciles à vérifier. D'une part, on a
\begin{equation}\label{herm-dual}
\forall p\in ]1,+\infty[ \quad \left\Vert\sqrt{e_d(n,\cdot)}\right\Vert_{L^p(\R^d)}\left\Vert\sqrt{e_d(n,\cdot)}\right\Vert_{L^{\frac{p}{p-1}}(\R^d)}  \lesssim_{d,p}  \dim(E_{d,n}).
\end{equation}
D'autre part, en utilisant \eqref{majo} et \eqref{hermLp}, on a
\begin{eqnarray}\nonumber
\forall n\in \N^\star \quad \Frac{\sqrt{e_d(n,x)}}{\left\Vert\sqrt{e_d(n,\cdot)}\right\Vert_{L^p(\R^d)}} & \lesssim_{d,p} & \Frac{1}{n^{\frac{d}{2p}}} \pun_{\mathbb{B}_d(0,\sqrt{2n+1})}(x) + \Frac{e^{-\frac{1}{2}\gamma |x|^2}}{n^{\frac{d}{2p}}} \\ \nonumber
\sup\limits_{n\geq 1}\Frac{\sqrt{e_d(n,x)}}{\left\Vert\sqrt{e_d(n,\cdot)}\right\Vert_{L^p(\R^d)}}& \lesssim_{d,p} & \Frac{1}{|x|^{\frac{d}{p}}}+e^{-\frac{1}{2}\gamma|x|^2} \\ \label{herm-lor}
\sup\limits_{n\geq 0}\Frac{\sqrt{e_d(n,x)}}{\left\Vert\sqrt{e_d(n,\cdot)}\right\Vert_{L^p(\R^d)}} & \in & L_x^{p,\infty}(\R^d) .
\end{eqnarray}

Ainsi, on sait par avance que les espaces
$\pl^p(\R^d,\oplus E_{d,n})$ sont stables par dualité.
Tout comme pour les variétés compactes, ces espaces vont s'identifier à des sous-espaces de distributions sur $\R^d$. Pour le voir, commençons par rappeler la définition des espaces de Sobolev naturellement associés à l'oscillateur harmonique.
En notant $\Pi_n:L^2(\R^d)\rightarrow L^2(\R^d)$ le projecteur orthogonal sur $E_{d,n}$, on a : 
\bemar   \forall s\geq 0 \quad \mathcal{H}^s(\R^d)& : =& \mbox{Dom}\left(\left(-\Delta+|x|^2\right)^{\frac{s}{2}}\right) \\
& =& \left\{u\in L^2(\R^d), \quad \Sum_{n\in \N} (1+n)^{s} \left\Vert \Pi_n(u) \right\Vert_{L^2(\R^d)}^2 <+\infty    \right\}. \enmar
Rappelons que cet espace abstrait admet, lorsque $s\in \N$, la caractérisation fonctionnelle suivante (voir la preuve de \cite[Lemma 2.4]{yaji04}) : 
\bema   \HH^s(\R^d)=\left\{u\in L^2(\R^d),   \forall (m_0,m_1,\dots,m_d)\in \N^d, \quad m_0+\dots+m_d\leq s \Rightarrow  |x|^{m_0} \partial_{x_1}^{m_1}\dots\partial_{x_d}^{m_d} u \in L^2(\R^d)  \right\}.             \enma
Toute fonction $\varphi$ de l'espace de Schwartz $\mathcal{S}(\R^d)=\bigcap\limits_{s\in \N} \mathcal{H}^s(\R^d)$ admet alors une décomposition 
\bema \varphi=\sum_{n\geq 0} \Pi_n(\varphi),\qquad \forall \alpha>0 \quad \forall n\in \N \quad \left\Vert \Pi_n(\phi)\right\Vert_{L^2(\R^d)}\leq \frac{C(\alpha)}{(1+n)^\alpha}.          \enma
On en déduit par dualité que toute distribution tempérée $u\in \mathcal{S}'(\R^d)$ admet une décomposition en série faiblement convergente pour la dualité $(\mathcal{S}(\R^d),\mathcal{S}'(\R^d))$ :
\bema u=\sum_{n\in \N} \Pi_n(u) ,\quad \exists \alpha>0\quad \forall n\in \N\quad \left\Vert \Pi_n(u) \right\Vert_{L^2(\R^d)}\leq C(\alpha)(1+n)^\alpha.    \enma
Cela nous amène à définir des espaces de Sobolev pour tout $s\in \R$
\bema   \forall s\in \R \quad \mathcal{H}^s(\R)=\left\{u\in \mathcal{S}'(\R^d), \quad \Sum_{n\in \N} (1+n)^{s} \left\Vert \Pi_n(u) \right\Vert_{L^2(\R^d)}^2 <+\infty    \right\}   .   \enma
Nous pouvons maintenant énoncer un lemme analogue 
au lemme \ref{identif}.

\begin{lemm}\label{identif2}
Pour tous $p\in [1,+\infty[$ et $(u_n)_{n\in \N} \in \pl^p(\R^d,\oplus E_{d,n})$, la série $\sum u_n$ converge pour la dualité $(\mathcal{S}(\R^d),\mathcal{S}'(\R^d))$ vers une distribution tempérée.
Par suite, on peut identifier $\pl^p(\R^d,\oplus E_{d,n})$ au sous-espace des distributions tempérées $u\in \mathcal{S}'(\R^d)$ vérifiant 
\bema \left\Vert u \right\Vert_{\pl^p(\R^d,\oplus E_n)}=\left\Vert \sqrt{\sum_{n\in \N} \left\Vert \Pi_n(u) \right\Vert_{L^2(\R^d)}^2 \frac{e_d(n,x)}{\dim(E_{d,n})} } \right\Vert_{L_x^p(\R^d)} <+\infty .       \enma
 \end{lemm}
\begin{demo}
On invoque l'inégalité triviale :
\bema \left\Vert u \right\Vert_{\pl^p(\R^d,\oplus E_{d,n}) }    \geq  \sup\limits_{n\in \N}  \Frac{\left\Vert u_n \right\Vert_{L^2(\R^d)}}{\sqrt{\dim(E_{d,n})}} \left\Vert\sqrt{e(n,\cdot)}\right\Vert_{L^p(\R^d)} . \enma
Les équivalents  $\dim(E_{d,n})\simeq_d n^{d-1}$ et \eqref{hermLp} assurent que $(\left\Vert u_n \right\Vert_{L^2(\R^d)})_{n\in \N}$ est à croissance polynomiale.
\end{demo}

La proposition \ref{ed-maj}, le théorème \ref{theo-dual} de dualité et le théorème \ref{theo-interpo} d'interpolation vont nous permettre de décrire complètement les sous-espaces $\hl$ et leurs propriétés de dualité et d'interpolation pour $p>1$.

\begin{theo}\label{HLpps}
Pour tout réel $p\in [1,+\infty[$ et toute distribution tempérée $u\in \mathcal{S}'(\R^d)$, on a 
\begin{equation}\label{harmo-cond}
\left\Vert u\right\Vert_{\pl^p(\R^d,\oplus E_{d,n})}\simeq_{d,p} \left\Vert\Pi_{0}(u)\right\Vert_{L^2(\R^d)}+\left[\Sum_{n\geq 1} n^{\frac{d}{2}-1}\left( \Sum_{k \geq n} \frac{ \left\Vert\Pi_k(u)\right\Vert_{L^2(\R^d)}^2 }{k^{\frac{d}{2}}} \right)^{\frac{p}{2}} \right]^{\frac{1}{p}}.\end{equation}
En outre, nous avons les propriétés de dualité et d'interpolation : 
\begin{enumerate}[i) ]
\item pour tout $p\in ]1,+\infty[$, on pose $q=\frac{p}{p-1}$ l'exposant conjugué. L'injection canonique 
\bema\Lambda_p:\pl^q(\R^d,\oplus E_{d,n})\rightarrow \hl'\enma
\item[ ] qui à un élément $w\in \pl^q(\R^d,\oplus E_{d,n})$ associe la forme linéaire 
\bema u\in \pl^p(\R^d,\oplus E_{d,n}) \mapsto  \sum_{n\geq 0} \langle \Pi_n(u),\Pi_n(w)\rangle_{L^2(\R^d)}     \enma
\item[ ] est bien définie et est un isomorphisme d'espaces de Banach.
\item les espaces $(\hl)_{p\in ]1,+\infty[}$ sont stables par interpolation complexe et réelle au sens du théorème \ref{theo-interpo}.
\end{enumerate}
\end{theo}
Le théorème précédent ressemble aux propositions \ref{beam-gaus} et \ref{prop-zon} mais sa preuve est bien plus simple car la concentration des fonctions $e_d(n,\cdot)$ est bien meilleure que celle des harmoniques sphériques $Y_{n}$ et $Z_n$ étudiées dans la partie \ref{harmo-sph}.
En effet, la proposition \ref{ed-maj} assure que l'on peut brutalement contrôler $e_d(n,x)$ par le terme gaussien $e^{-\gamma |x|^2}\ll 1$ à l'extérieur de la boule $\mathbb{B}_d(0,\sqrt{2n+1})$. En d'autres termes, l'approximation de $e_d(n,x)$ par $n^{\frac{d}{2}-1}\pun_{\mathbb{B}_d(0,\sqrt{2n+1})}(x)$ est si bonne que la preuve de  \eqref{harmo-cond} ne nécessite pas un argument d'interpolation.
Il s'agit d'une sorte d'effet régularisant de l'oscillateur harmonique.

L'équivalence de normes \eqref{harmo-cond} implique déjà quelques faits non triviaux : 
\begin{enumerate}[$\bullet$ ]
\item pour tout $p\in [1,+\infty[$ on a l'inclusion $\hl\subset \mathcal{H}^{-\frac{d}{2}}(\R^d)$. Cela signifie qu'il faut un minimum de régularité pour espérer arriver presque sûrement dans $L^p(\R^d)$.
Cela doit être comparé au lemme \ref{identif} sur les variétés compactes.
\item pour tous réel $p_1<p_2$ on a l'inclusion stricte $\pl^{p_1}(\R^d,\oplus E_{d,n})\subset \pl^{p_2}(\R^d,\oplus E_{d,n})$. D'une part, cela contraste fortement avec le fait que $L^{p_1}(\R^d)$ n'est pas inclus dans $L^{p_2}(\R^d)$.
D'autre part, le théorème de Paley-Zygmund est vérifié pour l'oscillateur harmonique.
\end{enumerate}

Il est temps à présent d'énoncer les injections de Sobolev probabilistes de l'oscillateur harmonique multidimensionnel qui s'avèrent bien meilleures que celles des variétés compactes.

\begin{theo}\label{inje-sobo}
Pour tout réel $p\in ]2,+\infty[$, on note $q=\frac{p}{p-1}\in ]1,2[$. Nous avons les inclusions strictes
\begin{equation*} \begin{array}{rcccl}
   \mathcal{H}^{-d\left(\frac{1}{2}-\frac{1}{p}\right)}(\R^d) &\subsetneq & \hl &  \subsetneq & \bigcap\limits_{\ep>0}     \mathcal{H}^{-d\left(\frac{1}{2}-\frac{1}{p}\right)-\ep}(\R^d) ,   \\
   \bigcup\limits_{\ep>0} \mathcal{H}^{d\left(\frac{1}{q}-\frac{1}{2}\right)+\ep}(\R^d) & \subsetneq &  \pl^q(\R^d,\oplus E_{d,n})  & \subsetneq &  \mathcal{H}^{d\left(\frac{1}{q}-\frac{1}{2}\right)}(\R^d). \end{array}\end{equation*}
\end{theo}

Remarquons que pour tout $q\in [1,2]$, l'espace $\pl^q(\R^d,\oplus E_{d,n})\subset L^2(\R^d)$ est un espace de fonctions (alors que pour une variété compacte $X$, on a vu que $\pl^q(X,\oplus E_n)$ est en général seulement un espace de distributions).
Dans le cas $p>2$, les injections de Sobolev déterministes s'écrivent ainsi : 
\bema  \mathcal{H}^{d\left(\frac{1}{2}-\frac{1}{p}\right)}(\R^d)
\subset H^{d\left(\frac{1}{2}-\frac{1}{p}\right)}(\R^d)\subset L^p(\R^d)  .\enma
En autorisant un aléa pour arriver dans $L^p(\R^d)$, 
l'injection de Sobolev probabiliste $ \mathcal{H}^{-d\left(\frac{1}{2}-\frac{1}{p}\right)}(\R^d)
\subset \pl^p(\R^d,\oplus E_{d,n})$ assure donc presque sûrement un gain de $2d\left( \frac{1}{2}-\frac{1}{p}\right)$ dérivées.

\section{Preuve du théorème \ref{mplp}, randomisation avec des matrices aléatoires}\label{part-mplp}

Comme nous manipulerons des matrices de tailles différentes, il sera commode de faire le raccourci suivant : à la place de 
``\textit{une suite $(b_n)_{n\in \N}$ de matrices carrées telles que pour tout $n\in \N$ la matrice $b_n$ soit de taille $d_n\times d_n$ à coefficients dans $B$}'', nous écrirons
``\textit{une suite de matrices carrées $b_n\in\mathcal{M}_{d_n}(B)$}".
Avec les changements évidents, on utilisera la même convention pour des matrices aléatoires.

Rappelons que $(\mathcal{E}_n)_n$ est une suite de matrices aléatoires indépendantes et que la matrice aléatoire $\mathcal{E}_n$ suit une loi uniforme dans le groupe orthogonal $O_{d_n}(\R)$ pour tout $n\in \N$ (de même pour $W_n$ qui suit un loi uniforme dans le groupe unitaire $U_{d_n}(\C)$).
Dans ce contexte, les inégalités de Kahane-Khintchine démontrées par Marcus et Pisier s'énoncent comme suit (voir \cite[Page 91, Corollary 2.12]{pisier1981}).

\begin{prop}
Pour tous réels $p>q\geq 1$, il existe une constante $K_{p,q}\geq 1$ telle que, pour tout espace de Banach $B$, pour tout entier $N\in \N$ et toute suite de matrices carrées $b_n\in\mathcal{M}_{d_n}(B)$, nous avons
\begin{equation}\label{KKMP}
\E \left[ \left\Vert \sum_{n=0}^N \tr (\mathcal{E}_n b_n) \right\Vert_B^q \right]^{\frac{1}{q}}\leq \E \left[\left\Vert  \sum_{n=0}^N \tr( \mathcal{E}_n b_n) \right\Vert_B^p \right]^{\frac{1}{p}} \leq K_{p,q} \E \left[ \left\Vert \sum_{n=0}^N \tr (\mathcal{E}_n b_n) \right\Vert_B^q \right]^{\frac{1}{q}}.
\end{equation}
De même que pour les inégalités de Kahane-Khintchine \eqref{KK}, il existe une constante numérique $K\geq 1$ telle que l'on a $K_{p,q}\leq K_{p,1}\leq K\sqrt{p}$.
Des inégalités similaires sont valides pour les matrices aléatoires $W_n$ à la place de $\mathcal{E}_n$.
\end{prop}

On se permet d'utiliser la même notation $K_{p,q}$ dans \eqref{KKMP} et dans \eqref{KK} car les inégalités \eqref{KKMP} impliquent les inégalités \eqref{KK}.
Les inégalités \eqref{KKMP} jouent un rôle clé pour montrer le théorème suivant \cite[page 92, Corollary 2.14]{pisier1981} où nous forçons l'apparition du terme multiplicatif $\sqrt{d_n}$ par cohérence avec la suite.

\begin{theo}\label{marcus-pis}
Considérons un espace de Banach complexe $B$ et une suite de matrices carrées $b_n\in \mathcal{M}_{d_n}(B)$.
Les propriétés suivantes sont équivalentes 
\begin{enumerate}[i) ]
\item la série $\sum_{n} \sqrt{d_n} \tr(\mathcal{E}_n b_n)$ converge dans $L^p(\Omega,B)$ pour \textbf{un} réel $p\in [1,+\infty[$,
\item la série $\sum_{n} \sqrt{d_n} \tr(\mathcal{E}_n b_n)$ converge dans $L^p(\Omega,B)$ pour \textbf{tout} réel $p\in [1,+\infty[$, 
\item la série $\sum_{n}\sqrt{d_n} \tr(\mathcal{E}_n b_n)$ converge presque sûrement dans $B$,
\item la série $\sum_{n} \sqrt{d_n} \tr(W_n b_n)$ converge dans $L^p(\Omega,B)$ pour \textbf{un}  réel $p\in [1,+\infty[$,
\item la série $\sum_{n} \sqrt{d_n} \tr(W_n b_n)$ converge dans $L^p(\Omega,B)$ pour \textbf{tout} réel $p\in [1,+\infty[$,
\item la série $\sum_{n}\sqrt{d_n} \tr(W_n b_n)$ converge presque sûrement dans $B$.
\end{enumerate}
\end{theo}

Sans hypothèse supplémentaire sur l'espace de Banach $B$, on ne peut pas en toute généralité ajouter d'autres lois matricielles, comme les gaussiennes $\frac{1}{\sqrt{d_n}}[g_{n,i,j}]$ (voir \eqref{defi-gau}), en plus des matrices aléatoires $(\mathcal{E}_n)$ et $(W_n)$.
En effet, ne serait-ce que pour le cas unidimensionnel, i.e. $d_n=1$ pour tout $n\in \N$, le théorème de Maurey-Pisier (Théorème \ref{theo-maupign}) énonce que l'équivalence suivante 
\bema  \sum \ep_n b_n \quad \mbox{converge presque sûrement} \quad\Leftrightarrow \quad\sum g_n b_n  \quad \mbox{converge presque sûrement} , \enma
où $(g_n)_{n\geq 0}$ est une suite i.i.d. de gaussiennes suivant une loi $\mathcal{N}_{\C}(0,1)$, caractérise les espaces de Banach de cotype fini.
Par exemple, $L^\infty(X)$ n'est généralement pas de cotype fini tandis que $L^p(X)$ est de cotype fini pour tout $p\in [1,+\infty[$.

Il n'est pas clair pour l'auteur que le théorème de randomisation universelle de Maurey-Pisier (Théorème \ref{theo-maupi}) soit applicable pour comprendre l'universalité de la randomisation dans un cadre multidimensionnel.
Cela dit, nous allons voir que dans le cas particulier $B=L^p(X)$, avec $p\in [1,+\infty[$, on peut considérablement préciser le théorème \ref{marcus-pis}.
Nous allons exploiter la preuve du théorème \ref{marcus-pis}, c'est-à-dire 
l'utilisation systématique des inégalités \eqref{KKMP} et d'un principe de contraction (Théorème \ref{contrac}).
L'annexe \ref{treillis} donnera un éclairage sur les propriétés géométriques de l'espace de Banach $L^p(X)$ utilisées dans notre argumentation.
 Commençons par le lemme élémentaire suivant.

\begin{lemm}\label{Lp-ccv}
Considérons un réel $p\in [1,+\infty[$, des matrices aléatoires $M_0\in \mathcal{M}_{d_1}(\R)$,...,$M_N\in \mathcal{M}_{d_N}(\R)$ et des matrices $b_0\in\mathcal{M}_{d_1}(L^p(X))$,...,$b_{N}\in \mathcal{M}_{d_N}(L^p(X))$. Alors on a l'inégalité : 
\bema  \forall q\in [p,+\infty[\qquad \E_\omega\left[ \left\Vert \sum_{n=0}^N \sqrt{d_n} \tr(M_n(\omega) b_n) \right\Vert_{L^p(X)}^{q} \right] ^{\frac{1}{q}}
\leq  \left\Vert \E_\omega\left[\left\vert  \sum_{n=0}^N \sqrt{d_n} \tr(M_n(\omega) b_n)\right\vert^{q}\right]^{\frac{1}{q}}\right\Vert_{L^p(X)} . \enma
\end{lemm}
\begin{demo}
Cela découle par interpolation entre $q=p$ et $q=+\infty$.
\end{demo}

Lorsque l'on choisit $B=\C$ (ou plus généralement un espace de Hilbert ou $L^p(X)$), les inégalités \eqref{KKMP} donnent un résultat simple et bien connu dans le cas unidimensionnel.
Pour exprimer ce résultat, on notera $\left\Vert \cdot\right\Vert_{HS}$ la norme de Hilbert-Schmidt d'une matrice carrée $A\in \mathcal{M}_{d_n}(\C)$ :
\bema \left\Vert A \right\Vert_{HS} := \sqrt{\Sum_{i,j=1}^{d_n} |a_{ij}|^2 } .\enma

\begin{lemm}
Pour toute suite de matrices $a_n\in \mathcal{M}_{d_n}(\C)$ et pour tout $N\in \N$, on a
\begin{eqnarray}\label{KKMP-R}   \E \left[ \left\vert \sum_{n=0}^N \sqrt{d_n} \tr (\mathcal{E}_n a_n) \right\vert^2 \right]= \sum_{n=0}^N \left\Vert a_n \right\Vert_{HS}^2 . \end{eqnarray}
Par suite, pour tous réels $p$ et $q$ appartenant à $[1,+\infty[$ et pour toute suite de matrices $b_n\in \mathcal{M}_{d_n}(L^p(X))$, on a 
\begin{equation}\label{KK-LP}
\E\left[\left\Vert \Sum_{n=0}^N \sqrt{d_n} \tr( \mathcal{E}_n b_n(x)) \right\Vert_{L_x^p(X)}^q \right]^{\frac{1}{q}} \simeq_{p,q} \left\Vert \sqrt{ \Sum_{n=0}^N \left\Vert b_n(x)\right\Vert_{HS}^2} \right\Vert_{L_x^p(X)} .
\end{equation}
La même conclusion est valide en remplaçant $(\mathcal{E}_n)_{n\geq 0}$ par $(W_n)_{n\geq 0}$.
\end{lemm}
\begin{demo}
Commençons par l'égalité \eqref{KKMP-R}.
Par indépendance des matrices aléatoires $\mathcal{E}_n$, il nous suffit de prouver les deux points suivants pour tout entier $d\in \N^\star$ et toute matrice $A\in \mathcal{M}_d(\C)$ : 
\begin{eqnarray} \nonumber
\Int_{O_d(\R)} \tr (PA) dP & = & 0 ,\\[3mm] \nonumber
\Int_{O_d(\R)}  |\tr (PA)|^2 dP & = & \Frac{\tr(\tt \overline{A} A)}{d} ,
\end{eqnarray}
où l'on note $dP$ la mesure de Haar normalisée du groupe compact $O_d(\R)$.
Le cas $d=1$ étant trivial, on suppose $d\geq 2$.
La nullité de la première intégrale est clair par le changement de variables $P\mapsto -P$. Pour la seconde intégrale, il suffit de traiter le cas où tous les coefficients de $A$ sont réels.
En effet, le cas complexe découle du cas réel en écrivant $A=A_x+iA_y$ avec $A_x,A_y\in \mathcal{M}_{d}(\R)$.
En notant la matrice symétrique $|A|=\sqrt{\tt AA}$, on peut factoriser $A$ sous la forme $P|A|=PQDQ^{-1}$ avec $P,Q\in O_d(\R)$ et $D$ matrice diagonale dont les valeurs propres $\mu_1,\dots,\mu_d$ sont les valeurs singulières de $A$. La propriété d'invariance de la trace et de la mesure de Haar $dP$ simplifie l'intégrale
\begin{equation*}
\begin{array}{rcl}
\Int_{O_d(\R)} \tr (PA)^2 dP  & =&  \Int_{O_d(\R)} \tr( PD)^2 dP\\[3mm]
& = &   \Int_{O_d(\R)} \left(p_{11} \mu_{1}+\dots+p_{dd} \mu_d \right)^2 dP    \\[3mm]
& = & \Sum_{i,j=1}^d \mu_i \mu_j \Int_{O_d(\R)} p_{ii} p_{jj} dP .
\end{array}
\end{equation*}
En effectuant les changements de coordonnées $P\mapsto E_{ij}P$ et $P\mapsto PE_{ij}$ où $E_{ij}$ est la matrice orthogonale associée à la permutation qui transpose $i$ et $j$ ainsi que les $d$ changements de coordonnées $P\mapsto \Delta_{k} P$ où $\Delta_{k}$ est la matrice diagonale $\mbox{Diag}(1,\dots,1,\underbrace{-1}_{k},1,\dots,1)$, on obtient 
\begin{equation*}
\begin{array}{rcl}
\forall i\neq j \quad \Int_{O_d(\R)} p_{ii}p_{jj} dP & =&  0 , \\
\forall i,j \qquad \Int_{O_d(\R)} p_{ij}^2 dP & = & \Int_{O_d(\R)} p_{11}^2 dP .
\end{array}
\end{equation*}
Par moyenne, la dernière intégrale précédente vaut 
\begin{equation*}
\frac{1}{d^2} \int_{O_d(\R)} \tr( \tt PP) d P =  \frac{1}{d}.
\end{equation*}
Cela prouve \eqref{KKMP-R} car $\tr(\tt A A)=\mu_1^2+\dots+\mu_d^2$.

Prouvons maintenant \eqref{KK-LP}. D'après les inégalités de  Kahane-Khintchine-Marcus-Pisier \eqref{KKMP}, il suffit de traiter le cas $q=p$. Le théorème de Fubini donne alors 
\bema \E\left[ \left\Vert \sum_{n=0}^N \sqrt{d_n} \tr (\mathcal{E}_n b_n(x)) \right\Vert_{L_x^p(X)}^p \right] =\int_{X} \E\left[ \left\vert \sum_{n=0}^N \sqrt{d_n} \tr(\mathcal{E}_n b_n(x))  \right\vert^p \right]d\mu(x)  .      \enma
De nouveau, les inégalités de Kahane-Khintchine-Marcus-Pisier sur $\C$ montrent que l'intégrale précédente est équivalente à la suivante
\bema    \int_X \E \left[ \left\vert \sum_{n=0}^N \sqrt{d_n} \tr( \mathcal{E}_n b_n(x))\right\vert^2 \right]^{\frac{p}{2}} d\mu(x)   .\enma
L'égalité \eqref{KKMP-R} achève la preuve.
\end{demo}

Nous avons à présent tous les moyens pour aborder randomisation de matrices aléatoires plus générales.
 
\begin{prop}\label{pietch-lem}
Fixons $p\in [2,+\infty[$ et considérons une suite de matrices aléatoires $M_{n}:\Omega \rightarrow \mathcal{M}_{d_n}(\R)$ indépendantes, orthogonalement invariantes et vérifiant l'hypothèse de bornitude des moments   
\bema   \sup\limits_{n\in \N} \E\left[ \left\Vert M_n\right\Vert_{op}^{p} \right] <+\infty  . \enma
Il existe une constante universelle $K\geq 1$ telle que, pour toute suite $a_n\in\mathcal{M}_{d_n}(\C)$, nous avons
\begin{equation}\label{pietch}\forall N\in \N\quad    \E\left[  \left\vert \sum_{n=0}^N \sqrt{d_n} \tr( M_n a_n) \right\vert^p \right]^{\frac{1}{p}} \leq    K\sqrt{p}\left( \sup\limits_{0 \leq n\leq N} \E\left[ \left\Vert M_n\right\Vert_{op}^p \right]\right)^{\frac{1}{p}} \sqrt{ \sum_{n=0}^N \left\Vert a_n \right\Vert_{HS}^2} .            \end{equation}
La même conclusion est valide en supposant que chaque matrice aléatoire $M_n$ est à valeurs dans $\mathcal{M}_{d_n}(\C)$ et est unitairement invariante.
\end{prop}

\begin{demo}
Fixons $\mathcal{E}_{0},\dots,\mathcal{E}_N$ des matrices aléatoires indépendantes deux à deux (et indépendantes des matrices aléatoires $M_n$) qui suivent des lois uniformes dans les groupes orthogonaux de tailles respectives $d_0,\dots,d_N$.
Remarquant que $M_n$ suit la même loi que $\mathcal{E}_n M_n$ et utilisant les inégalités de Kahane-Khintchine-Marcus-Pisier \eqref{KKMP} et \eqref{KKMP-R}, on peut écrire
\bemar 
 \E\left[  \left\vert \Sum_{n=0}^N \sqrt{d_n} \tr( M_n a_n) \right\vert^p \right] &= &  
\E_\omega \E_{\omega'} \left[  \left\vert \Sum_{n=0}^N \sqrt{d_n} \tr( \mathcal{E}_n(\omega') M_n(\omega) a_n) \right\vert^p \right]   \\
& \leq  & K^p p^{\frac{p}{2}} \E_{\omega}\left[  \E_{\omega'}\left[ \left\vert \Sum_{n=0}^N \sqrt{d_n} \tr( \mathcal{E}_n(\omega') M_n(\omega) a_n) \right\vert^2 \right]^{\frac{p}{2}}\right] \\
& \leq & K^p p^{\frac{p}{2}} \E_{\omega}\left[\left( \Sum_{n=0}^N \left\Vert M_n(\omega) a_n\right\Vert_{HS}^2\right)^{\frac{p}{2}} \right] \\
& \leq & K^p p^{\frac{p}{2}} \E_{\omega}\left[\left( \Sum_{n=0}^N \left\Vert M_n(\omega)\right\Vert_{op}^2 \left\Vert a_n\right\Vert_{HS}^2\right)^{\frac{p}{2}} \right] .
\enmar 
L'inégalité triangulaire dans $L^{\frac{p}{2}}(\Omega)$ nous amène alors à la conclusion
\bemar
 \E\left[  \left\vert \Sum_{n=0}^N \sqrt{d_n} \tr( M_n a_n) \right\vert^p \right] & \leq &
  K^p p^{\frac{p}{2}} \left(\Sum_{n=0}^N  \left\Vert a_n\right\Vert_{HS}^2 \E\left[\left\Vert M_n \right\Vert_{op}^p \right]^{\frac{2}{p}} \right)^{\frac{p}{2}} \\
  & \leq & K^p p^{\frac{p}{2}} \left( \sup\limits_{0\leq n\leq N} \E\left[\left\Vert M_n\right\Vert_{op}^p\right]\right)\left(\Sum_{n=0}^N \left\Vert a_n\right\Vert_{HS}^2 \right)^{\frac{p}{2}}.
\enmar

\end{demo}
\begin{rema}
Dans le cas particulier $d_n=1$, c'est-à-dire lorsque les matrices aléatoires $M_n:\Omega\rightarrow \R$ sont en fait des \textbf{variables aléatoires symétriques}, l'inégalité \eqref{pietch} s'écrit pour tout réel $p\in [2,+\infty[$ et tous réels $a_0,\dots,a_N$ :
\bema
\E\left[ \left\vert  \sum_{n=0}^N M_n a_n \right\vert^p \right]^{\frac{1}{p}}\leq K\sqrt{p} \left( \sup\limits_{0\leq n\leq N} \E[|M_n|^p] \right)^{\frac{1}{p}}\sqrt{\sum_{n=0}^N |a_n|^2 }.
\enma
Du fait que la dimension considérée est $1$, le principe de symétrisation montrerait que l'inégalité précédente est aussi vraie si les variables aléatoires $M_n$ sont seulement centrées (voir la preuve du lemme \ref{lata-commut}).
\end{rema}

On obtient facilement le résultat suivant.

\begin{coro}\label{part1}
Considérons $p\in [1,+\infty[$ et une suite de matrices aléatoires $M_n:\Omega \rightarrow \mathcal{M}_{d_n}(\R)$ indépendantes, orthogonalement invariantes et vérifiant
\bema \sup\limits_{n\in \N} \E\left[ \left\Vert M_n\right\Vert_{op}^{\max(2,p)} \right]  <+\infty .\enma
Il existe une constante universelle $K\geq 1$ telle que, pour toute suite de matrices $b_n\in \mathcal{M}_{d_n}(L^p(X))$ et tout entier $N\in \N$, on a 
\begin{equation}\label{max2p} \begin{array}{ll} \left\Vert \Sum_{n=0}^{N} \sqrt{d_n}\tr(M_n b_n )\right\Vert_{L^{\max(2,p)}(\Omega,L^p(X))} & \\[4mm]
\quad \quad \quad \leq  K\sqrt{p} \left( \sup\limits_{n\in \N} \E\left[ \left\Vert M_n\right\Vert_{op}^{\max(2,p)} \right] \right)^{\frac{1}{\max(2,p)}}\times &\left[ \Int_{X} \left(\Sum_{n=0}^N \left\Vert b_n(x)\right\Vert_{HS}^2\right)^{\frac{p}{2}}d\mu(x) \right]^{\frac{1}{p}}.         \end{array}          \end{equation}
En particulier, sous la condition
\bema\int_{X} \left(\Sum_{n\geq 0} \left\Vert b_n(x)\right\Vert_{HS}^2\right)^{\frac{p}{2}} <+\infty, \enma
la série aléatoire $\sum \sqrt{d_n} \tr (M_n b_n)$ converge dans $L^{\max(2,p)}(\Omega,L^p(X))$ et presque sûrement dans $L^p(X)$.
De nouveau, la même conclusion est valide en supposant que chaque matrice aléatoire $M_n$ est à valeurs dans $\mathcal{M}_{d_n}(\C)$ et est unitairement invariante.
\end{coro}
\begin{demo}
Le lemme \ref{Lp-ccv} avec $q=\max(2,p)$ fournit l'inégalité 
\bema\left\Vert \Sum_{n=0}^{N} \sqrt{d_n}\tr(M_n b_n )\right\Vert_{L^{\max(2,p)}(\Omega,L^p(X))}^p \leq \int_{X}  \E\left[ \left\vert \sum_{n=0}^N \sqrt{d_n}\tr(M_n b_n(x)) \right\vert^{\max(2,p)} \right]^{\frac{p}{\max(2,p)}}    d\mu(x)    . \enma
En remplaçant $p$ par $\max(2,p)$ dans l'inégalité \eqref{pietch} et remarquant que l'on a $\sqrt{\max(2,p)}\leq \sqrt{2p}$, on obtient \eqref{max2p}.
Pour obtenir la convergence dans $L^{\max(2,p)}(\Omega,L^p(X))$, il suffit d'appliquer \eqref{max2p} aux paquets de Cauchy de la série aléatoire de terme général $\sqrt{d_n} \tr(M_n b_n)$.
Ensuite, la convergence dans $L^{\max(2,p)}(\Omega,L^p(X))$ implique la convergence en probabilité.
Puisque les termes de la série aléatoire $\sum\sqrt{d_n} \tr(M_n b_n)$ sont indépendants deux à deux et sont à valeurs dans le sous-espace séparable de $L^p(X)$ engendré par les coefficients des matrices $b_n$, on en déduit sa convergence presque sûre (voir \cite[Theorem 6.1]{ledoux} ou \cite[Théorème II.3]{queff}).
\end{demo}
\begin{rema}
On peut reformuler l'inégalité \eqref{max2p} à l'aide de \eqref{KK-LP} 
\begin{equation}\label{pietch2} \begin{array}{lr} \left\Vert \Sum_{n=0}^{N} \sqrt{d_n}\tr(M_n b_n )\right\Vert_{L^{\max(2,p)}(\Omega,L^p(X))} &\\[6mm]
\qquad \leq  C(p) \left( \sup\limits_{n\in \N} \E\left[ \left\Vert M_n\right\Vert_{op}^{\max(2,p)} \right] \right)^{\frac{1}{\max(2,p)}} & \left\Vert \Sum_{n=0}^N \sqrt{d_n} \tr (\mathcal{E}_n b_n)\right\Vert_{L^{\max(2,p)}(\Omega,L^p(X))}.         \end{array}          \end{equation}
Cette inégalité a la même forme que celle qui apparaît dans \cite[page 69]{maurey-pisier76} suite à l'application du théorème de factorisation de Pietsch.
Il doit donc y avoir une interprétation de l'inégalité \eqref{pietch2} avec une notion de cotype adaptée.
\end{rema}

Rappelons maintenant un principe de contraction pour les variables $\mathcal{E}_n$ et $W_n$ obtenu par Marcus et Pisier.
\begin{theo}\label{contrac}
Fixons un espace de Banach complexe $B$, un réel $q\in [1,+\infty[$, une suite de matrices $b_n\in \mathcal{M}_{d_n}(B)$ et une suite de matrices aléatoires $M_n:\Omega \rightarrow \mathcal{M}_{d_n}(\R)$ indépendantes et orthogonalement invariantes.
Alors on a 
\begin{equation}\label{contrac2} \left(\inf\limits_{0\leq n\leq N} \sigma\left( \E\left[|M_n|\right] \right) \right)   \E\left[ \left\Vert \sum_{n=0}^N \sqrt{d_n} \tr(\mathcal{E}_n b_n) \right\Vert_B^q    \right]^{\frac{1}{q}} \leq 
    \E\left[ \left\Vert \sum_{n=0}^N \sqrt{d_n} \tr(M_n b_n) \right\Vert_B^q    \right]^{\frac{1}{q}}    \end{equation}
Une inégalité similaire est valide en remplaçant $\mathcal{E}_n$ par $W_n$, $\mathcal{M}_{d_n}(\R)$ par $\mathcal{M}_{d_n}(\C)$ et l'invariance orthogonale par l'invariance unitaire.
\end{theo}
\begin{demo}
Ce théorème est démontré dans \cite[Page 82, Proposition 2.1]{pisier1981} mais est énoncé avec le terme \bema \left\Vert \E[|M_n|]^{-1}\right\Vert_{op}^{-1} =\sigma(\E[|M_n|]) .\enma
En ce qui concerne les variables aléatoires $\mathcal{E}_n$, le théorème est énoncé si l'espace de Banach $B$ est réel mais la preuve fonctionne de façon similaire dans le cas complexe.
\end{demo}

Le résultat précédent semble suggérer qu'une minoration uniforme de la forme 
 $\inf\limits_{n\in \N} \sigma(\E[|M_n|])>0$ et
 la convergence presque sûre de la série aléatoire $\sum\sqrt{d_n}\tr (M_n b_n)$ impliquent celle de $\sum \sqrt{d_n} \tr(\mathcal{E}_n b_n)$.
En fait, il est aisé de construire un contre-exemple dans le cas plus simple, à savoir $d_n=1$ pour tout $n\in \N$ et $B=\R$.
Pour cela, considérons une suite indépendante de variables aléatoires symétriques $X_n:\Omega\rightarrow \R$, avec $n\geq 1$, qui se concentrent en $0$ au sens suivant : il existe une constante $C>1$ telle que 
\bema  \P[X_n=0]=1-\frac{1}{n^2}, \quad \P[|X_n|\geq 1]=\frac{1}{n^2},\quad  \E[|X_n|]\in \left[\frac{1}{C},C\right] . \enma
C'est par exemple le cas si la loi de $X_n$ est 
\bema \frac{\delta_{-n^2}}{2n^2} +\left(1-\frac{1}{n^2} \right)\delta_0+\frac{\delta_{n^2}}{2n^2}   . \enma
Le lemme de Borel-Cantelli assure que pour presque tout $\omega \in \Omega$,  la suite $(X_n(\omega))$ stationne en $0$, donc la série aléatoire
 $\sum X_n$ converge presque sûrement tandis que la série aléatoire $\sum \ep_n$ diverge toujours. Cet exemple élémentaire montre qu'il faut nécessairement imposer des conditions supplémentaires sur les matrices aléatoires $M_n$ pour obtenir la convergence presque sûre de la série aléatoire $\sum \tr \sqrt{d_n} \tr(\mathcal{E}_n b_n)$ (voir par exemple \cite[Part 5]{marcusjain} et \cite[Theorem 5.2]{randomh} dans le cas unidimensionnel).
Un défaut du contre-exemple précédent est l'explosion des moments d'ordre strictement plus grand que $1$.
Cela nous mène à la proposition suivante qui nous a été inspirée par la preuve de \cite[Page 55, Lemma 1.2]{pisier1981}.
Est examinée dans cette référence la situation pour un espace de Banach de la forme $\CC^0(K)$, où $K$ est compact d'un groupe localement compact, à la place de $L^p(X)$. L'idée de la preuve se résume simplement : nous allons majorer les moments d'ordre $2$ des sommes partielles de la série $\sum \sqrt{d_n} \tr(M_n b_n)$ par leurs moments d'ordre $1$ (ce qui inverse l'ordre naturel), puis la bornitude presque sûre impliquera que les moments d'ordre $1$ et $2$ sont uniformément bornés. Par comparaison avec les variables aléatoires  $\mathcal{E}_n$, nous obtiendrons la conclusion.

\begin{prop}\label{part2}
Fixons $p\in [1,+\infty[$ et une suite de matrices aléatoires $M_n:\Omega \rightarrow \mathcal{M}_{d_n}(\R)$ indépendantes, orthogonalement invariantes et vérifiant
\bema  \sup\limits_{n\in \N} \E \left[\left\Vert M_n\right\Vert_{op}^{\max(2,p)}\right]<+\infty \quad \mbox{et} \quad \inf\limits_{n\in \N}    \sigma(\E[|M_n|])>0 .\enma
Pour toute suite de matrices $b_n\in \mathcal{M}_{d_n}(L^p(X))$, si la série aléatoire $\sum \sqrt{d_n} \tr(M_n b_n)$ est bornée en probabilité dans $L^p(X)$, c'est-à-dire si l'on a 
\bema \lim\limits_{t\rightarrow +\infty} \sup\limits_{N\in \N} \P\left[\left\Vert \sum_{n=0}^N \sqrt{d_n} \tr(M_n b_n)  \right\Vert_{L^p(X)} >t \right] =0   ,   \enma
 alors 
\bema   \sqrt{ \sum_{n\in \N} \left\Vert b_n(x)\right\Vert_{HS}^2  } \in L_x^p(X)        . \enma
De nouveau, la même conclusion est valide en supposant que chaque matrice aléatoire $M_n$ est à valeurs dans $\mathcal{M}_{d_n}(\C)$ et est unitairement invariante.
\end{prop}
\begin{demo}
Posons pour tout entier $N\in\N$ les sommes partielles 
\bema\forall (x,\omega)\in X\times \Omega \qquad  S_N(x,\omega) = \Sum_{n=0}^N \sqrt{d_n} \tr (M_n(\omega) b_n(x)) . \enma
En outre, on considère $c>1$ de sorte que l'on a pour tout $n\in \N$
\bema \frac{1}{c} \leq \sigma(\E[|M_n|])\quad \mbox{et} \quad \E \left[\left\Vert M_n\right\Vert_{op}^{\max(2,p)}\right]\leq c . \enma
Le principe de contraction (Théorème \ref{contrac}) et les inégalités \eqref{KK-LP} montrent qu'il existe une constante $C(p,c)>1$ telle que 
\begin{eqnarray}\nonumber
\frac{1}{c} \E\left[ \left\Vert \Sum_{n=0}^N \sqrt{d_n} \tr(\mathcal{E}_n b_n(x)) \right\Vert_{L_x^p(X)} \right]    & \leq & \E\left[\left\Vert S_N\right\Vert_{L^p(X)}\right] \\
\label{prin-cont}   \frac{1}{C(p,c)} \left\Vert \sqrt{\sum_{n=0}^N \left\Vert b_n(x)\right\Vert_{HS}^2} \right\Vert_{L_x^p(X)}    & \leq & \E\left[\left\Vert S_N \right\Vert_{L^p(X)}\right] . \end{eqnarray}
Majorons maintenant les moments d'ordre $2$.
En invoquant \eqref{max2p} et quitte à augmenter $C(p,c)$, on a 
\begin{eqnarray}\nonumber  \E\left[\left\Vert S_N\right\Vert_{L^p(X)}^2\right]^{\frac{1}{2}} & \leq&   \E\left[\left\Vert S_N\right\Vert_{L^p(X)}^{\max(2,p)}\right]^{\frac{1}{\max(2,p)}}    \\ \nonumber
& \leq &  C(p,c) \left\Vert\sqrt{\sum_{n=0}^N \left\Vert b_n(x)\right\Vert_{HS}^2 }\right\Vert_{L_x^p(X)} \\[3mm] \nonumber
& \leq &    C(p,c)^2 \E\left[\left\Vert S_N\right\Vert_{L^p(X)} \right]   .  \end{eqnarray}
Il s'agit de l'inégalité d'inversion des moments que nous cherchions.
Utilisons l'inégalité de Paley-Zygmund 
\cite[Page 8, Inequality II]{kahane} pour obtenir 
\bema 
\P\left[\left\Vert S_N \right\Vert_{L^p(X)} \geq \frac{1}{2} \E\left[ \left\Vert S_N\right\Vert_{L^p(X)}\right] \right]  \geq   \Frac{\E\left[\left\Vert S_N\right\Vert_{L^p(X)}\right]^2 }{4 \E\left[ \left\Vert S_N\right\Vert_{L^p(X)}^2\right]}     
 \geq    \Frac{1}{4C(p,c)^4}.   \enma
Comme $(S_N(\cdot,\omega))_{N\geq 0}$ est une suite bornée en probabilité dans $L^p(X)$, il existe un réel $A>0$ indépendant de $N$ tel que 
\bema   \forall N\in \N \qquad \P\left[\left\Vert S_N\right\Vert_{L^p(X)}\geq A \right] < \frac{1}{4C(p,c)^4}.     \enma

Les moments $\E\left[\left\Vert S_N\right\Vert_{L^p(X)}\right]$ sont donc majorés par $2A$.
La conclusion découle en examinant \eqref{prin-cont}.
\end{demo}

On obtient le résultat suivant qui s'apparente à une version multidimensionnelle du théorème de randomisation universelle de Maurey-Pisier dans le cas très particulier de l'espace $L^p(X)$.

\begin{theo}\label{LPtr}
Fixons un réel $p\in [1,+\infty[$ et considérons une suite de matrices aléatoires $M_n:\Omega \rightarrow \mathcal{M}_{d_n}(\R)$ indépendantes, orthogonalement invariantes et vérifiant
\bema  \sup\limits_{n\in \N} \E \left[\left\Vert M_n \right\Vert_{op}^{\max(2,p)}\right]<+\infty \quad \mbox{et}  \quad \inf\limits_{n\in \N} \sigma(\E[|M_n|])>0. \enma
Pour toute suite de matrices $b_n\in\mathcal{M}_{d_n}(L^p(X))$, on a l'équivalence des propriétés suivantes : 
\begin{enumerate}[i) ]
\item la fonction $x\mapsto \sqrt{\Sum_{n\geq 0} \left\Vert b_n(x)\right\Vert_{HS}^2}$ appartient à $L^p(X)$,
\item la série aléatoire $\sum \sqrt{d_n} \tr (M_n b_n)$ converge presque sûrement dans $L^p(X)$,
\item la série aléatoire $\sum \sqrt{d_n} \tr (M_n b_n)$ converge dans $L^{\max(2,p)}(\Omega,L^p(X))$,
\item la série aléatoire $\sum \sqrt{d_n} \tr (M_n b_n)$ est bornée en probabilité dans $L^p(X)$.
\end{enumerate}
La même conclusion est valide en supposant que chaque matrice aléatoire $M_n:\Omega\rightarrow \mathcal{M}_{d_n}(\C)$ est unitairement invariante.
\end{theo}
\begin{demo}
En combinant le corollaire \ref{part1} et la proposition \ref{part2}, on vérifie les implications iv) $\Rightarrow$ i) et  i) $\Rightarrow$ iii).
L'implication iii) $\Rightarrow$ ii) a été vue au cours de la preuve du corollaire \ref{part1}.
Quant à l'implication ii) $\Rightarrow$ iv), elle est vraie en toute généralité.
\end{demo}

Nous pouvons maintenant démontrer très facilement le théorème \ref{mplp} à partir du théorème \ref{LPtr}. Pour tout entier $n\in \N$, on écrit 
\bema
\Sum_{i=1}^{d_n} \sum_{j=1}^{d_n}M_{n,i,j} \langle u_n,\phi_{n,j}\rangle \phi_{n,i} \\
 = \sqrt{d_n} \tr \left(M_n b_n \right),
\enma
avec 
\bema  \forall x\in X \quad    b_n(x)=\frac{1}{\sqrt{d_n}}\left[\langle u_n,\phi_{n,i} \rangle \phi_{n,j}(x)  \right]_{i,j}   \in \mathcal{M}_{d_n}(L_x^p(X)) . \enma
Grâce à \eqref{defi-en}, on a enfin
\bema \left\Vert b_n(x)\right\Vert_{HS}^2= \sum_{i,j=1}^{d_n} \frac{1}{d_n} |\langle u_n,\phi_{n,i} \rangle \phi_{n,j}(x)|^2  =\frac{1}{d_n} \left\Vert u_n \right\Vert_{L^2(X)}^2 \sum_{j=1}^{d_n} |\phi_{n,j}(x)|^2=\left\Vert u_n \right\Vert_{L^2(X)}^2\frac{e(n,x)}{d_n}  .  \enma

\begin{rema}\label{Lptr2}
D'après ce qui précède, les implications i) $\Rightarrow$ ii) et  i $\Rightarrow$ iii) des théorèmes \ref{LPtr} et \ref{mplp} n'utilisent que l'hypothèse de bornitude des moments $\sup\limits_{n\in \N} \E \left[\left\Vert M_n \right\Vert_{op}^{\max(2,p)}\right]<+\infty$.
En particulier, cette hypothèse et l'assertion 
\bema
\sqrt{\sum_{n\in \N} \left\Vert u_n \right\Vert_{L^2(X)}^2 \frac{e(n,x)}{d_n}} \in L_x^p(X)
\enma
impliquent la convergence presque sûre dans $L^p(X)$ et la convergence dans $L^{\max(2,p)}(\Omega,L^p(X))$ de la série aléatoire 
\bema
\sum_{n\in \N} \sum_{i=1}^{d_n} \left( \sum_{j=1}^{d_n} M_{n,i,j} \langle u_n,\phi_{n,j}\rangle \right)\phi_{n,i}.
\enma
\end{rema}

\section{Preuve de la proposition \ref{ex-mat}, inégalité de Latala précisée}

Notons $B_n=[\ep_{ij}]$ la matrice aléatoire de taille $n\times n$ et dont les coefficients $\ep_{ij}$ sont des variables aléatoires i.i.d. qui suivent une loi $\frac{1}{2}$-Bernoulli à valeurs dans $\{-1,+1\}$.
Les inégalités de Kahane-Khintchine \eqref{KK} dans l'espace de Banach $(\mathcal{M}_n(\R),\left\Vert \cdot\right\Vert_{op})$ nous donnent l'encadrement : 
\bema \forall q\in [1,+\infty[ \quad \forall n\in \N^\star \quad     \E\left[ \left\Vert B_n \right\Vert_{op}\right] \leq    \E\left[\left\Vert B_n \right\Vert_{op}^q\right]^{\frac{1}{q}} \leq K_{q,1}       \E\left[\left\Vert B_n \right\Vert_{op}\right].     \enma
Cela signifie que tous les moments $\E\left[\left\Vert  B_n \right\Vert_{op}^q\right]^{\frac{1}{q}}$ ont le même ordre de grandeur lorsque $n$ tend vers $+\infty$. 
Par ailleurs, la théorie des matrices aléatoires explique que le moment 
$\E\left[\left\Vert  B_n \right\Vert_{op}\right]$ est asymptotiquement de l'ordre de $\sqrt{n}$ (voir \cite[Part 2.3]{tao-rand-matrix}).
Nous allons exploiter cette idée pour démontrer l'inégalité \eqref{super-lata}.
Commençons par le lemme élémentaire suivant qui s'apparente à une version commutative de \eqref{super-lata}.

\begin{lemm}\label{lata-commut}
Considérons un réel $p\in [2,+\infty[$ ainsi que $N$ variables aléatoires $U_1,\dots,U_N$ réelles, centrées, i.i.d. et ayant un moment d'ordre $p$. 
Nous avons l'inégalité
\bema  \E\left[\left\vert \frac{U_1+\dots+U_N}{\sqrt{N}}\right\vert^p \right]\leq \left(C\sqrt{p}\right)^p \E[|U_{1}|^p].   \enma
\end{lemm}
\begin{demo}
L'idée se résume en deux points : on se ramène au cas où les variables $U_i$ sont symétriques et l'on invoque les inégalités de Kahane-Khintchine à l'aide du théorème de Fubini.
Si nous notons $\widetilde{U}_1,\dots,\widetilde{U}_N$ des copies indépendantes des variables $U_1,\dots,U_N$ alors l'inégalité de Jensen pour l'espérance en les variables $\widetilde{U}_1,\dots,\widetilde{U}_N$ donne 
\bema \E\left[\left\vert U_1+\dots+U_N\right\vert^p \right] \leq  
 \E\left[\left\vert U_1-\widetilde{U}_1+\dots+U_N-\widetilde{U}_N\right\vert^p \right]. \enma
Or les variables $U_{i}-\widetilde{U}_i$ et $\ep_i |U_{i}-\widetilde{U}_i|$ ont la même loi. 
Il s'agit maintenant d'utiliser les inégalités de Kahane-Khintchine \eqref{KK} et \eqref{KKMP-R} (avec $d_n=1$) pour obtenir 
\begin{eqnarray} \nonumber
 \E\left[\left\vert U_1+\dots+U_N\right\vert^p \right]& \leq &\E_{\omega} \E_{\omega'} \left[\left\vert \ep_1(\omega')|U_1(\omega)-\widetilde{U}_1(\omega)|+ \dots +\ep_N(\omega') |U_N(\omega)-\widetilde{U}_N(\omega)|\right\vert^p \right] \\ \label{typ1}
& \leq & K_{p,2}^p \E_{\omega} \left[ \E_{\omega'} \left[\left\vert \ep_1(\omega') |U_1(\omega)-\widetilde{U}_1(\omega)|+ \dots +\ep_{N}(\omega') |U_N(\omega)-\widetilde{U}_N(\omega)|\right\vert^2 \right]^{\frac{p}{2}} \right] \\ \label{typ2}
&\leq & \left(C\sqrt{p}\right)^p  \E_{\omega} \left[ \left(\Sum_{i=1}^N |U_i(\omega)-\widetilde{U}_i(\omega)|^2\right)^{\frac{p}{2}} \right] \\ \nonumber
& \leq & \left(C\sqrt{p}\right)^p  \left( \Sum_{i=1}^N \left\Vert (U_i-\widetilde{U}_i)^2\right\Vert_{L^\frac{p}{2}(\Omega)} \right)^{\frac{p}{2}}= \left(C\sqrt{pN}\right)^p \left\Vert (U_1-\widetilde{U}_1)^2\right\Vert_{L^{\frac{p}{2}(\Omega)}}^{\frac{p}{2}} \\ \nonumber 
& \leq & \left(C\sqrt{pN}\right)^{p} \left\Vert U_1\right\Vert_{L^p(\Omega)}^p .\end{eqnarray}
\end{demo}

La preuve précédente est similaire à celle de la proposition \ref{pietch-lem} avec $(d_n,a_n,M_n)=(1,1,U_n)$ à ceci près que nous pouvons utiliser en plus le principe de symétrisation.
Nous avons écrit la preuve précédente d'abord parce que nous aurons besoin plus loin de considérer des variables seulement centrées au lieu de symétriques, mais aussi pour des raisons pédagogiques. En effet, la disparition élémentaire de l'espérance $\E_{\omega'}$ de la ligne \eqref{typ1} à la ligne \eqref{typ2} peut être interprétée comme la propriété de type $2$ de l'espace de Banach $\R$ (voir la définition \ref{def-typ}).
Nous allons utiliser un substitut non-commutatif de cette propriété pour démontrer \eqref{super-lata}.
C'est l'objet de la proposition suivante dont la preuve est très technique.

\begin{prop}\cite[Theorem 1]{latala}
Il existe une constante universelle $C>1$ telle que pour tout entier $n\in \N^\star$, si l'on considère $n^2$ variables aléatoires i.i.d. $(g_{ij})_{1\leq i,j\leq n}$ qui suivent une loi normale $\mathcal{N}_\R(0,1)$ et une matrice $[a_{ij}]\in \mathcal{M}_{n}(\R)$ alors 
\begin{equation}\label{latala-gauss}
\E\left[\left\Vert  a_{ij} g_{ij}\right\Vert_{op} \right]\leq C \left(\sqrt[4]{\sum_{i,j=1}^n a_{ij}^4}+\max_{1\leq i\leq n} \sqrt{\sum_{j=1}^n a_{ij}^2}+\max_{1\leq j\leq n} \sqrt{\sum_{i=1}^n a_{ij}^2} \right).
\end{equation}
\end{prop}
Latala démontre la proposition précédente pour en déduire l'estimation \eqref{lata} qui se reformule
\bema   \E\left[\left\Vert X_{ij} \right\Vert_{op}\right] \leq C\sqrt{n} \E\left[|X_{1,1}|^4\right]^{\frac{1}{4}}.      \enma
On peut attaquer la preuve de \eqref{super-lata} et l'on commence comme dans \cite{latala}.
On note $\widetilde{X}_{ij}$ des copies indépendantes des variables aléatoires $X_{ij}$.
En particulier, on a $\E[\widetilde{X}_{ij}]=0$.
Pour tout réel $p\in [4,+\infty[$, l'argument classique de l'inégalité de Jensen en les variables $\widetilde{X}_{ij}$ donne 
\begin{equation}\label{symetriz}  \E\left[ \left\Vert X_{ij}\right\Vert_{op}^p \right] \leq   
 \E\left[\left\Vert X_{ij}-\widetilde{X}_{ij}\right\Vert_{op}^p \right] =\E\left[\left\Vert \ep_{ij} (X_{ij}-\widetilde{X}_{ij})\right\Vert_{op}^p \right] \leq 2^p \E\left[ \left\Vert \ep_{ij} X_{ij}\right\Vert_{op}^p \right] ,    \end{equation}
 où les $n^2$ variables aléatoires $\ep_{ij}$ sont indépendantes deux à deux entre elles et vis-à-vis des variables $X_{ij}$ et $\widetilde{X}_{ij}$.
On diffère maintenant de \cite{latala} en faisant appel aux inégalités de Kahane-Khintchine \eqref{KK} dans l'espace de Banach $(\mathcal{M}_n(\R),\left\Vert\cdot\right\Vert_{op})$ afin de récupérer le moment d'ordre $p$ :
\begin{eqnarray}   \E\left[ \left\Vert \ep_{ij} X_{ij}\right\Vert_{op}^p \right]&=& \E_{\omega} \E_{\omega'} \left[ \left\Vert \ep_{ij}(\omega') X_{ij}(\omega)\right\Vert_{op}^p \right] \nonumber \\ \nonumber
& \leq & K_{p,1}^p \E_{\omega} \left[  \E_{\omega'}\left[ \left\Vert \ep_{ij}(\omega') X_{ij}(\omega)\right\Vert_{op} \right]^p \right]\\ \label{fubin-ep} & \leq &  C(p)\E_{\omega} \left[  \E_{\omega'}\left[ \left\Vert \ep_{ij}(\omega') X_{ij}(\omega)\right\Vert_{op} \right]^p \right]        .   \end{eqnarray}
Et l'on reprend de nouveau l'argumentation de Latala. 
En supposant que toutes les variables aléatoires sont indépendantes deux à deux, \eqref{latala-gauss} et le principe de contraction (voir \eqref{contrac2}) donnent pour tout $\omega \in \Omega$
 \bemar
\E_{\omega'} \left[ \left\Vert \ep_{ij}(\omega') X_{ij}(\omega)\right\Vert_{op} \right] & \leq &  \Frac{1}{\E[|g_{11}|]}\E_{\omega'}\left[\left\Vert g_{ij}(\omega') X_{ij}(\omega)\right\Vert_{op} \right]=\Frac{\sqrt{\pi}}{\sqrt{2}}\E_{\omega'}\left[\left\Vert g_{ij}(\omega') X_{ij}(\omega)\right\Vert_{op} \right] \\
& \leq & C\left(\sqrt[4]{\Sum_{i,j=1}^n X_{ij}(\omega)^4}  +\max\limits_{1\leq i\leq n} \sqrt{\Sum_{j=1}^n X_{ij}(\omega)^2}+\max\limits_{1\leq j\leq n} \sqrt{\Sum_{i=1}^n X_{ij}(\omega)^2}  \right) .\enmar
Reprenant  \eqref{symetriz}, \eqref{fubin-ep} et tenant compte que les variables $X_{ij}$ sont i.i.d., on obtient

\bemar  \E\left[\left\Vert X_{ij}\right\Vert_{op}^p\right] & \leq & C(p)\E\left[\left(\Sum_{i,j=1}^n X_{ij}^4\right)^{\frac{p}{4}}+ \max\limits_{1\leq i\leq n} \left(\Sum_{j=1}^n X_{ij}^2\right)^{\frac{p}{2}}+\max\limits_{1\leq j\leq n} \left(\Sum_{i=1}^n X_{ij}^2\right)^{\frac{p}{2}} \right]             \\
& \leq & C(p)\E\left[ \left(\Sum_{i,j=1}^n X_{ij}^4\right)^{\frac{p}{4}}+\max\limits_{1\leq i\leq n} \left(\Sum_{j=1}^n X_{ij}^2\right)^{\frac{p}{2}} \right].\enmar
La fin de cette preuve est différente de \cite{latala} car on doit tenir compte de l'inégalité $p\geq 4$.
Le premier terme ne posera aucun problème tandis que le dernier est plus délicat, c'est pour cela que nous le forçons à être centré en majorant 
\bema \sum_{j=1}^n X_{ij}^2 \leq \sum_{j=1}^n \E\left[X_{ij}^2\right]+\left\vert \sum_{j=1}^n X_{ij}^2-\E\left[X_{ij}^2\right] \right\vert.      \enma
Cela nous donne 
\bemar 
 \E\left[\left\Vert X_{ij}\right\Vert_{op}^p\right] & \leq & C(p) \left[A(1)+A(2)+A(3)\right], \\
A(1) &:= & \E\left[\left(\Sum_{i,j=1}^n X_{ij}^4\right)^{\frac{p}{4}} \right] \\
A(2) &:= & \max\limits_{1\leq i\leq n} \left(\Sum_{j=1}^n \E\left[X_{ij}^2\right]\right)^{\frac{p}{2}}=n^{\frac{p}{2}} \E\left[X_{11}^2\right]^{\frac{p}{2}} \\
A(3) & :=& \E\left[ \max\limits_{1\leq i\leq n} \left\vert\Sum_{j=1}^n X_{ij}^2-\E\left[X_{ij}^2\right] \right\vert^{\frac{p}{2}}\right]      .
\enmar 
L'inégalité de Hölder donne d'abord $A(2)\leq n^{\frac{p}{2}} \E[|X_{11}|^p] $, puis 
\bema  A(1)\leq \E\left[(n^2)^{\frac{\frac{p}{4}-1}{\frac{p}{4}}\times \frac{p}{4}}\sum_{i,j=1}^n |X_{ij}|^p \right]  = n^{\frac{p}{2}}\E[|X_{11}|^p]. \enma
Concernant le terme $A(3)$, nous majorons grossièrement à l'aide du lemme \ref{lata-commut} et de l'inégalité $p\geq 4$

\bemar
A(3) & \leq& n \E\left[\left\vert\Sum_{j=1}^n X_{1j}^2-\E[X_{1j}^2] \right\vert^{\frac{p}{2}} \right] \\
& \leq & C(p) n^{1+\frac{p}{4}} \E\left[ \left\vert X_{11}^2 -\E\left[X_{11}^2\right]\right\vert^{\frac{p}{2}}\right] \\
& \leq & C(p) n^{1+\frac{p}{4}} \E[|X_{11}|^p] \\
& \leq & C(p) n^{\frac{p}{2}} \E[|X_{11}|^p] .\enmar

\begin{rema}\label{lata-optim}
Discutons de l'optimalité de l'inégalité \eqref{super-lata}.
La théorie des matrices aléatoires nous informe que l'on ne peut pas descendre en-dessous du moment d'ordre $4$ pour espérer contrôler $\E[\left\Vert M_n\right\Vert_{op}]$ et a fortiori les autres moments (\cite{bai88b} et \cite[Theorem 3.1]{bai88a}).
Cela se voit facilement avec une variable aléatoire $X_{11}$ qui vérifie 
\bema
\forall t\geq \sqrt[4]{e}\quad \P [|X_{11}|\geq t]=\frac{\ln(t)}{t^4}.\enma
D'une part, on a $\E[|X_{11}|^4]=+\infty$ et $\E[|X_{11}|^p]<\infty$ pour tout $p\in [1,4[$.
D'autre part, on peut minorer 
\bemar
\E\left[\left\Vert M_n\right\Vert_{op}\right] & \geq & \frac{1}{\sqrt{n}}\E\left[\sup\limits_{i,j} |X_{ij}| \right] \\
& \geq &  \Int_{0}^{+\infty} \P\left[ \sup\limits_{i,j} |X_{ij}|\geq t \sqrt{n} \right] dt  \\
& \geq & \Int_{0}^{+\infty} 1- \P\left[ |X_{11}| < t\sqrt{n}\right]^{n^2} dt   \\
& \geq & \Int_{\sqrt[4]{e}}^{+\infty} 1- \left( 1-\frac{\ln(t\sqrt{n} )}{t^4 n^2}\right)^{n^2} dt .
\enmar
Le lemme de Fatou implique que l'on a 
\bema \liminf\limits_{n\rightarrow +\infty} \E\left[\left\Vert M_n\right\Vert_{op}\right]\geq \int_{\sqrt[4]{e}}^{+\infty} 1 dt=+\infty. \enma
\end{rema}
\section{Preuve de la proposition \ref{ex-mat}, minoration de la plus petite valeur singulière}

On commence par un lemme dual au lemme \ref{lata-commut}. 

\begin{lemm}\label{lata-mino}
Considérons $N$ variables aléatoires $U_1,\dots,U_N$ réelles, centrées, i.i.d. et ayant un moment d'ordre $1$. 
Pour tout $(y_1,\dots,y_N)\in \R^N$, nous avons l'inégalité
\bema    \E[|U_1|] \sqrt{y_1^2+\dots+y_N^2} \leq C \E\left[ |y_1 U_1+\dots+y_N U_N| \right] .     \enma
\end{lemm}
\begin{demo}
Soient $\widetilde{U}_1,\dots,\widetilde{U}_N$ des copies indépendantes de $U_1,\dots,U_N$. On a 
\bema    \E\left[ |y_1 (U_1-\widetilde{U}_1)+\dots+y_N (U_N-\widetilde{U}_N)| \right]  \leq 2 \E\left[ |y_1 U_1+\dots+y_N U_N| \right].  \enma
Puisque les variables $U_i-\widetilde{U}_i$ sont symétriques, le principe de contraction (Théorème \ref{contrac}), l'inégalité de Khintchine et l'inégalité de Jensen donnent
\bemar
\E[|U_1|] \sqrt{y_1^2+\dots+y_N^2} & = & \E[|U_1|]\E\left[|\ep_1 y_1+\dots+\ep_n y_N|^2 \right]^{\frac{1}{2}} \\[2mm]
& \leq &  \E[|U_1-\widetilde{U}_1|] \times K_{2,1} \E\left[|\ep_1 y_1+\dots+\ep_n y_N| \right] \\
& \leq &  K_{2,1} \E\left[ |y_1 (U_1-\widetilde{U}_1)+\dots+y_N (U_N-\widetilde{U}_N)| \right] \\
& \leq & 2 K_{2,1} \E\left[ |y_1 U_1+\dots+y_N U_N| \right].
\enmar
\end{demo}

Passons à la preuve de \eqref{nucl}. Fixons $y=(y_1,\dots,y_n)\in \R^n$. Pour tout $\omega \in \Omega$, la diagonalisation de la matrice symétrique positive $|M_n(\omega)|=\sqrt{\tt M_n(\omega) M_n(\omega)}$ dans une base orthonormée fournit l'inégalité
\bema   \left\vert M_n(\omega) y\right\vert^2 :=\langle M_n(\omega)y,M_n(\omega)y \rangle =\tt y |M_n(\omega)|^2 y \leq \left\Vert M_n(\omega)\right\Vert_{op} \tt y |M_n(\omega)| y   . \enma
L'inégalité de Cauchy-Schwarz donne alors 
\bemar   \E_{\omega} [\left\vert M_n(\omega) y\right\vert] & \leq&   \E_\omega\left[\sqrt{ \left\Vert M_n(\omega)\right\Vert_{op}}\sqrt{ \tt y |M_n(\omega)| y }\right]   \\ 
 \E_{\omega} [\left\vert M_n(\omega) y\right\vert]^2 & \leq & \E_\omega\left[\left\Vert M_n(\omega)\right\Vert_{op}\right] \E_{\omega}\left[\tt y |M_n(\omega)| y\right] . \enmar
On invoque alors l'inégalité de Latala \eqref{lata} pour contrôler le moment d'ordre $1$ de $\left\Vert M_n\right\Vert_{op}$ :
\bema       \E_{\omega} [\left\vert M_n(\omega) y\right\vert]^2 \leq  C \E[|X_{11}|^4]^{\frac{1}{4}} \times  \tt y \E_{\omega}[|M_n(\omega)| ]y.     \enma
On va utiliser l'égalité
\bema
|M_n(\omega)y| = \Frac{1}{\sqrt{n}} \left(\Sum_{i=1}^n\left|\sum_{j=1}^n X_{ij}(\omega) y_j\right\vert^2\right)^{\frac{1}{2}}
\enma 
à l'aide de l'inégalité triangulaire entre $\E_\omega$ et la norme euclidienne $|\cdot|$ de $\R^n$ ainsi que le lemme \ref{lata-mino} :
\bemar
\E[|X_{11}|]^2(y_1^2+\dots+y_n^2) & \leq & \Frac{C}{n}\Sum_{i=1}^n \E_{\omega}\left[\left\vert \sum_{j=1}^n X_{ij}(\omega) y_j \right\vert \right]^2 \\
& \leq & C \E_{\omega}[|M_n(\omega)y|]^2.
\enmar 
On obtient alors \eqref{nucl} :
\bema \frac{\E[|X_{11}|]^2}{C\E[|X_{11}|^4]^{\frac{1}{4}}}(y_1^2+\dots+y_n^2)  \leq \tt y \E[|M_n|] y  .   \enma

\section{Preuves des théorèmes \ref{theo-interpo} et \ref{theo-dual}, partie I : rétracte d'un espace de Banach}\label{partI}

Afin d'alléger la rédaction, on conviendra que les éléments de $E_n\subset L^2(X)$ sont des fonctions de la variable $y$ et l'on préférera les deux écritures
\bema
L^2(X)\rightarrow L_y^2(X)\quad \mbox{et} \quad \pl^{p}(X,\oplus E_n)\rightarrow \pl_{y}^p(X,\oplus E_n).
\enma
Il sera aussi commode de définir l'espace de Hilbert abstrait $\oplus E_n$, c'est-à-dire que l'on pose
\bema
\forall (u_n)_{n\in \N}\in \bigoplus_{n\in \N} E_n \quad \left\Vert (u_n)_{n\in \N} \right\Vert_{\oplus E_n}:=\left(\int_{X}\sum_{n\in \N}  |u_n(y)|^2 d\mu(y) \right)^{\frac{1}{2}}.
\enma
On conseille au lecteur d'imaginer que les sous-espaces $E_n$ sont en somme directe orthogonale dans $L_y^2(X)$, si bien que l'espace abstrait $\oplus E_n$ s'identifie à un sous-espace fermé de $L_y^2(X)$. Mais il nous arrivera de considérer la situation générale dans certaines démonstrations.

Concernant l'espace de Bochner-Lebesgue, on conservera la lettre $x\in X$ pour écrire
$L_x^p(X,\oplus E_n)$ au lieu de $L^p(X,\oplus E_n)$.
La définition \ref{defi-cv} dit exactement que l'opérateur linéaire suivant est isométrique
\begin{equation}\label{defi-S} \begin{array}{rcl}
S_p:  \pl_{y}^p(X,\oplus E_n) & \rightarrow & L_x^p(X,\oplus E_n) \\
 (u_n(y))_{n\in \N} & \mapsto  &  \left( \frac{\sqrt{e(n,x)}}{\sqrt{d_n}} u_n(y)\right)_{n\in \N} .\end{array}        \end{equation}

Ainsi, $\pl_{y}^p(X,\oplus E_n)$ s'identifie à un sous-espace fermé de $L_x^p(X,\oplus E_n)$ et donc est complet.
Commençons par les points faciles.
\vspace*{0.5 cm}
\begin{center}
\textbf{Dualité}
\end{center}
\vspace*{0.5 cm}
Rappelons que l'on note $q=\frac{p}{p-1}$ l'exposant conjugué de $p$.
On montre facilement l'inégalité \eqref{dual-triv} pour tout $(u,w)\in \pl_{y}^p(X,\oplus E_n)\times \pl_{y}^q(X,\oplus E_n)$ :
\bemar
\Sum_{n\geq 0} |\langle u_n,w_n\rangle_{L^2(X)}| & = & \Int_{X} \sum_{n\in \N} \frac{e(n,x)}{d_n}|\langle u_n,w_n\rangle_{L^2(X)}|  d\mu(x)     \\
  & \leq &  \Int_X \sqrt{ \sum_{n\geq 0} \frac{e(n,x)}{d_n} \left\Vert u_n(y)\right\Vert_{L_y^2(X)}^2 }
\sqrt{ \sum_{n\geq 0} \frac{e(n,x)}{d_n} \left\Vert w_n(y)\right\Vert_{L_y^2(X)}^2 } d\mu(x) \\
  & \leq & \left\Vert u\right\Vert_{\pl^p(X,\oplus E_n)}\left\Vert w\right\Vert_{\pl^q(X,\oplus E_n)}.
\enmar
L'injection canonique $\Lambda_p:\pl_{y}^{q}(X,\oplus E_n)\rightarrow \pl_{y}^p(X,\oplus E_n)'$, définie dans l'énoncé du théorème \ref{theo-dual}, est donc continue.

\vspace*{0.5 cm}
\begin{center}
\textbf{Interpolation}
\end{center}
\vspace*{0.5 cm}
Rappelons que les espaces $\pl_{y}^{p}(X \oplus E_n)$ et leurs interpolés complexes et réels peuvent être vus comme des sous-espaces de $\prod_{n\in \N} E_n$.
En outre, on a le résultat suivant \cite[Part 1.18.4]{triebel}.
\begin{theo}\label{interpo-compl}
Considérons deux espaces de Banach complexes $B_1$ et $B_2$ ainsi que des réels $p_1\leq p\leq p_2$ appartenant à $[1,+\infty[$. En notant $\theta\in [0,1]$ un réel qui vérifie $\frac{1}{p}=\frac{1-\theta}{p_1}+\frac{\theta}{p_2}$, on a 
\bema
[L^{p_1}(X,B_1),L^{p_2}(X,B_2)]_{\theta}=L^{p}(X, [B_1,B_2]_{\theta}).
\enma
Si un opérateur linéaire $T$ borné de $L^{p_1}(X,B_1)$ dans lui-même et de $L^{p_2}(X,B_2)$ dans lui-même alors il est aussi borné de $L^p(X,[B_1,B_2]_{\theta})$ dans lui-même.
Le même énoncé est valide en remplaçant la méthode d'interpolation complexe $[\cdot,\cdot]_{\theta}$ par la méthode d'interpolation réelle $[\cdot,\cdot]_{\theta,p}$.
\end{theo}
Avec les notations du théorème précédent, on peut interpoler l'application \eqref{defi-S} et assurer que l'opérateur
\bemar
[\pl_{y}^{p_1}(X \oplus E_n),\pl_{y}^{p_2}(X,\oplus E_n)]_\theta  & \rightarrow & L_{x}^p(X,\oplus E_n) \\
(u_n(y))_{n\in \N} &\mapsto & \left(\frac{\sqrt{e(n,x)}}{\sqrt{d_n}}u_n(y) \right)_{n\in \N}
\enmar
est borné. Cela implique l'inclusion continue
\begin{equation}\label{inclu-inter}
[\pl_{y}^{p_1}(X \oplus E_n),\pl_{y}^{p_2}(X,\oplus E_n)]_\theta\subset \pl_{y}^p(X,\oplus E_n).
\end{equation}

\vspace*{0.5 cm}
\begin{center}
\textbf{Vers la notion de rétracte d'un espace de Banach}
\end{center}
\vspace*{0.5 cm}

Cependant, il ne paraît pas évident d'aller plus loin dans les deux analyses précédentes, c'est-à-dire de prouver que $\Lambda_p:\pl_{y}^{q}(X,\oplus E_n)\rightarrow \pl_{y}^p(X,\oplus E_n)'$ est surjective et que l'inclusion \eqref{inclu-inter} est une égalité.
La notion de rétracte d'un espace de Banach permet de reformuler la question.

\begin{defi}\label{defi-retra}
Soient $A$ et $B$ deux $\C$-espaces vectoriels normés, on dit que $A$ est un rétracte de $B$ s'il existe deux applications linéaires bornées 
$ S:A\rightarrow B $ et $R : B\rightarrow A     $ telles que $RS=\mbox{id}_{A}$.
\end{defi}
Dans la définition précédente, il faut imaginer que $B$ est un espace de référence qui est bien compris et $A$ un sous-espace que l'on veut analyser.
Le prototype de rétracte qu'il faut avoir à l'esprit est le cas où $A$ est un sous-espace complémenté de $B$, c'est-à-dire image d'un projecteur borné, l'application $S$ est alors l'application identité et $R$ un projecteur de $B$ sur $A$.
En effet, on vérifie facilement le résultat suivant.

\begin{prop}\label{prop-retra}
Avec les mêmes notations que dans la définition \ref{defi-retra}, on a 
\begin{enumerate}[i) ]
\item $SR:B\rightarrow B$ est un projecteur borné, son espace image est égal à $S(A)$ et est fermé.

\item $S$ est un isomorphisme d'espaces vectoriels normés de $A$ sur $S(A)\subset B$ : 
\bema
\forall a \in A \qquad    \frac{1}{\left\Vert R\right\Vert_{B\rightarrow A}}  \left\Vert a\right\Vert_A \leq \left\Vert S(a)\right\Vert_{B}\leq \left\Vert S \right\Vert_{A\rightarrow B} \left\Vert a \right\Vert_A.
\enma
\item si $B$ est complet alors $A$ l'est aussi.
\end{enumerate}
\end{prop}

En ce qui concerne la dualité, on a le corollaire facile.

\begin{coro}\label{coro-retra}
Avec les mêmes notations que dans la définition \ref{defi-retra}, si l'on note les applications duales (ou applications transposées) de $S$ et $R$ : 
\bema \begin{array}{rclcrclcl}
\tt S: B' & \rightarrow & A',&  & \tt R:A'& \rightarrow & B', &\\  
\phi & \mapsto & \phi\circ S & & \psi & \mapsto & \psi \circ R  
\end{array}     \enma
Alors on a $\tt S \tt R=\mbox{id}_{A'}$. Cela signifie que l'espace dual $A'$ est un rétracte de $B'$ et en particulier que $\tt R$ est un isomorphisme $A'$ sur l'image du projecteur $\tt R \tt S=\tt (SR):B'\rightarrow  B'$.
\end{coro}

En ce qui concerne l'interpolation réelle ou complexe, la réponse est donnée par le résultat suivant \cite[Page 22, Theorem 1.2.4]{triebel}. 

\begin{coro}\label{trieb-interpo}
On note $[\cdot,\cdot]$ une méthode d'interpolation.
Soient $(A_1,A_2)$ et $(B_1,B_2)$ deux couples d'interpolation d'espaces de Banach, on suppose qu'il existe un opérateur linéaire borné $S:A_1\rightarrow B_1$ 
et $S:A_2\rightarrow B_2$ et un opérateur linéaire borné $R:B_1\rightarrow A_1$ 
et $R:B_2\rightarrow A_2$ qui satisfont : 
\bema 
\forall a\in A_1 +A_2 \qquad  RS(a)=a.
\enma 
Alors $SR$ est un projecteur borné de l'espace interpolé $[B_1,B_2]$ et $S$ réalise un isomorphisme d'espaces de Banach de l'espace interpolé $[A_1,A_2]$ sur $SR([B_1,B_2])$.
\end{coro}
\begin{demo}
C'est immédiat puisque l'espace de Banach $[A_1,A_2]$ est un rétracte de $[B_1,B_2]$ par l'intermédiaire des opérateurs bornés $S:[A_1,A_2]\rightarrow [B_1,B_2]$ et $R:[B_1,B_2]\rightarrow [A_1,A_2]$.\end{demo}

Pour comprendre les espaces duaux et interpolés des espaces $\pl_{y}^p(X,\oplus E_n)$, il suffit d'étudier si $\pl_{y}^p(X,\oplus E_n)$ est un rétracte de l'espace de Banach $L_x^p(X,\oplus E_n)$ par le biais de l'application $S_p$ définie en \eqref{defi-S} lorsque $p$ parcourt $]p_1,p_2[$ pour les théorèmes \ref{theo-interpo} et \ref{theo-dual}. En d'autres termes, on cherche un opérateur borné 
\bema R_p:L_x^p(X, \oplus E_n) \rightarrow \pl_{y}^p(X,\oplus E_n)    \enma
tel que $R_p S_p= \mbox{id}_{\pl_{y}^p(X,\oplus E_n)}$.

\vspace*{0.5 cm}
\begin{center}
\textbf{Stratégie pour le théorème \ref{theo-dual} de dualité}
\end{center}
\vspace*{0.5 cm}

Puisque les fonctions $\frac{1}{d_n}e(n,\cdot)$ sont des densités de probabilité sur $X$, nous avons un candidat très naturel pour l'opérateur $R_p$. En effet, posons : 
\begin{equation}\label{defi-RRp}\begin{array}{rcl}  R_p: L_x^p(X,\oplus E_n) & \rightarrow & \pl_{y}^p(X,\oplus E_n) \\
\left( u_n(x,y)\right)_{n\geq 0} & \mapsto &     \left(\Int_X u_n(x',y) \frac{\sqrt{e(n,x')}}{\sqrt{d_n}} d\mu(x')  \right)_{n\geq 0}   .\end{array}\end{equation}
D'après la définition \eqref{defi-S}, on a bien formellement $R_p S_p=\mbox{id}_{\pl_{y}^p(X,\oplus E_n)}$.
Malheureusement, nous ne voyons aucune raison triviale assurant que $R_p$ arrive bien dans $\pl_{y}^p(X,\oplus E_n)$!
Puisque $S_p$ est une isométrie, la bornitude de $R_p$ équivaut à la bornitude du projecteur

\begin{equation}\label{SR}\begin{array}{rcl} S_p R_p : L_x^p(X,\oplus E_n) & \rightarrow &  L_x^p(X,\oplus E_n)   \\[3mm]
(u_n(x,y))_{n\in \N} & \mapsto &  \left( \Frac{\sqrt{e(n,x)}}{\sqrt{d_n}}\Int_X u_n(x',y) \frac{\sqrt{e(n,x')}}{\sqrt{d_n}}d\mu(x')\right)_{n\in\N}    .   \end{array} \end{equation}

Rappelons que sous les hypothèses du théorème \ref{theo-dual} de dualité, $p_1$ et $p_2$ sont supposés être deux exposants conjugués.
La bornitude du projecteur $S_p R_p$ sera prouvée dans la partie \ref{partIV} pour tout $p\in ]p_1,p_2[$.
La proposition suivante achèvera la preuve
du théorème \ref{theo-dual} de dualité
et explique pourquoi la notion de rétracte est la bonne notion pour aborder la dualité des espaces $\pl^{p}(X,\oplus E_n)$.

\begin{prop}\label{Lamba-surj}
Fixons $p\in ]1,+\infty[$ et posons $q=\frac{p}{p-1}$ l'exposant conjugué.
Pour tout $n\in \N$, on suppose que $\sqrt{e(n,\cdot)}$ appartient à $L^p(X)\cap L^q(X)$. Alors les propriétés suivantes sont équivalentes : 
\begin{enumerate}[i) ]
\item l'injection canonique $\Lambda_p:\pl_{y}^q(X,\oplus E_n)\rightarrow \pl_{y}^p(X,\oplus E_n)'$ est surjective,
\item l'injection canonique $\Lambda_q:\pl_{y}^p(X,\oplus E_n)\rightarrow \pl_{y}^q(X,\oplus E_n)'$ est surjective,
\item le projecteur $S_p R_p$ est borné sur $L_x^p(X,\oplus E_n)$,
\item le projecteur $S_q R_q$ est borné sur $L_x^q(X,\oplus E_n)$.
\end{enumerate}
Les assertions précédentes impliquent la condition 
\begin{equation}
\label{forcee}
\sup\limits_{n\in \N} \left\Vert\Frac{\sqrt{e(n,x)}}{\sqrt{d_n}} \right\Vert_{L_x^p(X)}\left\Vert \frac{\sqrt{e(n,x)}}{\sqrt{d_n}} \right\Vert_{L_x^q(X)} <+\infty   .\end{equation}
\end{prop}
\begin{demo}
\textbf{Preuve de \eqref{forcee}.} Pour tout $n\in \N$, on fixe $\phi \in E_n$ vérifiant $\left\Vert \phi\right\Vert_{L_y^2(X)}=1$ et l'on note
\bema u(x,y)=\Big( 0,\dots,0,\underbrace{\sqrt{e(n,x)}^{q-1} \phi(y)}_{n},0,\dots \Big) \in L_x^p(X,\oplus E_n). 
\enma
On peut alors écrire 
\bemar
\left\Vert S_p R_p u\right\Vert_{L^p(X,\oplus E_n)} & \leq & \left\Vert S_p R_p\right\Vert \times \left\Vert u\right\Vert_{L^p(X,\oplus E_n)} \\
\left\Vert \Frac{\sqrt{e(n,x)}}{\sqrt{d_n}}\right\Vert_{L_x^p(X)}\times \left\vert \Int_{X} \frac{\sqrt{e(n,x')}^q}{\sqrt{d_n}}d\mu(x') \right\vert & \leq & 
\left\Vert S_p R_p\right\Vert \times \left\Vert \sqrt{e(n,x)}^{q-1}\right\Vert_{L_x^p(X)}\\
\left\Vert \Frac{\sqrt{e(n,x)}}{\sqrt{d_n}}\right\Vert_{L_x^p(X)}\times \Frac{1}{\sqrt{d_n}}\left\Vert \sqrt{e(n,x)} \right\Vert_{L_x^q(X)}^q & \leq & \left\Vert S_p R_p\right\Vert\times \left\Vert \sqrt{e(n,x)}\right\Vert_{L_x^q(X)}^{q-1} \\
\left\Vert\Frac{\sqrt{e(n,x)}}{\sqrt{d_n}} \right\Vert_{L_x^p(X)}\left\Vert \Frac{\sqrt{e(n,x)}}{\sqrt{d_n}} \right\Vert_{L_x^q(X)}  & \leq & \left\Vert S_p R_p\right\Vert.
\enmar
On notera que le membre gauche est symétrique par permutation de $p$ et $q$. Cela suggère que l'on devrait avoir $\left\Vert S_p R_p\right\Vert=\left\Vert S_q R_q \right\Vert$. C'est effectivement le cas comme nous allons le voir maintenant.

\textbf{\'Equivalence des propriétés i), ii), iii) et iv).}
Puisque $\oplus E_n$ est un espace de Hilbert, les espaces de Bochner-Lebesgue $L_{x}^p(X,\oplus E_n)$ et $L_{x}^q(X,\oplus E_n)$ sont duaux l'un de l'autre et la dualité naturelle est donnée par 
\begin{equation}\label{dualit} \forall (u,w)\in    L_x^p(X,\oplus E_n)\times  L_x^q(X,\oplus E_n) \quad \langle u,w\rangle:= \int_{X} \left[\sum_{n\in \N} \int_{X} u_n(x,y) w_n(x,y)d\mu(y) \right] d\mu(x)     .\end{equation}
Or l'on a 
\bemar
\Int_{X} \left[\sum_{n\in \N} \int_{X} |u_n(x,y)| |w_n(x,y)|d\mu(y) \right] d\mu(x) & \leq & \Int_{X} \left[\sum_{n\in \N} \left\Vert u_n(x,y)\right\Vert_{L_y^2(X)}
\left\Vert w_n(x,y)\right\Vert_{L_y^2(X)} \right] d\mu(x) \\
& \leq & \left\Vert u \right\Vert_{L_x^p(X,\oplus E_n)}\left\Vert w \right\Vert_{L_x^q(X,\oplus E_n)} \\
& < & +\infty.
\enmar
Le théorème de Fubini assure que la dualité naturelle s'écrit aussi 
\bema \langle u,w\rangle= \sum_{n\in \N} \int_{X\times X} u_n(x,y) w_n(x,y)d\mu(x)d\mu(y).\enma
La définition \eqref{SR} nous permet alors d'écrire
\bemar \langle S_p R_p u,w\rangle & = &  \Sum_{n\in \N}  \frac{1}{d_n} \Int_{X\times X} \left[ \sqrt{e(n,x)}\int_{X} u_n(x',y)\sqrt{e(n,x')}d\mu(x')  \right]w_n(x,y)d\mu(x) d\mu(y)   .  \enmar
Remarquons maintenant, pour tout $n\in \N$, la majoration triviale de 
\bema \iiint_{X^3} \sqrt{e(n,x)} |w_n(x,y)| \sqrt{e(n,x')} |u_n(x',y)| d\mu(x) d\mu(x') d\mu(y) \enma
par 
\bemar
& & \displaystyle \left(\int_{X} \sqrt{e(n,x)} \left\Vert w_n(x,y)\right\Vert_{L_y^2(X)}d\mu(x)\right)\times \left(\int_{X} \sqrt{e(n,x')} \left\Vert u_n(x',y)\right\Vert_{L_y^2(X)}  d\mu(x') \right) \\
& \leq & \left\Vert \sqrt{e(n,x)}\right\Vert_{L_{x}^p(X)} \left\Vert w_n(x,y) \right\Vert_{L_{x}^q(X,L_y^2(X))} \times \left\Vert \sqrt{e(n,x')}\right\Vert_{L_{x'}^q(X)} \left\Vert u_n(x',y) \right\Vert_{L_{x'}^p(X,L_y^2(X))} \\
& \leq & \left\Vert \sqrt{e(n,x)}\right\Vert_{L_{x}^p(X)} \left\Vert w\right\Vert_{L_x^q(X,\oplus E_n)}\left\Vert \sqrt{e(n,x')}\right\Vert_{L_{x'}^q(X)}
\left\Vert u\right\Vert_{L_x^p(X,\oplus E_n)}\\
& < &+\infty.
\enmar
Le théorème de Fubini nous permet donc d'intervertir $x$ et $x'$ pour obtenir
\bemar \langle S_p R_p u,w\rangle & = &  \Sum_{n\in \N}  \frac{1}{d_n}\Int_{X\times X} \left[ \sqrt{e(n,x)}\int_{X} w_n(x',y)\sqrt{e(n,x')}d\mu(x')  \right]u_n(x,y)d\mu(x)d\mu(y)\\
& = &   \langle u, S_q R_q w\rangle   .     \enmar
Cela signifie que $S_p R_p$ et $S_q R_q$ sont adjoints l'un de l'autre.
L'équivalence iii) $\Leftrightarrow$ iv) est donc claire.

Avant d'aborder l'équivalence i) $\Leftrightarrow$ iv), on remarque que \eqref{defi-S} nous donne l'expression de $\tt S_p$ :
\bemar    \tt S_p : L_x^q(X,\oplus E_n)& \rightarrow & \pl_{y}^p(X,\oplus E_n)' \\
(w_n(x,y))_{n\in \N} & \mapsto & \left((u_n(y))_{n\geq 0} \mapsto \Int_{X\times X} \Sum_{n\in \N} w_n(x,y) u_n(y) \frac{\sqrt{e(n,x)}}{\sqrt{d_n}} d\mu(x)d\mu(y) \right) .\enmar

Montrons l'implication i) $\Rightarrow$ iv). 
Puisque l'application $\Lambda_p$ est une bijection linéaire continue entre deux espaces de Banach, elle est automatiquement un isomorphisme d'espaces de Banach.
\`A l'aide de \eqref{defi-RRp} et \eqref{dualit}, on vérifie la formule $R_q=\Lambda_p^{-1} \circ \tt S_p$.
A fortiori, $S_q R_q$ est un opérateur borné.

Passons à l'implication iv) $\Rightarrow$ i).
L'équivalence iv) $\Leftrightarrow$ iii) assure que $S_p R_p$ est borné.
Puisque $S_p$ et $S_q$ sont des isométries, $R_p$ et $R_q$ sont des opérateurs bornés et l'on a par construction
\bema
R_p S_p = \mbox{id}_{\pl_{y}^p(X,\oplus E_n)} \quad \mbox{et} \quad R_q S_q = \mbox{id}_{\pl_{y}^q(X,\oplus E_n)}.
\enma
L'idée est d'exprimer $\Lambda_p$ avec la formule suivante qui découle facilement de \eqref{defi-S} :
\bema
\Lambda_p= \tt S_p S_q.
\enma
Utilisant que $S_p R_p$ et $S_q R_q$ sont deux opérateurs duaux, on obtient
\begin{eqnarray} \nonumber
\tt R_p \Lambda_p & = & \tt R_p \tt S_p S_q \\ \nonumber
& = & \tt (S_p  R_p) S_q \\ \nonumber
& = & S_q  R_q  S_q \\ \label{form}
& = & S_q .
\end{eqnarray}
Le diagramme suivant est donc commutatif :
\begin{center}
\begin{picture}(8,4.5)
\put(0,0.75){\makebox(3,0.5)[c]{$\pl_{y}^q(X,\oplus E_n)$}}
\put(0,2.25){\makebox(3,0.5)[c]{$\Lambda_p(\pl_{y}^q(X,\oplus E_n))$}}
\put(1,1.5){\makebox(0.5,0.5)[c]{$\Lambda_p$}}
\put(1.5,1.25){\vector(0,1){1}}
\put(0,2.75){\makebox(3,0.5)[c]{$\cup$}}
\put(0,3.25){\makebox(3,0.5)[c]{$\pl_{y}^p(X,\oplus E_n)'$}}
\put(6,2.25){\makebox(2,0.5)[c]{$L_{x}^q(X,\oplus E_n)$}}
\put(3,2.5){\vector(1,0){3}}
\put(3,3.5){\vector(4,-1){3}}
\put(3,1.25){\vector(3,1){3}}
\put(4,1.75){\makebox(1,0.5)[c]{$S_q$}}
\put(4,2.5){\makebox(1,0.5)[c]{$\tt R_p$}}
\put(4,3.1){\makebox(1,0.5)[c]{$\tt R_p$}}
\end{picture}
\end{center}

On doit examiner le diagramme précédent en se rappelant que l'on a 
\bema
\tt S_p \tt R_p=\mbox{id}_{\pl_{y}^p(X,\oplus E_n)'} \quad \mbox{et} \quad R_q S_q = \mbox{id}_{\pl_{y}^q(X,\oplus E_n)}.
\enma
La proposition \ref{prop-retra} et le corollaire \ref{coro-retra} assurent
 que $S_q$ réalise un isomorphisme d'espaces de Banach de $\pl_{y}^{q}(X,\oplus E_n)$ sur l'image du projecteur $S_q R_q$ de $L_{x}^q(X,\oplus E_n)$ et que $\tt R_p$ réalise un isomorphisme d'espaces de Banach de $\pl_{y}^p(X,\oplus E_n)'$ sur l'image du même projecteur $\tt R_p \tt S_p=S_q R_q$ de $L_{x}^q(X,\oplus E_n)$.
Grâce à \eqref{form}, on vérifie que 
$ \Lambda_p ( \pl_{y}^q(X,\oplus E_n)     )$ et 
$\pl_{y}^p(X,\oplus E_n)'$ ont la même image par l'opérateur $\tt R_p$ : 
\bemar   \tt R_p ( \Lambda_p ( \pl_{y}^q(X,\oplus E_n)     )  )& =& S_q (    \pl^q(X,\oplus E_n)     ) = S_q R_q(L_{x}^q(X,\oplus E_n)),      \\
   \tt R_p( \pl_{y}^p(X,\oplus E_n)')& =& S_q R_q(L_{x}^q(X,\oplus E_n))   .  \enmar
Par application de $\tt S_p$, on obtient i) :
\bema  \Lambda_p(\pl_{y}^q(X,\oplus E_n))=\pl_{y}^p(X,\oplus E_n)'.      \enma
En permutant $p$ et $q$, on montre de même l'équivalence ii) $\Leftrightarrow$ iii).
\end{demo}

\vspace*{0.5 cm}
\begin{center}
\textbf{Stratégie pour le théorème \ref{theo-interpo}}
\end{center}
\vspace*{0.5 cm}

Il est naturel d'espérer utiliser une stratégie similaire pour prouver le théorème \ref{theo-interpo} d'interpolation en utilisant cette fois-ci le corollaire \ref{trieb-interpo}.
On suppose donc seulement que l'on a $\frac{1}{p_1}+\frac{1}{p_2}\leq 1$.
La proposition \ref{Lamba-surj} nous informe que la bornitude de $R_p$ pour tout $p\in ]p_1,p_2[$ implique l'assertion \eqref{forcee} :
\bema 
\forall p\in ]p_1,p_2[ \quad \sup\limits_{n\in \N} \left\Vert \frac{\sqrt{e(n,x)}}{\sqrt{d_n}}\right\Vert_{L_x^p(X)}
\left\Vert \frac{\sqrt{e(n,x)}}{\sqrt{d_n}}\right\Vert_{L_x^{\frac{p}{p-1}}(X)}<+\infty.
\enma
En raison de la symétrie entre $p$ et $\frac{p}{p-1}$, les inégalités précédentes sont aussi valides lorsque $p$ parcourt le plus petit ensemble contenant $]p_1,p_2[$ et stable par la fonction $p\mapsto \frac{p}{p-1}$.
Toutes ces inégalités sont trop violentes puisque l'on doit seulement se contenter de l'hypothèse \eqref{hypo-interpo}.
L'application $R_p$ paraît donc inutilisable.
Pour pallier ce problème, nous allons remplacer $R_p$ par une application de la forme 
\bema\begin{array}{rcl}  R_{p,\psi}: L_x^p(X,\oplus E_n) & \rightarrow & \pl_{y}^p(X,\oplus E_n) \\
(u_n(x,y))_{n\in \N} & \mapsto &     \left(\Int_X u_n(x',y) \psi_n(x') d\mu(x')  \right)_{n\geq 0}   ,\end{array}\enma
où $(\psi_n)_{n\in \N}$ est une suite de fonctions de $L^{\frac{p_{2}}{p_2-1}}(X)\cap L^{\frac{p_1}{p_1-1}}(X)$ vérifiant
\bema
\int_{X} \frac{\sqrt{e(n,x)}}{\sqrt{d_n}}\psi_n(x) d\mu(x)=1.
\enma
De nouveau, on a bien $R_{p,\psi} S_p=\mbox{id}_{\pl_{y}^p(X,\oplus E_n)}$ de manière formelle et l'on rencontre le même obstacle : il n'y a aucune raison pour que $R_{p,\psi}$ arrive bien dans $\pl_{y}^p(X,\oplus E_n)$.
Une nouvelle fois, puisque $S_p$ est une isométrie pour tout $p\in ]p_1,p_2[$, la bornitude de $R_{p,\psi}$ équivaut à la bornitude de l'opérateur
\bema\begin{array}{rcl} S_p R_{p,\psi} : L_x^p(X,\oplus E_n) & \rightarrow &  L_x^p(X,\oplus E_n)   \\[3mm]
(u_n(x,y))_{n\in \N} & \mapsto &  \left( \frac{\sqrt{e(n,x)}}{\sqrt{d_n}}\Int_X u_n(x',y) \psi_n(x') d\mu(x') \right)   .   \end{array} \enma
Sous les hypothèses \eqref{loren1} et \eqref{hypo-interpo}, l'existence d'une suite adéquate $(\psi_n)_{n\in \N}$ et
la bornitude de $S_p R_{p,\psi}$ seront établies dans la partie \ref{partVI}.

On peut maintenant expliquer la preuve du théorème \ref{theo-interpo} d'interpolation.
On remarque que les expressions de $S_p$ et $R_{p,\psi}$ sont indépendantes de $p$. On fixe alors $p_1'$ et $p_2'$ deux réels appartenant à $]p_1,p_2[$.
Pour tout $p\in ]p_1',p_2'[$ on note $\theta'\in [0,1]$ l'unique réel vérifiant $\frac{1}{p}=\frac{1-\theta'}{p_1'}+\frac{\theta'}{p_2'}$.
Le corollaire \ref{trieb-interpo} assure que l'image de l'opérateur  
\bema S_p: \left[ \pl_{y}^{p_1'}(X,\oplus E_n),\pl_{y}^{p_2'}(X,\oplus E_n)\right]_{\theta'} \rightarrow  \left[L_x^{p_1'}(X,\oplus E_n),L_x^{p_2'}(X,\oplus E_n)\right]_{\theta'}\enma
est $S_p R_{p,\psi} \left( \left[L_x^{p_1'}(X,\oplus E_n),L_x^{p_2'}(X,\oplus E_n)\right]_{\theta'} \right)$ et que $S_p$ induit un isomorphisme sur son image. 
Le théorème \ref{interpo-compl} et le point i) de la proposition \ref{prop-retra} donnent alors 
\bemar 
S_p\left( \left[ \pl_{y}^{p_1'}(X,\oplus E_n),\pl_{y}^{p_2'}(X,\oplus E_n)\right]_{\theta'}\right)& =& S_p R_{p,\psi}\left( \left[L_x^{p_1'}(X,\oplus E_n),L_x^{p_2'}(X,\oplus E_n)\right]_{\theta'} \right) \\
& =& S_p R_{p,\psi} \left( L_x^{p}(X,\oplus E_n) \right)\\
&=& S_p(\pl_{y}^p(X,\oplus E_n)).
\enmar
Se rappelant l'inclusion \eqref{inclu-inter}, on a donc le schéma
\begin{center}
\begin{picture}(8,2.5)
\put(0,0){\makebox(5,1)[c]{$\left[ \pl_{y}^{p_1'}(X,\oplus E_n),\pl_{y}^{p_2'}(X,\oplus E_n)\right]_{\theta'}$}}
\put(0,1){\makebox(5,0.5)[c]{$\cup$}}
\put(0,1.5){\makebox(5,0.5)[c]{$\pl_{y}^p(X,\oplus E_n)$}}
\put(7,1.5){\makebox(2,0.5)[c]{$L_{x}^p(X,\oplus E_n)$}}
\put(5,1){\vector(3,1){2}}
\put(4,1.75){\vector(1,0){3}}
\put(5,1.75){\makebox(1,0.5)[c]{$S_p$}}
\put(5.5,0.75){\makebox(1,0.5)[c]{$S_p$}}
\end{picture}
\end{center}

Or $S_p$ est isométrique sur $\pl_{y}^p(X,\oplus E_n)$ et donc injectif.
Cela nous amène à l'égalité
\bema
\left[ \pl_{y}^{p_1'}(X,\oplus E_n),\pl_{y}^{p_2'}(X,\oplus E_n)\right]_{\theta'}=\pl_{y}^p(X,\oplus E_n).
\enma
Enfin, les normes des espaces $\left[ \pl_{y}^{p_1'}(X,\oplus E_n),\pl_{y}^{p_2'}(X,\oplus E_n)\right]_{\theta'}$ et $\pl_{y}^p(X,\oplus E_n)$ sont équivalentes d'après le théorème du graphe fermé et l'inclusion continue \eqref{inclu-inter}. On aurait aussi pu invoquer le fait que $S_p$ est un isomorphisme d'espaces de Banach de $\left[ \pl_{y}^{p_1'}(X,\oplus E_n),\pl_{y}^{p_2'}(X,\oplus E_n)\right]_{\theta'}$ sur son image.
La même argumentation est valide en remplaçant la méthode d'interpolation complexe $[\cdot,\cdot]_{\theta'}$ par la méthode d'interpolation réelle $[\cdot,\cdot]_{\theta',p}$.

\section{Preuves des théorèmes \ref{theo-interpo} et \ref{theo-dual}, partie II : espaces de Lorentz}\label{partII}

On effectue quelques rappels sur les espaces de Lorentz $L^{p,\infty}(X)$, avec $p\in ]1,+\infty[$ (voir par exemple \cite[Chapter 1]{grafakos1}).
Afin d'exprimer l'inégalité de Hölder des espaces de Lorentz, il sera utile de remarquer la reformulation suivante de la quasi-norme $\left\Vert \cdot\right\Vert_{L^{p,\infty}(X)}$ : 
\bemar  \forall f\in L^{p,\infty}(X)\quad \left\Vert f\right\Vert_{L^{p,\infty}(X)}& :=&  \inf \left\{c>0,\quad \forall t>0\quad \mu\{x,\in X, |f(x)|>t\}\leq \Frac{c^p}{t^{p}} \right\} \\
& = &  \sup\limits_{T>0} T^{\frac{1}{p}} f^\star (T),  \enmar
où $f^\star:[0,+\infty[\rightarrow [0,+\infty[$ est le réarrangement décroissant de $f$ définie par la formule 
\bema  f^\star(T):=\inf\{ t>0, \quad \mu\{|f|>t \}\leq T \}         .\enma
En général, la quasi-norme $\left\Vert \cdot\right\Vert_{L^{p,\infty}(X)}$ ne vérifie pas l'inégalité triangulaire mais est toujours équivalente à la norme suivante dès lors que l'on a $p>1$ (\cite[Page 13 and Page 64]{grafakos1} ou \cite[Part V, Lemma 2.8]{garcia1985}) : 
\begin{equation}\label{loren-norm}
\forall f\in L^{p,\infty}(X) \qquad \tnorm{f}_{L^{p,\infty}(X)}:= \sup\limits_{\substack{A\in \mathcal{B}(X)\\ 0<\mu(A)<+\infty}} \frac{\left\Vert f \pun_{A}\right\Vert_{L^1(X)} }{ \mu(A)^{1-\frac{1}{p}} }  ,
\end{equation}
où $\mathcal{B}(X)$ désigne l'ensemble des parties mesurables de $X$.
Précisément, nous avons 
\bema      \left\Vert f \right\Vert_{L^{p,\infty}(X)}\leq \tnorm{f}_{L^{p,\infty}(X)}\leq  \frac{p}{p-1} \left\Vert f \right\Vert_{L^{p,\infty}(X)}.          \enma
En d'autres termes, quitte à perdre une constante multiplicative, on pourra utiliser l'inégalité triangulaire dans $L^{p,\infty}(X)$.
On aura aussi besoin de l'espace de Lorentz $L^{p,1}(X)$ : il s'agit de l'espace vectoriel des fonctions mesurables $g:X\rightarrow \C$ qui vérifient 
\bema \left\Vert g\right\Vert_{L^{p,1}(X)}:=\int_{0}^{+\infty} T^{\frac{1}{p}} g^\star (T) \frac{dT}{T}  <+\infty.   \enma
Venons-en maintenant aux inégalités de Hardy-Littlewood et de Hölder, on a pour toutes fonctions  $f\in L^{p,\infty}(X)$ et $g\in L^{\frac{p}{p-1},1}(X)$ 
\begin{eqnarray}  \nonumber \left\vert  \Int_{X} f(x)g(x)d\mu(x) \right\vert & \leq & \Int_{0}^{+\infty} f^\star(T)g^{\star}(T)dT =\Int_{0}^{+\infty} T^{\frac{1}{p}} f^\star(t) T^{\frac{p-1}{p}} g^{\star}(T) \frac{dT}{T}  \\ \nonumber 
& \leq & \left\Vert f\right\Vert_{L^{p,\infty}(X)} \left\Vert g\right\Vert_{L^{\frac{p}{p-1},1}(X)} \\ \label{hold-lore}
& \leq & \tnorm{f}_{L^{p,\infty}(X)} \left\Vert g\right\Vert_{L^{\frac{p}{p-1},1}(X)}.  \end{eqnarray}

On peut aussi définir des espaces de Lorentz $L^{p,r}(X,B)$, avec $r\in \{1,\infty\}$ (et même tout réel $r\geq 1$), à valeurs dans un espace de Banach complexe $B$. Il s'agit des fonctions mesurables $f:X\rightarrow B$ telles que $\left\Vert x\mapsto \left\Vert f(x)\right\Vert_{B}\right\Vert_{L^{p,r}(X)}<+\infty$.
La théorie de l'interpolation réelle fait jouer un rôle important aux espaces de Lorentz, en particulier on a le résultat suivant \cite[Part 1.18.6, Lemma, Theorem 2, line (16), $q=p$ and Part 1.18.7, Theorem 2]{triebel}.

\begin{theo}\label{interpo-loren}
Considérons un espace de Banach complexe $B$ et trois réels $p_1<p<p_2$ appartenant à $]1,+\infty[$.
Soit $\theta\in ]0,1[$ l'unique réel tel que $\frac{1}{p}=\frac{1-\theta}{p_1}+\frac{\theta}{p_2}$.
Pour tout $r\in \{1,+\infty\}$, on a pour la méthode d'interpolation réelle 
\bema 
[L^{p_1,r}(X,B),L^{p_2,r}(X,B)]_{\theta,p}=L^{p}(X,B).\enma
Si un opérateur linéaire $T$ est borné de $L^{p_1,1}(X,B)$ à valeurs dans $L^{p_1,\infty}(X,B)$ et de $L^{p_2,1}(X,B)$ à valeurs dans $L^{p_2,\infty}(X,B)$ alors il est borné de $L^p(X,B)$ à valeurs dans $L^p(X,B)$ et il existe une constante $K>0$ indépendante de $T$ telle que 
\bema
\left\Vert T\right\Vert_{L^p(X,B)\rightarrow L^p(X,B)}\leq K  \left\Vert T\right\Vert_{L^{p_1,1}(X,B)\rightarrow L^{p_1,\infty}(X,B)}^{\theta}   \left\Vert T\right\Vert_{L^{p_2,1}(X,B)\rightarrow L^{p_2,\infty}(X,B)}^{1-\theta}    .
\enma
\end{theo}

\section{Preuves des théorèmes \ref{theo-interpo} et \ref{theo-dual}, partie III : $R$-bornitude}\label{partIII}

\'Enonçons deux faits qui vont justifier que la bornitude d'un opérateur linéaire sur $L_x^p(X,\ell^2(\N))$ n'est généralement pas simple et qui vont motiver l'approche qui va suivre.

\textbf{Fait 1} : \textit{soit $(P_n)_{n\geq 0}$ une suite d'opérateurs uniformément bornés sur $L^2(X)$, c'est-à-dire que l'on a
\bema\sup\limits_{n\geq 0} \left\Vert P_n\right\Vert_{L^2(X)\rightarrow L^2(X)}<+\infty,\enma
 alors l'opérateur $\oplus P_n$ défini par l'expression suivante est borné sur $L^2(X,\ell^2(\N))$ }: 
\bema \forall (f_n)_{n\geq 0} \in L^2(X,\ell^2(\N))\quad (\oplus P_n) (f_n)=(P_n f_n).
 \enma
 En effet, cela découle immédiatement des formules : 
\bema  \left\Vert (f_n)\right\Vert_{L^2(X,\ell^2(\N))}=\sqrt{\int_{X}  \sum_{n\geq 0} |f_n(x)|^2 d\mu(x)   }=\sqrt{\sum_{n\geq 0} \left\Vert f_n\right\Vert_{L^2(X)}^2 }.  \enma

\textbf{Fait 2} : \textit{pour tout réel $p>2$ il existe une suite $(P_n)_{n\geq 0}$ d'opérateurs uniformément bornés de $L^p(\R)$ et tels que l'opérateur $\oplus P_n$ défini par 
\begin{equation}\label{oplusTn} \forall (f_n)_{n\geq 0} \in L^p(X,\ell^2(\N))\quad (\oplus P_n) (f_n)=(P_n f_n).   \end{equation}
n'est pas borné sur $L^p(\R,\ell^2(\N))$.}

L'exemple est élémentaire. On considère les isométries $P_n$ définies par
\bema\forall n\in \N \quad \forall f\in L^p(\R)\quad \forall x\in\R\quad (P_n f)(x)= f(x+n),     \enma
puis les fonctions $f_{N,n}$, paramétrées par $(n,N)\in \N^2$, définies par
\bema \forall x\in \R \quad f_{N,n}(x):=\left\{ \begin{array}{rcl} \pun_{[n,n+1[}(x) & \mbox{ si } & n< N \\
0& \mbox{ si } & n\geq N .\end{array}\right.   \enma
On a immédiatement 
$ \left\Vert P_n\right\Vert_{L^p(X)\rightarrow L^p(X)}=1   $. De même pour tout $N\in \N^\star$, on a 
\bema
\begin{array}{rcccl}
\left\Vert (\oplus P_n)(f_{N,n})\right\Vert_{L^p(\R,\ell^2(\N))}^p & = & \Int_{\R} \left(\sum_{n=0}^{N-1}\pun_{[0,1[}(x) \right)^{\frac{p}{2}}dx & = & N^{\frac{p}{2}}        ,    \\
\left\Vert (f_{N,n})\right\Vert_{L^p(\R,\ell^2(\N)}^p & =& \Int_{\R} \left( \sum_{n=0}^{N-1} \pun_{[n,n+1[}(x) \right)^{\frac{p}{2}} dx &=& N  \ll N^{\frac{p}{2}} .
\end{array}
\enma
En faisant tendre $N$ vers $+\infty$, on voit que l'opérateur $\oplus P_n$ n'est pas borné sur $L^p(\R,\ell^2(\N))$.
L'exemple précédent utilise le fait que $\R$ n'est pas compact afin de s'échapper vers l'infini.
En réalité, cela est un faux-semblant et l'on peut transférer la construction précédente sur $L^p(0,1)$ à l'aide d'une isométrie linéaire surjective adéquate de $L^p(0,1)$ sur $L^p(\R)$.
En outre, si l'on raisonne par dualité ou si l'on choisit $f_n=\pun_{[0,1[}$ à la place de $\pun_{[n,n+1[}$, la construction précédente s'adapte au cas $p<2$.  
Ainsi, le \textbf{Fait 2} contraste fortement avec le \textbf{Fait 1} et montre que la seule condition 
\bema\sup\limits_{n\in \N} \left\Vert P_n\right\Vert_{L^p(X)\rightarrow L^p(X)}<+\infty\enma
ne suffit pas pour assurer la bornitude de l'opérateur $\oplus P_n$ sur $L^p(X,\ell^2(\N))$ pour $p\neq  2$, même si $\mu(X)<+\infty$.
Un exemple qui illustre ce problème est la non continuité du multiplicateur de Fourier sur $L^p(\R^2)$ de l'indicatrice de la boule unité (voir la preuve de \cite{fefferman1971} dans laquelle les opérateurs $P_n$ sont des projecteurs de $L^p(\R^2)$).

 La notion mathématique qui s'est dégagée est la $\mathcal{R}$-bornitude.
\'Etant donné un espace de Banach $B$, une suite $(P_n)_{n\geq 0}$ d'opérateurs linéaires bornés de $B$ est dite $\mathcal{R}$-bornée s'il existe $K>0$ telle que pour toute suite $(f_n)_{n\in \N}$ de $B$ (nulle pour $n\gg 1$) on a un principe de contraction : 
\bema
\E\left[ \left\Vert \sum_{n\geq 0} \ep_n P_n(f_n) \right\Vert_{B} \right]\leq K \E\left[ \left\Vert \sum_{n\geq 0} \ep_n f_n \right\Vert_{B} \right].
\enma
Dans le cas $B=L^p(X)$, il est classique que les deux espérances précédentes s'estiment avec le théorème de Fubini et les inégalités de Kahane-Khintchine dans $L^p(X)$ et dans $\C$ (voir \eqref{KK-LP} dans le cas unidimensionnel).
L'inégalité précédente signifie alors que l'opérateur $\oplus P_n$, défini en \eqref{oplusTn}, est borné sur $L^p(X,\ell^2(\N))$. Quitte à changer $K>0$, cela équivaut à une estimation de la forme
\bema    \int_X \left( \sum_{n\geq 0} |(P_n f_n)(x)|^2 \right)^{\frac{p}{2}}d\mu (x) \leq K \int_X \left(\sum_{n\geq 0} |f_n(x)|^2 \right)^{\frac{p}{2}} d\mu(x)  .   \enma

Différents critères de $\mathcal{R}$-bornitude existent dans la littérature 
(voir par exemple \cite[Chapter V]{garcia1985}, \cite[Section 2]{rbounded2}, \cite{rbounded3} et les références indiquées).
Le lemme \ref{schur-1} donne un nouveau critère qui s'apparente au lemme de Schur et qui donne des conditions suffisantes pour obtenir la $\mathcal{R}$-bornitude.
Ce critère suffira pour la théorie de l'interpolation et de la dualité des espaces $\pl^{p}(X,\oplus E_n)$ dans le cas unidimensionnel $d_n=\dim(E_n)=1$.
Pour comprendre le cas $d_n\neq 1$, nous aurons besoin d'une version légèrement plus sophistiquée du lemme \ref{schur-1}, à savoir le lemme \ref{schur}, mais le coeur de l'idée est dans la démonstration suivante.

\begin{lemm}\label{schur-1}
Fixons des réels $p_1<p<p_2$ appartenant à $]1,+\infty[$.
Notons $q_1:=\frac{p_1}{p_1-1}$ et $q_2:=\frac{p_2}{p_2-1}$.
Considérons $K_n:X^2\rightarrow \C$ des fonctions mesurables, pour $n \in \N$, de sorte que l'on a
\begin{eqnarray} \label{schur1}
\sup\limits_{n\geq 0} \left\Vert K_n(x,x')\right\Vert_{L_{x'}^{q_1,\infty}(X)} &\in & L_x^{p_1,\infty}(X),\\ \label{schur2}
\sup\limits_{n\geq 0} \left\Vert K_n(x,x')\right\Vert_{L_{x'}^{q_2,\infty}(X)} &\in & L_x^{p_2,\infty}(X),\\ \label{schur3}
\sup\limits_{n\geq 0} \left\Vert K_n(x,x')\right\Vert_{L_x^{p_1,\infty}(X)} &\in & L_{x'}^{q_1,\infty}(X),\\ \label{schur4}
\sup\limits_{n\geq 0} \left\Vert K_n(x,x')\right\Vert_{L_x^{p_2,\infty}(X)} &\in & L_{x'}^{q_2,\infty}(X).
\end{eqnarray}
Les deux assertions suivantes sont vraies
\begin{enumerate}[a) ]
\item pour tout entier $n\in \N$ l'opérateur $P_n$ de noyau $K_n$ défini par 
\bema \forall f \in L_x^p(X) \quad \forall x\in X \quad  (P_n f)(x)=\int_{X} K_n(x,x') f(x') d\mu(x')     ,    \enma
est borné sur $L^p(X)$ et la borne supérieure $\sup\limits_{n\geq 0} \left\Vert P_n\right\Vert_{L^p(X)\rightarrow L^p(X)}$ est finie.
\item la famille d'opérateurs $(P_n)_{n\in \N}$ de $L^p(X)$ est $\mathcal{R}$-bornée. 
\end{enumerate}
\end{lemm}
\begin{demo}
Nous aurons besoin des propriétés des espaces de Lorentz rappelées dans la partie \ref{partII}.

a) On affirme que \eqref{schur1} implique la bornitude de l'opérateur 
\bema P_n : L_x^{p_1,1}(X)  \rightarrow  L_x^{p_1,\infty}(X). \enma
En effet, l'inégalité de Hölder \eqref{hold-lore} donne 
\bemar
|(P_n f)(x)| &\leq & \left\Vert K_n(x,x') \right\Vert_{L_{x'}^{q_1,\infty}(X)} \left\Vert f\right\Vert_{L^{p_1,1}(X)} \\
\left\Vert P_n f\right\Vert_{L^{p_1,\infty}(X)} & \leq & 
\left\Vert \left\Vert K_n(x,x') \right\Vert_{L_{x'}^{q_1,\infty}(X)} \right\Vert_{L_{x}^{p_1,\infty}(X)} \left\Vert f\right\Vert_{L^{p_1,1}(X)}.
\enmar

Afin de faciliter la lecture de la preuve de l'assertion b), nous donnons un autre argument qui utilise \eqref{schur3} et la norme $\tnorm{\cdot}_{L^{p_1,\infty}(X)}$ de $L^{p_1,\infty}(X)$, définie en \eqref{loren-norm}, car elle présente l'avantage de satisfaire l'inégalité triangulaire. En effet, on a 
\bemar
\tnorm{P_n f}_{L^{p_1,\infty}(X)}& \leq & \Int_{X} \tnorm{K_n(x,x')}_{L_x^{p_1,\infty}(X)} |f(x')| d\mu(x') \\[3mm]
& \leq & \left\Vert \tnorm{K_n(x,x')}_{L_x^{p_1,\infty}(X)}\right\Vert_{L_{x'}^{q_1,\infty}(X)} \left\Vert f\right\Vert_{L^{p_1,1}(X)}.
\enmar
Les deux arguments précédents sont en fait duaux.
En utilisant \eqref{schur2} ou \eqref{schur4}, on obtient de même la bornitude de l'opérateur 
\bema
P_n : L_x^{p_2,1}(X)  \rightarrow  L_x^{p_2,\infty}(X).
\enma
On conclut par interpolation réelle (à savoir le théorème \ref{interpo-loren} avec $B=\C$).

b) On affirme que l'opérateur suivant est borné 
\bema
\oplus P_n : L_x^{p_1,1}(X,\ell^\infty(\N)) \rightarrow L_x^{p_1,\infty}(X,\ell^\infty(\N)).
\enma
En effet, cela découle de l'inégalité de Hölder \eqref{hold-lore} et de \eqref{schur1} : 
\begin{equation}\nonumber
\begin{array}{rcl}
\sup\limits_{n\in \N} |(P_n f_n)(x)|&=& \sup\limits_{n\in \N}\left\vert \Int_{X} K_n(x,x')f_n(x')d\mu(x')    \right\vert \\
&\leq & \sup\limits_{n\in \N}   \Int_{X} |K_n(x,x')| \sup\limits_{m\in \N} |f_m(x')|d\mu(x')   \\
&\leq & \underbrace{\left( \sup\limits_{n\in \N} \left\Vert K_n(x,x')\right\Vert_{L_{x'}^{q_1,\infty}(X)} \right)}_{\in L_x^{p_1,\infty}(X)} \times \left\Vert \sup\limits_{m\in \N} |f_m(x')| \right\Vert_{L_{x'}^{p_1,1}(X)} .
\end{array} \end{equation}
On s'attelle maintenant à prouver la bornitude de l'opérateur 
\bema
\oplus P_n : L_{x}^{p_1,1}(X,\ell^1(\N)) \rightarrow L_{x}^{p_1,\infty}(X,\ell^1(\N)).
\enma
De nouveau, on utilise la norme $\tnorm{\cdot}_{L^{p_1,\infty}(X)}$ de $L^{p_1,\infty}(X)$.
L'hypothèse \eqref{schur3} et l'inégalité de Hölder \eqref{hold-lore} nous amènent aux estimations
\bemar
\tnorm{\Sum_{n\in \N} |(P_n f)(x)|}_{L_x^{p_1,\infty}(X)} & = & \tnorm{\Sum_{n\geq 0} \left\vert \Int_{X} K_n(x,x') f_n(x')d\mu(x') \right\vert }_{L_x^{p_1,\infty}(X)} \\
 & \leq & \Sum_{n\geq 0} \Int_{X} \tnorm{K_{n}(x,x')}_{L_{x}^{p_1,\infty}(X)}|f_n(x')|d\mu(x')      \\
 & \leq & \Sum_{n\geq 0} \Int_{X} \left(\sup\limits_{m\in \N} \tnorm{K_{m}(x,x')}_{L_{x}^{p_1,\infty}(X)}\right) |f_n(x')|d\mu(x')\\
 & \leq & \Int_{X} \left(\sup\limits_{m\in \N} \tnorm{K_{m}(x,x')}_{L_{x}^{p_1,\infty}(X)}\right) \left(\Sum_{n\geq 0} |f_n(x')| \right)d\mu(x') \\
 & \leq & \left\Vert \sup\limits_{m\in \N} \tnorm{K_{m}(x,x')}_{L_{x}^{p_1,\infty}(X)}\right\Vert_{L_{x'}^{q_1,\infty}(X)}\left\Vert \Sum_{n\geq 0} |f_n(x')|\right\Vert_{L_{x'}^{p_1,1}(X)}.
\enmar
La démonstration est évidemment similaire en remplaçant $p_1$ par $p_2$ et en utilisant \eqref{schur4} et \eqref{schur2}. On a ainsi obtenu la bornitude des quatre opérateurs 
\begin{displaymath}\begin{array}{lcl}   \oplus P_n : L_x^{p_1,1}(X,\ell^\infty(\N))& \rightarrow &  L_x^{p_1,\infty}(X,\ell^ \infty(\N))   , \\
\oplus P_n : L_x^{p_2,1}(X,\ell^\infty(\N))& \rightarrow &  L_x^{p_2,\infty}(X,\ell^\infty(\N)),\\ 
\oplus P_n : L_x^{p_1,1}(X,\ell^1(\N) & \rightarrow & L_x^{p_1,\infty}(X,\ell^1(\N)), \\
\oplus P_n : L_x^{p_2,1}(X,\ell^1(\N)) & \rightarrow & L_x^{p_2,\infty}(X,\ell^{1}(\N)) .\end{array}  \enma
Pour tout $p\in ]p_1,p_2[$, le théorème \ref{interpo-loren} assure la bornitude des deux opérateurs suivants par interpolation réelle
\bema\begin{array}{lcl}
\oplus P_n : L_x^{p}(X,\ell^\infty(\N))& \rightarrow &  L_x^{p}(X,\ell^\infty(\N))   ,\\
\oplus P_n : L_x^{p}(X,\ell^1(\N)) & \rightarrow & L_x^{p}(X,\ell^1(\N)).\end{array}
\enma
Par interpolation complexe (c'est-à-dire le théorème \ref{interpo-compl} avec $\theta=\frac{1}{2}$), on obtient la bornitude de l'opérateur 
\bema
\oplus P_n : L_x^{p}(X,\ell^2(\N))  \rightarrow  L_x^{p}(X,\ell^2(\N)). \enma
\end{demo}

Il s'agit maintenant d'écrire un résultat analogue au lemme \ref{schur-1} adaptée à la théorie générale avec $\dim(E_n)\neq 1$.
Pour cela, on aura besoin d'espaces analogues à $\ell^1(\N)$ et $\ell^\infty(\N)$ adaptés aux sous-espaces $E_n$, ce sera le rôle joué par les espaces $(\oplus E_n)_{\ell^r}$ dans la preuve du résultat suivant.

\begin{lemm}\label{schur}
Fixons des réels $p_1<p<p_2$ appartenant à $]1,+\infty[$ et considérons $K_n:X^2\rightarrow \C$ des fonctions mesurables, pour $n \in \N$, de sorte que 
\begin{eqnarray} \label{schur1-g}
\sup\limits_{n\geq 0} \left\Vert K_n(x,x')\right\Vert_{L_{x'}^{q_1,\infty}(X)} &\in & L_x^{p_1,\infty}(X),\\ \label{schur2-g}
\sup\limits_{n\geq 0} \left\Vert K_n(x,x')\right\Vert_{L_{x'}^{q_2,\infty}(X)} &\in & L_x^{p_2,\infty}(X),\\ \label{schur3-g}
\sup\limits_{n\geq 0} \left\Vert K_n(x,x')\right\Vert_{L_x^{p_1,\infty}(X)} &\in & L_{x'}^{q_1,\infty}(X),\\ \label{schur4-g}
\sup\limits_{n\geq 0} \left\Vert K_n(x,x')\right\Vert_{L_x^{p_2,\infty}(X)} &\in & L_{x'}^{q_2,\infty}(X).
\end{eqnarray}
Alors les deux assertions suivantes sont vraies 
\begin{enumerate}[a) ]
\item pour tout entier $n\in \N$, l'opérateur $P_n$ défini par 
\bema \forall u \in L_x^p(X, E_n) \quad \forall (x,y)\in X^2 \quad  (P_n u)(x,y)=\int_{X} K_n(x,x') u(x',y) d\mu(x')     ,    \enma
est borné sur $L_x^p(X, E_n)$ et la borne supérieure $\sup\limits_{n\geq 0} \left\Vert P_n\right\Vert$ est finie.
\item l'opérateur $\oplus P_n$, défini par l'expression suivante, est borné
\bemar \oplus P_n : L_x^p(X,\oplus E_n) & \rightarrow & L_x^p(X,\oplus E_n)  \\
   \left( u_n(x,y)\right)_{n\in \N} & \mapsto & \left( (P_n u_n)(x,y)\right)_{n\in \N}.
      \enmar 
\end{enumerate}
\end{lemm}
\begin{demo}
a) On raisonne comme pour le lemme \ref{schur-1}, c'est-à-dire que chaque opérateur $P_n$ est borné de $L_x^{p_1,1}(X, E_n)$ à valeurs dans $L_x^{p_1,\infty}(X, E_n)$ et de $L_x^{p_2,1}(X, E_n)$ à valeurs dans $L_x^{p_2,\infty}(X, E_n)$ avec des normes majorées indépendamment de $n$. On conclut par interpolation réelle (c'est-à-dire le théorème \ref{interpo-loren} avec $B=E_n$).

b) L'idée est d'interpoler l'espace de Hilbert $\oplus E_n$ entre les deux espaces 
\bemar
 \left(\oplus E_n \right)_{\ell^1}& := & \left\{(u_n)_{n\in \N},\quad \forall n\in \N, u_n\in E_n; \quad \Sum_{n\in \N} \left\Vert u_n \right\Vert_{L_y^2(X)}<+\infty  \right\} ,\\
  \left(\oplus E_n \right)_{\ell^\infty} &:= &\{(u_n)_{n\in \N},\quad \forall n\in \N, u_n\in E_n; \quad \sup\limits_{n\in \N} \left\Vert u_n \right\Vert_{L_y^2(X)}<+\infty  \} ,
\enmar
qui sont munis de leurs normes naturelles.
 On vérifie aisément que les espaces $(\oplus E_n)_{\ell^1}$ et $(\oplus E_n)_{\ell^\infty}$ sont des rétractes, au sens de la définition \ref{defi-retra}, des espaces $\ell^1(\N,\oplus E_n)$ et $\ell^\infty(\N,\oplus E_n)$.
Par suite, le corollaire \ref{trieb-interpo} et la théorie de l'interpolation complexe des espaces $\ell^r(\N,\oplus E_n)$ (c'est-à-dire le théorème \ref{interpo-compl}) montrent que l'on a l'égalité
\bema [\left(\oplus E_n \right)_{\ell^1},\left(\oplus E_n \right)_{\ell^\infty}]_{\frac{1}{2}} = \oplus E_n,\enma
avec équivalence de normes.
Nous pouvons donc reprendre l'argumentation du lemme \ref{schur-1} et il nous suffit manifestement de prouver la bornitude des quatre opérateurs 
\begin{displaymath}\begin{array}{lcl}   
\oplus P_n : L_x^{p_1,1}(X,(\oplus E_n)_{\ell^1}) & \rightarrow & L_x^{p_1,\infty}(X,(\oplus E_n)_{\ell^1}), \\
\oplus P_n : L_x^{p_2,1}(X,(\oplus E_n)_{\ell^1}) & \rightarrow & L_x^{p_2,\infty}(X,(\oplus E_n)_{\ell^1}) ,\\
\oplus P_n : L_x^{p_1,1}(X,(\oplus E_n)_{\ell^\infty})& \rightarrow &  L_x^{p_1,\infty}(X,(\oplus E_n)_{\ell^\infty})   , \\
\oplus P_n : L_x^{p_2,1}(X,(\oplus E_n)_{\ell^\infty})& \rightarrow &  L_x^{p_2,\infty}(X,(\oplus E_n)_{\ell^\infty}).\end{array}  \enma
Par symétrie évidente, il nous suffit de justifier la bornitude des deux opérateurs 
\begin{displaymath}\begin{array}{lcl}  \oplus P_n : L_x^{p_1,1}(X,(\oplus E_n)_{\ell^1}) & \rightarrow & L_x^{p_1,\infty}(X,(\oplus E_n)_{\ell^1}),\\
 \oplus P_n : L_x^{p_1,1}(X,(\oplus E_n)_{\ell^\infty})& \rightarrow &  L_x^{p_1,\infty}(X,(\oplus E_n)_{\ell^\infty})   .\end{array}  \enma
L'hypothèse \eqref{schur3-g} et l'inégalité de Hölder \eqref{hold-lore} nous donnent les majorations
\bema \begin{array}{l}
 \tnorm{\Sum_{n\geq 0} \left\Vert (P_n u_n)(x,y) \right\Vert_{L_y^2(X)} }_{L_{x}^{p_1,\infty}(X)} \\
\hspace*{1 cm} \leq   \tnorm{\Sum_{n\geq 0} \left\Vert \Int_{X} K_n(x,x') u_n(x',y)d\mu(x') \right\Vert_{L_y^2(X)}}_{L_{x}^{p_1,\infty}(X)} \\
\hspace*{1 cm} \leq   \Sum_{n\geq 0} \Int_{X} \tnorm{K_n(x,x')}_{L_{x}^{p_1,\infty}(X)} \left\Vert u_n(x',y)\right\Vert_{L_y^2(X)} d\mu(x') \\
\hspace*{1 cm} \leq     \Int_{X} \left( \sup\limits_{m\in \N} \tnorm{K_m(x,x')}_{L_{x}^{p_1,\infty}(X)} \right)\left( \Sum_{n\geq 0}\left\Vert u_n(x',y)\right\Vert_{L_y^2(X)} \right)d\mu(x')            \\
 \hspace*{1 cm} \leq   \left\Vert \sup\limits_{m\in \N} \tnorm{K_m(x,x')}_{L_x^{p_1,\infty}(X)}\right\Vert_{L_{x'}^{q_1,\infty}(X)} \left\Vert \Sum_{n\geq 0} \left\Vert u_n(x',y)\right\Vert_{L_y^2(X)}\right\Vert_{L_{x'}^{p_1,1}(X)}.
\end{array} \enma
Enfin, \eqref{schur1-g} et l'inégalité de Hölder \eqref{hold-lore} nous permettent d'écrire
\bemar
\sup\limits_{n\in \N} \left\Vert (P_n u_n)(x,y)\right\Vert_{L_y^2(X)}     
&=& \sup\limits_{n\in \N}\left\Vert \Int_{X} K_n(x,x') u_n(x',y)d\mu(x')    \right\Vert_{L_y^2(X)} \\
&\leq & \sup\limits_{n\in \N}   \Int_{X} |K_n(x,x')| \sup\limits_{m\in \N} \left\Vert  u_m(x',y)\right\Vert_{L_y^2(X)} d\mu(x')   \\
&\leq &\underbrace{\left(  \sup\limits_{n\in \N} \left\Vert K_n(x,x')\right\Vert_{L_{x'}^{q_1,\infty}(X)} \right)}_{\in L_{x}^{p_1,\infty}(X)} \times \left\Vert \sup\limits_{m\in \N} \left\Vert  u_m (x',y)\right\Vert_{L_y^2(X)} \right\Vert_{L_{x'}^{p_1,1}(X)} .
\enmar
\end{demo}

\section{Preuves des théorèmes \ref{theo-interpo} et \ref{theo-dual}, partie IV : bornitude de $S_p R_p$}\label{partIV}

Sous les hypothèses du théorème \ref{theo-dual} de dualité, les exposants $p_1$ et $p_2$ sont conjugués et l'on souhaite montrer que $S_p R_p$ est un projecteur borné sur $L_x^p(X,\oplus E_n)$ pour tout $p\in ]p_1,p_2[$.
Avec les notations usuelles, on a $q_1=\frac{p_1}{p_1-1}=p_2$ et $q_2=\frac{p_2}{p_2-1}=p_1$.
D'après la forme du projecteur \eqref{SR}, il s'agit d'appliquer le lemme \ref{schur} avec la suite de noyaux $K_n$ définis par
\bema
\forall (x,x')\in X^2 \quad K_n(x,x')=\frac{1}{d_n} \sqrt{e(n,x) e(n,x')}.
\enma
Il est facile de constater que \eqref{loren1} et \eqref{hol1} impliquent l'assertion \eqref{schur1-g} (qui est identique à \eqref{schur4-g}). En effet, on a pour tout $(x,n)\in X\times \N$ 
\bemar
\left\Vert \Frac{1}{d_n} \sqrt{e(n,x) e(n,x')}\right\Vert_{L_{x'}^{q_1,\infty}(X)} &\leq &
\Frac{\sqrt{e(n,x)}}{d_n} \left\Vert  \sqrt{e(n,x')}\right\Vert_{L_{x'}^{q_1}(X)} \\
& \leq & \Frac{\sqrt{e(n,x)}}{\left\Vert \sqrt{e(n,x')}\right\Vert_{L_{x'}^{p_1}(X) }} \times \Frac{\left\Vert \sqrt{e(n,x')}\right\Vert_{L_{x'}^{p_1}(X) }\left\Vert \sqrt{e(n,x')}\right\Vert_{L_{x'}^{q_1}(X) }}{d_n}.
\enmar

Les deux autres hypothèses \eqref{schur2-g} et \eqref{schur3-g} 
du lemme \ref{schur} se vérifient de la même façon.
La preuve du théorème \ref{theo-dual} de dualité est finie.

\section{Preuves des théorèmes \ref{theo-interpo} et \ref{theo-dual}, partie V :  défaut d'interpolation}\label{partV}

Nous allons définir la notion de ``défaut d'interpolation'' afin d'aborder l'interpolation des espaces $\pl^p(X,\oplus E_n)$.
Commençons par introduire quelques notations.
Pour tous réels $p_1,p$ et $p_2$ appartenant à $[1,+\infty[$ et vérifiant $p_1\leq p\leq p_2$, on définit les nombres $\theta_1(p_1,p,p_2)$ et $ \theta_2(p_1,p,p_2)$ tels que 
\bema    \frac{\theta_1(p_1,p,p_2)}{p_1}+\frac{\theta_2(p_1,p,p_2)}{p_2}=\frac{1}{p} \qquad \mbox{et} \qquad \theta_1(p_1,p,p_2)+\theta_2(p_1,p,p_2)=1.\enma
De façon précise, on a les formules 
\bema
\theta_1(p_1,p,p_2)=\frac{\frac{1}{p}-\frac{1}{p_2}}{\frac{1}{p_1}-\frac{1}{p_2}} \qquad \mbox{et} \qquad \theta_2(p_1,p,p_2)=\frac{\frac{1}{p_1}-\frac{1}{p}}{\frac{1}{p_1}-\frac{1}{p_2}} .
\enma

On peut alors poser la définition suivante. 

\begin{defi}
Considérons $p_1$ et $p_2$ appartenant à $[1,+\infty[$ et vérifiant $p_1<p_2$ ainsi qu'une fonction non nulle $\phi\in L^{p_1}(X)\cap L^{p_2}(X)$.
Nous définissons $Q(\phi,[p_1,p_2])\in [1,+\infty[$ le défaut d'interpolation de $\phi$ sur $[p_1,p_2]$ par la formule 
\begin{equation}\label{defi-Q}  Q(\phi,[p_1,p_2]):=\sup\limits_{p\in [p_1,p_2]} \frac{\left\Vert \phi\right\Vert_{L^{p_1}(X)}^{\theta_1(p_1,p,p_2)}\left\Vert \phi\right\Vert_{L^{p_2}(X)}^{\theta_2(p_1,p,p_2)}}{\left\Vert \phi\right\Vert_{L^p(X)}} . \end{equation}
\end{defi}

L'inégalité $Q(\phi,[p_1,p_2])<+\infty$ peut évidemment se voir par continuité et compacité de l'intervalle $[p_1,p_2]$ mais, dans notre contexte, il est plus intéressant d'utiliser un argument de convexité (voir le lemme \ref{majo-Q} ci-après).
Quant à l'inégalité $1\leq Q(\phi,[p_1,p_2])$, elle découle de l'inégalité de Hölder avec les exposants conjugués $\frac{p_1}{p\theta_1(p_1,p,p_2)}$ et $\frac{p_2}{p\theta_2(p_1,p,p_2)}$ : 
\begin{eqnarray}\nonumber
\left(\Int_{X} |\phi(x)|^p d\mu(x) \right)^{\frac{1}{p}}& = & \left(\Int_{X} |\phi(x)|^{p\theta_1(p_1,p,p_2)}|\phi(x)|^{p\theta_2(p_2,p,p_2)}d\mu(x) \right)^{\frac{1}{p}}\\[3mm] \label{cvxit}
& \leq & \left( \Int_{X} |\phi(x)|^{p_1} d\mu(x)\right)^{\frac{\theta_1(p_1,p,p_2)}{p_1}}\left( \Int_{X} |\phi(x)|^{p_2} d\mu(x)\right)^{\frac{\theta_2(p_1,p,p_2)}{p_2}} .
\end{eqnarray}
La condition d'égalité de l'inégalité de Hölder montre alors l'équivalence 
\bema      Q(\phi,[p_1,p_2])=1 \quad \Leftrightarrow \quad  \exists a>0 \quad \exists A\in \mathcal{B}(X), \quad 0<\mu(A)<+\infty, \quad |\phi|=a\pun_{A}, \enma
où $\mathcal{B}(X)$ est l'ensemble des parties mesurables de $X$.
Ainsi, le défaut d'interpolation de la fonction $\phi$ permet de tester si elle se concentre complètement sur une même partie de l'espace mesuré $X$.
Par comparaison avec \eqref{cvxit}, 
le défaut d'interpolation permet de \textbf{minorer} $\left\Vert \phi\right\Vert_{L^p(X)}$ lorsque l'on connaît $\left\Vert \phi\right\Vert_{L^{p_1}(X)}$
et $\left\Vert \phi\right\Vert_{L^{p_2}(X)}$ : 
\begin{equation}\label{defa-mino}
\forall p\in ]p_1,p_2[ \qquad     \frac{\left\Vert \phi\right\Vert_{L^{p_1}(X)}^{\theta_1(p_1,p,p_2)}\left\Vert \phi\right\Vert_{L^p_2(X)}^{\theta_2(p_1,p,p_2)}}{Q(\phi,[p_1,p_2])}   \leq \left\Vert \phi\right\Vert_{L^p(X)}.
\end{equation}
Le lemme suivant montre qu'il suffit d'examiner un seul point de l'intervalle $]p_1,p_2[$ pour contrôler $Q(\phi,[p_1,p_2])$.

\begin{lemm}\label{majo-Q}
Fixons des réels $p_1<p<p_2$ appartenant à $[1,+\infty[$ et une suite de fonctions non nulles $(\phi_n)_{n\in \N}$ de $L^{p_1}(X)\cap L^{p_2}(X)$, alors on a l'équivalence :
 \bema \sup\limits_{n\in \N} 
\frac{\left\Vert \phi_n \right\Vert_{L^{p_1}(X)}^{\theta_1(p_1,p,p_2)}\left\Vert \phi_n\right\Vert_{L^{p_2}(X)}^{\theta_2(p_1,p,p_2)}  }{\left\Vert \phi_n\right\Vert_{L^p(X)}} <+\infty \quad \Leftrightarrow \quad \sup\limits_{n\in \N}  Q(\phi_n,[p_1,p_2])<+\infty.    \enma
\end{lemm}
\begin{demo}
Pour toute fonction non nulle $\phi\in L^{p_1}(X)\cap L^{p_2}(X)$, il est bien connu que la fonction 
\bema
\Phi:\wp\in [p_1,p_2]\mapsto \ln \left( \left\Vert \phi \right\Vert_{L^{\wp }(X)} \right)
\enma
est convexe par rapport à $\frac{1}{\wp }$ (voir les inégalités \eqref{cvxit}). Introduisons alors la fonction
\bema \widetilde{\Phi}: \wp \in [p_1,p_2]\mapsto \theta_1(p_1,\wp,p_2) \Phi(p_1)+ \theta_2(p_1,\wp,p_2) \Phi(p_2)-\Phi(\wp). \enma
La fonction $\widetilde{\Phi}$ s'annule en $p_1$ et $p_2$ et est concave par rapport à $\frac{1}{\wp}$.

\begin{center}
\begin{picture}(7,5)(-1.75,-1)
\put(-0.1,0){\line(1,0){4.2}}
\put(0,-0.1){\line(0,1){3.3}}
\put(4,-0.1){\line(0,1){3.3}}
\put(4,0){\line(-2,1){4}}
\put(0,0){\line(4,3){4}}
\put(4,3){\line(1,0){0.1}}
\put(0,2){\line(-1,0){0.1}}
\put(1.6,-0.1){\line(0,1){1.3}}
\put(2.8,-0.1){\line(0,1){1.7}}
\multiput(2.8,1.6)(0,0.2){3}{\line(0,1){0.1}}
\put(-1.75,1.5){\makebox(1.75,1)[c]{$\frac{\widetilde{\Phi}(p)}{\theta_2(p_1,p,p_2)}$}}
\multiput(1.6,1.2)(0.2,0){13}{\line(1,0){0.1}}
\put(1.1,-1){\makebox(1,1)[c]{$\frac{1}{p}$}}
\put(4,0.9){\makebox(1.3,0.75)[c]{$\scriptstyle\widetilde{\Phi}(p)$}}
\multiput(2.8,1.6)(0.2,0){7}{\line(1,0){0.1}}
\put(2.3,-1){\makebox(1,1)[c]{$\frac{1}{\wp}$}}
\put(4,1.3){\makebox(1.3,0.75)[c]{$\scriptstyle\widetilde{\Phi}(\wp)$}}
\put(-0.5,-1){\makebox(1,1)[c]{$\frac{1}{p_2}$}}
\put(3.5,-1){\makebox(1,1)[c]{$\frac{1}{p_1}$}}
\put(4,2.5){\makebox(1.75,1)[c]{$\frac{\widetilde{\Phi}(p)}{\theta_1(p_1,p,p_2)}$}}
\linethickness{0.25 mm}
\qbezier(0,0)(3.6,3.2)(4,0)
\end{picture}
\end{center}
Par concavité et application du théorème de Thalès, nous avons les estimations 
\bema 
\forall \wp\in [p_1,p_2]\qquad \widetilde{\Phi}(\wp) \leq \widetilde{\Phi}(p)\max \left(\frac{\frac{1}{p_1}-\frac{1}{p_2}}{\frac{1}{p}-\frac{1}{p_2}} ,\frac{\frac{1}{p_1}-\frac{1}{p_2}}{\frac{1}{p_1}-\frac{1}{p}} \right).
\enma
En passant à l'exponentielle, on obtient les estimations suivantes qui donnent la conclusion :
\bema
\frac{\left\Vert \phi \right\Vert_{L^{p_1}(X)}^{\theta_1(p_1,p,p_2)}\left\Vert \phi\right\Vert_{L^{p_2}(X)}^{\theta_2(p_1,p,p_2)}  }{\left\Vert \phi\right\Vert_{L^p(X)}} \leq Q(\phi, [p_1,p_2]) \leq  \left( \frac{\left\Vert \phi \right\Vert_{L^{p_1}(X)}^{\theta_1(p_1,p,p_2)}\left\Vert \phi\right\Vert_{L^{p_2}(X)}^{\theta_2(p_1,p,p_2)}  }{\left\Vert \phi\right\Vert_{L^p(X)}}\right)^{C(p_1,p,p_2)}.
\enma
\end{demo}

On sait que pour tout réel $p\in ]1,+\infty[$ et toute fonction $\phi \in L^p(X)$, il existe une fonction $\psi\in L^{\frac{p}{p-1}}(X)$ telle que 
\bema   \int_{X} \phi(x)\psi(x)d\mu(x)=\left\Vert \phi\right\Vert_{L^p(X)}\left\Vert \psi\right\Vert_{L^{\frac{p}{p-1}}(X)}.    \enma
Par exemple, si $\phi$ est positive alors on peut choisir $\psi(x)=\phi(x)^{p-1}$.
Le défaut d'interpolation permet de formuler des propriétés analogues lorsque $\phi$ appartient à deux espaces de Lebesgue $L^{p_1}(X)$ et $L^{p_2}(X)$.
Dans l'énoncé de la proposition suivante, on utilise les notations habituelles $q_1=\frac{p_1}{p_1-1}$ et $q_2=\frac{p_2}{p_2-1}$.

\begin{prop}\label{hahn}
Fixons $p_1$ et $p_2$ appartenant à $]1,+\infty[$ et vérifiant
$p_1 <p_2$ et $\frac{1}{p_1}+\frac{1}{p_2}\leq 1$.
Il existe un réel $r=r(p_1,p_2)>1$ tel que pour toute fonction non nulle $\phi\in L^{p_1}(X)\cap L^{p_2}(X)$ on peut trouver une fonction $\psi\in L^{q_2}(X)\cap L^{q_1}(X)$ qui vérifie 
\begin{eqnarray} \nonumber
\int_{X} \phi(x)\psi(x)d\mu(x)& =& 1, \\ \nonumber
\left\Vert \phi\right\Vert_{L^{p_1}(X)} \left\Vert\psi\right\Vert_{L^{q_1}(X)} & \leq & Q(\phi,[p_1,p_2])^r ,          \\ \nonumber
          \left\Vert \phi\right\Vert_{L^{p_2}(X)} \left\Vert\psi\right\Vert_{L^{q_2}(X)} & \leq & Q(\phi,[p_1,p_2])^r ,          \\ \nonumber
\Frac{|\psi|^{q_1}}{\int_X |\psi(x)|^{q_1} d\mu(x)} & \leq &  Q(\phi,[p_1,p_2])^{(r-1)q_1}\left[ \Frac{|\phi|^{p_1}}{\int_X |\phi(x)|^{p_1}d\mu(x)}+ \frac{|\phi|^{p_2}}{\int_X |\phi(x)|^{p_2}d\mu(x)}  \right], \\[3mm] \nonumber
\Frac{|\psi|^{q_2}}{\int_X |\psi(x)|^{q_2} d\mu(x)} & \leq &  Q(\phi,[p_1,p_2])^{(r-1)q_2}\left[ \Frac{|\phi|^{p_1}}{\int_X |\phi(x)|^{p_1}d\mu(x)}+ \frac{|\phi|^{p_2}}{\int_X |\phi(x)|^{p_2}d\mu(x)}  \right].
\end{eqnarray}
\end{prop}
\begin{demo}
Par hypothèse, le nombre $(p_2-p_1)(p_1 p_2-p_1 -p_2)$ est positif. Or on peut le factoriser  
\bema 
(p_2-p_1)(p_1 p_2-p_1 -p_2)=p_2^2(p_1-1) -p_1^2(p_2-1)=(p_2-1)(p_1-1)(p_2 q_2-p_1 q_1).
\enma 
Par conséquent on a $p_1 q_1 \leq p_2 q_2$ et l'on peut alors choisir un réel $r \in [1+\frac{p_1}{q_2},1+\frac{p_2}{q_1}]$. En particulier, on a 
\begin{equation}\label{hahn-b}
p_1 \leq (r-1)q_2\leq (r-1)q_1 \leq p_2.
\end{equation}
On peut maintenant considérer la fonction $\psi$ définie par 
\bema 
\forall x\in X \quad \psi(x)=\left\{\begin{array}{cc} \Frac{1}{\int_X |\phi(x')|^r d\mu(x')} \times \frac{|\phi(x)|^{r}}{\phi(x)} & \mbox{ si } \phi(x)\neq 0, \\ 0 & \mbox{ si }  \phi(x)=0. \end{array} \right.
\enma
Il est clair que la fonction $\phi \psi$ est positive est d'intégrale égale à $1$.
D'après \eqref{cvxit}, la fonction $\phi$ appartient à $L^p(X)$ pour tout $p\in [p_1,p_2]$.
Les inégalités \eqref{hahn-b} impliquent donc que $\psi$ appartient à $L^{q_2}(X)\cap L^{q_1}(X)$. 
Grâce à \eqref{cvxit} et \eqref{defa-mino}, il vient 
\begin{equation}\label{hanh-casp}   \left\Vert \phi\right\Vert_{L^{p_1}(X)} \left\Vert \psi\right\Vert_{L^{q_1}(X)} = \Frac{\left\Vert \phi\right\Vert_{L^{p_1}(X)}\left\Vert \phi\right\Vert_{L^{(r-1)q_1}(X)}^{r-1} }{\left\Vert \phi\right\Vert_{L^r(X)}^r} \leq  Q(\phi,[p_1,p_2])^r \left\Vert \phi\right\Vert_{L^{p_1}(X)}^{\alpha}  \left\Vert \phi\right\Vert_{L^{p_2}(X)}^{\beta} , \end{equation}
 où l'on a noté
 \bemar
\alpha & :=& 1+(r-1)\theta_1(p_1,(r-1)q_1,p_2)-r\theta_1(p_1,r,p_2),\\
\beta & :=& (r-1)\theta_2(p_1,(r-1)q_1,p_2)-r\theta_2(p_1,r,p_2). \enmar
L'égalité $\theta_1+\theta_2=1$ implique que l'on a $\alpha+\beta=0$. 
En fait, il s'avère que $\alpha=0$. Cela peut se voir par calcul, mais puisque $\alpha$ ne dépend que de $p_1,p_2,r$, il suffit de traiter le cas particulier $X=\R$ muni de la mesure de Lebesgue.
Pour tout $t>0$, si l'on pose $\phi_t=\pun_{[0,t]}$ alors on a 
$\psi_t=\frac{1}{t}\pun_{[0,t]}$ et 
$\left\Vert \phi_t\right\Vert_{L^{p}} =t^{\frac{1}{p}}$  pour tout $p\in [1,+\infty[$ et donc 
\bema  \left\Vert \phi_t \right\Vert_{L^{p_1}(X)} \left\Vert \psi_t\right\Vert_{L^{q_1}(X)}=t^{\frac{1}{p_1}+\frac{1}{q_1}-1}=1     \quad    \mbox{et} \quad    Q(\phi_t,[p_1,p_2])=1    .\enma
L'inégalité \eqref{hanh-casp} devient
\bema \forall t>0 \qquad 1\leq t^{\alpha\left( \frac{1}{p_1}-\frac{1}{p_2}\right)}. \enma
Cela force les égalités $\alpha=\beta=0$.
Un argument similaire permet d'estimer $\left\Vert \phi \right\Vert_{L^{p_2}(X)}
\left\Vert \psi \right\Vert_{L^{q_2}(X)}$.

Passons aux estimations de $|\psi|^{q}$ avec $q\in \{q_1,q_2\}$. 
On se permet de noter $\theta_1=\theta_1(p_1,(r-1)q,p_2)$
et $\theta_2=\theta_2(p_1,(r-1)q,p_2)$. On obtient alors pour tout $x\in X$
\bema
\frac{|\psi(x)|^q }{\left\Vert \psi \right\Vert_{L^q(X)}^q}=\Frac{|\phi(x)|^{(r-1)q}}{ \left\Vert \phi\right\Vert_{L^{(r-1)q}(X)}^{(r-1)q} } \leq  Q(\phi,[p_1,p_2])^{(r-1)q}\Frac{|\phi(x)|^{(r-1)q}}{\left\Vert \phi\right\Vert_{L^{p_1}(X)}^{(r-1)q\theta_1}
\left\Vert \phi\right\Vert_{L^{p_2}(X)}^{(r-1)q\theta_2}} .\enma
L'égalité $\frac{(r-1)q\theta_1}{p_1}+\frac{(r-1)q\theta_2}{p_2}=1$ et la convexité de la fonction exponentielle nous donnent la conclusion
\bemar \Frac{|\psi(x)|^q }{\left\Vert \psi \right\Vert_{L^q(X)}^q} &  \leq & Q(\phi,[p_1,p_2])^{(r-1)q} \left( \Frac{|\phi(x)|^{p_1}}{\left\Vert \phi\right\Vert_{L^{p_1}(X)}^{p_1}}\right)^{\frac{(r-1)q\theta_1}{p_1}}\left( \Frac{|\phi(x)|^{p_2}}{\left\Vert \phi\right\Vert_{L^{p_2}(X)}^{p_2}}\right)^{\frac{(r-1)q\theta_2}{p_2}} \\
& \leq & Q(\phi,[p_1,p_2])^{(r-1)q} \left( \Frac{(r-1)q\theta_1}{p_1} \Frac{|\phi(x)|^{p_1}}{\left\Vert \phi\right\Vert_{L^{p_1}(X)}^{p_1}}+  \frac{(r-1)q\theta_2}{p_2}\Frac{|\phi(x)|^{p_2}}{\left\Vert \phi\right\Vert_{L^{p_2}(X)}^{p_2}}\right) \\
& \leq & Q(\phi,[p_1,p_2])^{(r-1)q} \left(  \Frac{|\phi(x)|^{p_1}}{\left\Vert \phi\right\Vert_{L^{p_1}(X)}^{p_1}}+ \Frac{|\phi(x)|^{p_2}}{\left\Vert \phi \right\Vert_{L^{p_2}(X)}^{p_2}}\right).
\enmar
\end{demo}

Passons maintenant à la condition nécessaire d'interpolation des espaces $\pl^p(X,\oplus E_n)$. 
On rappelle que les espaces $\pl^{p}(X,\oplus E_n)$ et leurs espaces interpolés sont inclus dans $\prod_{n\in \N} E_n$.

\begin{prop}\label{neces-interpo}
Fixons des réels $p_1<p<p_2$ appartenant à $[1,+\infty[$ et supposons que l'on a l'égalité d'espaces vectoriels 
\bema \pl^p(X,\oplus E_n)=[\pl^{p_1}(X,\oplus E_n),\pl^{p_2}(X,\oplus E_n)]_{\theta_1(p_1,p,p_2)} \enma
et que les normes des deux précédents espaces sont équivalentes. Alors 
\begin{equation}\label{interpo-nec}
\sup\limits_{n\in \N} Q\left( \sqrt{e(n,\cdot)},[p_1,p_2] \right)<+\infty.
\end{equation}
La même conclusion est valide en remplaçant la méthode d'interpolation complexe $[\cdot,\cdot]_{\theta_1(p_1,p,p_2)}$ par 
la méthode d'interpolation réelle $[\cdot,\cdot]_{\theta_1(p_1,p,p_2),p}$.
\end{prop}
\begin{demo}
On considère pour tout entier $k\in \N$ le ``projecteur sur $E_k$'' défini par 
\bemar
\Lambda_k : \Prod_{n\in \N} E_n & \rightarrow & L^2(X) \\
(u_n)_{n\in \N}& \mapsto & u_k .
\enmar
Le calcul de la norme d'opérateur de $\Lambda_k:\pl^p(X,\oplus E_n)\rightarrow L^2(X)$ est immédiat
\bema
\left\Vert \Lambda_k\right\Vert_{\pl^{p}(X,\oplus E_n) \rightarrow L^2(X)}=
\sup\limits_{(u_n)\neq 0} \Frac{\left\Vert u_k\right\Vert_{L^2(X)}}{\left\Vert
 \sqrt{\Sum_{n\in \N} \left\Vert u_n \right\Vert_{L^2(X)}^2 \frac{e(n,\cdot)}{d_n}}
\right\Vert_{L^p(X)}}=\frac{\sqrt{d_k}}{\left\Vert\sqrt{e(k,\cdot)}\right\Vert_{L^p(X)}}.
\enma 
Par interpolation complexe, il existe une constante $K>0$ telle que pour tout $k\in \N$ on a
\bema
\left\Vert \Lambda_k\right\Vert_{\pl^{p}(X,\oplus E_n) \rightarrow L^2(X)}  \leq  
K \left\Vert \Lambda_k\right\Vert_{\pl^{p_1}(X,\oplus E_n) \rightarrow L^2(X)}^{\theta_1(p_1,p,p_2)}
\left\Vert \Lambda_k\right\Vert_{\pl^{p_2}(X,\oplus E_n) \rightarrow L^2(X)}^{\theta_2(p_1,p,p_2)} .\enma
\`A noter que $K=1$ ne convient pas forcément car la norme de l'espace $\pl^p(X,\oplus \C \phi_n)$ n'est a priori pas identique à celle de l'espace interpolé.
On obtient donc
\bema
\Frac{\sqrt{d_k}}{\left\Vert\sqrt{e(k,\cdot)}\right\Vert_{L^p(X)}}  \leq   K 
\Frac{\sqrt{d_k}}{\left\Vert\sqrt{e(k,\cdot)}\right\Vert_{L^{p_1}(X)}^{\theta_1(p_1,p,p_2) }\left\Vert\sqrt{e(k,\cdot)}\right\Vert_{L^{p_2}(X)}^{\theta_2(p_1,p,p_2) }}.
\enma 
Le lemme \ref{majo-Q} implique \eqref{interpo-nec}.
Le raisonnement est similaire pour la méthode d'interpolation réelle.
\end{demo}

\section{Preuves des théorèmes \ref{theo-interpo} et \ref{theo-dual}, partie VI : bornitude de $S_p R_{p,\psi}$}\label{partVI}

Nous avons maintenant les moyens d'achever 
la preuve, développée dans la partie \ref{partI}, du théorème \ref{theo-interpo}.
Sous les hypothèses \eqref{loren1} et \eqref{hypo-interpo}, il s'agit de justifier l'existence d'une suite $(\psi_n)_{n\in \N}$ de $L^{\frac{p_2}{p_2-1}}(X)\cap L^{\frac{p_1}{p_1-1}}(X)$ telle que 
\begin{equation}\label{proj}
\forall n\in \N \quad \Int_{X} \frac{\sqrt{e(n,x)}}{\sqrt{d_n}} \psi_n(x) d\mu(x)=1,
\end{equation}
et que l'opérateur $S_{p} R_{p,\psi}$ soit borné sur $L_{x}^p(X,\oplus E_n)$.
L'hypothèse \eqref{hypo-interpo} et le lemme \ref{majo-Q} nous apprennent que la suite des défauts d'interpolation $(Q(\sqrt{e(n,\cdot)},[p_1,p_2]))_{n\in \N}$ est bornée.
Par suite, la proposition \ref{hahn} nous assure l'existence d'une constante $K>0$ et d'une suite de fonctions $(\psi_n)_{n\in \N}$ de $L^{\frac{p_2}{p_2-1}}(X)\cap L^{\frac{p_1}{p_1-1}}(X)$ qui vérifient \eqref{proj} et 
\begin{eqnarray} \nonumber
\left\Vert \frac{\sqrt{e(n,\cdot)}}{\sqrt{d_n}} \right\Vert_{L^{p_1}(X)} \left\Vert \psi_n\right\Vert_{L^{q_1}(X)} & \leq & K,\\ \nonumber
\left\Vert \frac{\sqrt{e(n,\cdot)}}{\sqrt{d_n}} \right\Vert_{L^{p_2}(X)} \left\Vert \psi_n\right\Vert_{L^{q_2}(X)} & \leq & K,
\end{eqnarray}
ainsi que les inégalités suivantes pour tout $x\in X$
\begin{eqnarray}\label{majo-psi}
\Frac{|\psi_n(x)|^{q_1}}{ \left\Vert \psi_n \right\Vert_{L^q(X)}^{q_1}} & \leq & K \left(\Frac{\sqrt{e(n,x)}^{p_1}}{\left\Vert \sqrt{e(n,\cdot)}\right\Vert_{L^{p_1}(X)}^{p_1}}+ \Frac{\sqrt{e(n,x)}^{p_2}}{\left\Vert \sqrt{e(n,\cdot)}\right\Vert_{L^{p_2}(X)}^{p_2}}\right), \\ \nonumber
\Frac{|\psi_n(x)|^{q_2}}{ \left\Vert \psi_n \right\Vert_{L^q(X)}^{q_2}} & \leq & K \left(\Frac{\sqrt{e(n,x)}^{p_1}}{\left\Vert \sqrt{e(n,\cdot)}\right\Vert_{L^{p_1}(X)}^{p_1}}+ \Frac{\sqrt{e(n,x)}^{p_2}}{\left\Vert \sqrt{e(n,\cdot)}\right\Vert_{L^{p_2}(X)}^{p_2}}\right).
\end{eqnarray}
On va appliquer le lemme \ref{schur} avec les noyaux $K_n$ définis par 
\bema
\forall n\in \N \quad \forall (x,x')\in X^2 \quad K_n(x,x'):=\frac{\sqrt{e(n,x)}}{\sqrt{d_n}}\psi_n(x').
\enma
Les propriétés \eqref{schur1-g} et \eqref{schur2-g} se traitent comme pour la bornitude de $S_p R_p$. Par exemple, l'inclusion continue $L^{q_1}(X)\subset L^{q_1,\infty}(X)$ nous permet d'écrire pour tout $x\in X$
\bemar
\left\Vert K_n(x,x')\right\Vert_{L_{x'}^{q_1,\infty}(X)} & \leq & \Frac{\sqrt{e(n,x)}}{\sqrt{d_n}} \left\Vert \psi_n \right\Vert_{L^{q_1}(X)} \\
& \leq & K\Frac{\sqrt{e(n,x)}}{\left\Vert \sqrt{e(n,\cdot)}\right\Vert_{L^{p_1}(X)}}.
\enmar
L'hypothèse \eqref{loren1} implique alors \eqref{schur1-g}.
On montre de même \eqref{schur2-g}.
L'intérêt de la proposition \ref{hahn} apparaît pour démontrer \eqref{schur3-g} et \eqref{schur4-g}.
En effet, on a 
\bemar
\forall x'\in X \qquad \left\Vert K_n(x,x')\right\Vert_{L_{x}^{p_1,\infty}(X)} & \leq &  |\psi_n(x')| \left\Vert \frac{\sqrt{e(n,\cdot)}}{\sqrt{d_n}}\right\Vert_{L^{p_1}(X)} \\
& \leq & K \Frac{|\psi_n(x')|}{\left\Vert \psi_n \right\Vert_{L^{q_1}(X)} }.
\enmar
D'après \eqref{loren1} et \eqref{majo-psi}, on obtient
\bema
\left(\sup\limits_{n\in \N} \left\Vert K_n(x,x')\right\Vert_{L_{x}^{p_1,\infty}(X)}\right)^{q_1} \in L_x^{1,\infty}(X),
 \enma
c'est-à-dire 
\bema
\sup\limits_{n\in \N} \left\Vert K_n(x,x')\right\Vert_{L_{x}^{p_1,\infty}(X)} \in L_x^{q_1,\infty}(X).
 \enma
On a ainsi obtenu \eqref{schur3-g}. Et \eqref{schur4-g} se traite de même.
On a validé toutes les hypothèses du lemme \ref{schur} et l'on peut donc conclure que $S_{p} R_{p,\psi}$ est borné sur $L_x^p(X,\oplus E_n)$.
Cela achève la preuve du théorème \ref{theo-interpo} d'interpolation.

\section{Preuve du théorème \ref{pz-linf}, théorème de Paley-Zygmund sur une variété}

Commençons par un lemme qui fait intervenir la classique perte multiplicative $\sqrt{\ln(\alpha)}$ dans l'estimation de l'espérance d'un maximum de $\alpha$ variables aléatoires.
\begin{lemm}\label{zoll-linf}
Soit $(\alpha,N)\in  \N^\star\times \N$ et considérons $\alpha (N+1)$ matrices 
\bema
b_{(0,\beta)} \in \mathcal{M}_{d_0}(\C),\dots, b_{(N,\beta)}\in \mathcal{M}_{d_N}(\C),
\enma
où $\beta$ parcourt l'ensemble des entiers compris entre $1$ et $\alpha$.
Alors on a 
\bema \E\left[ \sup\limits_{1\leq \beta\leq \alpha} \left\vert\sum_{n=0}^N \sqrt{d_n}\tr(\mathcal{E}_{n} b_{(n,\beta)})\right\vert \right]\leq C\sqrt{\ln(\alpha+1) \sup\limits_{1\leq \beta\leq \alpha} \left(\sum_{n=0}^N \left\Vert b_{(n,\beta)} \right\Vert_{HS}^2\right) }.\enma
La même conclusion est valide en remplaçant les matrices orthogogonales aléatoires $\mathcal{E}_n$ par les matrices unitaires aléatoires $W_n$.
\end{lemm}
\begin{demo}
Pour tout réel $p\geq 1$, on peut écrire grâce à l'inégalité de Hölder et aux inégalités de Kahane-Khintchine-Marcus-Pisier \eqref{KKMP} et \eqref{KKMP-R} :
\bemar
\E\left[ \sup\limits_{1\leq \beta\leq \alpha} \left\vert\Sum_{n=0}^N \sqrt{d_n}\tr(\mathcal{E}_n b_{(n,\beta)})\right\vert \right] & \leq &
\E\left[ \left( \Sum_{\beta=1}^{\alpha} \left\vert\Sum_{n=0}^N \sqrt{d_n}\tr(\mathcal{E}_n  b_{(n,\beta)} ) \right\vert^p\right)^{\frac{1}{p}} \right] \\[3mm]
& \leq & 
\E\left[ \Sum_{\beta=1}^{\alpha} \left\vert\Sum_{n=0}^N \sqrt{d_n}\tr(\mathcal{E}_n b_{(n,\beta)})  \right\vert^p \right]^{\frac{1}{p}} \\[3mm]
& \leq &  \left(\Sum_{\beta=1}^{\alpha} \E\left[\left\vert\Sum_{n=0}^N \sqrt{d_n}\tr(\mathcal{E}_n  b_{(n,\beta)})  \right\vert^p \right]\right)^{\frac{1}{p}}         \\
&  \leq & \alpha^{\frac{1}{p}} \sup\limits_{1\leq \beta \leq \alpha} \E\left[\left\vert\Sum_{n=0}^N \sqrt{d_n}\tr(\mathcal{E}_n b_{(n,\beta)})  \right\vert^p \right]^{\frac{1}{p}}   \\[3mm]
 & \leq &  C\alpha^{\frac{1}{p}}\sqrt{p} \sup\limits_{1\leq \beta \leq \alpha} \E\left[\left\vert\Sum_{n=0}^N \sqrt{d_n}\tr(\mathcal{E}_n b_{(n,\beta)})  \right\vert^2 \right]^{\frac{1}{2}}   \\[3mm]
 & \leq &  C\alpha^{\frac{1}{p}}\sqrt{p} \sup\limits_{1\leq \beta \leq \alpha} \sqrt{\Sum_{n=0}^N \left\Vert b_{(n,\beta)} \right\Vert_{HS}^2}.
\enmar
On obtient la conclusion voulue en choisissant $p=\ln(2+\alpha)\geq 1$.
Une démarche similaire est valide pour les matrices aléatoires $W_n:\Omega\rightarrow U_{d_n}(\C)$.
\end{demo}

On aura besoin d'une généralisation de l'inégalité de Salem-Zygmund sur une variété riemannienne compacte sans bord.
Ce type d'inégalité est généralement présenté sous forme probabiliste (\cite[Chapter 6]{kahane}, \cite[Theorem 2.5]{randomh} ou \cite[Théorème 2]{burq-lebeau}) mais nous avons déjà souligné que les versions intégrales sont plus efficaces dans notre contexte.
En outre, nous avons besoin d'une forme matricielle de l'inégalité de Salem-Zygmund (voir la raison plus loin à la ligne \eqref{utili-horm}).
La preuve repose sur une idée classique : on utilise une inégalité de Bernstein et un maillage de la variété compacte $X$.

\begin{prop}
Fixons un entier $N\geq 2$ et considérons $N$ matrices 
\bema
b_0\in \mathcal{M}_{d_0}(L^2(X)), \dots ,b_N\in \mathcal{M}_{d_N}(L^2(X)),
\enma
dont les coefficients sont des fonctions spectralement localisées au sens suivant :
\bema 
\forall n\in [0,N] \quad \forall (i,j)\in [1,d_n]^2 \qquad b_{n,i,j}\in \bigoplus_{k=0}^N E_k.
\enma
Alors on a une majoration à la Salem-Zygmund :
\begin{equation}\label{salem}
\E\left[ \left\Vert 
\sum_{n=0}^N \sqrt{d_n}\tr(\mathcal{E}_n b_n(x)) \right\Vert_{L_x^\infty(X)} \right] \leq C(X)\sqrt{\ln(N)} \left\Vert \sqrt{\sum_{n=0}^N \left\Vert b_n(x)\right\Vert_{HS}^2} \right\Vert_{L_x^\infty(X)}.
\end{equation}
\end{prop}

\begin{demo}
Commençons par considérer une fonction spectralement localisée 
\bema
u:=\sum_{n=0}^N \sum_{i=1}^{d_n} \langle u,\phi_{n,i}\rangle \phi_{n,i}.
\enma
Pour tout $(x,y)\in X^2$, on peut écrire
\bema
|u(x)-u(y)|\leq \sum_{n=0}^N \sum_{i=1}^{d_n} \left\Vert u\right\Vert_{L^\infty(X)} \left\Vert \phi_{n,i}\right\Vert_{L^1(X)} |\phi_{n,i}(x)-\phi_{n,i}(y)|. 
\enma
D'après l'inégalité \eqref{sogge} de \horm et la majoration $d_n\leq C(X) (n+1)^{d-1}$ (voir le lemme \ref{en-var}), on obtient 
\bema
|u(x)-u(y)|\leq c_1(X) N^{c_2(X)} \left\Vert u \right\Vert_{L^\infty(X)} \sum_{n=0}^N \sum_{i=1}^{d_n}  |\phi_{n,i}(x)-\phi_{n,i}(y)|. 
\enma
Pour un réel fixé $\tau>1+\frac{d}{2}$, l'injection de Sobolev $H^{\tau}(X)\subset W^{1,\infty}(X)$ nous permet de contrôler la norme de Lipschitz des fonctions $\phi_{n,i}$. Quitte à augmenter $c_2(X)$, cela donne une inégalité de Bernstein :
\bema
|u(x)-u(y)|\leq c_1(X) N^{c_2(X)} \left\Vert u\right\Vert_{L^\infty(X)}\mbox{dist}(x,y).
\enma
Par compacité de $X$, il existe un ensemble fini $\{x_1,\dots,x_\alpha\}\subset X$ maximal pour la propriété suivante : les boules ouvertes de centres $x_1,\dots,x_\alpha$ et de rayon 
$\frac{1}{4c_1(X) N^{c_2(X)}}$ sont disjointes deux à deux.
En examinant la somme du volume de ces boules, on vérifie aisément l'inégalité
\begin{equation}\label{majo-alpha} \alpha\leq C(X) N^{c_2(X)\dim(X)}. \end{equation}
L'inégalité de Bernstein et la maximalité de $\{x_1,\dots,x_\alpha\}$ assurent, pour tout $x\in X$, l'existence d'un point $x_\beta$ tel que l'on a 
\bema
|u(x)-u(x_\beta)| \leq  c_1(X)N^{c_2(X)} \left\Vert u\right\Vert_{L^\infty(X)} \mbox{dist}(x,x_\beta)\leq \frac{1}{2} \left\Vert u\right\Vert_{L^\infty(X)}.
\enma
On déduit par conséquent que l'on a 
\bema
\left\Vert u \right\Vert_{L^\infty(X)} \leq 2 \sup\limits_{1\leq\beta \leq \alpha} |u(x_\beta)|.
\enma
Appliquant l'inégalité précédente à une fonction aléatoire à la place de $u$, on obtient
\bema
\E\left[\left\Vert \sum_{n=0}^N \sqrt{d_n} \tr(\mathcal{E}_n b_n(x)) \right\Vert_{L_x^\infty(X)} \right]\leq 2 
\E\left[\sup\limits_{1\leq \beta \leq \alpha} \left\vert \sum_{n=0}^N \sqrt{d_n} \tr(\mathcal{E}_n b_n(x_\beta)) \right\vert \right].
\enma
On conclut avec le lemme \ref{zoll-linf} et l'estimation \eqref{majo-alpha}.
\end{demo}

\begin{rema}
Dans la théorie $L^p$, le théorème de Fubini nous permettait de travailler sur un espace mesuré $\sigma$-fini $X$. Par contre, la démonstration précédente suggère que la théorie $L^\infty$ nécessite que les espaces mesurés soient en fait des espaces métriques (munis de mesures boréliennes).
Cela n'est guère surprenant puisque la continuité est une notion métrique!
\end{rema}
\begin{rema}
Nous avons choisi d'appeler \eqref{salem} une inégalité de Salem-Zygmund car elle implique que pour toute fonction $u\in L^2(X)$ et pour tour entier $N\geq 2$ l'on a
\bema
\E\left[ \left\Vert 
\sum_{n=0}^N \sum_{i=1}^{d_n} \sum_{j=1}^{d_n} \mathcal{E}_{n,i,j} \langle u,\phi_{n,j}\rangle \phi_{n,i}(x)\right\Vert_{L_x^\infty(X)}\right]\leq C(X) \sqrt{\ln(N)} \left\Vert u \right\Vert_{L^2(X)} .\enma
En effet, comme dans la théorie $L^p$, il s'agit d'appliquer \eqref{salem} et 
le lemme \ref{en-var} avec les matrices $b_n$ dont les coefficients sont donnés par 
\bema
b_{n,i,j}(x)=\frac{1}{\sqrt{d_n}} \langle u,\phi_{n,i}\rangle \phi_{n,j}(x),\qquad \left\Vert b_n(x)\right\Vert_{HS}^2 = \frac{e(n,x)}{d_n} \left\Vert \Pi_n(u)\right\Vert_{L^2(X)}^2 .
\enma
Dans le cas $X=\mathbb{T}$, il s'agit de la classique inégalité de Salem-Zygmund.
\end{rema}
On peut maintenant démontrer le théorème \ref{pz-linf}.
Comme déjà utilisé, la convergence dans $L^p(\Omega,L^\infty(X))$ des séries aléatoires considérées implique leur convergence presque sûre dans $L^\infty(X)$ et il nous suffit de prouver l'estimation \eqref{zoll-pal}.
Il sera commode de noter $b_n\in \mathcal{M}_{d_n}(L^\infty(X))$ la matrice dont l'élément $(i,j)$ est donné par 
\bema
\forall x\in X \quad b_{n,i,j}(x)=\frac{1}{\sqrt{d_n}}\langle u,\phi_{n,i}\rangle \phi_{n,j}(x).
\enma 
Cela permet de reformuler \eqref{zoll-pal} en 
\bema
\forall N\in \N \quad \E\left[\left\Vert \sum_{n=0}^N \sqrt{d_n} \tr(M_n b_n(x))\right\Vert_{L_x^\infty(X)}^p\right]\leq C(X,p,s)\left(\sup\limits_{n\in \N} \E\left[\left\Vert M_n\right\Vert_{op}^p \right] \right)\left\Vert u\right\Vert_{H^s(X)}^p.
\enma
Comme déjà effectué dans la théorie $L^p$, on peut remplacer $M_n$ par $\mathcal{E}_n M_n$ où les matrices aléatoires $\mathcal{E}_n:\Omega\rightarrow O_{d_n}(\R)$ sont indépendantes entre elles et vis-à-vis des matrices aléatoires $M_n$. Les inégalités \eqref{KKMP} de Kahane-Khintchine-Marcus-Pisier nous amènent à l'estimation 
\bemar
\quad \E\left[\left\Vert \Sum_{n=0}^N \sqrt{d_n} \tr(M_n b_n(x))\right\Vert_{L_x^\infty(X)}^p\right] & =& \quad \E_{\omega'}\E_{\omega}\left[\left\Vert \Sum_{n=0}^N \sqrt{d_n} \tr(\mathcal{E}_n(\omega)M_n(\omega') b_n(x))\right\Vert_{L_x^\infty(X)}^p\right] \\
& \leq & C(p) \E_{\omega'}\left[\E_{\omega}\left[\left\Vert \Sum_{n=0}^N \sqrt{d_n}\tr(\mathcal{E}_n(\omega)M_n(\omega') b_n(x))\right\Vert_{L_x^\infty(X)}\right]^p\right].
\enmar
On va momentanément geler $\omega'$ et s'atteler à démontrer l'inégalité suivante : 
\begin{equation}\label{pz-p1}  \E_{\omega}\left[\left\Vert \Sum_{n=0}^N \sqrt{d_n}\tr(\mathcal{E}_n(\omega)M_n(\omega') b_n(x))\right\Vert_{L_x^\infty(X)}\right] \leq C(X,s) \sqrt{\sum_{n\geq 0} \left\Vert M_n(\omega')\right\Vert_{op}^2  (1+n)^{2s} \left\Vert \Pi_n(u)\right\Vert_{L^2}^2 }    .    \end{equation}
L'inégalité triangulaire dans $L^{\frac{p}{2}}(\Omega)$ (car on a $p\geq 2$) et \eqref{pz-p1} impliqueront alors la conclusion : 
\bemar
\E\left[\left\Vert \Sum_{n=0}^N \sqrt{d_n} \tr(M_n b_n)\right\Vert_{L^\infty(X)}^p \right] & \leq & C(X,p,s)\E_{\omega'} \left[   \left( \Sum_{n=0}^N \left\Vert M_n(\omega')\right\Vert_{op}^2  (1+n)^{2s} \left\Vert \Pi_n(u)\right\Vert_{L^2(X)}^2 \right)^{\frac{p}{2}}   \right]  \\
& \leq &  C(X,p,s) \left(\Sum_{n=0}^N \E\left[ \left\Vert M_n \right\Vert_{op}^p\right]^{\frac{2}{p}}  (1+n)^{2s} \left\Vert \Pi_n(u)\right\Vert_{L^2(X)}^2
 \right)^{\frac{p}{2}} \\
 & \leq & C(X,p,s)\left(\sup\limits_{n\geq 0} \E\left[\left\Vert M_n \right\Vert_{op}^p\right] \right) \left\Vert u \right\Vert_{H^s(X)}^p.
\enmar
Il s'agit précisément de l'analogue multidimensionnel de l'argument que nous avons esquissé pour prouver l'implication \eqref{univ-linf}. 
Démontrons maintenant \eqref{pz-p1}.
La première chose à remarquer est la croissance en $N$ :
\begin{equation}\label{zoll-crois}
\E_{\omega}\left[ \left\Vert 
\sum_{n=0}^N \sqrt{d_n} \tr\left(\mathcal{E}_n(\omega) M_n(\omega')b_n(x)\right)\right\Vert_{L_x^\infty(X)}\right]\leq 
\E_{\omega}\left[ \left\Vert 
\sum_{n=0}^{N+1} \sqrt{d_n}\tr\left( \mathcal{E}_n(\omega) M_n(\omega') b_n(x)\right) \right\Vert_{L_x^\infty(X)} \right] .\end{equation}
En effet, cela découle de l'indépendance des matrices aléatoires $\mathcal{E}_0,\dots,\mathcal{E}_{N+1}$ et de l'inégalité de Jensen en intégrant selon les $d_{N+1}^2$ variables de la matrice aléatoire $\mathcal{E}_{N+1}$ (on utilisera l'égalité $\E[\mathcal{E}_{N+1,i,j}]=0$).
L'inégalité \eqref{zoll-crois} nous assure qu'il nous suffit de prouver \eqref{pz-p1} pour une infinité d'entiers $N\in \N$.
On choisit alors $N$ de la forme $2^{2^k}$ avec $k\in \N$ et l'on invoque l'inégalité de Salem-Zygmund \eqref{salem} :
\bema
\E_\omega \left[\left\Vert  \sum_{2^{2^k}<n\leq 2^{2^{k+1}}}\sqrt{d_n} \tr(\mathcal{E}_n(\omega) M_n(\omega')b_n(x)) \right\Vert_{L_x^\infty(X)} \right]\leq C(X) 2^{\frac{k}{2}} \left\Vert \sqrt{\sum_{2^{2^k}<n\leq 2^{2^{k+1}}} \left\Vert M_n(\omega') b_n(x)\right\Vert_{HS}^2 } \right\Vert_{L_x^\infty(X)}.
\enma
Pour tout $x\in X$, le lemme \ref{en-var} nous permet de majorer :
\begin{eqnarray} \label{utili-horm}
\left\Vert M_n(\omega') b_n(x)\right\Vert_{HS}^2& \leq & \left\Vert M_n(\omega')\right\Vert_{op}^2 \left\Vert b_n(x)\right\Vert_{HS}^2 \\  \nonumber
  & \leq & \left\Vert M_n(\omega')\right\Vert_{op}^2\Frac{e(n,x)}{d_n} \left\Vert \Pi_n(u)\right\Vert_{L^2(X)}^2 \\[2mm]  \nonumber
  &  \leq &  C(X) \left\Vert M_n(\omega')\right\Vert_{op}^2\left\Vert \Pi_n(u)\right\Vert_{L^2(X)}^2   .   
\end{eqnarray}

Comme la fréquence de $\Pi_n(u)$ est de l'ordre de $n$, on peut forcer l'apparition des paquets de la norme de Sobolev $\left\Vert u\right\Vert_{H^s(X)}$ : 
\bema
2^{\frac{k}{2}} \left(\sum_{2^{2^k}<n\leq 2^{2^{k+1}}} \left\Vert M_n(\omega')b_n(x)\right\Vert_{HS}^2\right)^{\frac{1}{2}}\leq \frac{C(X)}{2^{2^k s  -\frac{k}{2}}}   \left(\sum_{2^{2^k}<n\leq 2^{2^{k+1}}} \left\Vert M_n(\omega')\right\Vert_{op}^2 n^{2s}\left\Vert \Pi_n(u)\right\Vert_{L^2(X)}^2\right)^{\frac{1}{2}}  .
\enma
La condition $s>0$ nous permet d'utiliser l'inégalité de Cauchy-Schwarz en sommant sur $k$ : 
\bema
\sum_{k\in \N} \E\left[\left\Vert  \sum_{2^{2^k}<n\leq 2^{2^{k+1}}}\sqrt{d_n} \tr(\mathcal{E}_{n}(\omega)M_n(\omega')b_n(x))\right\Vert_{L_x^\infty(X)} \right]\leq C(X,s)\sqrt{\sum_{n>2} \left\Vert M_n(\omega')\right\Vert_{op}^2 n^{2s}\left\Vert \Pi_n(u)\right\Vert_{L^2(X)}^2}.
\enma
Il manque a priori les entiers $n=0,1$ et $2$. Mais il est clair que l'on a toujours 
\bemar
\forall n\in \{0,1,2\} \qquad \E_{\omega}\left[ \left\Vert \tr (\mathcal{E}_n(\omega) M_n(\omega') b_n(x) \right\Vert_{L_x^\infty(X)}\right] &\leq &   C(n)\left\Vert M_n(\omega')\right\Vert_{op} \sup\limits_{x\in X}\left\Vert b_n(x)\right\Vert_{HS}    \\
& \leq &  C(X,n)\left\Vert M_n(\omega')\right\Vert_{op} \left\Vert \Pi_n(u)\right\Vert_{L^2(X)}.
\enmar
L'inégalité triangulaire implique \eqref{pz-p1} sans difficulté.
La preuve du théorème \ref{pz-linf} est complète.

\section{Preuve du théorème \ref{edp-zoll}, équation cubique des ondes}

Les propositions 2.1 et 2.2 de \cite{burq-tz-jems2011} montrent en fait le résultat suivant : 
\begin{theo}\label{zoll-tz}[Burq-Tzvetkov] Considérons un réel $s>0$, une variété riemannienne compacte $X$ sans bord de dimension $3$ et une fonction
\begin{eqnarray}\nonumber
\Omega & \rightarrow & H^s(X)\times H^{s-1}(X),\\
\omega & \mapsto & (v_0^\omega,v_1^\omega), \label{zoll-random}
\end{eqnarray}
qui vérifie la propriété suivante pour presque tout $\omega\in \Omega$ :
\begin{equation}\label{flot-lin}
\cos(t\sqrt{-\Delta})v_0^\omega(x)+\frac{\sin(t\sqrt{-\Delta})}{\sqrt{-\Delta}}v_1^\omega(x) \in L_{t,loc}^3(\R,L_x^6(X))\cap L_{t,loc}^1(\R,L_x^\infty(X))       .
\end{equation}
Alors pour presque tout $\omega\in \Omega$, l'équation cubique des ondes sur $X$
\bema
(\partial_t^2-\Delta)v+v^3=0,\quad v(0,\cdot)=v_0^\omega,\quad \dot{v}(0,\cdot)=v_1^\omega,
\enma
admet une unique solution globale $v$ qui vérifie 
\bema v(t)-\cos(t\sqrt{-\Delta})v_0^\omega-\frac{\sin(t\sqrt{-\Delta})}{\sqrt{-\Delta}}v_1^\omega \in \CC_t^0(\R,H^1(X))\cap \CC_t^1(\R,L^2(X)).
\enma 
\end{theo}
La dimension de la variété intervient dans le résultat d'existence locale \cite[Proposition 2.1]{burq-tz-jems2011}, grâce à l'injection de Sobolev $H^1(X)\subset L^6(X)$, afin de gérer la non-linéarité $v^3$.
La preuve du théorème \ref{zoll-tz} est \textbf{déterministe} et l'intérêt des probabilités intervient seulement lorsqu'il s'agit de trouver un procédé de randomisation \eqref{zoll-random} qui vérifie presque sûrement la propriété \eqref{flot-lin}.
Cette dernière est démontrée dans \cite[Appendix A]{burq-tz-jems2011} pour la randomisation unidimensionnelle du tore $\mathbb{T}^3$.
Dans \cite{suzzo-S3}, la preuve a été généralisée par de Suzzoni pour la variété $X=\S^3$ qui s'avère admettre, comme le tore, une base de fonctions propres uniformément bornées dans tous les espaces $L^p(X)$ (voir \cite{burq-lebeau}).

Discutons maintenant de la vérification de 
l'appartenance presque sûre dans $L_{t,loc}^1(\R,L_x^\infty(X))$ du flot linéaire \eqref{flot-lin} par Burq et Tzvetkov. 
Cela est démontré grâce à l'injection de Sobolev $W^{s,p}(X)\subset L^\infty(X)$ (valide pour tout $p>\frac{3}{s}$).
Cela n'est pas gênant car la randomisation \eqref{zoll-random} utilisée dans \cite{burq-tz-jems2011} fait intervenir des variables aléatoires sous-gaussiennes donc qui ont des moments finis de tout ordre (en particulier les moments d'ordre $\frac{3}{s}$).
Dans le cas du tore $\mathbb{T}^3$, il est possible d'améliorer cet argument en remplaçant l'injection de Sobolev $W^{s,p}(\mathbb{T}^3)\subset L^\infty(\mathbb{T}^3)$, avec $p\gg 1$, par le théorème de Paley-Zygmund dans $L^\infty(\mathbb{T}^3)$ sous sa forme universelle, c'est-à-dire l'implication 
\eqref{univ-linf}.
Cela permet de se ramener presque sûrement à $W^{s,2}(\mathbb{T}^3)=H^s(\mathbb{T}^3)$ et donc de n'utiliser que très peu de moments des variables aléatoires en jeu.
En d'autres termes, l'idée est tout simplement de considérer le théorème de Paley-Zygmund comme une injection de Sobolev probabiliste!

En exprimant $\cos(t\sqrt{-\Delta})$ et $\sin(t\sqrt{-\Delta})$ en fonction de $e^{it\sqrt{-\Delta}}$ et $e^{-it\sqrt{-\Delta}}$, nous voyons que le théorème \ref{edp-zoll} est une conséquence du théorème \ref{zoll-tz} et des deux résultats suivants.

\begin{prop}\label{bien-surc}
Considérons une suite de matrices aléatoires $M_n:\Omega\rightarrow \mathcal{M}_{d_n}(\R)$ indépendantes, orthogonalement invariantes et vérifiant
$\sup\limits_{n\in \N} \E[\left\Vert M_n \right\Vert_{op}^2]<+\infty$.
Pour toute distribution $u\in \mathcal{D}'(X)$ et pour tout $s\in \R$, la propriété i) implique la propriété ii) ci-après : 
\begin{enumerate}[i) ]
\item $u$ appartient à $H^s(X)$,
\item la série aléatoire suivante converge dans $H^{s}(X)$ pour presque tout $\omega$ :
\bema
\sum_{n\in \N} \sum_{i=1}^{d_n} \left(\sum_{j=1}^{d_n}M_{n,i,j}(\omega) \langle u,\phi_{n,j}\rangle   \right) \phi_{n,i}.\enma
\end{enumerate}
Si l'on suppose en outre que l'on a $\inf\limits_{n\in \N} \sigma\left( \E[|M_n|]\right)>0$, alors i) et ii) sont deux assertions équivalentes.
Une conclusion similaire est valide si l'on suppose que les matrices $M_n:\Omega\rightarrow \mathcal{M}_{d_n}(\C)$ sont unitairement invariantes.
\end{prop}
\begin{demo}
\textbf{Cas $s=0$.} L'implication i) $\Rightarrow$ ii) découle de la remarque \ref{Lptr2}.
L'équivalence i) $\Leftrightarrow$ ii) découle du théorème \ref{mplp} et du théorème \ref{zoll-Lp} avec $p=2$.

\textbf{Cas $s\neq 0$.} On perturbe l'opérateur de Laplace-Beltrami afin de le rendre diagonal en posant
\bemar
\widetilde{\Delta} :  \quad \mathcal{D}'(X)& \rightarrow & \mathcal{D}'(X) \\ \nonumber
\Sum_{n,i} \alpha_{n,i} \phi_{n,i} & \mapsto & \Sum_{n,i} (n+1)^2 \alpha_{n,i} \phi_{n,i}.
\enmar
Comme les fonctions $\phi_{n,i}$ sont spectralement localisées dans $ [\kappa n,\kappa n+\kappa]$, il est clair que $\widetilde{\Delta}$ se comporte de la même manière que l'opérateur de Laplace-Beltrami : 
\bema
u\in H^{s}(X)   \quad  \Leftrightarrow  \quad \widetilde{\Delta}^{\frac{s}{2}} u \in L^2(X).
\enma
Comme $\widetilde{\Delta}$ est diagonalisé par les fonctions $\phi_{n,i}$, on se ramène aisément au cas $s=0$.
\end{demo}

\begin{prop}\label{zoll-verif}
Soit une suite de matrices aléatoires $M_n:\Omega\rightarrow \mathcal{M}_n(\R)$ indépendantes, orthogonalement invariantes et vérifiant $     \sup\limits_{n\in \N} \E\left[\left\Vert M_n \right\Vert_{op}^{3} \right]<+\infty  $.
On considère de plus un réel $s>0$ et une fonction $u\in H^s(X)$. Alors pour presque tout $\omega \in \Omega$, nous avons
\bema e^{i t\sqrt{-\Delta}} \sum_{n\in \N} \sum_{i=1}^{d_n} \left( \sum_{j=1}^{d_n} M_{n,i,j}(\omega) \langle u,\phi_{n,j}\rangle \right)\phi_{n,i}(x) \in L_{t,loc}^3(\R,L_x^\infty(X)).\enma
Un résultat similaire est vrai si les matrices aléatoires $M_n$ sont complexes et unitairement invariantes.
\end{prop}
\begin{demo}
D'après la proposition \ref{bien-surc}, on peut bien considérer la somme infinie.
L'opérateur $e^{it\sqrt{-\Delta}}$ présente l'avantage d'être unitaire et de conserver la localisation spectrale :
\bema  \forall t\in \R \qquad e^{it\sqrt{-\Delta}} E_n=E_n.       \enma
Pour tout $T\in \N^\star$ et $N\in \N$, il nous suffit de majorer 
\bema
\E\left[\left\Vert e^{i t\sqrt{-\Delta}} \Sum_{n=0}^N \sum_{i=1}^{d_n} \left( \sum_{j=1}^{d_n} M_{n,i,j}(\omega) \langle u,\phi_{n,j}\rangle \right)\phi_{n,i}(x)\right\Vert_{L_{t}^3([-T,T],L_x^\infty(X))}^3 \right]
\enma
par 
\bema C(T,X,s)\left(\sup\limits_{n\in \N}\E\left[\left\Vert M_n\right\Vert_{op}^3 \right] \right)\left\Vert u\right\Vert_{H^s(X)}^3 .\enma
Quitte à examiner des paquets de Cauchy, la majoration précédente montrera bien la convergence de la série aléatoire considérée dans l'espace de Banach $L_\omega^3(\Omega,L_{t}^3([-T,T],L_x^\infty(X)))$.
Comme déjà utilisé plusieurs fois dans cet article, cela implique la convergence en probabilité et a fortiori la convergence presque sûre de la série aléatoire considérée dans $L_{t}^3([-T,T],L_x^\infty(X))$. Comme $T$ parcourt $\N^\star$, on en déduira 
la convergence presque sûre dans $L_{t,loc}^3(\R,L_x^\infty(X))$.
Une nouvelle fois, le théorème de Fubini-Tonelli nous ramène à l'étude de l'intégrale 
\bema  
\int_{-T}^{T}\E \left[\left\Vert \sum_{n=0}^N \sum_{i=1}^{d_n}\left(\sum_{j=1}^{d_n} M_{n,i,j}(\omega)\langle u,\phi_{n,j}\rangle \right)e^{it\sqrt{-\Delta}}\phi_{n,i}(x)\right\Vert_{L_x^\infty(X)}^3 \right] dt.
\enma
\`A $t$ fixé, on peut écrire $\langle u,\phi_{n,j}\rangle=\langle e^{it\sqrt{-\Delta}}u,e^{it\sqrt{-\Delta}}\phi_{n,j}\rangle$.
Puisque $e^{it\sqrt{-\Delta}}$ envoie $\phi_{n,1},\dots,\phi_{n,d_n}$ sur une base hilbertienne de $E_n$ et que nos estimations sont toutes indépendantes de la base hilbertienne choisie, on peut appliquer le théorème \ref{pz-linf} pour obtenir l'estimation 
\bema \begin{array}{lcl}
\E \left[\left\Vert \Sum_{n=0}^N \sum_{i=1}^{d_n}\left(\sum_{j=1}^{d_n} M_{n,i,j}(\omega)\langle u,\phi_{n,j}\rangle \right)e^{it\sqrt{-\Delta}}\phi_{n,i}(x)\right\Vert_{L_x^\infty(X)}^3 \right]& &  \\
\hspace{4 cm} \leq C(X,s)\left(\sup\limits_{n\in \N}\E\left[\left\Vert M_n\right\Vert_{op}^3 \right] \right)\left\Vert e^{it\sqrt{-\Delta}} u\right\Vert_{H^s(X)}^3 .& & 
\end{array}
\enma

Cela achève notre preuve car $e^{it\sqrt{-\Delta}}$ est une isométrie de $H^s(X)$.
\end{demo}

\section{Preuve de la proposition \ref{beam-gaus}, randomisation des fonctions $Y_n$}

Rappelons en quel sens la fonction $Y_n$ se concentre de façon gaussienne sur une bande de largeur $\frac{1}{\sqrt{n}}$ autour de la géodésique $\{x_1^2+x_2^2=1\}$.
On part des inégalités
\bema
\forall \delta \in \left[0,\frac{\pi}{2}\right]\quad 1-\frac{\delta^2}{2}\leq \cos(\delta)\leq e^{-\frac{1}{2}\delta^2}.\enma 
Le nombre $\delta=\arccos\left(\sqrt{x_1^2+x_2^2}\right)$ désigne la distance géodésique de $x$ à $\{x_1^2+x_2^2=1\}$, puis la définition $Y_n(x)=c_{d,n} (x_1+ix_2)^n$ et l'équivalent $c_{d,n}\simeq_{d} n^{\frac{d-1}{4}}$ assurent qu'il existe $C(d)>1$ de sorte que 
\begin{equation}\label{circu-conc}\forall x\in \S^d \quad \forall n\in \N^\star \quad \frac{1}{C(d)}  \pun_{\left\{ \arccos\left(\sqrt{x_1^2+x_2^2}\right)\leq \frac{1}{\sqrt{n}}	 \right\}}  \leq
\frac{  |Y_{n}(x)|}{n^{\frac{d-1}{4}}} \leq C(d)  e^{-\frac{1}{2} n \left[\arccos\left(\sqrt{x_1^2+x_2^2}\right)\right]^2}.  \end{equation}
L'idée que l'on doit garder à l'esprit est que l'on peut sous certaines conditions assimiler $|Y_n|$ à la fonction $\widetilde{Y}_n$ définie comme suit
  \bema 
\widetilde{Y}_n (x):= n^{\frac{d-1}{4}}\pun_{\left\{\arccos\left(\sqrt{x_1^2+x_2^2}\right)\leq \frac{1}{\sqrt{n}}\right\}}.
\enma
Pour tout $p>1$, on peut alors vérifier les équivalences 
\begin{equation}\label{Yn-normLp} \left\Vert \widetilde{Y}_n\right\Vert_{L^p(\S^d)}\simeq_{d,p} 
\left\Vert Y_n\right\Vert_{L^p(\S^d)}\simeq_{d,p} n^{\frac{d-1}{2}\left(\frac{1}{2}-\frac{1}{p} \right)} , \end{equation}
ce qui signifie que la concentration autour de la géodésique $\{x_1^2+x_2^2=1\}$ est significative dans $L^p(\S^d)$.

Venons-en maintenant aux estimées multilinéaires des fonctions $Y_{n}$ qui font ``disparaître" la plus grande fréquence.
Le résultat suivant énonce que les estimées multilinéaires des fonctions $Y_n$ et $\widetilde{Y}_n$ sont équivalentes.

\begin{lemm}\label{ineg-mult}
Pour tout entier $\alpha \geq 2$ et pour tous entiers naturels $n_1\geq \dots\geq n_\alpha \geq 1$, on a 
\bemar
\Int_{\S^d} |Y_{n_1}(x) \cdots Y_{n_\alpha}(x)|^2 d\mu_d(x) & \simeq_{d,\alpha} &(n_2 \cdots n_\alpha)^\frac{d-1}{2}, \\[3mm]
\Int_{\S^d} |\widetilde{Y}_{n_1}(x) \cdots \widetilde{Y}_{n_\alpha}(x)|^2 d\mu_d(x)
& \simeq_{d,\alpha} &(n_2 \cdots n_\alpha)^\frac{d-1}{2}.
\enmar
\end{lemm}
\begin{demo}
Concernant les estimées des fonctions $Y_n$, 
il est possible d'obtenir une preuve qui n'utilise que l'encadrement \eqref{circu-conc} (voir la remarque \ref{rema-conc}).
On va cependant invoquer l'argument algébrique de \cite[pages 5-8]{burq2005multi}, on a pour tout $x\in \S^d$ : 
\bema
Y_{n_1}(x)\dots Y_{n_\alpha}(x)   = c_{d,n_1} \dots c_{d,n_\alpha} (x_1+ix_2)^{n_1+\dots+n_\alpha} = \Frac{c_{d,n_1} \dots c_{d,n_\alpha} }{c_{d,n_1+\dots+n_\alpha}} Y_{n_1+\dots +n_\alpha}(x).
\enma
Et donc 
\bema
\Int_{\S^d} |Y_{n_1}(x)\dots Y_{n_\alpha}(x)|^2 d\mu_d(x)  =  \left( \Frac{c_{d,n_1} \dots c_{d,n_\alpha} }{c_{d,n_1+\dots+n_\alpha}}\right)^2
 \simeq_{d}  \left(\Frac{n_1 \dots n_\alpha}{n_1+\dots+n_\alpha} \right)^{\frac{d-1}{2}}.
\enma 
On conclut en invoquant les inégalités $n_1\leq n_1+\dots+n_\alpha\leq \alpha n_1$.

Les intégrales multilinéaires des fonctions $\widetilde{Y}_n$ sont faciles à calculer l'aide d'une formule de changement de variables (voir plus loin \eqref{chg-var}) en tenant compte que $n_1$ est le plus grand entier parmi $n_2,\dots,n_\alpha$ :
\begin{eqnarray}\nonumber
\Frac{\Int_{\S^d} |\widetilde{Y}_{n_1}(x)\dots \widetilde{Y}_{n_\alpha}(x)|^2 d\mu_d(x) }{(n_1 \cdots n_\alpha )^{\frac{d-1}{2}}}& =& \Int_{\S^d} \pun_{\left\{\arccos\left(\sqrt{x_1^2+x_2^2}\right)\leq \frac{1}{\sqrt{n_1}} \right\}} d\mu_d(x)  \\ \nonumber
 &= &  C(d)\Int_{x_1^2+x_2^2<1}
(1-x_1^2-x_2^2)^{\frac{d-3}{2}}  \pun_{\left\{\arccos\left(\sqrt{x_1^2+x_2^2}\right)\leq \frac{1}{\sqrt{n_1}} \right\}} dx_1 dx_2 \\ \nonumber
& = & C(d) \Int_{\cos\left(\frac{1}{\sqrt{n_1}} \right)}^1 (1-r^2)^{\frac{d-3}{2}} rdr \\ \label{holomo}
& = & C(d)\Int_{0}^{\frac{1}{\sqrt{n_1}}} \sin^{d-2}(u) \cos(u)  du = \Frac{1}{n_1^{\frac{d-1}{2}}}\sum_{k\geq 0}\Frac{C(k,d)}{n_1^{\frac{k}{2}}} ,
\end{eqnarray}
où l'on a intégré le développement en série entière en $\xi=0$ de la fonction holomorphe $\xi\mapsto \sin^{d-2}(\xi)\cos(\xi)$.
\end{demo}

\begin{prop}\label{YN-form}
Pour tout réel $p>1$ et pour toute suite complexe $(a_n)_{n\geq 1}$, on a 
\bema
\left\Vert \sqrt{\sum_{n\geq 1} |a_n \widetilde{Y}_n|^2} \right\Vert_{L^p(\S^d)}\simeq_{d,p} \left[ \Sum_{n\geq 1} \frac{1}{n^{\frac{d+1}{2}}}\left(\sum_{k=1}^n k^{\frac{d-1}{2}}|a_k|^2 \right)^{\frac{p}{2}}\right]^{\frac{1}{p}}.
\enma
\end{prop}
\begin{demo}
Il s'agit de décomposer en supports disjoints
\bemar
\Sum_{n\geq 1} |a_n \widetilde{Y}_n(x)|^2 & =&  \Sum_{n\geq 1} n^{\frac{d-1}{2}} |a_n|^2 \pun_{\left\{\arccos(x)\leq\frac{1}{\sqrt{n}}\right\}} \\[3mm]
& = & \Sum_{n\geq 1}  \left(\sum_{k=1}^n k^{\frac{d-1}{2}} |a_k|^2\right)\pun_{\left\{\frac{1}{\sqrt{n+1}}< \arccos(x)\leq\frac{1}{\sqrt{n}}\right\}},
\enmar
puis l'on intègre et l'on fait usage du développement \eqref{holomo}.
\end{demo}

Le résultat suivant dit précisément que l'approximation de $|Y_n|$ par $\widetilde{Y}_n$ est légitime dans la théorie $L^p$ presque sûre des modes propres $Y_n$.
\begin{prop}\label{YN-equiv}
Pour tout réel $p\in ]1,+\infty[$ et pour toute suite complexe $(a_n)_{n\geq 1}$, on a 
\begin{equation}\label{YN-approx}
\frac{1}{C(p,d)} \left\Vert \sqrt{\sum_{n\geq 1} |a_n \widetilde{Y}_n|^2} \right\Vert_{L^p(\S^d)} \leq 
\left\Vert \sqrt{\sum_{n\geq 1} |a_n Y_n|^2} \right\Vert_{L^p(\S^d)} \leq C(p,d) \left\Vert \sqrt{\sum_{n\geq 1} |a_n \widetilde{Y}_n|^2} \right\Vert_{L^p(\S^d)}.
\end{equation}
\end{prop}
Faisons une remarque importante avant d'entamer la preuve.
L'idée consiste à vérifier \eqref{YN-approx} pour $p\in 2\N^\star$ puis à appliquer les résultats d'interpolation des espaces $\pl^p(\S^d,\oplus \C Y_{n})$.
Il convient de remarquer que cela ne semble pas découler gratuitement de la théorie de l'interpolation, bien connue, des espaces $L^p(\S^d,\ell^2(\N))$. La minoration dans \eqref{YN-approx} est évidente. Seule la majoration est difficile et équivaut à la bornitude de l'opérateur 
\bemar
\mathcal{U}:\pl^p(\S^d,\oplus \C \widetilde{Y}_n) & \rightarrow & L_x^p(\S^d,\ell^2(\N^\star)) \\
(a_n \widetilde{Y}_n)_{n\geq 1} & \mapsto & (a_n Y_n(x))_{n\geq 1}.
\enmar
D'après \eqref{defi-S}, l'espace $\pl^p(\S^d,\oplus \C\widetilde{Y}_n)$ s'identifie isométriquement à un sous-espace de $L^p(\S^d,\ell^2(\N^\star))$.
Si l'on veut faire de l'interpolation à partir des espaces $L^p(\S^d,\ell^2(\N^\star))$, alors il est tentant d'essayer de prolonger $\mathcal{U}$ sur tout l'espace $L^p(\S^d,\ell^2(\N^\star))$ dans lui-même (et de plus avec une expression indépendante de $p$).
Comme il n'y a pas de théorème de prolongement de Hahn-Banach pour les opérateurs, on pourrait essayer de justifier que $\pl^{p}(\S^d,\oplus \C \widetilde{Y}_n)$ (ou plutôt son image dans $L^p(\S^d,\ell^2(\N^\star))$) est un sous-espace complémenté de $L^p(\S^d,\ell^2(\N^\star))$.
Cela signifie que l'on doit trouver un projecteur borné de $L^p(\S^d,\ell^2(\N^\star))$ dont l'image est $\pl^p(\S^d,\oplus \C \widetilde{Y}_n)$.
Cela est parfaitement cohérent avec la démarche expliquée dans la partie \ref{partI} concernant la notion de rétracte d'espaces de Banach et explique pourquoi cela revenait à justifier la bornitude d'un projecteur!
Passons à la preuve de la  proposition \ref{YN-equiv}.

\begin{demo}
Si $p=2\alpha$ est un entier naturel non nul, \eqref{YN-approx} découle du lemme \ref{ineg-mult} et des deux formules
\bemar
\left\Vert \sqrt{\Sum_{n\geq 1} |a_n Y_n|^2} \right\Vert_{L^p(\S^d)}^p & = & \Sum_{n_1,\dots,n_\alpha} |a_{n_1}\dots a_{n_\alpha}|^2 \Int_{\S^d} |Y_{n_1}(x)\dots Y_{n_\alpha}(x)|^2 d\mu_d(x) \\
\left\Vert \sqrt{\Sum_{n\geq 1} |a_n \widetilde{Y}_n|^2} \right\Vert_{L^p(\S^d)}^p & = & \Sum_{n_1,\dots,n_\alpha} |a_{n_1}\dots a_{n_\alpha}|^2 \Int_{\S^d} |\widetilde{Y}_{n_1}(x)\dots \widetilde{Y}_{n_\alpha}(x)|^2 d\mu_d(x) .
\enmar 
Pour assurer que l'équivalence \eqref{YN-approx} est encore valide pour tout $p>1$, il nous suffit de prouver que les espaces $\pl^p(\S^d,\oplus\C Y_n)$ et 
$\pl^p(\S^d,\oplus\C \widetilde{Y}_n)$ sont stables par interpolation et dualité au sens des théorèmes \ref{theo-interpo} et \ref{theo-dual} avec $\sqrt{e(n,x)}=|Y_n(x)|$ et $d_n=1$.
En examinant \eqref{Yn-normLp} et en se rappelant que 
$\frac{1}{2}-\frac{1}{p}$ change de signe en remplaçant $p$ par son exposant conjugué $\frac{p}{p-1}$, on comprend que l'hypothèse \eqref{hol1} est vérifiée.

Quant à l'hypothèse \eqref{loren1}, elle va s'avérer être une conséquence de la concentration gaussienne des fonctions $Y_n$.
En effet, grâce à \eqref{circu-conc} et \eqref{Yn-normLp}, nous avons pour tout $n\in \N^\star$ et $x\in \S^d$ 
\bema
\frac{|Y_n(x)|}{\left\Vert Y_n\right\Vert_{L^p(\S^d)}}\leq C(d,p) n^{\frac{d-1}{2p}} e^{-\frac{n}{2} \left[\arccos\left( \sqrt{x_1^2+x_2^2}\right)\right]^2}\leq \frac{C(d,p)}{\arccos\left( \sqrt{x_1^2+x_2^2}\right)^{\frac{d-1}{p}}}.
\enma
Et il se trouve que $x\mapsto \arccos\left( \sqrt{x_1^2+x_2^2}\right)^{-\frac{(d-1)}{p}}$ appartient à l'espace de Lorentz $L^{p,\infty}(\S^d)$, ou ce qui revient au même que 
$x\mapsto \arccos\left( \sqrt{x_1^2+x_2^2}\right)^{-(d-1)}$ appartient à $L^{1,\infty}(\S^d)$ : la formule de changement de variables \eqref{chg-var} donne pour tout $t>1$ 
\begin{eqnarray}\nonumber 
&  & \mu_d\left\{x\in \S^d, \arccos\left( \sqrt{x_1^2+x_2^2}\right)^{-(d-1)}>t \right\} \\  \nonumber
& & \hspace{2 cm} =  C(d)\Int_{x_1^2+x_2^2<1} (1-x_1^2-x_2^2)^{\frac{d-3}{2}} \pun_{\left\{\arccos\left(\sqrt{x_1^2+x_2^2} \right)^{-(d-1)}>t\right\}}(x_1,x_2) dx_1 dx_2\\[4mm]\nonumber
& & \hspace{2 cm} =  C(d) \Int_{\cos\left(t^{\frac{-1}{d-1}} \right)}^{1} (1-r^2)^{\frac{d-3}{2}} rdr \\  \nonumber 
& & \hspace{2 cm} \simeq_d  \Int_{\cos\left(t^{\frac{-1}{d-1}} \right)}^{1} (1-r)^{\frac{d-3}{2}} dr \\ \nonumber
& & \hspace{2 cm}\simeq_d   \left[1-\cos\left(t^{\frac{-1}{d-1}} \right)\right]^{\frac{d-1}{2}} \\  \label{non-integ}
& & \hspace{2 cm} \lesssim_d   \Frac{1}{t}.
\end{eqnarray}
Comme $\mu_d(\S^d)$ est fini, l'estimation précédente est aussi valide si $t$ appartient à $]0,1]$.
\end{demo}
\begin{rema}\label{loren-perti}
Concernant l'hypothèse \eqref{loren1}, l'intérêt des espaces de Lorentz est désormais flagrant.
En effet, en utilisant la minoration de \eqref{circu-conc}, nous arrivons à 
\bemar
\forall x\in \S^d\quad \sup\limits_{n\geq 1} \Frac{|Y_n(x)|^p}{\left\Vert Y_n\right\Vert_{L^p(\S^d)}^p}& \geq & C(d,p)
\sup\limits_{n\geq 1} n^{\frac{d-1}{2}} \pun_{\left\{ \arccos\left( \sqrt{x_1^2+x_2^2}\right)\leq \frac{1}{\sqrt{n}} \right\}}(x)\\
&\geq & C(d,p)\arccos\left( \sqrt{x_1^2+x_2^2}\right)^{-(d-1)}.
\enmar
Or \eqref{non-integ} est une équivalence si $t$ tend vers $+\infty$.
Cela implique que $x\mapsto \arccos\left( \sqrt{x_1^2+x_2^2}\right)^{-(d-1)}$ n'est pas intégrable sur $\S^d$.
Il s'ensuit que $\sup\limits_{n\geq 1} \frac{|Y_n|}{\left\Vert Y_n\right\Vert_{L^p(\S^d)}}$ n'appartient pas à $L^p(\S^d)$.
\end{rema}

La proposition \ref{beam-gaus} est alors une conséquence des propositions \ref{YN-form} et \ref{YN-equiv} et des théorèmes \ref{theo-interpo} et \ref{theo-dual}.
On conclut avec une formule standard de changement de variables.

\begin{prop}
Considérons un entier naturel $k\in [1,d]$ et une fonction intégrable $f:\S^d\rightarrow \C$ qui ne dépend que des $k$ premières coordonnées $y=(x_1,x_2,\dots,x_k)$. En notant $\widetilde{f}(y)$ la valeur commune des nombres $f(y,z)$ pour $(y,z)\in \S^d$, on a 
\begin{equation}\label{chg-var}  \int_{\S^d } f(x) d\mu_d(x)  = \mu_{d-k}(\S^{d-k})\int_{\mathbb{B}_k(0,1)} \widetilde{f}(y) (1-|y|^2)^{\frac{d-k-1}{2}}d y,\end{equation}
où l'on a noté 
$\mathbb{B}_k(0,1):=\{y\in \R^k, |y| <1\}$. Dans le cas $k=d$, on convient que 
$\mu_{0}(\S^{0})=2$.
\end{prop}
\begin{demo}
On donne les principales lignes la preuve du cas $k\leq d-1$.
Le cas $k=d$ se traite de même.
On introduit le changement de variables 
\bemar  \Psi:\mathbb{B}_k(0,1)\backslash\{0\} \times \S^{d-k}& \rightarrow & \S^{d}  \\
(y,u)& \mapsto &   \left( y, \left(1-|y|^2\right)^{\frac{1}{2}} u\right) \enmar
dont la différentielle au point $(y,u)$ est l'application linéaire
\bemar
D_{(y,u)}\Psi:\R^k \times T_u \S^{d-k}& \rightarrow &    T_{\Psi(y,u)} \S^d \\
(\eta,w) & \mapsto & \left(\eta, \Frac{-\langle y,\eta \rangle }{\sqrt{1-|y|^2}}u+\left(1-|y|^2\right)^{\frac{1}{2}} w\right).
\enmar
L'injectivité de $\Psi$ et l'inversibilité de $D_{(y,u)}\Psi$ pour tout $(y,u)$ sont immédiates de sorte que le théorème d'inversion globale sur les variétés montre que $\Psi$ est $\CC^1$-difféomorphisme sur son image (qui est de mesure pleine dans $\S^d$).
Il nous reste à étudier le transport des formes volumes.
Considérons d'abord dans $\R^k$ une base orthonormée de la forme $\frac{y}{|y|},\overrightarrow{\xi_2},\dots,\overrightarrow{\xi_{k}}$, puis dans 
$T_u \S^{d-k}$ une base orthonormée quelconque $\overrightarrow{\vartheta_1},\dots,\overrightarrow{\vartheta_{d-k}}$.
Ainsi, 
\bema \left(\frac{y}{|y|},0\right),(\overrightarrow{\xi_2},0),\dots,(\overrightarrow{\xi_{k}},0),(0,\overrightarrow{\vartheta_1}),\dots,(0,\overrightarrow{\vartheta_{d-k}}) \enma
 est une base orthonormée de $\R^k\times T_u \S^{d-k}$.
 Il s'avère que cette base orthonormée est envoyée sur une famille orthogonale de l'espace tangent $T_{\Psi(y,u)}\S^d$ : 
 \bemar
D_{(y,u)}\Psi\left(\frac{y}{|y|},0\right)&=&\frac{1}{\sqrt{1-|y|^2}}\left(\frac{y}{|y|}\sqrt{1-|y|^2},-|y|u \right)  \\
2\leq i \leq k, \qquad D_{(y,u)}\Psi( \overrightarrow{\xi_i},0 ) & =&( \overrightarrow{\xi_i},0 ) \\
1\leq j\leq d-k,\qquad D_{(y,u)}\Psi(0,\overrightarrow{\vartheta_j}) & = & \sqrt{1-|y|^2}(0,\overrightarrow{\vartheta_j}).
 \enmar
Par conséquent, $D_{(y,u)}\Psi$ multiplie les volumes par $(1-|y|^2)^{\frac{d-k-1}{2}}$.
Le changement de variables $x=\Psi(y,u)$ nous amène à la formule 
\bemar  \Int_{\S^d } f(x) d\mu_d(x) &  =&  \Int_{\mathbb{B}_k(0,1)\times \S^{d-k}} f\left( y, u\sqrt{1-|y|^2}\right) (1-|y|^2)^{\frac{d-k-1}{2}} d y d\mu_{d-k}(u) \\
& = & \mu_{d-k}(\S^{d-k}) \Int_{\mathbb{B}_k(0,1)} \widetilde{f}(y) (1-|y|^2)^{\frac{d-k-1}{2}} d y .
\enmar 
\end{demo}

\section{Preuve du corollaire \ref{Yn-injesobo}, injections de Sobolev probabilistes des fonctions $Y_n$}

Nous avons vu dans la preuve de la proposition \ref{YN-equiv} que les espaces $\pl^p(\S^d,\oplus \C Y_n)$ sont stables par dualité.
Comme il en est de même des espaces de Sobolev, on peut se restreindre au cas $p>2$.

\textbf{L'inclusion stricte $H^{\frac{d-1}{2}\left(\frac{1}{2}-\frac{1}{p}\right)}(\S^d)  \subsetneq  \pl^p(\S^d,\oplus \C Y_n)$.}
L'inégalité triangulaire dans $L^{\frac{p}{2}}(\S^d)$ donne immédiatement :
\bema
\left\Vert \sqrt{\Sum_{n\geq 1} |a_n Y_n|^2} \right\Vert_{L^{p}(\S^d)}  =\left\Vert \Sum_{n\geq 1} |a_n Y_n|^2 \right\Vert_{L^{\frac{p}{2}}(\S^d)}^{\frac{1}{2}}  \leq  \sqrt{  \Sum_{n\geq 1} |a_n|^2 \left\Vert Y_n\right\Vert_{L^p(\S^d)}^2}. \enma
L'inclusion voulue découle alors des estimations \eqref{circul-normlp} 
des normes $\left\Vert Y_n\right\Vert_{L^p(\S^d)}$.
Nous allons donner deux preuves de l'inclusion stricte. L'une des deux preuves utilise la proposition suivante.

\begin{prop}\label{basincoLp}
Considérons $p\in ]1,+\infty[$, un espace métrique séparable $X$ muni d'une mesure borélienne $\sigma$-finie.
Si une suite $(u_n)_{n\in \N}$ de $L^p(X)$ converge faiblement vers $0$ et vérifie $\inf\limits_{n\in \N} \left\Vert u_n \right\Vert_{L^p(X)}>0$,
alors il existe un ensemble infini $I\subset \N$ et une constante $K \geq 1$ tels que pour toute suite complexe $(\alpha_n)_{n\in I}$ à support fini, on a 
\bema
\frac{1}{K} \left( \sum_{n\in I} |\alpha_n|^{\max(2,p)} \right)^{\frac{1}{\max(2,p)}} \leq  \left\Vert \sum_{n\in I} \alpha_n u_n \right\Vert_{L^p(X)} \leq K  \left( \sum_{n\in I} |\alpha_n|^{\min(2,p)} \right)^{\frac{1}{\min(2,p)}}.
\enma
\end{prop}

Nous n'utiliserons que la majoration précédente mais la minoration ne présente aucune difficulté supplémentaire.
La proposition précédente est facile à démontrer dans le cas $p=2$ mais nous ne connaissons pas de preuve élémentaire pour $p\neq 2$.
On proposera une preuve dans l'appendice \ref{incond} avec les notions de suite basique inconditionnelle, de type, de cotype et le théorème de \bess.

Le fait que l'inclusion $H^{\frac{d-1}{2}\left(\frac{1}{2}-\frac{1}{p}\right)}(\S^d)  \subset  \pl^p(\S^d,\oplus \C Y_n)$ soit stricte est une conséquence facile du résultat suivant.

\begin{prop}\label{lemm-incon}
Pour tout réel $p>2$, il existe un réel $r(p)>1$ et un sous-ensemble infini $I\subset \N^\star$ tel que pour toute suite $(a_n)_{n\geq 1}$ à croissance polynomiale et à support dans $I$, on a l'implication 
\bema
\sum_{n\in I} \left[n^{\frac{d-1}{2}\left(\frac{1}{2}-\frac{1}{p} \right)} |a_n| \right]^{2r(p)} <+\infty \quad \Rightarrow \quad \sum_{n\in I} a_n Y_n \in \pl^p(\S^d,\oplus \C Y_n).
\enma
\end{prop}
\begin{demo}
Commençons par une preuve simple mais un peu technique.
On va utiliser la condition \eqref{circu-cond} qui décrit les espaces $\pl^p(\S^d,\oplus \C Y_n)$. 
 Pour tout ensemble $I\subset \N^\star$ et toute suite $(a_n)_{n\in I}$ on estime à l'aide de l'inégalité de Hölder avec exposants $\left( \frac{p}{2},\frac{p}{p-2}\right)$
\bemar
\Sum_{n\geq 1} \Frac{1}{n^{\frac{d+1}{2}}} \left(\sum_{k=1}^n k^{\frac{d-1}{2}} |a_k|^2 \right)^{\frac{p}{2}} & \leq & \Sum_{n\geq 1} \Frac{1}{n^{\frac{d+1}{2}}} \left(\sum_{\ell=1}^n \ell^{\frac{(d-1)p}{2(p-2)}} \pun_{I}(\ell) \right)^{\frac{p-2}{2}}\left( \sum_{k=1}^n |a_k|^p\right) \\
& \leq & \Sum_{k\geq 1} |a_k|^p \sum_{n\geq k} \frac{1}{n^{\frac{d+1}{2}}} \left(\sum_{\ell=1}^n \ell^{\frac{(d-1)p}{2(p-2)}} \pun_{I}(\ell) \right)^{\frac{p-2}{2}}.
\enmar
Puisque l'on a $\frac{(d-1)p}{2(p-2)}>0$, on peut construire un ensemble $I\subset \N^\star$ tel que 
\bema
\forall n\in \N^\star \quad  \sum_{\ell=1}^n \ell^{\frac{(d-1)p}{2(p-2)}} \pun_{I}(\ell) \leq 2 n^{\frac{(d-1)p}{2(p-2)}}.
\enma
La condition $p>2$ implique que l'on a $\frac{p(d-1)}{4}-\frac{d}{2}-\frac{1}{2}>-1$.
Cela nous amène à la conclusion en choisissant $r(p)=\frac{p}{2}$ :
\bemar
\Sum_{n\geq 1} \Frac{1}{n^{\frac{d+1}{2}}} \left(\sum_{k=1}^n k^{\frac{d-1}{2}} |a_k|^2 \right)^{\frac{p}{2}} & \leq & 2^{\frac{p-2}{p}}\Sum_{k\geq 1} |a_k|^p \sum_{n\geq k} n^{\frac{p(d-1)}{4}-\frac{d}{2}-\frac{1}{2}} \\
& \lesssim_{d,p} & \Sum_{k\geq 1} |a_k|^p  k^{\frac{p(d-1)}{4}-\frac{d}{2}+\frac{1}{2}}, \\
& \lesssim_{d,p} & \Sum_{k\geq 1} |a_k|^p  k^{\frac{(d-1)}{2}\left(\frac{p}{2}-1\right)}.
\enmar 
On peut donner une seconde preuve plus abstraite à l'aide de la proposition \ref{basincoLp} et qui n'utilise que les estimations des normes $\left\Vert Y_n\right\Vert_{L^{p}(\S^d)}$.
On considère la suite de fonctions $L^{\frac{p}{2}}(\S^d)$-normalisées :
\bema
\frac{|Y_n|^2}{\left\Vert Y_n^2 \right\Vert_{L^{\frac{p}{2}}(\S^d)}}.
\enma
Par réflexivité de $L^{\frac{p}{2}}(\S^d)$, on peut en extraire une sous-suite qui converge faiblement vers une fonction $Y_\infty\in L^{\frac{p}{2}}(\S^d)$ au sens de la dualité avec $L^{\frac{p}{p-2}}(\S^d)$.
 Il est clair que $Y_{\infty}$ est une fonction positive.
Or nous avons d'après \eqref{circul-normlp} 
\bema  \int_{\S^d} Y_{\infty}(x) d\mu_d(x) = \lim\limits_{n\rightarrow +\infty}\int_{\S^d} \frac{|Y_n(x)|^2}{\left\Vert Y_n^2 \right\Vert_{L^{\frac{p}{2}}(\S^d)}} d\mu_d(x)= \lim\limits_{n\rightarrow+\infty} \frac{1}{\left\Vert Y_n \right\Vert_{L^{p}(\S^d)}^2}=0  .  \enma
Cela implique que l'on a $Y_{\infty}=0$.
La proposition \ref{basincoLp} assure que la suite de terme général $\frac{|Y_n|^2}{\left\Vert Y_n^2 \right\Vert_{L^{\frac{p}{2}}(\S^d)}}$
admet une sous-suite, indexée par un ensemble infini $I$, et une constante $K>0$ telle que pour toute suite complexe $(\alpha_n)_{n\in I}$ à support fini on a
\bema
\left\Vert \sum_{n\in I} \alpha_n \frac{|Y_n|^2}{\left\Vert Y_n^2 \right\Vert_{L^{\frac{p}{2}}(\S^d)}}\right\Vert_{L^{\frac{p}{2}}(\S^d)} \leq K \left( \sum_{n\in I} |\alpha_n|^{\min(2,\frac{p}{2})}\right)^{\frac{1}{\min(2,\frac{p}{2})}}.
\enma
On conclut aisément en choisissant la suite positive de terme général
\bema \alpha_n= |a_n|^2 \left\Vert Y_n^2 \right\Vert_{L^{\frac{p}{2}}(\S^d)}\simeq_{d,p} |a_n|^2 n^{(d-1)\left(\frac{1}{2}-\frac{1}{p} \right)}. \enma
\end{demo}
\textbf{L'inclusion stricte $\pl^p(\S^d,\oplus \C Y_n)\subsetneq \bigcap\limits_{\ep>0} H^{\frac{d-1}{2}\left(\frac{1}{2}-\frac{1}{p}\right)-\ep.}(\S^d)$.}
Pour tout $\ep>0$ et toute distribution $u=\sum_{n\geq 1} a_n Y_n$, il s'agit de vérifier l'inégalité suivante 
\begin{equation}\label{Yn-inclust}
\left\Vert u\right\Vert_{H^{\frac{d-1}{2}\left(\frac{1}{2}-\frac{1}{p}\right)-\ep}(\S^d)}+\left(\sum_{n\geq 1} \left\vert n^{\frac{(d-1)}{2}\left(\frac{1}{2}-\frac{1}{p} \right)} a_n\right\vert^p\right)^{\frac{1}{p}}\lesssim_{d,p,\ep} \left\Vert u\right\Vert_{\pl^p(\S^d,\oplus \C Y_n)}.
\end{equation}
Cette inégalité implique évidemment l'inclusion voulue. 
En outre, si l'on pose $J:=\{2^{k},k\in \N^\star\}$ et l'on définit $(a_n)_{n\geq 1}$ de sorte que 
\bemar n\not\in  J & \Rightarrow & a_n=0, \\
n\in J & \Rightarrow & a_n= n^{-\frac{(d-1)}{2}\left(\frac{1}{2}-\frac{1}{p}\right)},
\enmar
alors \eqref{Yn-inclust} assure que la distribution $\sum_{n\in J} a_n Y_n$ n'appartient pas à 
$\pl^p(\S^d,\oplus \C Y_n)$ mais appartient à tous les espaces de Sobolev $H^{\frac{d-1}{2}\left(\frac{1}{2}-\frac{1}{p}\right)-\ep}(\S^d)$.

Prouvons maintenant \eqref{Yn-inclust}.
On pose $S_0=0$ et $S_n=\sum_{k=1}^n k^{\frac{d-1}{2}} |a_k|^2$ pour tout $n\geq 1$. En utilisant une transformation d'Abel et une inégalité de Hölder, on a 
\bemar
\left\Vert u\right\Vert_{H^{\frac{d-1}{2}\left(\frac{1}{2}-\frac{1}{p}\right)-\ep}(\S^d)}^2 & = & \Sum_{n\geq 1} n^{(d-1)\left(\frac{1}{2}-\frac{1}{p}\right)-\ep}|a_n|^2 \\
& = & \Sum_{n\geq 1} n^{(d-1)\left(\frac{1}{2}-\frac{1}{p}\right)-\ep-\left(\frac{d-1}{2}\right)} \left[S_{n}-S_{n-1}\right] \\
& = &  \Sum_{n\geq 1} n^{-\frac{(d-1)}{p}-\ep} \left[S_{n}-S_{n-1}\right]\\
& \lesssim_{d,p,\ep} & \Sum_{n\geq 1} \Frac{1}{n^{\frac{(d-1)}{p}+1+\ep }} S_n \\
& \lesssim_{d,p,\ep} & \Sum_{n\geq 1} \Frac{1}{n^{1-\frac{2}{p}+\ep }} \times \frac{S_n}{n^{\frac{d+1}{p}}} \\
& \lesssim_{d,p,\ep} & \left(\Sum_{n\geq 1} \frac{S_n^{\frac{p}{2}}}{n^{\frac{d+1}{2}}} \right)^{\frac{2}{p}}.
\enmar
Pour achever la preuve de l'inégalité \eqref{Yn-inclust}, on utilise \eqref{circul-normlp} et le fait que $L^p(\S^d)$ est un espace de Banach de cotype $\max(2,p)=p$ (voir \eqref{Lp-type}) :
\bema
\left( \sum_{n\geq 1} \left[|a_n| \left\Vert Y_n \right\Vert_{L^p(\S^d)} \right]^p \right)^{\frac{1}{p}} \leq \left\Vert u\right\Vert_{\pl^p(\S^d,\oplus \C Y_n)}.
\enma

\section{Preuve de la proposition \ref{prop-zon}, randomisation des fonctions zonales $Z_n$}

La description des espaces $\pl^p(\S^d,\oplus \C Z_n)$ est plus délicate que celle des espaces $\pl^p(\S^d,\oplus \C Y_n)$ car on ne peut pas raisonner par interpolation en faisant parcourir $p$ dans $2\N$.
Cela dit, une fois cette description obtenue, les injections de Sobolev probabilistes \eqref{Zn-inclus} des fonctions $Z_n$ se démontrent de façon rigoureusement semblable à celle des fonctions $Y_n$ du corollaire \ref{Yn-injesobo}.
Avant d'expliquer plus en détail la difficulté rencontrée dans cette preuve par rapport à celle de la proposition \ref{beam-gaus}, commençons par rappeler les estimations précises des polynômes de Jacobi.

\begin{lemm}
Pour tout $\alpha>-1$, il existe des constantes $c=c(\alpha)\in ]0,\frac{\pi}{2}]$ et $C(\alpha)\geq 1$ de sorte que pour tous $n\in \N^\star$ l'on a 
\begin{equation}\label{jaco-es1}
\Theta \in \left[0,\frac{c}{n}\right]\cup \left[ \pi-\frac{c}{n},\pi\right] \quad \Rightarrow \quad \frac{n^\alpha}{C(\alpha)}  \leq P_n^{(\alpha,\alpha)}(\cos(\Theta)) \leq C(\alpha) n^{\alpha}.
\end{equation}
On note ensuite $N=n+\alpha+\frac{1}{2}$ et $\varrho =\frac{\pi}{2}\left(\alpha+\frac{1}{2}\right)$.
Si $\Theta$ appartient à $[\frac{c}{n},\pi-\frac{c}{n}]$ alors 
\begin{equation}\label{jaco-es2}
P_n^{(\alpha,\alpha)}(\cos(\Theta))=\frac{2^{\alpha+\frac{1}{2}}}{\sqrt{\pi n}(\sin \Theta)^{\alpha+\frac{1}{2}}}     \left[\cos\left(N\Theta-\varrho \right) +\frac{\mathcal{O}_{\alpha}(1)}{n\sin(\Theta)}\right],
\end{equation}
où le terme $\mathcal{O}_{\alpha}(1)$ vérifie $ 
|\mathcal{O}_{\alpha}(1)|\leq C(\alpha) $.
\end{lemm}
\begin{demo}
Ces estimées découlent des formules (4.1.3),(4.21.7), (7.32.5) et (8.21.18) du livre \cite{szeg}. 
La formule (4.1.3) nous donne $P_n^{(\alpha,\alpha)}(-x_1)=(-1)^{n} P_n^{(\alpha,\alpha)}(x_1)$, ce qui nous ramène au cas $\Theta\in [0,\frac{\pi}{2}]$.
D'une part, on a toujours 
\bema P_n^{(\alpha,\alpha)}(1)=\comb{n+\alpha}{n}\geq \frac{n^\alpha }{C(\alpha)} . \enma
Pour tout $x_1\in [1-\frac{1}{n^2},1]$, on peut estimer grâce aux formules 
(4.21.7) page 63 et (7.32.5) page 169 : 
\bema
\left\vert \frac{d}{dx_1} P_n^{(\alpha,\alpha)}(x_1) \right\vert  = \frac{1}{2}\left\vert(n+2\alpha+1)P_{n-1}^{(\alpha+1,\alpha+1)}(x_1)\right\vert  \leq  C(\alpha) n^{\alpha+2}.
\enma
Choisissons $c(\alpha)=\frac{1}{2C(\alpha)^2}$ de sorte que $\frac{1}{C(\alpha)}-c(\alpha) C(\alpha)= \frac{1}{2C(\alpha)}$. On a 
\bema
1-\frac{c(\alpha)}{n^2}\leq x_1 \leq 1 \quad \Rightarrow \quad P_n^{(\alpha,\alpha)}(x_1)\geq \frac{n^\alpha}{2 C(\alpha)}  .
\enma
De nouveau, d'après la formule (7.32.5) de la page 169 et quitte à augmenter $C(\alpha)>1$, on a aussi $P_n^{(\alpha,\alpha)}(x_1)\leq C(\alpha) n^\alpha$.
Cela nous donne \eqref{jaco-es1}.
Quant à \eqref{jaco-es2}, c'est la formule (8.21.18) de la page 198.
\end{demo}

Dans la suite, on notera 
\bema \Theta:=\arccos(x_1) \in [0,\pi] \enma
 la distance géodésique d'un point $x\in \S^d$ au pôle $(1,0,\dots,0)$.
D'après \eqref{defi-Zn}, \eqref{jaco-es1}, \eqref{jaco-es2}, on a 
\begin{eqnarray}  \label{legend1}
\Theta\in \left[ 0,\frac{c}{n}\right]\cup \left[\pi-\frac{c}{n},\pi \right] & \Rightarrow & \frac{n^{\frac{d-1}{2}}}{C(d)} \leq |Z_{n}(x)| \leq C(d) n^{\frac{d-1}{2}},\\[3mm] \label{legend2}
\Theta \in \left] \frac{c}{n},\pi-\frac{c}{n}\right[& \Rightarrow & |Z_{n}(x)| \leq \frac{C(d)}{\sin(\Theta)^{\frac{d-1}{2}}}.
\end{eqnarray}
C'est d'ailleurs avec ces estimations que l'on peut obtenir les estimations \eqref{born-Lp-Yn0} des normes dans $L^p(\S^d)$ des fonctions $Z_{n}$.
Ces dernières disent que seule la concentration au voisinage des pôles est significative dans l'échelle des espaces $L^p(\S^d)$, avec $p>\frac{2d}{d-1}$.
Il est donc naturel de comparer $Z_n$ à sa restriction $\widetilde{Z}_n$ au voisinage du pôle $(1,0,\dots,0)$ : 
\bema
\widetilde{Z}_n(x):= \pun_{[0,\frac{c}{n}]}(\Theta) \times Z_{n}(x).
\enma 
La fonction $\widetilde{Z}_n$ se concentre sur une boule de centre $(1,0,\dots,0)$ de rayon $\frac{1}{n}$ et avec une amplitude d'ordre $n^{\frac{d-1}{2}}$.
La preuve de la proposition \ref{beam-gaus} consistait à comparer $|Y_{n}|$ à sa restriction $\widetilde{Y}_n$ autour d'une géodésique.
Le lemme \ref{ineg-mult} assurait alors que les fonctions $|Y_{n}|$ ont les mêmes estimations multilinéaires que les fonctions $\widetilde{Y}_n$.
Malheureusement, il est illusoire de refaire le même argument en approchant les fonctions $Z_{n}$ par les fonctions $\widetilde{Z}_n$ par exemple pour étudier l'espace $\pl^6(\S^d,\oplus \C Z_{n})$.
Utilisant l'équivalent $\left\Vert Z_n \right\Vert_{L^2(\S^d)}\simeq _{d} 1$, \eqref{legend1} et \cite[line (2.3)]{burq2005multi}, nous avons en effet pour tout entier $n\in \N^\star$
\bema
\Int_{\S^d} |\widetilde{Z}_{1}(x)^2 \widetilde{Z}_n(x) |^2 d\mu_d (x)  \leq  C(d) \Int_{\S^d} |\widetilde{Z}_n(x) |^2 d\mu_d (x)  \leq  \frac{C(d)}{n}  \ll  \Int_{\S^d} |Z_1(x)^2 Z_{n}(x)|^2 d\mu_d (x).
\enma
L'estimation précédente est une manifestation de la mauvaise qualité de l'approximation de $Z_{n}$ par $\widetilde{Z}_n$ dans $L^2(\S^d)$, phénomène qui ne se produit pas pour les fonctions $Y_n$.
Malgré cela, la proposition \ref{mult-YN0} montrera que la concentration polaire des fonctions $Z_n$ est tout de même significative pour écrire leur théorie $L^p$ presque sûre pour $p>\frac{2d}{d-1}$.
Commençons par le résultat suivant.

\begin{prop}\label{YN0-LPinf}
Pour tout réel $p>\frac{2d}{d-1}$, il existe une constante $C(p,d)\geq 1$ telle que pour toute suite complexe $(a_n)_{n\geq 1}$ et tout $x\in \S^d$ vérifiant $\Theta(x)\in [0,\frac{\pi}{2}]$, on a
\begin{equation}\label{Z-majo} \sqrt{\Sum_{n\geq 1} |a_n Z_n(x)|^2} \leq C(p,d)\sqrt{\Sum_{n\geq 1} |a_n \widetilde{Z}_n(x)|^2}+\Frac{C(p,d)}{\sin^{\frac{d}{p}}(\Theta)}\left\Vert \sqrt{\Sum_{n\geq 1} |a_n \widetilde{Z}_n|^2} \right\Vert_{L^p(\S^d)} . \end{equation}
Par conséquent, on a 
\begin{equation}\label{Z-equiva}
\left\Vert \sqrt{\sum_{n\geq 1} |a_n Z_{n}|^2 }\right\Vert_{L^{p,\infty}(\S^d)}
\leq C(p,d)\left\Vert \sqrt{\sum_{n\geq 1} |a_n \widetilde{Z}_n|^2} \right\Vert_{L^{p}(\S^d)}.
\end{equation}
\end{prop}
\begin{demo}
Sans perte de généralité, on suppose que la suite $(a_n)$ n'a qu'un nombre fini de termes non nuls.
L'idée consiste essentiellement à décomposer les différentes fonctions en jeu en somme de fonctions à supports disjoints deux à deux. La formule de changement de variables \eqref{chg-var} donne 
\begin{eqnarray} \nonumber
\Int_{\S^d} \pun_{]\frac{c}{n+1},\frac{c}{n}]}(\Theta) d\mu_d(x) & = & \mu_{d-1}(\S^{d-1})\Int_{-1}^{1} \pun_{]\frac{c}{n+1},\frac{c}{n}]}(\Theta) (1-x_1^2)^{\frac{d-2}{2}}dx_1 \\ \nonumber 
& = & \mu_{d-1}(\S^{d-1}) \Int_{0}^{\pi} \pun_{]\frac{c}{n+1},\frac{c}{n}]}(\Theta) \sin(\Theta)^{d-1}d \Theta \\ \nonumber
& \simeq_d &  \Frac{1}{n^{d+1}}  .
\end{eqnarray}

En posant $S_n=\Sum_{k=1}^n  k^{d-1}|a_k|^2$ pour tout $n\in \N^\star$, nous avons grâce à \eqref{legend1} :
\begin{eqnarray} \nonumber
\sqrt{\Sum_{n\geq 1} |a_n \widetilde{Z}_n(x)|^2} & \simeq_d &   \sqrt{\Sum_{n\geq 1} n^{d-1} |a_n|^2 \pun_{[0,\frac{c}{n}]}(\Theta)} \\ \label{form-fn}
& \simeq_d & \Sum_{n\geq 1} \sqrt{S_n} \pun_{]\frac{c}{n+1},\frac{c}{n}]}(\Theta), \end{eqnarray}
et donc 
\begin{eqnarray} \nonumber
\left\Vert \sqrt{\sum_{n\geq 1} |a_n \widetilde{Z}_n|^2 }\right\Vert_{L^p(\S^d)} &\simeq_d & \left( \sum_{n\geq  1} S_n^{\frac{p}{2}} \int_{\S^d} \pun_{]\frac{c}{n+1},\frac{c}{n}]}(\Theta)d\mu_d(x) \right)^{\frac{1}{p}} \\ \label{zon-nor}
& \simeq_{d,p} & \left( \sum_{n\geq  1} \frac{S_n^{\frac{p}{2}}}{n^{d+1}} \right)^{\frac{1}{p}}   .
\end{eqnarray}
On peut maintenant estimer $Z_n(x)-\widetilde{Z}_n(x)$.
D'après \eqref{legend2}, nous avons 

\bemar
\sqrt{\Sum_{n\geq 1}|a_n|^2 |Z_{n}(x)-\widetilde{Z}_n(x)|^2 } & \lesssim_d & \left(\Sum_{k\geq 1} |a_k|^2\right)^{\frac{1}{2}}
 \Frac{\pun_{]c,\frac{\pi}{2}]}(\Theta)}{\sin(\Theta)^{\frac{d-1}{2}}}
+\Sum_{n\geq 1} \left(\Sum_{k> n} |a_k|^2 \right)^{\frac{1}{2}} \frac{\pun_{]\frac{c}{n+1},\frac{c}{n}]}(\Theta)}{\sin(\Theta)^{\frac{d-1}{2}}} \\
& \lesssim_d &  \underbrace{\left(\Sum_{k\geq 1} |a_k|^2\right)^{\frac{1}{2}}
 \pun_{]c,\frac{\pi}{2}]}(\Theta)}_{:=A_1(\Theta)}
+\underbrace{\Sum_{n\geq 1} \left(n^{d-1} \Sum_{k> n} |a_k|^2 \right)^{\frac{1}{2}} \pun_{]\frac{c}{n+1},\frac{c}{n}]}(\Theta)  }_{:=A_2(\Theta)}     .   \\
\enmar 

On va faire quelques calculs avant d'attaquer l'estimation des termes $A_1(\Theta)$ et $A_2(\Theta)$.
En convenant que $S_{0}=0$, nous pouvons effectuer une transformation d'Abel pour tout $n\in\N$
\begin{equation}\label{Abel}
\sum_{k> n} |a_k|^2= 
\sum_{k> n} \frac{k^{d-1}|a_k|^2}{k^{d-1}}=\frac{-S_{n}}{(n+1)^{d-1}}+\sum_{k> n} \left(\frac{1}{k^{d-1}}-\frac{1}{(k+1)^{d-1}} \right)S_k.
\end{equation}
Puisque la suite $(a_n)_{n\geq 1}$ n'a qu'un nombre fini de termes non nuls, la suite $(S_n)_{n\geq 0}$ est bornée et la dernière série converge bien .
Remarquons maintenant l'égalité triviale 
\bema
\left(d-\frac{2(d+1)}{p}\right)\frac{p}{p-2}=d-\frac{2}{p-2}.
\enma
Exploitons alors que l'on a $p>\frac{2d}{d-1}$, ou encore $d-\frac{2}{p-2}>1$, 
 à l'aide de l'inégalité de Hölder avec les exposants conjugués $\frac{p}{p-2},\frac{p}{2}$ :
\begin{eqnarray} \nonumber
(n+1)^{d-1}\Sum_{k> n} \frac{S_k}{k^d} & \leq  & (n+1)^{d-1}\Sum_{k> n}\overbrace{\frac{1}{k^{d-\frac{2(d+1)}{p}}}}^{\in \ell^{\frac{p}{p-2}}(\N)} \times \frac{S_k}{k^\frac{2(d+1)}{p}} \\ \nonumber
& \leq & 
(n+1)^{d-1}\left(\Sum_{k> n} \frac{1}{k^{d-\frac{2}{p-2}}} \right)^{\frac{p-2}{p}} \left( \Sum_{k> n} \frac{S_k^{\frac{p}{2}}}{k^{d+1}}\right)^{\frac{2}{p}} \\ \nonumber
& \lesssim_{d,p} & \Frac{(n+1)^{d-1}}{(n+1)^{\frac{p-2}{p}\left(d-\frac{2}{p-2}-1\right)}}\left( \Sum_{k> n} \frac{S_k^{\frac{p}{2}}}{k^{d+1}}\right)^{\frac{2}{p}} \\ \label{Z-majo-2}
& \lesssim_{d,p} & (n+1)^{\frac{2d}{p}}\left( \Sum_{k> n} \frac{S_k^{\frac{p}{2}}}{k^{d+1}}\right)^{\frac{2}{p}} \\ \label{Z-majo-3}
& \lesssim_{d,p} &  (n+1)^{\frac{2d}{p}} \left\Vert \sqrt{\Sum_{n\geq 1} |a_n \widetilde{Z}_n|^2 }\right\Vert_{L^p(\S^d)}^2      ,
\end{eqnarray}
où l'on a utilisé \eqref{zon-nor}.
On peut alors contrôler $A_1(\Theta)$ avec \eqref{Abel} et \eqref{Z-majo-3} pour $n=0$ : 

\bema
A_1(\Theta)=\sqrt{\Sum_{k\geq 1} |a_k|^2} \times \pun_{]c,\frac{\pi}{2}[}(\Theta)  \lesssim_{d,p}  
  \left\Vert \sqrt{\Sum_{n\geq 1} |a_n \widetilde{Z}_n|^2 }\right\Vert_{L^p(\S^d)} \Frac{1}{\sin^{\frac{d}{p}}(\Theta)} .
\enma
Pour contrôler $A_2(\theta)$, on utilise \eqref{form-fn}, \eqref{Abel} et \eqref{Z-majo-3}   : 
\begin{eqnarray} \nonumber   
A_2(\Theta) & \lesssim_{d,p} & \Sum_{n\geq 1} \sqrt{S_{n}} \pun_{]\frac{c}{n+1},\frac{c}{n}]}(\Theta) +\Sum_{n\geq 1} n^{\frac{d}{p}} \left\Vert \sqrt{\Sum_{n\geq 1}|a_n \widetilde{Z}_n|^2}\right\Vert_{L^p(\S^d)} \pun_{]\frac{c}{n+1},\frac{c}{n}]}(\Theta) \\ \nonumber
& \lesssim_{d,p} &  \sqrt{\Sum_{n\geq 1} |a_n \widetilde{Z}_n(x)|^2}+\left\Vert \sqrt{\Sum_{n\geq 1}|a_n \widetilde{Z}_n|^2}\right\Vert_{L^p(\S^d)} \Frac{1}{\sin^{\frac{d}{p}}(\Theta)}.
\end{eqnarray}
On a donc obtenu \eqref{Z-majo}.
Passons à l'inégalité \eqref{Z-equiva}.
Comme les polynômes de Jacobi sont pairs ou impairs, les fonctions $x\mapsto |Z_n(x)|$ sont invariantes par la transformation $\Theta(x) \rightarrow \pi-\Theta(x)$.
Se rappelant que $\left\Vert \cdot\right\Vert_{L^{p,\infty}(\S^d)}$ n'est pas une norme, on sait qu'il existe néanmoins une constante universelle $C>1$ telle que
\bema \left\Vert \sqrt{\sum_{n\geq 1} |a_n Z_{n}|^2 }\right\Vert_{L^{p,\infty}(\S^d)}\leq 
C \left\Vert \sqrt{\sum_{n\geq 1} |a_n Z_{n}|^2 } \times \pun_{[0,\frac{\pi}{2}]}(\Theta)\right\Vert_{L^{p,\infty}(\S^d)}. \enma
Pour obtenir l'estimation \eqref{Z-majo} à partir de \eqref{Z-equiva}, on remarque d'une part l'inégalité triviale
\bema
\left\Vert \sqrt{\sum_{n\geq 1} |a_n \widetilde{Z}_n|^2}\right\Vert_{L^{p,\infty}(\S^d)}\leq \left\Vert \sqrt{\sum_{n\geq 1} |a_n \widetilde{Z}_n|^2}\right\Vert_{L^{p}(\S^d)} ,
\enma
d'autre part que la fonction $x\mapsto \sin(\Theta(x))^{-\frac{d}{p}}$ appartient à $L^{p,\infty}(\S^d)$. Cela peut se vérifier avec la formule de changement de variable \eqref{chg-var} mais il est plus simple de remarquer que sur un voisinage $\mathcal{V}$ du pôle $( 1,0,\dots 0)\in \S^d$ nous avons l'équivalent $\sin(\Theta)^{-\frac{d}{p} }	\sim \Theta^{-\frac{d}{p}}$ et que la mesure de $\S^d$ sur $\mathcal{V}$ est comparable à la mesure de Lebesgue d'un voisinage de l'origine de $\R^d$.
\end{demo}

Si l'on essaie d'estimer directement la norme dans $L^p(\S^d)$ de $\sqrt{\sum_{n\geq 0} |a_n Z_n|^2}$ en comparant $Z_n$ avec $\widetilde{Z}_n$ à l'aide de \eqref{legend1} et \eqref{legend2}, alors la preuve précédente montre que l'on commet une perte avec l'inégalité de Hölder.
De façon précise, après application de l'inégalité de Hölder, la formule \eqref{Z-majo-2} et le contrôle de $\left\Vert A_2(\Theta)\right\Vert_{L^p(\S^d)}$ conduisent aux inégalités 
\bemar
\left\Vert \Sum_{n\geq 1} n^{\frac{d}{p}} \left( \Sum_{k> n} \frac{S_k^{\frac{p}{2}}}{k^{d+1}}\right)^{\frac{1}{p}} \pun_{]\frac{c}{n+1},\frac{c}{n}]}(\Theta)\right\Vert_{L^p(\S^d)}^p & \simeq_d & \Sum_{n\geq 1} \frac{n^{\frac{d}{p}p}}{n^{d+1}}\sum_{k> n} \frac{S_k^{\frac{p}{2}}}{k^{d+1}} \\
& \simeq_d & \Sum_{n\geq 1} \frac{1}{n}\sum_{k> n} \frac{S_k^{\frac{p}{2}}}{k^{d+1}}\\
& \simeq_d & \Sum_{k\geq 1} \frac{\ln(k+1)}{k^{d+1}} S_k^{\frac{p}{2}}.
\enmar
Par comparaison avec \eqref{zon-nor}, l'estimation précédente est mauvaise en raison du terme logarithmique.
Un argument d'interpolation réelle bien connu va nous permettre de corriger le facteur logarithmique.
C'est maintenant que le théorème \ref{theo-interpo} intervient via le lemme suivant.

\begin{lemm}\label{lemm-interp}
La famille d'espaces de Banach $(\pl^p(\S^d,\oplus \C \widetilde{Z}_n))_{p\in ]1,+\infty[}$ est stable par interpolation réelle au sens du théorème \ref{theo-interpo}.
\end{lemm}
\begin{demo}
Il s'agit d'appliquer le théorème \ref{theo-interpo} avec $E_n=\C \widetilde{Z}_n$
 et  $\sqrt{e(n,x)}=|\widetilde{Z}_n(x)|$.
D'une part, on a le comportement asymptotique 
\bema \forall p>1 \quad \left\Vert \widetilde{Z}_n\right\Vert_{L^p(\S^d)}  \simeq_{d,p}
n^{\frac{d-1}{2}-\frac{d}{p}}, \enma
ce qui nous assure la validité de \eqref{hypo-interpo}.
Quant à l'hypothèse \eqref{loren1}, on peut écrire d'après \eqref{legend1} :
\bema
\forall x\in \S^d \qquad \sup\limits_{n\geq 1} \frac{|\widetilde{Z}_n(x)|}{\left\Vert \widetilde{Z}_n\right\Vert_{L^p(\S^d)}} \leq C(d,p)
\sup\limits_{n\geq 1} \frac{|\widetilde{Z}_n(x)|}{n^{\frac{d-1}{2}-\frac{d}{p}}}= C(d,p)
\sup\limits_{n\geq 1} n^{\frac{d}{p}}\pun_{[0,\frac{c}{n}]}(\Theta)\leq \frac{C(d,p)}{\sin(\Theta)^{\frac{d}{p}}}.
\enma 
Or l'on a remarqué à la fin de la proposition \ref{YN0-LPinf} que la fonction $x\mapsto \sin(\Theta)^{-\frac{d}{p}}$ appartient à $L^{p,\infty}(\S^d)$.
\end{demo}

Nous sommes en mesure d'obtenir la partie de la proposition \ref{prop-zon} qui décrit $\pl^p(\S^d,\oplus \C Z_n)$.

\begin{prop}\label{mult-YN0}
Pour tout réel $p>\frac{2d}{d-1}$, et pour toute suite complexe $(a_n)_{n\geq 1}$ on a
\begin{equation}\label{equiv-YN0}
\left\Vert \sqrt{\sum_{n\geq 1} |a_n \widetilde{Z}_n|^2} \right\Vert_{L^{p}(\S^d)}\simeq_{d,p} \left\Vert \sqrt{\sum_{n\geq 1} |a_n Z_n|^2 }\right\Vert_{L^{p}(\S^d)}.
\end{equation}
En outre, on a
\begin{equation}\label{YN0-carac} \left\Vert \sqrt{\Sum_{n\geq 1} |a_n Z_n|^2 }\right\Vert_{L^{p}(\S^d)} \simeq_{d,p}
\left[\Sum_{n\geq 1}\frac{1}{n^{d+1}} \left(\sum_{k=1}^n k^{d-1}|a_k|^2 \right)^{\frac{p}{2}}  \right]^{\frac{1}{p}} .\end{equation}
\end{prop}
\begin{demo}
L'équivalence \eqref{YN0-carac} découlera de \eqref{zon-nor} et \eqref{equiv-YN0}.
Montrons donc \eqref{equiv-YN0}.
L'inégalité $|\widetilde{Z}_n|\leq |Z_{n}|$ donne un sens de l'équivalence \eqref{equiv-YN0}.
L'autre sens équivaut à la bornitude de l'opérateur : 
\bemar
\pl^p(\S^d,\oplus \C \widetilde{Z}_n) & \rightarrow & L_x^p(\S^d,\ell^2(\N^\star)) \\
(a_n \widetilde{Z}_n)_{n\geq 1} & \mapsto & (a_n Z_{n}(x))_{n\geq 1}.
\enmar
L'inégalité \eqref{Z-equiva} de la proposition \ref{YN0-LPinf} montre, pour tout $p>\frac{2d}{d-1}$, la continuité de l'opérateur 
\bemar
\pl^p(\S^d,\oplus \C \widetilde{Z}_n) & \rightarrow & L_x^{p,\infty}(\S^d,\ell^2(\N^\star)) \\
(a_n \widetilde{Z}_n)_{n\geq 1} & \mapsto & (a_n Z_n(x))_{n\geq 1}.
\enmar
\`A présent, il s'agit de raisonner par interpolation réelle. Fixons deux réels $p_1'$ et $p_2'$ tels que 
\bema \frac{2d}{d-1}<p_1'<p<p_2'.\enma
Considérons de plus l'unique réel $\theta' \in [0,1]$ tel que 
 \bema
 \frac{1}{p}=\frac{1-\theta'}{p_1'}+\frac{\theta'}{p_2'}.
 \enma
Par application de la méthode d'interpolation réelle $[\cdot,\cdot]_{\theta',p}$, l'opérateur suivant est borné : 
\bemar
[\pl^{p_1'}(\S^d,\oplus \C \widetilde{Z}_n),\pl^{p_2'}(\S^d,\oplus \C \widetilde{Z}_n)]_{\theta',p}& \rightarrow & \left[L_x^{p_1',\infty}(\S^d,\ell^2(\N^\star)),L_x^{p_2',\infty}(\S^d,\ell^2(\N^\star))\right]_{\theta',p} \\
(a_n \widetilde{Z}_n)_{n\geq 1} & \mapsto & (a_n Z_{n}(x))_{n\geq 1}.
\enmar
L'espace de départ s'identifie à $\pl^p(\S^d,\oplus \C \widetilde{Z}_n)$ d'après le lemme \ref{lemm-interp}.
Quant à l'espace d'arrivée, il s'identifie à $L_x^{p}(\S^d,\ell^2(\N^\star))$
d'après le théorème \ref{interpo-loren}.
Cela prouve \eqref{equiv-YN0}.
\end{demo}
\begin{rema}
Dans cette preuve, le point difficile est la bornitude de l'opérateur 
\bemar
\pl^p(\S^d,\oplus \C \widetilde{Z}_n) & \rightarrow & \pl^{p}(\S^d, \oplus\C Z_n) \\
(a_n \widetilde{Z}_n)_{n\geq 1} & \mapsto & (a_n Z_{n})_{n\geq 1}.
\enmar
En ce qui concerne les espaces d'arrivée, on a remplacé les espaces 
$\pl^{p}(\S^d, \oplus\C Z_n)$ par les espaces plus grands $L^p(\S^d,\ell^2(\N^\star))$ qui admettent une théorie de l'interpolation connue depuis longtemps.
Mais, souhaitant raisonner par interpolation, nous sommes obligés d'écrire une théorie d'interpolation adaptée aux espaces de départ $\pl^p(\S^d,\oplus \C \widetilde{Z}_n)$ (voir aussi la discussion sur la complémentation des espaces de Lebesgue probabilistes après la proposition \ref{YN-equiv})
\end{rema}

En apparence, on n'a pas montré l'interpolation des espaces $\pl^p(\S^d,\oplus \C Z_n)$.
La preuve de la proposition suivante montre en fait que cela est inclus dans l'inégalité \eqref{Z-majo}.

\begin{prop}
Les espaces $\pl^p(\S^d,\oplus \C Z_n)$ sont stables par interpolation réelle et complexe pour $p$ parcourant $]\frac{2d}{d-1},+\infty[$ au sens du théorème \ref{theo-interpo}.
\end{prop}
\begin{demo}
On peut donner deux preuves. La proposition \ref{mult-YN0} montre que les espaces $\pl^p(\S^d,\oplus \C Z_n)$ et $\pl^p(\S^d,\oplus \C \widetilde{Z}_n)$ sont isomorphes pour tout $p>\frac{2d}{d-1}$.
Le lemme \ref{lemm-interp}, qui contient aussi dans sa preuve l'interpolation complexe, donne alors la conclusion.

La seconde preuve utilise aussi la démonstration du lemme \ref{lemm-interp} avec l'argument additionnel suivant.
D'après \eqref{Z-majo}, on a pour tout $x\in \S^d$ vérifiant $\Theta(x)\in [0,\frac{\pi}{2}]$ :
\bema
|Z_n(x)|\leq C(d,p) |\widetilde{Z}_n(x)|+\frac{C(d,p)}{\sin^{\frac{d}{p}}(\Theta)}\left\Vert \widetilde{Z}_n\right\Vert_{L^p(\S^d)} .
\enma
La symétrie $\Theta(x)\leftrightarrow \pi-\Theta(x)$ et le fait que $\sin^{-d}(\Theta)$ appartient à $L^{1,\infty}(\S^d)$ nous donnent 
\bema
\sup\limits_{n\in \N^\star} \frac{|Z_n(x)|}{\left\Vert Z_n \right\Vert_{L^p(\S^d)}}\in L_x^{p,\infty}(\S^d).
\enma
On conclut par application du théorème \ref{theo-interpo} d'interpolation pour les espaces $\pl^p(\S^d,\oplus \C Z_n)$.
\end{demo}

\section{Rappels sur les fonctions de Hermite}

On reprend les notations de la partie \ref{Lp-oscillo}, $(h_n)_{n\in \N}$ désigne donc la suite des fonctions de Hermite $L^2(\R)$-normalisées et la fonction spectrale $e_d(n,\cdot)$ s'exprime à l'aide de \eqref{h-sumcar}. Par récurrence, on a immédiatement 
\bemar
e_d(n,x)& =& \Sum_{\substack{(i_1,\dots,i_d)\in \N^d\\ i_1+\dots+i_d=n }} h_{i_1}(x_1)^2\dots h_{i_d}(x_d)^2 \\ 
& = & \Sum_{k=0}^n e_{d-1}(k,\underbrace{(x_1,\dots,x_{d-1})}_{\in \R^{d-1}}) h_{n-k}(x_d)^2.
\enmar
Cela dit, nous allons utiliser d'autres formules plus maniables.
Remarquons que toute rotation linéaire de $\R^d$ commute avec l'oscillateur harmonique $-\Delta+|x|^2$ et laisse donc invariant son sous-espace propre $E_{d,n}$ associé à la valeur propre $d+2n$.
On en déduit que la fonction spectrale $e_{d}(n,\cdot)$ de $E_{d,n}$ est invariante par rotation linéaire. En d'autres termes, $e_d(n,x)$ ne dépend que de la norme euclidienne $|x|=\sqrt{x_1^2+\dots+x_d^2}$ et l'on peut écrire 
\begin{eqnarray}\label{edn-1}
e_d(n,x) & =& \Sum_{k=0}^n e_{d-1}(k,\overbrace{(|x|,0,\dots,0)}^{\in \R^{d-1}})h_{n-k}(0)^2 \\ \label{edn-2}
& = & \Sum_{k=0}^n e_{d-1}(k,0)h_{n-k}(|x|)^2.
\end{eqnarray}

Le lemme suivant nous donne des estimations faciles
concernant les fonctions de Hermite et les fonctions spectrales en $0$.

\begin{prop}\label{estim-hn}
Pour tout entier $k\in \N$, on a $h_{2k+1}(0)=0$ et l'équivalence
$(-1)^k h_{2k}(0)\simeq \frac{1}{(k+1)^{1/4}}$.
Si $n$ est impair alors on a $e_d(n,0)=0$. Par contre, si $n$ est pair alors on a  $e_d(n,0)\simeq_d n^{\frac{d}{2}-1}$.
\end{prop}
\begin{demo}
La première équivalence découle de la formule (5.5.5) de \cite{szeg} : 
\bema (-1)^k h_{2k}(0)=
(-1)^k\frac{H_{2k}(0)}{\sqrt{(2k)! 4^k \sqrt{\pi}}} =\frac{\sqrt{(2k)!}}{k!2^k \pi^{1/4}} \simeq \frac{1}{(k+1)^{\frac{1}{4}}}.  \enma
Pour le cas $n$ impair, on a déjà remarqué l'égalité $e_d(n,0)=0$ après l'énoncé de la proposition \ref{ed-maj}.
Pour le cas $n$ pair, on peut écrire :
\bemar   e_d(n,0)&= &\Sum_{2i_1+\dots+2i_d=n} h_{2i_1}(0)^2 \cdots h_{2i_d}(0)^2 \\[3mm]
&  \simeq_d &   \Sum_{2i_1+\dots+2i_d=n} \frac{1}{\sqrt{(1+i_1)\cdots(1+i_d)}} \\
& \gtrsim_d &  \Sum_{2i_1+\dots+2i_d=n}n^{-\frac{d}{2}}\\
& \gtrsim_d & n^{\frac{d}{2}-1}.  
\enmar
La majoration $e_d(n,0)\lesssim_{d} n^{\frac{d}{2}-1}$ peut se démontrer par récurrence sur $d$ en séparant les sommes suivantes selon que $\ell$ est plus petit ou plus grand que $\frac{n}{4}$ :
\bema
 e_{d}(n,0)=\sum_{\ell=0}^{\frac{n}{2}}e_{d-1}(n-2\ell,0)h_{2\ell}(0)^2 \simeq 
\sum_{\ell=0}^{\frac{n}{2}}\frac{e_{d-1}(n-2\ell,0)}{\sqrt{1+\ell}}.
\enma
\end{demo}

Nous aurons besoin d'estimations précises concernant le comportement des fonctions de Hermite. Pour le résultat suivant, on fait référence à \cite[Lemma 1.5.1 and line 3.2.13]{thanga}.

\begin{prop}\label{hermi-majo}
Il existe deux constantes universelles $C>1$ et $\gamma>0$ telles que pour tous $x\in \R$ et $n\in \N$ on a 
\bemar
|x|\leq \sqrt{2n+1}& \Rightarrow & |h_n(x)|\leq \Frac{1}{\sqrt[4]{2n+2-x^2}}, \\[3mm]
|x|> \sqrt{2n+1} & \Rightarrow & |h_n(x)|\leq C e^{-\frac{\gamma}{2}x^2}.
\enmar
\end{prop}

Comme les fonctions $h_n$ oscillent et s'annulent plusieurs fois, les estimations de la proposition \ref{hermi-majo} sont inutilisables pour minorer les fonctions $|h_n|$.
Pour remédier à cela, on fait appel à un résultat d'approximation des fonctions de Hermite dû à Muckenhoupt 
et prouvé dans \cite[line 2.4]{muckenI} à partir de \cite[page 23]{erdelyi} (voir aussi une utilisation de cette formule dans \cite[line 7.2]{muckenII}).

\begin{prop}\label{muck}
Introduisons la fonction croissante $\Phi :[0,1]\rightarrow [0,\frac{\pi}{4}]$ définie par 
\bema \Phi(u)=\int_{0}^u \sqrt{1-s^2} ds=\frac{1}{2}\arcsin(u)+\frac{1}{2}u\sqrt{1-u^2}.    \enma
Il existe une constante universelle $C>1$ telle que pour tous $n\in \N$ et $x\in [0,\sqrt{2n+1}-(2n+1)^{-\frac{1}{6}}]$ on a 
\bema\left\vert h_n(x)-   \frac{\sqrt{2}}{\sqrt{\pi}(2n+1-x^2)^{\frac{1}{4}}}\cos\left[(2n+1)\Phi \left(\frac{x}{\sqrt{2n+1}} \right)-\frac{n\pi}{2}\right] \right\vert \leq  \frac{C\sqrt{2n+1}}{(2n+1-x^2)^{\frac{7}{4}}}.    \enma
\end{prop}
En fait, le terme principal est exprimé dans \cite[line 2.4]{muckenI} sous la forme différente  mais équivalente 
\bema\frac{\sqrt{2}}{\sqrt{\pi}(2n+1-x^2)^{\frac{1}{4}}}\cos\left[\frac{1}{4}(2n+1)[2\theta-\sin(2\theta)]-\frac{\pi}{4}\right],\enma
avec $\theta=\frac{\pi}{2}-\arcsin\left(\frac{x}{\sqrt{2n+1}} \right)$ et donc $\sin(2\theta)=2\sin(\theta)\cos(\theta)=2\frac{x}{\sqrt{2n+1}}\sqrt{1-\frac{x^2}{2n+1}}$.
\`A noter que sur l'intervalle $[0,\sqrt{2n+1}-(2n+1)^{-\frac{1}{6}}]$, le reste est contrôlé par l'amplitude du cosinus 
\bema    \frac{\sqrt{2n+1}}{(2n+1-x^2)^{\frac{7}{4}}}\leq \frac{C}{(2n+1-x^2)^{\frac{1}{4}}}  .         \enma
Pour tout $\beta \in ]0,1[$ et $x\in [0,\beta\sqrt{2n+1}]$, 
les propositions \ref{hermi-majo} et \ref{muck} nous amènent à 
\begin{equation}\label{herm-car}   \left\vert h_n(x)^2 -\frac{2}{\pi \sqrt{2n+1-x^2}}\cos^2\left[(2n+1)\Phi \left(\frac{x}{\sqrt{2n+1}} \right)-\frac{n\pi}{2}\right] \right\vert \leq \frac{C(\beta)}{(2n+1)^{\frac{3}{2}}}  .   \end{equation}
C'est la formule précédente et la proposition \ref{hermi-majo} que nous allons utiliser pour comprendre la localisation spatiale de la fonction spectrale $e_d(n,\cdot)$ de l'oscillateur harmonique multidimensionnel.

\section{Preuve de la proposition \ref{ed-maj}, majoration \eqref{majo} de la fonction spectrale}

Il s'agit juste de suivre correctement les constantes dans \cite[Lemma 3.2.2]{thanga}
et d'appliquer la proposition \ref{hermi-majo}. Nous ne traitons que le cas $n$ pair
car le cas $n\in 2\N+1$ se traite de la même façon. La formule \eqref{edn-2} et la proposition \ref{estim-hn} donnent
\bema     e_d(n,x)=\sum_{k=0}^{\frac{n}{2}} h_{2k}(|x|)^2 e_{d-1}(n-2k,0)\simeq \sum_{k=0}^{ \frac{n}{2} } h_{2k}(|x|)^2 (1+n-2k)^{\frac{d}{2}-\frac{3}{2}}    .\enma
\textbf{Premier cas $|x|> \sqrt{2n+1}$.} On a 
\bema e_{d}(n,x)\leq C  e^{-\gamma x^2} \sum_{k=0}^{\frac{n}{2}}(1+n-2k)^{\frac{d}{2}-\frac{3}{2}} \leq C n^{\frac{d}{2}-1} e^{-\gamma x^2}.    \enma
\textbf{Deuxième cas $|x|\leq \sqrt{2n+1}$.}
Notant $k_x \in [0,\frac{n}{2}]$ le plus petit entier tel que $|x|\leq \sqrt{4k_x+1}$, on a 
\bemar \Sum_{k=k_x}^{\frac{n}{2} } h_{2k}(|x|)^2 (1+n-2k)^{\frac{d}{2}-\frac{3}{2}} & \lesssim & \Sum_{k=k_x}^{ \frac{n}{2}  } \frac{(1+n-2k)^{\frac{d}{2}-\frac{3}{2}}}{(4k+2-x^2)^{\frac{1}{2}}}  \\
& \lesssim &  \Sum_{k=k_x}^{\frac{n}{2}} \frac{(1+n-2k)^{\frac{d}{2}-\frac{3}{2}}}{(1+k-k_x)^{\frac{1}{2}}} \\ 
 & \lesssim & \Sum_{\ell=0}^{\frac{n}{2}-k_x} \frac{(n-2k_x-2\ell+1)^{\frac{d}{2}-\frac{3}{2}}}{\sqrt{1+\ell}} . \enmar 
 En distinguant selon que $\ell$ est plus petit ou plus grand que $\frac{1}{2}\left( \frac{n}{2}-k_x\right)$, on obtient 
 \bema     \Sum_{k=k_x}^{\frac{n}{2} } h_{2k}(|x|)^2 (1+n-2k)^{\frac{d}{2}-\frac{3}{2}}   \lesssim  (n-2k_x+1)^{\frac{d}{2}-1}\lesssim n^{\frac{d}{2}-1}. \enma
Pour l'autre somme, on écrit
\bema \sum_{k=0}^{k_x-1} h_{2k}(|x|)^2 (1+n-2k)^{\frac{d}{2}-\frac{3}{2}}\lesssim      e^{-\gamma |x|^2} \sum_{k=0}^{k_x-1}   (1+n-2k)^{\frac{d}{2}-\frac{3}{2}} \lesssim e^{-\gamma |x|^2} n^{\frac{d}{2}-\frac{1}{2}}.\enma

\section{Preuve de la proposition \ref{ed-maj}, minoration \eqref{mino} de la fonction spectrale}

Pour tout entier $d\geq 2$, on veut démontrer l'hypothèse de récurrence 
\bema
H(d): \quad \forall \alpha\in \left]0,\sin\left( \frac{1}{4}\right)\right[ \quad \exists C(d,\alpha)>1 \quad  \quad\liminf\limits_{n\rightarrow +\infty} \left( \inf\limits_{\frac{C(d,\alpha)}{\sqrt{2n+1}}\leq |x|\leq \alpha\sqrt{2n+1}} \frac{e_d(n,x)}{n^{\frac{d}{2}-1}}\right)>0.
\enma

Dans toute cette preuve, on utilisera la notation suivante 
\begin{equation}\label{defi-beta} \beta:=\frac{1}{2}\left(\alpha+\sin\left(\frac{1}{4}\right) \right).    \end{equation}
En outre, $o_{\alpha}(1)$ désignera une fonction de $(n,\alpha)$ qui vérifie 
$\lim\limits_{n\rightarrow +\infty} o_{\alpha}(1)=0$ pour tout $\alpha$ fixé à l'avance.

\textbf{Premier cas : $H(2)$ avec $n\in 2\N$.}
Il s'agit du cas le plus technique. Tout d'abord, l'invariance radiale de $e_2(n,\cdot)$ et la proposition \ref{estim-hn} donnent pour tout $x\in \R^2$
\bema    e_2(n,x)=\sum_{k=0}^{\frac{n}{2}} h_{2k}(|x|)^2 h_{n-2k}(0)^2  \simeq  \sum_{k=0}^{\frac{n}{2}} \frac{h_{2k}(|x|)^2 }{\sqrt{1+n-2k}}.   \enma
On va sélectionner les indices $k$ tels que $\frac{\alpha^2}{2\beta^2}n\leq k\leq \frac{n}{2}$. Puisque l'on a $\beta>\alpha$, on déduit que l'on a 
\bemar
n &\leq &  \frac{2\beta^2}{\alpha^2} k+  \frac{1}{2}\left(\frac{\beta^2}{\alpha^2}-1\right)\\[3mm]
2n+1& \leq &  \frac{\beta^2}{\alpha^2}(4k+1) \\[3mm]
\alpha\sqrt{2n+1} & \leq & \beta \sqrt{4k+1}.
\enmar
En utilisant \eqref{herm-car} et en imposant $|x|\leq \alpha\sqrt{2n+1}$, il existe une constante $C(\alpha)>1$ telle que
\bema e_2(n,x)  \gtrsim    \Sum_{\frac{\alpha^2}{2\beta^2}n \leq k \leq \frac{n}{2}} \Frac{1}{\sqrt{(1+n-2k)k}} \cos^2 
\left[(4k+1)\Phi \left(\Frac{|x|}{\sqrt{4k+1}} \right)\right] -\Sum_{\frac{\alpha^2}{2\beta^2}n \leq k \leq \frac{n}{2} } \Frac{C(\alpha)}{(4k+1)^{\frac{3}{2}}(1+n-2k)^{\frac{1}{2}}} , \enma
ce qui se minore grossièrement par 
\bema \Frac{1}{Cn}\Sum_{ \frac{\alpha^2}{2\beta^2}n \leq  k \leq \frac{n}{2}} \cos^2 
\left[(4k+1)\Phi \left(\Frac{|x|}{\sqrt{4k+1}} \right)\right]    -\Sum_{\frac{\alpha^2}{2\beta^2}n \leq k \leq \frac{n}{2}} \Frac{C(\alpha)}{(4k+1)^{\frac{3}{2}}(1+n-2k)^{\frac{1}{2}}}   .
 \enma
Or on a immédiatement 
\bema \Sum_{\frac{\alpha^2}{2\beta^2}n \leq k \leq \frac{n}{2}} \Frac{C(\alpha)}{(4k+1)^{\frac{3}{2}}(1+n-2k)^{\frac{1}{2}}} \leq C(\alpha)\frac{\sqrt{n}}{n^{3/2}}=\frac{C(\alpha)}{n} .  \enma
Notons à présent $\lceil \frac{\alpha^2}{2\beta^2}n \rceil$ le plus petit entier supérieur ou égal à $\frac{\alpha^2}{2\beta^2}n$ et $\lfloor \frac{n}{2}\rfloor$ le plus grand entier inférieur ou égal à $\frac{n}{2}$.
Si l'on définit
\bema S(n,x,\alpha,\beta) :=  \Sum_{\lceil \frac{\alpha^2}{2\beta^2}n \rceil \leq  k \leq \lfloor \frac{n}{2} \rfloor} \cos\left[2(4k+1)\Phi \left(\Frac{|x|}{\sqrt{4k+1}} \right)\right],\enma
alors notre minoration se reformule en 
\begin{equation}\label{minor-Sn} 
e_2(n,x) \gtrsim_\alpha  \Frac{1}{4}\left( 1-\frac{\alpha^2}{\beta^2}\right)+\Frac{1}{2n} S(n,x,\alpha,\beta)+o_\alpha(1).
\end{equation}
Pour gérer $S(n,x,\alpha,\beta)$, on fait appel à la formule d'Euler-Maclaurin au rang $0$ en posant
\bemar
\forall t\in \left[\frac{\alpha^2}{2\beta^2} n,\frac{n}{2} \right]\qquad    a_{x}(t)& :=& 2(4t+1)\Phi \left( \Frac{|x|}{\sqrt{4t+1}} \right)   =2(4t+1) \Int_{0}^{\frac{|x|}{\sqrt{4t+1}}} \sqrt{1-s^2}ds ,\\[3mm]
   a_{x}'(t) & = & 4 \arcsin\left(\Frac{|x|}{\sqrt{4t+1}} \right) \in [0,4\arcsin(\beta)] . \enmar
Cela nous permet de reformuler $S(n,x,\alpha,\beta)$ en
\bema \frac{\cos\left[a_x\left( \left\lceil \frac{\alpha^2}{2\beta^2}n \right\rceil \right)\right]+\cos\left[a_x\left( \left\lfloor \frac{n}{2} \right\rfloor\right)\right]}{2}+\int_{\lceil \frac{\alpha^2}{2\beta^2}n \rceil }^{\lfloor \frac{n}{2}\rfloor } \cos(a_{x}(t))dt -\int_{\lceil \frac{\alpha^2}{2\beta^2}n\rceil }^{\lfloor \frac{n}{2}\rfloor}
  a_{x}'(t) \sin(a_{x}(t)) \left[t-\lfloor t \rfloor -\frac{1}{2} \right] dt.
\enma
 La deuxième intégrale est facile à contrôler : 
\bema \left\vert \int_{\lceil \frac{\alpha^2}{2\beta^2}n\rceil }^{\lfloor \frac{n}{2}\rfloor}
  a_{x}'(t) \sin(a_{x}(t)) \left[t-\lfloor t \rfloor -\frac{1}{2} \right] dt  \right\vert   \leq \frac{4\arcsin(\beta)}{2}\left(\left\lfloor \frac{n}{2} \right\rfloor -\left\lceil \frac{\alpha^2}{2\beta^2}n \right\rceil\right) \leq n\left(1-\frac{\alpha^2 }{\beta^2} \right) \arcsin(\beta)                    .  \enma
En ce qui concerne la première intégrale, on utilise une intégration par parties 
et la concavité de $t\mapsto a_x(t)$ : 
\bemar
\Int_{\lceil \frac{\alpha^2}{2\beta^2}n \rceil }^{\lfloor \frac{n}{2}\rfloor } \cos(a_{x}(t))dt & = & \left[\Frac{\sin(a_x(t)) }{a_x'(t)} \right]_{\lceil \frac{\alpha^2}{2\beta^2}n \rceil}^{\lfloor\frac{n}{2} \rfloor} -\Int_{\lceil \frac{\alpha^2}{2\beta^2}n \rceil }^{\lfloor \frac{n}{2}\rfloor} \left(\frac{1}{a_x'(t)} \right)'\sin(a_x(t))dt\\[3mm]
\left\vert \Int_{\lceil \frac{\alpha^2}{2\beta^2}n \rceil }^{\lfloor \frac{n}{2}\rfloor } \cos(a_{x}(t))dt \right\vert & \leq & \left\vert\Frac{1}{a_x'(\lfloor\frac{n}{2}\rfloor)}+\frac{1}{a_x'\left( \lceil \frac{\alpha^2 n}{2\beta^2}\rceil \right)} \right\vert+ \left\vert 
\Frac{1}{a_x'( \lfloor \frac{n}{2} \rfloor )}-\frac{1}{a_x'\left( \lceil\frac{\alpha^2 n}{2\beta^2}\rceil\right)}\right\vert \\[4mm]
&   \leq & \Frac{2}{a_x'(\lfloor\frac{n}{2}\rfloor)} .
\enmar 
Et l'on a donc
\bema  \left\vert \int_{\lceil \frac{\alpha^2}{2\beta^2}n\rceil }^{\lfloor\frac{n}{2}\rfloor} \cos(a_x(t))dt\right\vert \leq       \frac{1}{2\arcsin \left(\frac{|x|}{\sqrt{2n+1}} \right)} \leq \frac{\sqrt{2n+1}}{2|x|}   . \enma
La définition \eqref{defi-beta} de $\beta$ nous autorise à majorer pour tout $n\in \N^\star$
\bemar    |S(n,x,\alpha,\beta)| &\leq & 1+ n\left(1-\Frac{\alpha^2}{\beta^2}\right)\arcsin(\beta)   +\Frac{\sqrt{2n+1}}{2|x|}  \\
 &  \leq & 1+ \Frac{n}{4}\left(1-\Frac{\alpha^2}{\beta^2}\right)  +\frac{2n}{|x|\sqrt{2n+1}}  . \enmar
On peut donc choisir $C(2,\alpha)>1$ de sorte que l'inégalité $\frac{C(2,\alpha)}{\sqrt{2n+1}}\leq |x|$ implique 
\bema
|S(n,x,\alpha,\beta)|\leq  1+\frac{n}{3}\left(1-\frac{\alpha^2}{\beta^2}\right)   .
\enma

On conclut à l'aide de \eqref{minor-Sn}.

\textbf{Deuxième cas : $H(2)$ et $n\in 2\N+1$.}
Par invariance radiale et par imparité des fonctions de Hermite $h_{n-2k}$, on obtient
\bema  e_2(n,x)=\sum_{k=1}^{\frac{n+1}{2} } h_{2k-1}(|x|)^2 h_{n-2k+1}(0)^2  \simeq   \sum_{k=1}^{ \frac{n+1}{2}  } h_{2k-1}(|x|)^2 \frac{1}{\sqrt{n-2k+2}}   .     \enma
En notant $K$ l'ensemble des indices $k$ vérifiant $ \frac{\alpha^2}{2\beta^2}n+\frac{\alpha^2+\beta^2}{4\beta^2}\leq k   \leq \frac{n+1}{2}$, on a 
\bemar
\forall k\in K \quad \frac{\alpha^2}{\beta^2}\left(\frac{n}{2}+\frac{1}{4} \right) & \leq & k-\frac{1}{4} \\ 
\alpha \sqrt{2n+1}& <& \beta \sqrt{4k-1}.
\enmar
En utilisant \eqref{herm-car} et la même argumentation que celle du premier cas, nous arrivons à 
 \bemar e_2(n,x)& \geq   &   \Frac{1}{C(\alpha)\sqrt{n}} \sum_{K}  \frac{1}{\sqrt{n-2k+2}}\sin^2\left[(4k-1)\Phi \left(\frac{x}{\sqrt{4k-1}} \right)\right]    +o_\alpha (1)             \\
& \geq & \Frac{1}{C(\alpha) n} \Sum_{K }  \sin^2\left[(4k-1)\Phi \left(\frac{x}{\sqrt{4k-1}} \right)\right]               +o_\alpha(1)     \\
& \gtrsim_\alpha & \Frac{1}{4}\left(1-\frac{\alpha^2}{\beta^2} \right)-\frac{1}{2n} \Sum_{K }  \cos\left[2(4k-1)\Phi \left(\frac{x}{\sqrt{4k-1}} \right) \right]   +o_\alpha(1)     .        \enmar
On finit exactement comme dans le premier cas.

\textbf{Troisième cas : $H(d)$ avec $d>2$ et $n\in 2\N$.}
On utilise le lemme \ref{estim-hn} et l'invariance radiale de $e_d$ exprimée par la formule \eqref{edn-1} pour obtenir
\bemar        e_{d}(n,x)& = & \Sum_{k=0}^{\frac{n}{2}} e_{d-1} (2k,\overbrace{(|x|,0,\dots,0)}^{\in \R^{d-1}}) h_{n-2k}(0)^2 \\
&   \simeq & \Sum_{k=0}^{\frac{n}{2}} \frac{e_{d-1} (2k,(|x|,0,\dots,0)) }{\sqrt{1+n-2k}}  \\
& \gtrsim &  \Sum_{\frac{\alpha^2}{2\beta^2}n\leq k\leq \frac{n}{2}} \frac{e_{d-1} (2k,(|x|,0,\dots,0)) }{\sqrt{1+n-2k}} .   \enmar
Comme dans le premier cas, les indices $k$ sélectionnés vérifient l'inégalité $\alpha\sqrt{2n+1}\leq \beta \sqrt{2(2k)+1}$.
Puisque l'on a $\beta<\sin\left( \frac{1}{4}\right)$, l'hypothèse de récurrence $H(d-1)$ nous fournit un nombre $C(d-1,\beta)>1$.
En imposant la condition $\frac{\beta C(d-1,\beta)}{\alpha\sqrt{2n+1}}\leq |x|\leq \alpha \sqrt{2n+1}$, on a 
 \bema
\frac{C(d-1,\beta)}{\sqrt{2(2k)+1}}\leq  |x|\leq \beta \sqrt{2(2k)+1}.
\enma
Cela nous mène à 
\bemar  e_{d}(n,x)& \gtrsim_{d,\beta} &  \Sum_{\frac{\alpha^2}{2\beta^2}n \leq k\leq \frac{n}{2}} \frac{k^{\frac{d-1}{2}-1}}{\sqrt{1+n-2k}} \\
&\gtrsim_{d,\beta} &  n^{\frac{d-1}{2}-1} \Sum_{\frac{\alpha^2}{2\beta^2}n \leq k\leq \frac{n}{2}}
\frac{1}{\sqrt{1+n-2k}} \\
& \gtrsim_{d,\beta} &  n^{\frac{d-1}{2}-1} \times \sqrt{n}=n^{\frac{d}{2}-1}.
     \enmar

\textbf{Dernier cas : $H(d)$ avec $d>2$ et $n\in 2\N+1$.} On se ramène à $H(d-1)$ comme dans le troisième cas.

\section{Preuve du théorème \ref{HLpps}, randomisation des fonctions de Hermite}

Nous avons déjà vérifié les hypothèses des théorème \ref{theo-dual} de dualité 
et théorème \ref{theo-interpo} d'interpolation (voir 
\eqref{herm-dual} et \eqref{herm-lor}). On a donc les points i) et ii).
Décomposant $u=\sum_{n\in \N} \Pi_n(u)\in \pl^p(\R^d,\oplus E_{d,n})$, on souhaite maintenant montrer que la norme $\left\Vert \cdot \right\Vert_{\pl^p(\R^d,\oplus E_{d,n}) }$ est équivalente à la suivante 
\bema
N(u)  :=  \left\Vert\Pi_{0}(u)\right\Vert_{L^2(\R^d)}+\left[\Sum_{n\geq 1} n^{\frac{d}{2}-1}\left( \Sum_{k \geq n} \frac{ \left\Vert\Pi_k(u)\right\Vert_{L^2(\R^d)}^2 }{k^{\frac{d}{2}}} \right)^{\frac{p}{2}} \right]^{\frac{1}{p}}.\enma

On va appliquer la proposition \ref{ed-maj}.
Quitte à augmenter l'entier $n(d)$, on peut supposer $C(d)\leq \sqrt{2n_d+1}$.
On vérifie facilement que la norme $\left\Vert \cdot\right\Vert_{\pl^p(\R^d,\oplus E_{d,n})}$ domine la norme $N$ :

\bemar
\left\Vert u\right\Vert_{\pl^p(\R^d,\oplus E_{d,n})}^p& \gtrsim_{d,p} & 
\Int_{\R^d} \left( \sum_{n\geq n(d)} (1+n)^{-(d-1)} \left\Vert \Pi_n(u) \right\Vert_{L^2(\R^d)}^2 e_d(n,x)\right)^{\frac{p}{2}}dx \\[3mm]
& \gtrsim_{d,p} &  \Int_{\R^d} \left( \sum_{n\geq n(d)} (1+n)^{-\frac{d}{2}} \left\Vert \Pi_n(u) \right\Vert_{L^2(\R^d)}^2 \pun_{1\leq |x|\leq \alpha_0 \sqrt{2n+1}} \right)^{\frac{p}{2}} dx  \\[3mm]
& \gtrsim_{d,p} & \Int_{\R^d} \sum_{n\geq n(d)} \left(\sum_{k> n} (1+k)^{-\frac{d}{2}} \left\Vert\Pi_k(u)\right\Vert_{L^2(\R^d)}^2 \right)^{\frac{p}{2}} \pun_{\sqrt{2n+1}<\frac{|x|}{\alpha_0}\leq \sqrt{2n+3}} dx \\
& \gtrsim_{d,p} & \Sum_{n> n(d)} n^{\frac{d}{2}-1}  \left(\Sum_{k \geq  n} (1+k)^{-\frac{d}{2}}\left\Vert\Pi_k(u)\right\Vert_{L^2(\R^d)}^2 \right)^{\frac{p}{2}}.
\enmar
Pour récupérer les premiers termes d'indice $n\leq n(d)$, il s'agit de remarquer les inégalités triviales
\bema  \forall n\in \N\cap [0,n(d)] \qquad     \left\Vert u\right\Vert_{\pl^p(\R^d,\oplus E_{d,n})} \gtrsim_{d,p}  (1+n)^{-\frac{(d-1)}{2}} \left\Vert \Pi_n(u) \right\Vert_{L^2(\R^d)} \left\Vert\sqrt{ e_d(n,\cdot)} \right\Vert_{L^{p}(\R^d)} .       \enma
On obtient alors facilement l'estimation $\left\Vert u\right\Vert_{\pl^p(\R^d,\oplus E_{d,n})} \gtrsim_{d,p} N(u)$.
Montrons maintenant l'estimation réciproque $N(u)\gtrsim_{d,p} \left\Vert u\right\Vert_{\pl^p(\R^d,\oplus E_{d,n})}$ à l'aide des deux fonctions
\bemar  A(u,x)& :=& \Sum_{n\geq 0} (1+n)^{-\frac{d}{2}}\left\Vert \Pi_n(u) \right\Vert_{L^2(\R^d)}^2 \pun_{0\leq |x|\leq \sqrt{2n+1}}   ,   \\
 B(u,x)& :=& e^{-\gamma |x|^2} \Sum_{n\geq 0} (1+n)^{-\frac{d}{2}}\left\Vert \Pi_n(u) \right\Vert_{L^2(\R^d)}^2  \pun_{\sqrt{2n+1}\leq |x|} . \enmar
Par une argumentation similaire à celle que nous venons d'employer, on vérifie l'équivalence
\bema       \int_{\R^d} A(u,x)^{\frac{p}{2}}dx \simeq_{d,p} N(u)^p   .  \enma
Nous affirmons que nous avons aussi
\bema\int_{\R^d} B(u,x)^{\frac{p}{2}}dx  \lesssim_{d,p} N(u)^p.\enma
En effet, nous pouvons écrire
\bemar    \Int_{\R^d} B(u,x)^{\frac{p}{2}}dx & = & \Int_{\R^d} 
e^{-\frac{p\gamma}{2}|x|^2} \sum_{n\geq 0} \left(\sum_{k=0}^n (1+k)^{-\frac{d}{2}} \left\Vert \Pi_k(u)\right\Vert_{L^2(\R^d)}^2 \right)^{\frac{p}{2}} \pun_{\sqrt{2n+1}\leq |x|\leq \sqrt{2n+3}}dx           \\
 & \lesssim_{d,p} & \Sum_{n\geq 0} (1+n)^{\frac{d}{2}-1} e^{-p\gamma n} \left(\sum_{k=0}^n (1+k)^{-\frac{d}{2}} \left\Vert \Pi_k(u)\right\Vert_{L^2(\R^d)}^2 \right)^{\frac{p}{2}} \\
&    \lesssim_{d,p} & \left\Vert u \right\Vert_{\mathcal{H}^{-\frac{d}{2}}(\R^d)}^{p}\Sum_{n\geq 0} (1+n)^{\frac{d}{2}-1}e^{-p\gamma n}\\
& \lesssim_{d,p} & \left\Vert u \right\Vert_{\mathcal{H}^{-\frac{d}{2}}(\R^d)}^{p}\\
&  \lesssim_{d,p} & N(u)^p. \enmar
D'après la proposition \ref{ed-maj}, nous pouvons conclure 
\bema  \left\Vert u\right\Vert_{\pl^p(\R^d,\oplus E_{d,n})}^p\lesssim_{d,p} \int_{\R^d} [A(u,x)+B(u,x)]^{\frac{p}{2}}dx  
\lesssim_{d,p} \int_{\R^d} A(u,x)^{\frac{p}{2}}+B(u,x)^{\frac{p}{2}} dx \lesssim_{d,p} N(u)^p.            \enma

\section{Preuve du théorème \ref{inje-sobo}, injections de Sobolev probabilistes hermitiennes}

Le théorème \ref{HLpps} permet par dualité de se ramener au cas $p\in [2,+\infty[$.

\subsection{L'inclusion stricte $\mathcal{H}^{-d\left(\frac{1}{2}-\frac{1}{p}\right) }(\R^d)\subset \hl$}

On invoque l'inégalité triangulaire dans $L^{\frac{p}{2}}(\R^d)$ et \eqref{hermLp} pour obtenir pour tout $u\in \mathcal{H}^{ -d\left(\frac{1}{2}-\frac{1}{p}\right) }(\R^d)$ :
\bemar
 \left\Vert \Sum_{n\geq 0}  \left\Vert \Pi_n(u) \right\Vert_{L^2(\R^d)}^2 \Frac{e_d(n,\cdot)}{(1+n)^{d-1}} \right\Vert_{L^{\frac{p}{2}}(\R^d)} & 
  \leq & 
\Sum_{n\geq 0}   \left\Vert \Pi_n(u) \right\Vert_{L^2(\R^d)}^2  \Frac{
\left\Vert e_d(n,\cdot) \right\Vert_{L^{\frac{p}{2}}(\R^d)}}{(1+n)^{d-1}}\\[3mm]
 & \leq & \Sum_{n\geq 0} \Frac{\left\Vert \Pi_n(u) \right\Vert_{L^2(\R^d)}^2}{(1+n)^{d\left(\frac{1}{2}-\frac{1}{p}\right)}} = \left\Vert u\right\Vert_{\mathcal{H}^{ -d\left(\frac{1}{2}-\frac{1}{p}\right) }(\R^d)}^2  . \enmar

Pour montrer que l'inclusion est stricte, on adopte la même stratégie que pour le corollaire \ref{Yn-injesobo} et l'on écrit un résultat analogue à la proposition \ref{lemm-incon}.

\begin{prop}\label{Lpcond}
Pour tout réel $p>2$ il existe un réel $r=r(p)>1$ et un ensemble infini $I\subset \N$ tels que, pour toute distribution tempérée $u=\Sum_{n\in I} \Pi_n(u)$ à support dans $I$, on a 
\bema   \sum_{n\in I} \left[\Frac{ \left\Vert \Pi_n(u) \right\Vert_{L^2(\R^d)}^2}{(1+n)^{d\left(\frac{1}{2}-\frac{1}{p} \right)}} \right]^{r}    <+\infty \quad \Rightarrow \quad u\in \hl.\enma
\end{prop}

\begin{demo}
Comme pour la proposition \ref{lemm-incon}, on donne deux preuves.
La première est duale de celle de la proposition \ref{lemm-incon} : on raisonne avec des restes de séries au lieu de sommes partielles et l'on part de l'équivalence de normes \eqref{harmo-cond}.
 Pour tout sous-ensemble infini $I\subset \N^\star$ et toute distribution tempérée $u$ on a

\bemar
\Sum_{n\geq 1} n^{\frac{d}{2}-1} \left(\sum_{k\geq n} k^{-\frac{d}{2}} \left\Vert \Pi_{k}(u)\right\Vert_{L^2(\R^d)}^2 \right)^{\frac{p}{2}}&
 \leq &
\Sum_{n\geq 1} n^{\frac{d}{2}-1} \left( \sum_{k\geq n}  \left\Vert\Pi_k(u)\right\Vert_{L^2(\R^d)}^p \right) \left(\sum_{\ell\geq n} \ell^{-\frac{dp}{2(p-2)}} \pun_{I}(\ell)\right)^{\frac{p-2}{2}}  \\[3mm]
& \leq & \Sum_{k\geq 1} \left\Vert \Pi_k(u)\right\Vert_{L^2(\R^d)}^p \sum_{n=1}^k n^{\frac{d}{2}-1}\left(\sum_{\ell\geq n} \ell^{-\frac{dp}{2(p-2)}} \pun_{I}(\ell)\right)^{\frac{p-2}{2}} .
\enmar 
Par récurrence, on peut construire un ensemble lacunaire infini $I\subset \N^\star$ tel que 
\bema    \forall n \in I \quad \sum_{\ell\geq n} \ell^{-\frac{dp}{2(p-2)}} \pun_{I}(\ell) \leq 2 n^{-\frac{dp}{2(p-2)}}   .  \enma
Pour un tel ensemble $I$, on a 
\bemar \Sum_{n\geq 1} n^{\frac{d}{2}-1} \left(\sum_{k\geq n} k^{-\frac{d}{2}} \left\Vert\Pi_{k}(u)\right\Vert_{L^2(\R^d)}^2 \right)^{\frac{p}{2}} & \leq & 2  \Sum_{k\geq 0} \left\Vert\Pi_k(u)\right\Vert_{L^2(\R^d)}^p \Sum_{n=1}^k n^{\frac{d}{2}-1-\frac{dp}{4}} \\
& \lesssim & \Sum_{k\geq 1} \left[k^{-d\left( \frac{1}{2}-\frac{1}{p}\right)} \left\Vert \Pi_k(u) \right\Vert_{L^2(\R^d)}^2 \right]^{\frac{p}{2}} .\enmar

Cela achève notre première preuve avec $r(p)=\frac{p}{2}$.

Notre seconde preuve est plus abstraite et utilise la norme originelle des espaces $\hl$ qui apparaît dans le lemme d'identification \ref{identif2}, à savoir 
\bema
\left\Vert u \right\Vert_{\pl^p(\R^d,\oplus E_{d,n})} = \left( \int_{\R^d} \left[ \sum_{n\in \N}\frac{ \left\Vert \Pi_n(u) \right\Vert_{L^2(\R^d)}^2}{\dim(E_{d,n})} e_d(n,x) dx\right]^{\frac{p}{2}}\right)^{\frac{1}{p}}.
\enma
L'expression précédente suggère d'analyser les sous-suites de $(e_d(n,\cdot))_{n\in \N}$ dans l'espace $L^{\frac{p}{2}}(\R^d)$.
Introduisons pour tout $n\in \N$ la fonction normalisée dans $L^{\frac{p}{2}}(\R^d)$
\bema f_n:=\frac{e_d(n,\cdot)}{\left\Vert e_d(n,\cdot)\right\Vert_{L^{\frac{p}{2}}(\R^d)}}.\enma
Par réflexivité de $L^{\frac{p}{2}}(\R^d)$, on peut extraire une sous-suite faiblement convergente de $(f_n)_{n\geq 0}$ au sens de la dualité  avec $L^{\frac{p}{p-2}}(\R^d)$.
Au niveau de l'analyse spectrale, on aura seulement besoin de l'estimation $\left\Vert e_d(n,\cdot) \right\Vert_{L^{\frac{p}{2}}(\R^d)} \simeq_{d,p} n^{\frac{d}{2}-1+\frac{d}{p}}$ (voir \eqref{hermLp}).
Contrairement au cas des variétés compactes, l'équivalent précédent est décroissant en $p$ et l'on préfère donc une estimation avec $\left\Vert e_d(n,\cdot)\right\Vert_{L^\infty(\R^d)}$ :
\bema
\Int_{\R^d} e_d(n,x)  e^{-|x|^2}dx   \leq   \left\Vert e_d(n,\cdot)\right\Vert_{L^\infty(\R^d)} \Int_{\R^d} e^{-|x|^2}dx  = o\left(\left\Vert e_d(n,\cdot) \right\Vert_{L^{\frac{p}{2}}(\R^d)} \right).
\enma
En notant $f_{\infty}\in L^{\frac{p}{2}}(\R^d)$ la limite faible de la sous-suite considérée de $(f_n)_{n\in \N}$, on a donc obtenu 
\bema 
0=\lim\limits_{n\rightarrow +\infty} \int_{\R^d} f_n(x)e^{-|x|^2} dx = \int_{\R^d} f_\infty(x)e^{-|x|^2}dx.
\enma
Puisque les fonctions $f_n$ sont à valeurs positives, il en est de même de $f_{\infty}$ et l'on a donc $f_\infty=0$.
D'après la proposition \ref{basincoLp}, on peut extraire de la suite $(f_n)_{n\in \N}$ une sous-suite $(f_n)_{n\in I}$ dans $L^{\frac{p}{2}}(\R^d)$ qui vérifie pour toute suite complexe $(\alpha_n)_{n\in I}$ à support fini
\bema
\left\Vert \sum_{n\in I} \alpha_n f_n \right\Vert_{L^{\frac{p}{2}}(\R^d)} \leq K \left( \sum_{n\in I} |\alpha_n|^{\min(2,\frac{p}{2})}\right)^{\frac{1}{\min(2,\frac{p}{2})} },   \enma
où $K\geq 1$ est indépendant de la suite $(\alpha_n)_{n\in I}$.
L'inégalité précédente est encore valide si la suite $(\alpha_n)_{n\in I}$ n'est pas à support fini mais constituée de termes positifs. 
On conclut alors avec 
\bema \alpha_n=\Frac{\left\Vert \Pi_n(u) \right\Vert_{L^2(\R^d)}^2 \left\Vert e_d(n,\cdot)\right\Vert_{L^{\frac{p}{2}}(\R^d)}}{\dim(E_{d,n})} \simeq_{d,p} \Frac{\left\Vert \Pi_n(u) \right\Vert_{L^2(\R^d)}^2}{(1+n)^{d\left( \frac{1}{2}-\frac{1}{p}\right)} }. \enma
\end{demo}

\subsection{L'inclusion stricte $\pl^p(\R^d,\oplus E_{d,n})\subsetneq \mathcal{H}^{-d(\frac{1}{2}-\frac{1}{p})-\ep}(\R^d)$}

Le coeur de la preuve réside dans l'inégalité suivante valide pour tout $u\in \pl^p(\R^d,\oplus E_{d,n}) $
\begin{equation}\label{inclu-stri1} 
\left\Vert u\right\Vert_{\HH^{-d\left( \frac{1}{2}-\frac{1}{p}\right)-\ep}(\R^d)} + \left( \sum_{n\geq 0} \left[ \Frac{\left\Vert \Pi_n(u) \right\Vert_{L^2(\R^d)}^{2}}{(1+n)^{d\left(\frac{1}{2}-\frac{1}{p}\right)}}\right]^{\frac{p}{2}}\right)^{\frac{1}{p}}  \leq C(p,d,\ep) \left\Vert u\right\Vert_{\pl^p(\R^d,\oplus E_{d,n})}.
\end{equation}
Sans perte de généralité, on peut supposer que l'on a $\Pi_n(u)=0$ pour $n\gg 1$.
La première partie de \eqref{inclu-stri1} est facile. Posons à cet effet 
\bema
\forall n\in \N \quad R_n:=\sum_{k\geq n}(1+k)^{-\frac{d}{2}} \left\Vert \Pi_k(u)\right\Vert_{L^2(\R^d)}^2.
\enma
On a alors
\bemar
\left\Vert u\right\Vert_{\HH^{-d\left(\frac{1}{2}-\frac{1}{p}\right)-\ep}(\R^d)}^2 & =& \Sum_{n\in \N} (1+n)^{-d\left(\frac{1}{2}-\frac{1}{p}\right)-\ep} \left\Vert \Pi_n(u) \right\Vert_{L^2(\R^d)}^2 \\
& = & \Sum_{n\in \N} (1+n)^{\frac{d}{p}-\ep} (R_n-R_{n+1}) \\
& \leq  & R_0+C(d,p,\ep)\Sum_{n\in \N}(1+n)^{\frac{d}{p}-\ep-1} R_n \\
& \leq & R_0+C(d,p,\ep)\Sum_{n\in \N} \Frac{1}{(1+n)^{\ep+\frac{p-2}{p}}}\times (1+n)^{\frac{d-2}{p}}R_n \\
& \leq & R_0+C(d,p,\ep) \left(\Sum_{n\in \N} (1+n)^{\frac{d}{2}-1}R_n^{\frac{p}{2}}\right)^{\frac{2}{p}} \\
& \leq & C(d,p,\ep) \left(\Sum_{n\in \N} (1+n)^{\frac{d}{2}-1}R_n^{\frac{p}{2}}\right)^{\frac{2}{p}}. 
\enmar 
On conclut avec \eqref{harmo-cond}.

Le deuxième terme est le plus simple à voir car exprime le fait que $L^p(\R^d)$ est un espace de Banach de cotype $\max(2,p)=p$ (voir \eqref{Lp-type}) : 
\bemar
\left\Vert u\right\Vert_{\pl^p(\R^d,\oplus E_{d,n})} &  = & \left\Vert \sqrt{\Sum_{n\in \N} \left\Vert \Pi_n(u) \right\Vert_{L^2(\R^d)}^2 \frac{e_d(n,x)}{\dim(E_{d,n})}} \right\Vert_{L_x^p(\R^d)} \\
& \geq & \left(  \Sum_{n\in \N}  \left\Vert \left\Vert \Pi_n(u) \right\Vert_{L^2(\R^d)} \frac{\sqrt{e_d(n,x)}}{\sqrt{\dim(E_{d,n})}}\right\Vert_{L_x^p(\R^d)}^p\right)^{\frac{1}{p}}.
\enmar
On obtient alors \eqref{inclu-stri1} grâce à \eqref{herm-dn} et \eqref{hermLp} :
\bema
\left\Vert u\right\Vert_{\pl^p(\R^d,\oplus E_{d,n})} \geq C(d,p)  \left( \sum_{n\geq 0} \left[ \Frac{\left\Vert \Pi_n(u) \right\Vert_{L^2(\R^d)}^{2}}{(1+n)^{d\left(\frac{1}{2}-\frac{1}{p}\right)}}\right]^{\frac{p}{2}}\right)^{\frac{1}{p}} .
\enma
On déduit aisément de l'inégalité précédente le fait que l'inclusion 
$\pl^p(\R^d,\oplus E_{d,n})\subset \mathcal{H}^{-d(\frac{1}{2}-\frac{1}{p})-\ep}(\R^d)$ est stricte.
Considérons par exemple un sous-ensemble infini $J\subset \N\backslash \{0,1\}$ tel que $\sum_{n\in J} \frac{1}{\ln(n)}$ converge.
A fortiori, pour tout $\ep>0$ la série $\sum_{n\in J}n^{-\ep}$ converge aussi.
On choisit une distribution $v\in \mathcal{S}'(\R^d)$ de manière lacunaire, à savoir 
$n^{-d\left( \frac{1}{2}-\frac{1}{p} \right)} \left\Vert \Pi_n v\right\Vert_{L^2(\R^d)}^2=1$ pour $n\in J$ et $0$ ailleurs. Il est clair que $v$ appartient à $\mathcal{H}^{-d\left( \frac{1}{2}-\frac{1}{p} \right)-\ep}(\R^d)$ mais n'appartient pas à $\pl^p(\R^d,\oplus E_{d,n})$.

\appendix 

\section{Preuve du théorème \ref{maurey-pi} : conséquence du théorème \ref{mplp}}\label{preuve-mp1}

Il s'agit d'examiner le théorème \ref{mplp} avec $d_n=1$ pour tout $n\in \N$.
La seule chose à justifier est l'équivalence ii) $\Leftrightarrow$ iii) dans le cas où les variables aléatoires $X_n:\Omega \rightarrow \R$ sont indépendantes et centrées (au lieu de symétriques).
On va invoquer des arguments classiques de symétrisation pour se ramener au cas où chaque variable aléatoire $X_n$ est symétrique.
Pour cela, on notera $(\widetilde{X}_n)_n$ une suite i.i.d. de copies indépendantes de la suite $(X_n)_n$. Ainsi, les variables aléatoires $X_n-\widetilde{X}_n$ sont symétriques et vérifient 
\bemar
\sup\limits_{n\in \N} \left\Vert X_n-\widetilde{X}_n\right\Vert_{L^{\max(2,p)}(\Omega)}& \leq & 2 \sup\limits_{n\in \N}\left\Vert X_n\right\Vert_{L^{\max(2,p)}(\Omega)}  <+\infty, \\
\inf\limits_{n\in \N} \E[|X_n-\widetilde{X}_n|] & \geq & \inf\limits_{n\in \N} \E[|X_n|]>0,
\enmar
où l'on a utilisé l'inégalité triangulaire (en intégrant en $\widetilde{X}_n$) pour la seconde inégalité.

Supposons iii). En particulier la série aléatoire $\sum \widetilde{X}_n u_n$ converge presque sûrement, et il en est donc de même de la série aléatoire $\sum (X_n- \widetilde{X}_n)u_n$. 
Le théorème \ref{maurey-pi}, avec $d_n=1$ pour tout $n\in \N$, assure que la série aléatoire $\sum \ep_n u_n$ converge presque sûrement dans $L^p(X)$. C'est-à-dire ii).

Supposons ii). De nouveau, le théorème \ref{maurey-pi} assure que la série aléatoire $\sum (X_n-\widetilde{X}_n)u_n$ converge dans l'espace $L^{\max(2,p)}(\Omega,L^p(X))$ et a fortiori dans $L^1(\Omega,L^p(X))$.
L'inégalité triangulaire donne à nouveau
\bema
\forall p>q\in \N \qquad \E\left[ \left\Vert \sum_{n=q}^p (X_n-\widetilde{X}_n) u_n \right\Vert_{L^p(X)} \right] \geq  \E\left[ \left\Vert \sum_{n=q}^p X_n u_n \right\Vert_{L^p(X)}\right].
\enma
L'inégalité de Markov nous permet d'affirmer que la série aléatoire $\sum X_n u_n$ converge en probabilité et donc presque sûrement (voir \cite[Theorem 6.1]{ledoux} ou \cite[Théorème II.3]{queff}), c'est-à-dire ii).

\section{Preuve du théorème \ref{maurey-pi} : optimalité de l'exposant $\max(2,p)$}\label{preuve-mp2}

Examinons l'espace mesuré $X=\N\backslash\{0,1\}$ muni de la mesure de comptage si bien que l'on a $L^p(X)=\ell^p(\N\backslash\{0,1\})$.
Pour tout réel $p\in [2,\infty[$, on note $(X_{n,p})_{n\geq 2}$ une suite i.i.d. de variables aléatoires symétriques, réelles et vérifiant
\bema  \forall t\gg 1 \quad  \P[|X_{n,p}|\geq t] = \frac{\ln(t)}{t^p}.   \enma
Il est clair que $X_{n,p}$ n'a pas de moment d'ordre $p$ et a des moments d'ordre $q\in [1,p[$. On a aussi
\begin{equation}\label{Xnps}
\mbox{p.s. }  \qquad \sup\limits_{n\geq 2}\frac{|X_{n,p}|}{n^{\frac{1}{p}} \ln^{\frac{2}{p}}(n)}=+\infty.
\end{equation}
En effet, il s'agit de remarquer que pour tout entier $K\in \N$, on a 
\bema
\sum_{n\geq 2} \P\left[|X_{n,p}|\geq K n^{\frac{1}{p}} \ln^{\frac{2}{p}}(n) \right] =+\infty.
\enma
Par indépendance des variables $X_{n,p}$, il existe presque sûrement une infinité d'entiers $n\geq 2$ tels que $|X_{n,p}|\geq Kn^{\frac{1}{p}} \ln^{\frac{2}{p}}(n)$. 
On en déduit facilement \eqref{Xnps}.
Revenons à l'optimalité de l'exposant $\max(2,p)$.

\textbf{Cas $p\in [1,2]$.} On fixe $u$ une suite non nulle appartenant à $\ell^p(\N\backslash\{0,1\})$ et l'on examine les deux séries aléatoires dans $\ell^p(\N\backslash\{0,1\})$
\begin{equation}\label{contrex}
\sum \frac{\ep_n }{\sqrt{n}\ln(n)}u     \quad \mbox{et} \quad \sum \frac{X_{n,2}}{\sqrt{n}\ln(n)} u,
\end{equation}
La série aléatoire $\sum \frac{\ep_n}{\sqrt{n}\ln(n)}$ converge presque sûrement dans $\R$.
Il s'ensuit que la première série aléatoire dans \eqref{contrex} converge presque sûrement dans $\ell^p(\N\backslash\{0,1\})$.
La divergence presque sûre de la seconde série aléatoire dans \eqref{contrex} découle de \eqref{Xnps}.

\textbf{Cas $p\in [2,+\infty[$.} Une démarche similaire est valide en examinant les deux séries aléatoires 
\bema 
\Sum \frac{\ep_n}{n^{\frac{1}{p} }\ln^{\frac{2}{p}}(n)} w_n \quad \mbox{et} \quad \Sum \frac{X_{n,p}}{n^{\frac{1}{p} }\ln^{\frac{2}{p}}(n)} w_n,
\enma
où $w_n=(0,\dots,0,1,0,\dots)\in \ell^p(\N\backslash\{0,1\})$ est la suite qui admet $1$ à la position $n$ et $0$ ailleurs.
Il est clair que la première série converge de manière déterministe dans $\ell^p(\N\backslash\{0,1\})$.
De nouveau, \eqref{Xnps} implique la divergence presque sûre de la seconde série aléatoire.

\section{Type, cotype et théorème de Maurey-Pisier}\label{cotype}

Nous rappelons des résultats standards et un énoncé essentiellement optimal du théorème de randomisation universelle de Maurey-Pisier.
On consultera principalement \cite[Chapter e]{linden1} et \cite[Chapitre 4]{queff} pour la théorie générale ainsi que \cite[corollaire 1.3]{maurey-pisier76} et \cite[annex]{randomh} pour le théorème de Maurey-Pisier. Par commodité, on considérera des espaces de Banach complexes.

\begin{defi}\label{def-typ}
Un espace de Banach $B$ est de type $p\in [1,2]$ s'il existe $c>0$ et $r\in [1,+\infty[$ tel que pour toute suite finie $(u_n)$ de $B$ on a 
\bema
\E\left[\left\Vert \Sum_{n} \ep_n u_n \right\Vert_{B}^r \right]^{\frac{1}{r}}\leq c \left(\sum_{n} \left\Vert u_n \right\Vert_{B}^p    \right)^{\frac{1}{p}}  .
\enma
L'espace de Banach $B$ est de cotype $q\in [2,+\infty[$ s'il existe $c>0$ et $r\in [1,+\infty[$ tel que pour toute suite finie $(u_n)$ de $B$ on a 
 \bema
\left(\sum_{n} \left\Vert u_n \right\Vert_{B}^q  \right)^{\frac{1}{q}} \leq c \E\left[\left\Vert \Sum_{n} \ep_n u_n \right\Vert_{B}^r \right]^{\frac{1}{r}}.
\enma
\end{defi}
La définition précédente est usuellement énoncée pour $r=1$ ou $r=2$ mais les inégalités de Kahane-Khintchine \eqref{KK} assurent que l'on peut, de façon équivalente, choisir tout réel $r\in [1,+\infty[$.
Tout espace de Banach est à la fois de type $1$ et de cotype $+\infty$ (avec le changement évident) et tout espace de Hilbert est à la fois de type $2$ et de cotype $2$.
En ce qui concerne l'espace $L^p(X)$ avec $p\in [1,+\infty[$, le théorème de Fubini et les inégalités de Kahane-Khintchine \eqref{KK} permettent de montrer l'équivalence
\bema
\forall (u_0,\dots,u_N)\in (L^p(X))^{N+1} \qquad \E\left[ \left\Vert \sum_{n=0}^N  \ep_n u_n \right\Vert_{L^p(X)} \right]\simeq_p \left\Vert \sqrt{\sum_{n=0}^N |u_n|^2 }\right\Vert_{L^p(X)}.
\enma
En distinguant selon que $p$ est plus petit ou plus grand que $2$, on peut démontrer de façon élémentaire que l'espace $L^p(X)$ est de type $\min(2,p)$ et de cotype $\max(2,p)$ : 
\begin{equation}\label{Lp-type}
\left(\sum_{n=0}^N \left\Vert u_n \right\Vert_{L^p(X)}^{\max(2,p)} \right)^{\frac{1}{\max(2,p)}}\leq \left\Vert \sqrt{\sum_{n=0}^N |u_n(x)|^2 }\right\Vert_{L_x^p(X)}\leq 
\left(\sum_{n=0}^N \left\Vert u_n \right\Vert_{L^p(X)}^{\min(2,p)} \right)^{\frac{1}{\min(2,p)}}.
\end{equation}

Rappelons que $(g_n)_{n\in \N}$ désigne une suite i.i.d. de variables aléatoires suivant une loi $\mathcal{N}_\C(0,1)$. Pour tout espace de Banach $B$ et pour toute suite $(u_n)_{n\in \N}$ de $B$, on a l'implication 
\begin{equation}\label{gn-prin-con}
\sum g_n u_n \quad \mbox{converge p.s. dans } B \quad \Rightarrow \quad 
\sum \ep_n u_n \quad \mbox{converge p.s. dans } B.
\end{equation}
Il s'agit là d'un fait non trivial qui découle par exemple des propriétés d'intégrabilité des variables gaussiennes, d'un principe de contraction de Hoffman-Jorgensen (voir \cite[annex]{randomh}) ou encore d'un principe de contraction de Jain-Marcus \cite[Theorem 5.1]{marcusjain}.

On dit qu'un espace de Banach $B$ est de cotype fini s'il existe $q\in [2,+\infty[$ tel que $B$ soit de cotype $q$.
L'intérêt des espaces de cotype fini est que ce sont précisément les bons espaces de Banach pour lesquels la randomisation est universelle \cite[corollaire 1.3]{maurey-pisier76}.

\begin{theo}[Maurey-Pisier]\label{theo-maupign}
Un espace de Banach $B$ est de cotype fini si et seulement si pour toute suite $(u_n)_{n\in \N}$ de $B$, la convergence de la série aléatoire $\sum \ep_n u_n$ équivaut à celle de la série aléatoire $\sum g_n u_n$.
\end{theo}

Examinant la preuve du théorème précédent, on obtient la version quantitative suivante.

\begin{theo}[Maurey-Pisier]\label{theo-maupi}
Soit $B$ un espace de Banach de cotype $q\in [2,+\infty[$, considérons de plus une suite $(X_n)_{n\in \N}$ de variables aléatoires réelles, centrées, indépendantes et vérifiant 
\begin{enumerate}[H1) ]
\item $\inf\limits_{n\in \N} \E[|X_n|]>0$,
\item il existe $\ep>0$ tel que $\sup\limits_{n\in \N} \E[|X_n|^{q+\ep}]<+\infty$.
\end{enumerate}
Alors pour toute suite $(u_n)_{n\in \N}$ de $B$ les deux propriétés suivantes sont équivalentes 
\begin{enumerate}[i) ]
\item la série aléatoire $\sum \ep_n u_n$ converge presque sûrement dans $B$,
\item la série aléatoire $\sum X_n u_n$ converge presque sûrement dans $B$.
\end{enumerate}
\end{theo}
\begin{demo}
Le sens i) $\Rightarrow$ ii) est prouvé dans \cite[corollaire 1.3]{maurey-pisier76} et n'utilise que l'hypothèse H2.
De façon précise, la preuve de \cite[corollaire 1.3]{maurey-pisier76} donne une estimation de la forme 
\begin{equation}\label{mp-ineg}
\exists C>0\quad \forall k>\ell\in \N \quad \E\left[ \left\Vert\sum_{n=\ell}^k X_n u_n\right\Vert_{B}^2 \right] \leq C \left(\sup\limits_{N\in \N} \E[|X_N|^{q+\ep}] \right)\E\left[ \left\Vert\sum_{n=\ell}^k \ep_n u_n\right\Vert_{B}^2\right].
\end{equation}
Par suite, on invoque le fait que la convergence presque sûre de la série de Rademacher $\sum \ep_n u_n$ équivaut à sa convergence dans $L^2(\Omega,B)$ (il s'agit d'une conséquence des inégalités de Kahane-Khintchine \eqref{KK}, voir par exemple \cite[page 128]{queff}).
L'inégalité \eqref{mp-ineg} implique la convergence dans $L^2(\Omega,B)$, et a fortiori en probabilité, de la série aléatoire $\sum X_n u_n$.
Cela implique aussi la convergence presque sûre de $\sum X_n u_n$ (voir \cite[page 119, théorème II.3]{queff} ou \cite[Theorem 6.1]{ledoux}).

En ce qui concerne le sens ii) $\Rightarrow$ i), on pourrait refaire le même argument que celui de la preuve de la proposition \ref{part2} en exploitant l'hypothèse H2. 
De même que pour l'implication \eqref{gn-prin-con}, on préfère invoquer un principe de contraction dû à Jain et Marcus.
D'abord, quitte à remplacer $X_n$ par sa variable symétrisée, on peut supposer que chaque variable aléatoire $X_n$ est elle-même symétrique. L'inégalité de Paley-Zygmund donne
\bema   \inf\limits_{n\in \N}\P(|X_n|\geq c) \geq  \inf\limits_{n\in \N} \frac{\E[|X_n|]^2}{4\E[|X_n|^2]}   >0,\qquad c:=\frac{1}{2}\inf\limits_{k\in \N} \E[|X_k|] ,\enma
puis d'après \cite[Theorem 5.1]{marcusjain}, on sait que l'assertion  $\inf\limits_{n\in \N}\P(|X_n|\geq c) >0$  force l'implication ii) $\Rightarrow$ i).
\end{demo}

\section{Randomisation dans un treillis de Banach de cotype fini}\label{treillis}

On explique comment les arguments développés dans la partie \ref{part-mplp} permettent d'interpréter l'universalité de la randomisation multidimensionnelle en termes de treillis complet de cotype fini.
Pour la notion de treillis de Banach, on se réfère au livre \cite{linden2}.
On notera respectivement $\leq$, $\vee$ et $|\cdot|$ la relation d'ordre, la borne supérieure et la valeur absolue.
Pour tout $N\in \N^\star$, on notera aussi $\overline{\mathcal{H}_N}$ l'espace vectoriel réel des fonctions 
continues et $1$-homogènes sur $\R^N$.
En particulier, $\phi_i:(t_1,\dots,t_N)\mapsto t_i$ désignera la projection de $\R^N$ sur la $i$-ème coordonnée.
Le théorème suivant assure l'existence d'un calcul fonctionnel dans un treillis de Banach basé sur $\overline{\mathcal{H}_n}$ (voir \cite[Theorem 1d1]{linden2}).

\begin{theo}\label{fonc-ho}
Soit $B$ un treillis de Banach et fixons $f_1,\dots,f_N$ des éléments de $B$. Il existe une unique application linéaire
\bemar
\overline{\mathcal{H}_N} & \rightarrow & B \\
F & \mapsto & F(f_1,\dots,f_N)
\enmar
qui vérifie les deux propriétés suivantes : 
\begin{enumerate}[i) ]
\item $\phi_i(f_1,\dots,f_N)=f_i$ pour tout entier $i\in [1,N]$,
\item pour tout $(F_1,F_2)\in \overline{\mathcal{H}_N}^2$ on a $\max(F_1,F_2)(f_1,\dots,f_N)=F_1(f_1,\dots,f_N)\vee F_2(f_1,\dots,f_N)$.
\end{enumerate}
Dans ce cas, en notant $f=|f_1|\vee \dots \vee |f_N|$, on a l'estimation de continuité
\begin{equation}\label{conc-calc}
\left\Vert F(f_1,\dots,f_N)\right\Vert_B \leq \left\Vert  f \right\Vert_B
\sup\limits_{\left\Vert t\right\Vert_{\infty}\leq 1} |F(t_1,\dots,t_N)|. 
\end{equation}
\end{theo}

Rappelons maintenant la définition suivante.
\begin{defi}
Un treillis de Banach $B$ est $q$-concave, avec $q\in [1,+\infty[$, s'il existe un réel $M_{(q)}(B)>0$ tel que l'on a l'inégalité suivante 
pour tout entier $N\in \N^\star$ et tout élément $(f_1,\dots,f_N)\in B^N$ :
\bema   \left(\sum_{i=1}^N \left\Vert f_i \right\Vert_B^q \right)^{\frac{1}{q}}   \leq M_{(q)}(B)\left\Vert  \left(\sum_{i=1}^N |f_i|^q \right)^{\frac{1}{q}}\right\Vert_B  , \enma
où 
le terme $\left(\sum_{i=1}^N |f_i|^q \right)^{\frac{1}{q}}\in B$ est défini par calcul fonctionnel.
\end{defi}

La notion de $q$-concavité est reliée à celle de cotype comme le montre le résultat suivant \cite[Proposition 1f3, Corollary 1f9, page 100]{linden2}.

\begin{prop}
Soient $B$ un treillis de Banach et un réel $q\in [2,+\infty[$, on a
\begin{enumerate}[i) ]
\item si $B$ est $q$-concave, alors $B$ est de cotype $q$,
\item si $B$ est de cotype $q\in [2,+\infty[$ alors $B$ est $r$-concave pour tout $r\in]q,+\infty[$.
\end{enumerate}
Par conséquent, on a 
\begin{equation}\label{born-inf}
\inf\{ q\geq 2, \quad B \mbox{ est } q\mbox{-concave}\}=\inf\{ q\geq 2,\quad  B \mbox{ est de cotype } q\}.
\end{equation}
\end{prop}

Pour tous réels $q\geq p\geq 1$, on vérifie que $L^p(\R)$ est $q$-concave (il s'agit d'interpoler l'injection canonique $L^p(\R,\ell^q(\N))\rightarrow \ell^q(\N,L^p(\R))$ entre $q=p$ et $q=\infty$) et la borne inférieure \eqref{born-inf} vaut $\max(2,p)$.
Un autre exemple intéressant est fourni par les espaces de Lorentz $L^{p,1}(\R)$ qui sont de cotype $\max(2,p)$ pour $p\in ]1,2[\cup]2,+\infty[$ (voir les calculs exacts dans \cite{creekmore1981type}).
Afin de prouver un théorème de type Maurey-Pisier dans la catégorie des treillis de Banach, nous aurons besoin du lemme suivant.
\begin{lemm}
Considérons des variables aléatoires $X_1,\dots,X_N$ appartenant à $L^1(\Omega)$
 et des éléments $f_1,\dots,f_N$ d'un treillis de Banach $B$.
Alors on a ``l'inégalité triangulaire'' :
\begin{equation}\label{sans-p-con}
\left\Vert \E_\omega\left[ \left\vert \sum_{n=1}^N  X_n(\omega) f_n\right\vert\right] \right\Vert_{B}\leq \E_\omega \left[ \left\Vert \sum_{n=1}^N  X_n(\omega) f_n \right\Vert_{B}\right].
\end{equation}
S'il existe un réel $q\in [1,+\infty[$ tel que $B$ soit $q$-concave et que les variables aléatoires $X_1,\dots,X_N$ appartiennent à $L^q(\Omega)$, alors on a aussi
\begin{equation}\label{p-conc-esp}
\E_{\omega} \left[\left\Vert\sum_{n=1}^N X_n(\omega) f_n \right\Vert_{B}^q\right]^{\frac{1}{q}} \leq M_{(q)}(B)\left\Vert \E_\omega \left[ \left\vert \sum_{n=1}^N X_n(\omega) f_n \right\vert^q\right]^{\frac{1}{q}} \right\Vert_{B},
\end{equation}
où les espérances dans le membre gauche de \eqref{sans-p-con} et dans le membre droit de \eqref{p-conc-esp} sont définies par calcul fonctionnel sur les $N$ variables $f_1,\dots,f_N$.
\end{lemm}

\begin{demo}
On commence par \eqref{p-conc-esp}.
Pour tout entier $n\in [1,N]$, on note $(X_{n,k})_{k\in \N}$ une suite de variables aléatoires qui prend un nombre fini de valeurs et qui converge dans $L^q(\Omega)$ vers $X_n$.
Pour tout $k\in \N$, il existe donc une partition finie 
\bema
\Omega=\bigsqcup_{\ell=1}^{L} \Omega_{k,\ell}
\enma
en parties mesurables telle que chaque variable aléatoire $X_{n,k}$ prend une valeur fixe, disons $x_{n,k,\ell}$, sur $\Omega_{k,\ell}$. La $q$-concavité de $B$ nous donne alors
\bemar
\E_{\omega} \left[\left\Vert\Sum_{n=1}^N X_{n,k}(\omega) f_n \right\Vert_{B}^q\right]^{\frac{1}{q}} & = & \left(\Sum_{\ell=1}^{L} \P(\Omega_{k,\ell}) \left\Vert \Sum_{n=1}^N x_{n,k,\ell} f_n \right\Vert_{B}^q\right)^{\frac{1}{q}}, \\[5mm]
& \leq & M_{(q)}(B) \left\Vert  \left( \Sum_{\ell=1}^L \P(\Omega_{k,\ell}) \left\vert \sum_{n=1}^N x_{n,k,\ell} f_n \right\vert^q \right)^{\frac{1}{q}}\right\Vert_{B}   .
\enmar
Or on a évidemment pour tout $(t_1,\dots,t_N)\in \R^N$
\begin{equation}\label{desn}
\left( \Sum_{\ell=1}^L \P(\Omega_{k,\ell}) \left\vert \sum_{n=1}^N x_{n,k,\ell} t_n \right\vert^q \right)^{\frac{1}{q}} = \E\left[\left\vert \Sum_{n=1}^N X_{n,k} t_n\right\vert^q \right]^{\frac{1}{q}}.
\end{equation}
Si $k$ tend vers $+\infty$, l'expression précédente converge uniformément sur chaque compact de $\R^N$ vers 
\bema 
\E\left[\left\vert \Sum_{n=1}^N X_{n} t_n\right\vert^q \right]^{\frac{1}{q}}.
\enma
Le calcul fonctionnel du théorème \ref{fonc-ho} assure que l'égalité \eqref{desn} est encore valide en substituant $(f_1,\dots,f_N)$ à $(t_1,\dots,t_N)$, ce qui nous donne
\begin{equation}\label{densit}
\E \left[\left\Vert\Sum_{n=1}^N X_{n,k} f_n \right\Vert_{B}^q\right]^{\frac{1}{q}}  \leq    M_{(q)}(B)\left\Vert \E\left[\left\vert \sum_{n=1}^N X_{n,k} f_n \right\vert^q\right]^{\frac{1}{q}}\right\Vert_{B}   .
\end{equation}
L'estimation \eqref{conc-calc} de continuité du calcul fonctionnel assure aussi que le membre droit de \eqref{densit} tend vers 
\bema 
M_{(q)}(B)\left\Vert\E\left[\left\vert \Sum_{n=1}^N X_{n} f_n\right\vert^q \right]^{\frac{1}{q}}\right\Vert_{B}.
\enma
Enfin, il est clair que le membre gauche de \eqref{densit} tend vers le membre gauche de \eqref{p-conc-esp}.
Pour démontrer \eqref{sans-p-con}, on refait la même démarche de densité dans $L^1(\Omega)$ à l'aide de l'inégalité triangulaire dans $B$ : 
\bemar
\E_{\omega} \left[\left\Vert\Sum_{n=1}^N X_{n,k}(\omega) f_n \right\Vert_{B}\right] & = & \Sum_{\ell=1}^{L} \P(\Omega_{k,\ell}) \left\Vert \Sum_{n=1}^N x_{n,k,\ell} f_n \right\Vert_{B}, \\[5mm]
& \geq & \left\Vert  \Sum_{\ell=1}^L \P(\Omega_{k,\ell}) \left\vert \sum_{n=1}^N x_{n,k,\ell} f_n \right\vert \right\Vert_{B}   .
\enmar
\end{demo}

\begin{coro}\label{treill-coro}
Considérons des matrices aléatoires $M_0\in \mathcal{M}_{d_1}(\R)$,...,$M_N\in \mathcal{M}_{d_N}(\R)$ à coefficients dans $L^1(\Omega)$ et des matrices $b_0\in\mathcal{M}_{d_1}(B)$,...,$b_{N}\in \mathcal{M}_{d_N}(B)$ dont les coefficients appartiennent à un treillis de Banach $B$.
Alors on a l'inégalité : 
\bema   \left\Vert \E\left[\left\vert  \sum_{n=0}^N \sqrt{d_n} \tr(M_n b_n)\right\vert\right]\right\Vert_{B} \leq 
\E\left[ \left\Vert \sum_{n=0}^N \sqrt{d_n} \tr(M_n b_n) \right\Vert_{B} \right]. \enma

S'il existe de plus un réel $q\in [1,+\infty[$ tel que $B$ soit $q$-concave et tel que les coefficients des matrices aléatoires $M_0,\dots,M_N$ appartiennent à $L^q(\Omega)$, alors on a 

\bema  \E\left[ \left\Vert \sum_{n=0}^N \sqrt{d_n} \tr(M_n b_n) \right\Vert_{B}^{q} \right] ^{\frac{1}{q}}
\leq M_{(q)}(B) \left\Vert \E\left[\left\vert  \sum_{n=0}^N \sqrt{d_n} \tr(M_n b_n)\right\vert^{q}\right]^{\frac{1}{q}}\right\Vert_{B} . \enma
\end{coro}

Il nous reste à remarquer que les inégalités de Kahane-Khintchine-Marcus-Pisier permettent d'étendre le théorème de Maurey \cite[Theorem 1.d.6]{linden2} au cas multidimensionnel.
\begin{theo}\label{treill-maur}
Soit $B$ un treillis de Banach $q$-concave, avec $q\in [1,+\infty[$ et considérons une suite de matrices $b_n \in \mathcal{M}_{d_n}(B)$. Alors, pour tout entier $N\in \N$, on a 
\begin{equation}\label{maurey-K}
\frac{1}{C} \left\Vert \sqrt{\sum_{n=0}^N     \left\Vert b_n \right\Vert_{HS}^2 } \right\Vert_B
\leq  \E\left[ \left\Vert \sum_{n=0}^N \sqrt{d_n} \tr(\mathcal{E}_n b_n) \right\Vert_B \right]\leq C(B,q)\left\Vert \sqrt{\sum_{n=0}^N     \left\Vert b_n \right\Vert_{HS}^2 } \right\Vert_B ,\end{equation}
où l'élément $\sqrt{\sum_{n=0}^N     \left\Vert b_n \right\Vert_{HS}^2 }\in B$ est défini par calcul fonctionnel.
\end{theo}
\begin{demo}
La minoration est vraie sans hypothèse de $q$-concavité. 
La version scalaire des inégalités de Kahane-Khintchine-Marcus-Pisier (voir \eqref{KKMP} et \eqref{KKMP-R}) rend triviale l'inégalité suivante lorsque les matrices $b_n$ sont à coefficients réels
\bema
\sqrt{\sum_{n=0}^N \left\Vert b_n \right\Vert_{HS}^2 } \leq  K_{2,1}\E\left[ \left\vert \sum_{n=0}^N \sqrt{d_n} \tr(\mathcal{E}_n b_n) \right\vert\right]  . 
\enma
L'inégalité précédente s'étend par calcul fonctionnel lorsque les coefficients des matrices $b_n$ sont à coefficients dans $B$.
La minoration de \eqref{maurey-K} découle alors du corollaire \ref{treill-coro}.
On raisonne de même pour la majoration.
\end{demo}

Le corollaire \ref{treill-coro} et le théorème \ref{treill-maur} nous permettent de reprendre mutatis mutandis les arguments de la partie \ref{part-mplp} afin de prouver le résultat suivant qui éclaire le théorème \ref{LPtr}.

\begin{theo}
Fixons un réel $q\in [2,+\infty[$, un treillis de Banach $B$ $q$-concave et considérons une suite de matrices aléatoires $M_n:\Omega \rightarrow \mathcal{M}_{d_n}(\R)$ indépendantes, orthogonalement invariantes et vérifiant
\bema  \sup\limits_{n\in \N} \E \left[\left\Vert M_n \right\Vert_{op}^{q}\right]<+\infty \quad \mbox{et}  \quad \inf\limits_{n\in \N} \sigma(\E[|M_n|])>0. \enma
Pour toute suite de matrices $b_n\in\mathcal{M}_{d_n}(B)$, on a l'équivalence des propriétés suivantes : 
\begin{enumerate}[i) ]
\item les normes $\left\Vert\sqrt{\sum_{n=0}^N \left\Vert b_n \right\Vert_{HS}^2}\right\Vert_{B}$ sont majorées indépendamment de $N\in \N$,
\item la série aléatoire $\sum \sqrt{d_n} \tr (M_n b_n)$ converge presque sûrement dans $B$,
\item la série aléatoire $\sum \sqrt{d_n} \tr (M_n b_n)$ converge dans $L^{q}(\Omega,B)$,
\item la série aléatoire $\sum \sqrt{d_n} \tr (M_n b_n)$ est bornée en probabilité dans $B$.
\end{enumerate}
\end{theo}

Tout comme à la fin de la partie \ref{part-mplp}, le théorème précédent implique une version du théorème \ref{mplp} où $L^p(X)$ est remplacé par un espace fonctionnel muni d'une structure de treillis de Banach de cotype fini, par exemple l'espace de Lorentz $L^{p,1}(X)$.
\section{Suites basiques inconditionnelles et preuve de la proposition \ref{basincoLp}}\label{incond}

On se réfère par exemple à \cite[parties 1.I, 2.I et 2.II]{queff} ou \cite[Chapter 1]{linden1}.
On notera $B$ un espace de Banach complexe et l'on pose la définition :
\begin{defi}
Une suite $(u_n)_{n\in \N}$ de $B$ est semi-normalisée si elle vérifie
\bema
0<\inf\limits_{n\in \N} \left\Vert u_n \right\Vert_{B} \leq \sup\limits_{n\in \N} \left\Vert u_n \right\Vert_B <+\infty.
\enma
\end{defi}

On va rappeler les propriétés importantes des sous-suites semi-normalisées qui convergent faiblement vers $0$ dans $B$ (``semi-normalized weakly null sequences").
Remarquons que $B$ est nécessairement de dimension infinie pour que cette notion ait un sens.
Cela va notamment se traduire par le fait que ces suites ont des sous-suites infinies constituées de termes linéairement indépendants (et même beaucoup plus).

\begin{defi}
Une suite $(u_n)_{n\in \N}$ de $B$ est une base de Schauder, si pour tout élément $b\in B$ il existe une \textbf{unique} suite numérique notée $(u_n^\star(b))_{n\in \N}$
telle que la série $\sum u_n^\star(b) u_n$ converge vers $b$.
La suite $(u_n)_{n\in \N}$ est dite basique si elle est une base de Schauder du sous-espace fermé de $B$ engendré par les éléments $u_n$.
\end{defi}

Venons-en maintenant à la notion de convergence inconditionnelle.

\begin{prop}\label{prop-inco}
Soit $(u_n)_{n\in \N}$ une suite de $B$, les propriétés suivantes sont équivalentes  :
\begin{enumerate}[i) ]
\item pour toute bijection $\tau$ de $\N$, la série $\sum u_{\tau(n)}$ est convergente,
\item pour toute suite $(\vartheta_n)_{n \in \N}$ à valeurs dans $\{-1,+1\}$, la série $\sum \vartheta_n u_n$ converge,
\item la famille $\sum u_n$ est sommable.
\end{enumerate}
On dit alors que la série $\sum u_n$ converge inconditionnellement.
Dans le point i), les sommes des séries $\sum u_{\tau(n)}$ sont égales et on dit que $\sum u_n$ converge inconditionnellement vers cette valeur.
\end{prop}

On définit alors de façon évidente la notion de suite basique inconditionnelle.
\begin{defi}
Une suite $(u_n)_{n\in \N}$ de $B$ est basique inconditionnelle si elle est basique et si pour tout $b\in \overline{\mbox{Vect}\{u_n,n\in \N\}}$ la série $\sum u_n^{\star}(b) u_n$ converge inconditionnellement vers $b$. 
\end{defi}

Le lien avec la proposition \ref{basincoLp} est fait par le théorème de sélection de \bess, on pourra consulter \cite[Theorem 3]{bessaga1958} ou \cite[Théorème II.6, page 57]{queff}.

\begin{theo}[\bess]
On suppose que $B$ admet une base de Schauder inconditionnelle, alors toute suite semi-normalisée $(u_n)_{n\in \N}$ de $B$ qui converge faiblement vers $0$ admet une sous-suite basique inconditionnelle.
\end{theo}

Nous aurons besoin du lemme élémentaire suivant qui quantifie la notion
de suite basique inconditionnelle.

\begin{lemm}\label{inco-lemm}
Soit $(u_n)_{n\in \N}$ une suite basique inconditionnelle de $B$, alors il existe une constante $K\geq 1$ telle que pour toute suite complexe $(\alpha_n)_{n\in \N}$ et toute suite $(\vartheta_n)_{n\in \N}$ à valeurs dans $\{-1,+1\}$ on a 
\bema 
\forall N\in \N\qquad  \frac{1}{K}\left\Vert \sum_{n=0}^{N}  \alpha_n u_n \right\Vert_{B}\leq \left\Vert \sum_{n=0}^{N} \vartheta_n \alpha_n u_n \right\Vert_{B} \leq K\left\Vert \sum_{n=0}^{N}  \alpha_n u_n \right\Vert_{B}.
\enma
\end{lemm}
\begin{demo}
Remarquons que, quitte à remplacer $\alpha_n$ par $\vartheta_n \alpha_n$, il nous suffit de prouver la majoration.
Sans perte de généralité, on peut supposer que $(u_n)_{n\in \N}$ est une base de Schauder de $B$ et se restreindre aux suites $(\alpha_n)_{n\in \N}$ telles que la série
$\sum \alpha_n u_n$ converge dans $B$.
En particulier, la série $\sum \alpha_n u_n$ est de la forme $\sum u_n^\star(b) u_n$, avec $b\in B$, et converge inconditionnellement. On en déduit les deux points :
\begin{enumerate}[a) ]
\item pour toute suite $(\vartheta_n)_{n \in \N}$ à valeurs dans $\{-1,+1\}$, la série $\sum \vartheta_n \alpha_n u_n$ converge,
\item pour tout $u^\star \in B'$, la série $\sum u^{\star}(\alpha_n u_n)$ converge inconditionnellement dans $\R$.
\end{enumerate}
D'après la proposition \ref{prop-inco}, on a $\sum |\alpha_n u^\star(u_n)|<+\infty$.
Le théorème du graphe fermé montre alors que l'opérateur linéaire
\bemar
T_{(\alpha_n)} : B' & \rightarrow & \ell^1(\N) \\
u^\star & \mapsto & (\alpha_n u^\star (u_n))_{n\in \N}
\enmar
est borné. Pour tout entier $N\in \N$, on déduit que l'on a
\bema
\left\Vert \Sum_{n=0}^{N} \vartheta_n \alpha_n u_n \right\Vert_{B}=
\sup\limits_{\substack{u^\star \in B' \\ \left\Vert u'\right\Vert \leq 1}}\left\vert \Sum_{n=0}^{N} \vartheta_n \alpha_n u^\star(u_n) \right\vert \leq \left\Vert T_{(\alpha_n)}\right\Vert_{B'\rightarrow \ell^1(\N)}.
\enma
Comme la majoration précédente est indépendante de la suite $(\vartheta_n)$, le théorème de Banach-Steinhauss avec $(\alpha_n)=(u_n^\star(b))$ nous assure que les opérateurs linéaires
\bemar
M_{(\vartheta_n,N)} : B & \rightarrow & B \\
b& \mapsto &    \Sum_{n=0}^N \vartheta_n u_n^\star(b) u_n
\enmar 
sont uniformément bornés. Cela achève notre preuve.
\end{demo}

Lorsque $B$ est admet un type ou un cotype non trivial, on déduit de ce qui précède le résultat suivant.

\begin{prop}\label{basinco}
On suppose que $B$ a une base de Schauder inconditionnelle, est de cotype $q\in [2,+\infty[$ et de type $p\in [1,2]$.
On considère une suite $(u_n)_{n\in \N}$ de $B$ convergeant faiblement vers $0$ et vérifiant $\inf\limits_{n\in \N} \left\Vert u_n\right\Vert_{B}>0$.
Alors il existe une constante $K\geq 1$ et un sous-ensemble infini $I\subset \N$ tels que, pour toute suite complexe $(\alpha_n)_{n\in I}$ à support fini, on a 
\bema   \frac{1}{K} \left( \sum_{n\in I} |\alpha_n|^q \right)^{\frac{1}{q}}\leq  \left\Vert \sum_{n\in I} \alpha_n u_n\right\Vert_B  \leq  K  \left( \sum_{n\in I} |\alpha_n|^p \right)^{\frac{1}{p}}.  \enma
\end{prop}
\begin{demo}
On notera que la suite $(u_n)$ est semi-normalisée car la convergence faible implique la bornitude.
On considère une sous-suite basique inconditionnelle $(u_n)_{n\in I}$ fournie
par le théorème de \bess.
On applique alors le lemme \ref{inco-lemm}.
Cela nous donne une constante $K\geq 1$ tel que pour tout choix de signe $\vartheta_n=\pm 1$ et pour toute suite $(\alpha_n)_{n\in I}$ à support fini on a 
\bema 
\frac{1}{K}\left\Vert \sum_{n\in I} \vartheta_n \alpha_n u_n \right\Vert_{B}\leq \left\Vert \sum_{n\in I}  \alpha_n u_n \right\Vert_{B} \leq K\left\Vert \sum_{n\in I}  \vartheta_n \alpha_n u_n \right\Vert_{B}.
\enma
En considérant maintenant une suite $(\ep_n)_{n\in I}$ 
 de variables aléatoires indépendantes suivant une loi $\frac{1}{2}$-Bernoulli, on obtient immédiatement
\bema 
\frac{1}{K}\E\left[\left\Vert \sum_{n\in I}  \ep_n \alpha_n u_n \right\Vert_{B}\right]\leq \left\Vert \sum_{n\in I}  \alpha_n u_n \right\Vert_{B}  \leq K\E\left[\left\Vert \sum_{n\in I} \ep_n \alpha_n u_n \right\Vert_{B}\right].
\enma
Quitte à changer $K$, la définition \ref{def-typ} de type et de cotype achève la preuve.
\end{demo}

Il s'agit maintenant d'examiner l'espace de Banach $L^p(X)$.

\begin{prop}
Soit $X$ un espace métrique séparable muni d'une mesure borélienne $\sigma$-finie, on note $I\subset X$ l'ensemble des atomes de $\mu$.
Pour tout $p\in[1,+\infty[$, l'espace de Banach $L^p(X)$ est linéairement isométrique à $L^p(0,1)\times \ell^p(I)$ ou à $\ell^p(I)$.
\end{prop}
\begin{demo}
En écrivant $X=I\sqcup (X\backslash I)$, il est immédiat que $L^p(X)$ est isométrique à $L^p(I)\times L^p(X\backslash I)$ où $I$ et $X\backslash I$ sont munis de la mesure, notée $\mu$, de $X$. Comme $I$ est au plus dénombrable, on vérifie que $L^p(I)$ est isométrique à $\ell^p(I)$.
En ce qui concerne $L^p(X\backslash I)$, si $I=X$ il n'y a rien à dire.
Supposons $I\neq X$ et l'on considère une densité de probabilité $g\mu$ sur $X\backslash I$ où $g$ est une fonction intégrable strictement positive. L'application linéaire
\bemar
L^p(X\backslash I,\mu) & \rightarrow & L^p(X\backslash I, g \mu) \\
f & \mapsto &  g^{-\frac{1}{p}} f
\enmar
est isométrique et surjective. Il s'agit donc d'étudier $L^p(X\backslash I, g \mu)$.
Rappelons que la densité de probabilité $g\mu$ est régulière \cite[Proposition 15.11, page 407]{royden}. Si on note $(x_n)_{n\in \N}$ une suite dense de $X\backslash I$, on déduit que tout borélien $\mathcal{B}\subset X\backslash I$ peut être approché en mesure par une réunion finie de boules $B(x_n,r_n)$ avec $r_n$ rayon rationnel.
D'après \cite[Theorem 15.4, page 399]{royden}, il existe un isomorphisme de tribus $\Phi$ de la tribu des boréliens de $X\backslash I$ sur la tribu des boréliens de $[0,1]$.
En outre, $\Phi$ transporte la densité de probabilité $g\mu$ sur la mesure de Lebesgue de $[0,1]$.
Pour toute famille finie $(\mathcal{B}_n)_{n}$ de boréliens de $X\backslash I$ disjoints deux à deux, on peut définir 
\bema \widetilde{\Phi}\left( \sum_{n} \alpha_n \pun_{\mathcal{B}_n} \right)=\sum_{n} \alpha_n \pun_{\Phi(\mathcal{B}_n)}  .  \enma 
On vérifie aisément que $\widetilde{\Phi}$ se prolonge en une isométrie surjective de $L^p(X\backslash I)$ dans $L^p(0,1)$.
\end{demo}

\begin{theo}[Paley]
Soit $X$ un espace métrique séparable muni d'une mesure borélienne $\sigma$-finie, alors $L^p(X)$ admet une base de Schauder inconditionnelle pour tout $p\in ]1,+\infty[$.
\end{theo}
\begin{demo}
D'après la proposition précédente, il suffit d'examiner le cas $L^p(0,1)\times \ell^p(I)$ avec $I$ fini ou dénombrable.
La base de Haar $(f_n)_{n\in \N}$ est une base de Schauder inconditionnelle de $L^p(0,1)$ (\cite[Chapitre 6, Théorème III.8, page 261]{queff} ou \cite[Theorem 2.c.5, page 155]{linden2}).
Le point ii) de la définition \ref{prop-inco} montre d'abord que la base canonique $(e_i)_{i\in I}$ de $\ell^p(I)$ est une base de Schauder inconditionnelle, puis que $(f_n,e_i)_{(n,i)}$ est une base de Schauder inconditionnelle de $L^p(0,1)\times \ell^p(I)$.
\end{demo}

La proposition \ref{basincoLp} découle immédiatement de la proposition \ref{basinco}, du théorème de Paley et du fait que l'espace de Banach $L^p(X)$ est de type $\min(2,p)$ et de cotype $\max(2,p)$ (voir \eqref{Lp-type}).

\bibliographystyle{plain}
\bibliography{random-space}

\end{document}